\documentclass{amsart}
\usepackage{amsmath}
\usepackage{amsthm}
\usepackage{amssymb}
\usepackage{color}
\usepackage{graphicx}
\usepackage{apptools}
\usepackage{chngcntr}
\usepackage{yfonts}
\usepackage{enumitem}
\usepackage{tikz}
\usepackage[numbers,sort&compress]{natbib}
\usepackage{hyperref}
\hypersetup{
    colorlinks=true,
    linkcolor=blue,   
    citecolor=blue,   
    urlcolor=blue     
}

\usetikzlibrary{calc,arrows.meta,positioning}
\usetikzlibrary{calc,positioning,arrows.meta,shapes.geometric}

\usepackage[utf8]{inputenc}

\newtheorem{thm}{Theorem}\newtheorem{lemma}[thm]{Lemma}

\newtheorem{corollary}[thm]{Corollary}
\newtheorem{defi}[thm]{Definition}
\newtheorem{prop}[thm]{Proposition}
\newtheorem{rk}[thm]{Remark}

\AtAppendix{\counterwithin{thm}{section}}
\AtAppendix{\counterwithin{equation}{section}}

\newcommand{\rr}{{\mathbb{R}}}

\newcommand{\nn}{{\mathbb{N}}}
\newcommand{\zz}{{\mathbb{Z}}}
\newcommand{\sS}{{\mathbb{S}}}
\newcommand{\E}{{\mathbb{E}}}

\newcommand{\mM}{{\mathbb{M}}}
\newcommand{\PP}{{\mathbb{P}}}
\newcommand{\QQ}{{\mathbb{Q}}}
\newcommand{\e}{\varepsilon}
\newcommand{\vip}{\vskip.2cm}
\newcommand{\indiq}{\hbox{\rm 1}{\hskip -2.8 pt}\hbox{\rm I}}
\newcommand{\dd}{{\rm d}}

\newcommand{\cG}{{\mathcal{G}}}

\newcommand{\cZ}{{\mathcal{Z}}}
\newcommand{\cS}{{\mathcal{S}}}

\newcommand{\cU}{{\mathcal{U}}}
\newcommand{\cL}{{\mathcal{L}}}
\newcommand{\cK}{{\mathcal{K}}}
\newcommand{\cM}{{\mathcal{M}}}
\newcommand{\cN}{{\mathcal{N}}}
\newcommand{\cF}{{\mathcal{F}}}
\newcommand{\cP}{{\mathcal{P}}}

\newcommand{\cB}{{\mathcal{B}}}

\newcommand{\cC}{{\mathcal{C}}}
\newcommand{\cD}{{\mathcal{D}}}
\newcommand{\cE}{{\mathcal{E}}}
\newcommand{\cH}{{\mathcal{H}}}
\newcommand{\cI}{{\mathcal{I}}}

\newcommand{\xX}{{\mathbb{X}}}

\newcommand{\gK}{{\mathfrak{K}}}
\newcommand{\gM}{{\mathfrak{M}}}
\newcommand{\gb}{{\mathfrak{b}}}
\newcommand{\gn}{{\mathfrak{n}}}
\newcommand{\gge}{{\mathfrak{e}}}

\newcommand{\gu}{{\mathfrak{u}}}

\newcommand{\ig}{[\![}
\newcommand{\id}{]\!]}

\newcommand{\bla}{\color{black}}

\begin{document}

\title[Keller-Segel equation]
{Existence, uniqueness, Feller property and Dirichlet form of the Keller-Segel particle system up to blow up}

\author{Yoan Tardy}
\address{Polytechnique, CMAP-UMR 7641, 91120 Palaiseau, France.}
\email{yoan.tardy@polytechnique.edu}


\begin{abstract}
The Keller--Segel particle system consists of $N$ planar Brownian motions interacting through the kernel $\theta K(x)=-\theta x/\|x\|^2$, $\theta>0$. Depending on $\theta$, particles collide, and the drift is then so singular that giving a pathwise meaning to the corresponding SDE beyond collisions, and proving uniqueness in law of its solutions for every initial condition and up to blow up, were open questions: the particle system was defined through the Dirichlet form theory, and analytic methods only provide a semigroup, an approximation uniqueness and heat kernel bounds for attraction strengths far from the blow-up regime. We introduce a notion of solution of this SDE which does not rely on Dirichlet forms, and we prove existence and uniqueness in law of such solutions up to blow up. The construction relies on an excursion decomposition based on the barycenter--radius--angle decomposition of the system, and on the fact that after a collision, the angular part of the process becomes stationary. The drift term, which is not of finite variation, is then defined through the compensations induced by the symmetry of the angular dynamics. This approach also gives a fine description of the process: we prove that it is strong Markov and Feller, that it coincides with the process associated with the Dirichlet form, and that its law is explicit when all the particles start from the same position. The method seems to apply to other singular SDEs, such as Dyson-type models or signed Coulomb particle systems.
\end{abstract}

\maketitle

\tableofcontents

\section{Introduction}
\subsection{Informal introduction of the model}
Let $\theta >0$ and consider the PDE
\begin{align}\label{eq:PDE}
    \partial_t f_t(x) = \frac12\Delta f_t(x) - \theta \mathrm{div}_x (f_t(x)K\star f_t(x)),
\end{align}
with for all $x\in \rr^2$, 
\begin{align*}
    K(x) = -\indiq_{x\ne 0} \frac{x}{\|x\|^2}.
\end{align*}
This model, introduced by \cite{ks}, describes the motion of bacteria diffusing on the plane and attracting one another. A striking feature of this model is the emergence of particle clusters, depending on the value of $\theta$, due to the singularity of the kernel $K$ :
\begin{itemize}
    \item if $\theta <2$, then there exists a global-in-time solution of \eqref{eq:PDE} (see \cite{bdp}),
    \item if $\theta =2$ and $f_0$ has a finite second moment, then there exists a global in time solution of \eqref{eq:PDE} and $f_t$ converges to a Dirac with total mass as $t\to \infty$,
    \item if $\theta >2$, then there is no global in time solution in the classical sense due to the formation of a Dirac mass in finite time, see \cite{bdp}.
\end{itemize}
Regarding the last point, the work of \cite{DolbeaultSchmeiser} gives a weak notion of solution using a defect measure and a precise heuristic description of a potential solution, however the question of the existence of a well-posed notion of solution for \eqref{eq:PDE} or a rigorous description of a weak solution remains  open in the literature.  

\vip

A natural approximation of this model is the Keller--Segel particle system which can be informally written as follows: fix $N\ge 1$,
\begin{align}\label{eq:EDS}
    \mbox{ for all } i\in \ig 1,N\id,\quad \dd X^i_t = \dd B^i_t +\frac{\theta}{N}\sum_{j\ne i} K( X^i_t - X^j_t)\dd t,
\end{align}
where $(B_t)_{t\ge 0} = (B^1_t,\dots, B^N_t)_{t\ge 0}$ is a $2N$-dimensional Brownian motion. Observe that if $N=1$, this is simply a $2$-dimensional Brownian motion. In \cite[Subsection 3.13 page 4099]{ft}, it is proven that (2) is not well-posed in the classical sense at times when too many particles collide (i.e. too many particles overlap). Two approaches have been developed so far to define the particle system beyond collisions. The Dirichlet form theory provides a process for every $\theta>0$ up to blow up \cite{ft}, but for quasi-every initial condition only, and without an a priori meaning for the SDE. Analytic methods provide a semigroup, an approximation uniqueness and heat kernel bounds (see \cite{DamirKinzebulatov},\cite{BoutiahKinzebulatov}), but only for attraction strengths far from the blow-up regime. Giving a pathwise meaning to (2) beyond collisions, and proving uniqueness in law of its solutions for every initial condition and in the whole range of parameters up to blow up, were thus open questions.

\vip

In this paper, we settle these questions: under the mild assumptions $N>3\theta$ and $\mathrm{(H)}$ (see Section \ref{subsec:collisionsdim}) give a meaning to (2) beyond collisions without relying on Dirichlet forms, and we prove uniqueness in law of its solutions. The key new idea, described in Section~1.4, is that after a collision the angular part of the process becomes stationary, which allows us to define the drift through compensations, although it is not of finite variation. This also gives a fine description of the process (Feller property, explicit law when starting from a total collision), and the method seems to apply to other singular SDEs, such as Dyson-type models \cite{cepalepingle}, \cite{dyson} or signed Coulomb particle systems \cite{vanmeurs2026strong}. Our main result is Theorem \ref{mainresult}, stated in Section \ref{sec:mainresult} together with the precise notion of solution.

\subsection{Bessel processes and principal value}\label{subsec:bessel}
Let us recall (see \cite{ry} for more details) that a squared Bessel process $(V_t)_{t\ge 0}$ with dimension $\delta \in \rr$ is a nonnegative process solving the SDE
\begin{align*}
    \dd V_t = 2\sqrt{V_t}\dd W_t + \delta \dd t,
\end{align*}
where $(W_t)_{t\ge 0}$ is a standard Brownian motion. The following properties hold
\begin{itemize}
    \item if $\delta \ge 2$, then a.s. $V_t >0$ for all $t>0$,
    \item if $\delta \in (0,2)$, then a.s. $V_t$ is reflected infinitely often at $0$,
    \item if $\delta \le 0$, then a.s., $V_t$ hits $0$ and is stopped. 
\end{itemize}
We define the Bessel process with dimension $\delta$ as $(Z_t)_{t\ge 0} := (\sqrt{V_t})_{t\ge 0}$. The following properties hold, see for example \cite{bertoin(1990)}:
\begin{itemize}
    \item if $\delta >1$, then $(Z_t)_{t\ge 0}$ solves the SDE
    \begin{align}\label{eq:eqbessel}
     Z_t = Z_0 + B_t + \int_0^t\frac{\delta-1}{2 Z_s}\dd s,
    \end{align}
   in the classical sense, where $(B_t)_{t\ge 0}$ is a standard Brownian motion,
   \item if $\delta =1$, then $(Z_t)_{t\ge 0}$ is a one-dimensional reflected Brownian motion,
   \item if $\delta \in (0,1)$, then $(Z_t)_{t\ge 0}$ solves
   \begin{align}\label{eq:eqbesselpv}
     Z_t = Z_0 + B_t + \mathrm{p.v.}\int_0^t\frac{\delta-1}{2 Z_s}\dd s,
    \end{align}
    \item if $\delta \le 0$, then $(Z_t)_{t\ge 0}$ solves \eqref{eq:eqbessel} until the first hitting time of $0$.
\end{itemize}
Equations (2) and (3) share the same structure, namely the competition between a Brownian noise and a $1/r$ drift, and the collisions of (2) play the role of the zeros of the Bessel process. Our construction can be seen as a multidimensional analogue of the principal value (4), where the cancellations come from the symmetry of the angular part of the process instead of the symmetry of the excursions of a one-dimensional process.

\subsection{Collisions and dimensions}\label{subsec:collisionsdim}
We set for all nonempty $K\subset \ig 1,N\id$, all $x = (x^1,\dots,x^N) \in (\rr^2)^N$, 
$$
S_K(x) = \frac{1}{|K|}\sum_{k\in K} x^k \qquad \mbox{ and } \qquad R_K(x) = \sum_{i\in K} \|x^i-S_K(x)\|^2= \frac{1}{2|K|}\sum_{i,j\in K}\|x^i-x^j\|^2,
$$
with the convention that $R_K = 0$ whenever $K$ is a singleton. Observe that $S_K$ is the empirical mean of the particles indexed by $K$, while $R_K$ is, up to a normalizing constant, the empirical variance of these particles.

\vip

As already shown in \cite{ft} and as we will prove again, for any non empty $K\subset \ig 1,N\id$, $(R_K(X_t))_{t\ge 0}$ behaves like a squared Bessel process with dimension $d_{\theta,N}(|K|)$ where
\begin{align}\label{eq:dimension}
    d_{\theta,N}(|K|) = (|K|-1)\Big(2-\frac{\theta|K|}{N}\Big).
\end{align}
We state a slight variant of \cite[Lemma 1]{ft} which we prove in Appendix \ref{app:preuveLemma1}:
\begin{lemma}\label{lem:dimension}
Fix $\theta > 0$ and $N \geq 2$ such that $N>3\theta$. For
\[
k_0 = \left\lceil \frac{2N}{\theta} \right\rceil ,
\]
we have
\begin{align}\label{eq:dimpositive}
d_{\theta,N}(k) > 0
\quad\text{if } k \in \ig 2,k_0-1\id,
\qquad
d_{\theta,N}(k) \leq 0
\quad\text{if } k \geq k_0.
\end{align}
We also define
\[
k_1 = k_0-1,
\qquad
k_2 =
\begin{cases}
k_0-2, & \text{if } d_{\theta,N}(k_0-2) < 2,\\
k_0-1, & \text{if } d_{\theta,N}(k_0-2) \geq 2.
\end{cases}
\]
If $\theta \geq 2$, then
\begin{itemize}
\item if $k_2\ge 3$, then $d_{\theta,N}(2) \in (1,2)$,
    \item $d_{\theta,N}(k) \geq 2$ if $k \in \ig 3,k_2-1\id$,
    \item $d_{\theta,N}(k) \in (0,2)$ if $k \in \{k_2,k_1\}$,
    \item if $N/\theta \ne 2$, then $d_{\theta,N}(k_1-1)>1$,
    \item $d_{\theta,N}(k_1) = 1$ if and only if there exists $n\in \nn\setminus \{0,1,2,3\}$ such that $2N/\theta = n - (n-4)/(2n-5)$.
\end{itemize}
\end{lemma}
It is proven in \cite{ft} that, for their notion of solution, for all $K\subset \ig 1,N\id$, 
\begin{itemize}
    \item if $d_{\theta,N}(|K|)\ge 2$, no $K$-collisions (meaning the particles indexed by $K$ and only them overlap) occurs,
    \item if $d_{\theta,N}(|K|)\in (0,2)$, then infinitely many reflected $K$-collisions occur,
    \item if $d_{\theta,N}(|K|)\le 0$, then the process is killed as soon as a $K$-collision happens. Moreover if there exists such a $K$, then there exists a (random) $L\subset \ig 1,N\id$ with $d_{\theta,N}(|L|)\le 0$ such that a $L$-collision occurs. 
\end{itemize}
Our focus will be on defining the particle system. As in the Bessel case (see Subsection \ref{subsec:bessel}), the main difficulty lies in understanding $K$-collisions with $K\in \cB^N_\theta$ where
\begin{align}\label{eq:btheta}
\cB^N_\theta := \Big\{ K\subset \ig 1,N\id : d_{\theta,N}(|K|)\in (0,1) \Big\}.
\end{align}
According to Lemma \ref{lem:dimension}, only $k_1$-collisions can have this property (however, for some values of the parameters $\theta$ and $N$ one can have $d_{\theta,N}(k_1)>1$). Throughout all the paper, we work under the assumption 
$$
\mathrm{(H)} : \quad 2N/\theta \notin \{ n - (n-4)/(2n-5), n\in \nn \setminus \{ 0,1,2,3\} \}.
$$
This restriction is inessential and spares us treating the case where the dimension is $1$. The technique would be different in this case, since it would have to be compared with a reflected standard Brownian motion : it should involve local times in the SDE in a way that takes into account the geometry of the collisions. We have chosen not to pursue this case further in this already lengthy paper. However we are optimistic about the fact that this should not be more difficult than the present work. 

\subsection{Methodological originality}
The proof relies heavily on the ''barycenter-radius-angle'' decomposition proved in \cite{ft}, see Theorem \ref{thm:decompo}. Roughly speaking, if $(X_t)_{t\ge 0}$ solves \eqref{eq:EDS}, then, as suggested in Subsection \ref{subsec:collisionsdim}, the process $(R_{\ig 1,N\id}(X_t))_{t\ge 0}$ is a squared Bessel process with dimension $d_{\theta,N}(N)$. On every excursion $(a,b)$ of $(R_{\ig 1,N\id}(X_t))_{t\ge 0}$, we can write $X_t = ( M_t ,\dots, M_t)+ D_t^{1/2} U_{A_t}$ where $A_t = \int_{(a+b)/2}^t\dd s/D_s$, $(M_t)_{t\ge 0}$ is a $2$-dimensional Brownian motion with diffusion parameter $N^{-1/2}$ and the processes $(M_t)_{t\ge 0}$, $(D_t)_{t\ge 0}$ and $(U_t)_{t\in \rr}$ are independent. The process $(U_t)_{t\in \rr}$ is what we will call an eternal spherical Keller-Segel process, and since this process is indexed by $\rr$ and is ergodic (see Proposition \ref{prop:ergodicity}), then its marginals are distributed according to its invariant measure, which is crucially even. Moreover, the angles associated with different excursions are independent because of a loss of memory after a collision, due to ergodicity. This phenomenon was already investigated by \cite{fj} in the case $N=2$. Thus, one can informally write if $T_1$ and $T_2$ are two ordered zeros of $(D_t)_{t\ge 0}$,
$$
\int_{T_1}^{T_2}b(X_s)\dd s = \sum_{i\ge 0 : a_i\in [T_1,T_2] } \int_{a_i}^{b_i} \frac{b(U^i_{A_s})}{D^{1/2}_s} \dd s = \sum_{i\ge 0 :  a_i\in [T_1,T_2] } \int_{-\infty}^{\infty} Z^i_vb(U^i_v)\dd v,
$$
where we have performed the change of variables $A_s = v$ in the integral and $(Z^i_v)_{v\in \rr}$ is a $(D_t)_{t\in (a_i,b_i)}$-measurable process. Observe that this computation is highly informal since the excursions of a squared Bessel process have a more intricate structure, however this provides a good intuition of the phenomenon. Now we see that there is a chance that the integral has meaning since we end up with a sum of i.i.d. random variables with $0$ mean, since $b$ is odd, the law of $U^i_v$ is even for all $i\ge 0$, all $v\in \rr$, the excursions of the squared Bessel process are independent from one another and the squared Bessel process is independent of $\{(U^i_t)_{t\in \rr},i\ge 0\}$. Rigorously, we study a Lévy process which plays the role of the infinite sum of i.i.d. variable in this informal computation (see Lemmas \ref{lem:drifttendzero} and \ref{lem:EDSpv}). The integral of the drift is thus not absolutely convergent, but convergent thanks to cancellations due to the angular component.

\vip

The difficulty is to construct the processes $(U^i_t)_{t\in \rr}$ which involve the same singularity as the one appearing in $(X_t)_{t\ge 0}$, but with the additional constraint of living on $\mathbb S$, which helps. As mentioned in Section 1.1, this pathwise approach developed here to give meaning to the drift in an {\it a priori} ill-posed SDE seems new and adaptable to many other SDEs. Moreover, this construction gives a good understanding of the underlying process. As an example, we will derive the explicit law of the process when starting from an initial configuration where all the particles overlap. 

\subsection{The Dirichlet form theory}
Although this is not the main focus of our work, we will also deal with the Dirichlet form theory, see \cite{f}. We will in particular show that the process we build is the one considered in \cite{ft} by showing that it is associated with the same Dirichlet space in the sense of \cite{f}. Moreover, the use of this theory is crucial in our approach to identify the invariant measure of $SKS$-process rigorously.

\vip

As a by-product, the results of \cite{ft} apply to the process constructed in this paper, and conversely, the process of \cite{ft}, which was only defined for quasi-every starting point, can be started from every point and is Feller. We will frequently refer to \cite[Appendix B]{ft} for the main ingredients of the theory used throughout this paper, as they are presented there in a particularly concise and convenient form. These results build on the works \cite{cz}, \cite{f}, and \cite{ly}.

\subsection{Relation between weak uniqueness for SDEs and independence of processes}
A related question, which is of independent interest, is addressed in Appendix \ref{sec:unicitefaible}. We show that two solutions of an SDE with the uniqueness in law property are independent whenever they are driven by independent noises, under the assumptions of Proposition \ref{prop:indepuniciteloi}. This question is related to the classical results of Cherny \cite{cherny2002uniqueness}, who established that weak uniqueness for an SDE implies uniqueness in law of the joint pair consisting of the solution and its driving Brownian motion. Our result provides a useful complementary statement in the setting considered here. It is a key ingredient in establishing the independence of two Keller--Segel processes driven by independent noises. Indeed, since strong existence for \eqref{eq:EDS} is only known for very small values of $\theta$ relative to $N$ (see \cite{DamirKinzebulatov}), one cannot assume that a solution is measurable with respect to its driving noise. Proposition \ref{prop:indepuniciteloi}, together with the uniqueness in law established in this paper, therefore provides a way to establish independence without relying on strong existence.

\subsection{Comparison with other results}
To our knowledge, the first attempts to define solutions of such SDE are given by \cite{fj} and \cite{cp}. In the work of \cite{fj}, the authors succeed in constructing a solution up to the first time at which three particles or more collide, meaning that their positions coincide. Moreover it is shown that for all $\theta \le 2 \frac{N-2}{N-1}$, a $k$-collision (which means a collision between exactly $k$ particles) with $k\ge 3$ almost surely never occurs. Moreover, if $\theta > 2\frac{N-2}{N-1}$, such collisions occur almost surely in finite time. They also prove that $2$-collisions occur with positive probability and show uniqueness in law for the system with $N=2$ particles.

\vip

The work of \cite{cp} is the first to use the Dirichlet form theory for defining this particular model, for roughly the same range of values of $\theta$. Finally, in \cite{ft}, the use of Dirichlet forms is pushed further and gives existence of a solution for every $\theta>0$ until the first time of blow up. However, although the result is sharp in the sense that existence is established precisely in the cases where intuition suggests it should hold, the use of this theory is not entirely satisfactory. Indeed it does not provide an {\it a priori} meaning for the SDE solved by the process. In addition, the process is defined for {\it quasi-every} starting point, with no concrete description of the non-authorized initial configurations. 

\vip

Another viewpoint is the analytic one developed in \cite{DamirKinzebulatov} and \cite{BoutiahKinzebulatov}. In \cite{DamirKinzebulatov}, the authors establish weak existence and approximation uniqueness, that is uniqueness of the limits of regularized systems, for general singular drifts, and cover almost the whole admissible range of attraction strengths in dimensions $d\ge3$. In dimension two, the kernel $K$ is not form-bounded, and they use weakly form-bounded drifts. In \cite{BoutiahKinzebulatov}, a many-particle fractional Hardy inequality yields a semigroup and upper heat kernel bounds in dimension two, for small attraction strengths. These methods apply to general kernels, but do not give a pathwise meaning to the SDE nor uniqueness in law among all solutions, and they do not reach the blow-up regime. On the contrary, we take advantage of the special structure and symmetries of the kernel $K$, which allows us to obtain a pathwise construction, uniqueness in law for every initial condition, and the optimal range of parameters up to blow up. 

\vip

Our work can be compared to the one of \cite{erickson1990continuous} where it is explained how a Markov process with independent radius and angle can be extended after hitting $0$ and under which condition this extension is unique. Their techniques rely on spectral theory through the use of resolvents, which does not give {\it a priori} pathwise information or a definition of the SDE solved by the process, contrary to our work. Another difference is that we need to build the angle process, which is given in \cite{erickson1990continuous}. 
\vip

The results of \cite{bethencourt2025reflected}, which construct a reflected fractional diffusion in a convex domain, provide an interesting point of comparison with the present work. At a conceptual level, the two constructions share a common philosophy: both rely on a decomposition into excursions and on understanding how these excursions can be consistently pieced together. The main difference lies in the nature of the technical challenges involved. In \cite{bethencourt2025reflected}, the analysis addresses a different structural setting, in which the absence of a Langevin-type structure presents its own difficulties. For example it is not clear {\it a priori} what an excursion looks like and how to connect them together. In the present work, the Langevin structure of \eqref{eq:EDS} provides additional structure that makes this question easier, but the singular nature of the dynamic requires finer integrability control near collision times in order to justify the rigorous definition of an excursion.

\vip

Let us also mention the work of \cite{BJW} which, even if it is not the main purpose of their work, succeeds in providing a notion of {\it entropic solution} to the Liouville equation that should be solved by $\cL(X^1,\dots,X^N)$ on the torus for finite free energy initial conditions. However a rigorous connection with any particle system is not clear at the moment.  

\section{Main results}\label{sec:mainresult}
We fix $N\ge 1$ and $\theta >0$ for all this section. We set for all $\eta>0$, all $k\in \nn$, the open sets of configurations

$$
E^\eta_k = \Big\{ x\in B(0,\eta^{-1/k}) :  \min_{\substack{K\subset \ig 1,N\id\\ |K|=k}}R_K(x) >\eta^{1/k} \Big\}\quad \mbox{ and } \quad E_k = \cup_{\eta >0} E_k^\eta.
$$
The set $E^\eta_k$ is intuitively the set of configurations of $N$ particles such that there is no subset of $k$ particles or more which are ''$\eta$-close''. The fact that particles are not going too far is more a technical restriction to simplify the analysis.  
\begin{rk}\label{rk:equivalencerelation}
    One example of a situation where considering $\eta^{1/k}$ instead of $\eta$ in the above setting seems more adapted is the following. Consider a configuration $(x^1,x^2,x^3)\in(\mathbb{R}^2)^3$. If $|x^1-x^2|\leq\eta^{1/2}$ and $|x^1-x^3|\leq\eta^{1/2}$, this does not guarantee that $|x^2-x^3|\leq\eta^{1/2}$, and thus the relation “being at distance at most $\eta^{1/2}$” is not an equivalence relation. However, if $\eta$ is chosen sufficiently small, the latter condition implies that $|x^2-x^3|\leq\eta^{1/3}$. Thus, the notion of being “$\eta$-close” is closer in spirit to an equivalence relation and provides a more natural way of defining clusters of particles.
\end{rk}

Notice that if $k>N$, then $E_k = (\rr^2)^N$.

\vip

 We define for all $\eta>0$, all $(\rr^2)^N\cup\{\triangle\}$-valued process $(Z_t)_{t\ge 0}$, 

\begin{gather}\label{eq:tempsarret}
\tau^\eta_k (Z) := \inf \Big\{t\ge 0: \mbox{ there exists } k'\ge k : Z_t \notin E^\eta_{k'} \Big\},\\
\tau_k(Z) := \lim_{\eta\to 0} \tau^\eta_k(Z),
\end{gather}
where $\triangle$ is a cemetery point. Observe that the limit is well defined thanks to monotonicity. We will abusively write $\tau^\eta_k$ and $\tau_k$ when the underlying process is obvious according to the context.

\vip

Since the processes we consider will be killed and sent to a cemetery point $\triangle$, if $E$ is a locally compact separable metrizable space, we need to put a topology on $E_\triangle := E\cup \{ \triangle \}$. We endow it with the one-point compactification topology consisting of all the open sets of $E$ and of all sets of the form $K^c\cup \{ \triangle\}$ with $K$ compact in $E$, see \cite[page 69]{f}.
\begin{rk}
It will be implicit in all this paper that when considering a continuous $E_\triangle$-valued process $(Z_t)_{t\ge 0}$, we have $Z_t = \triangle$ for all $t\ge \inf\{ s\ge 0 : Z_s = \triangle\}$.
\end{rk}
 Before introducing our main definition, we need some preparation. Let 
 $$
 H=H_N := \Big\{ x = (x^1,\dots,x^N)\in (\rr^2)^N : \sum_{i=1}^N x^i =0\Big\}
 $$ 
 and $\sS_N$ the sphere of radius $1$ of $H_N$. We will omit the dependence of $H_N$ and $\sS_N$ with respect to $N$ when the context makes it obvious to understand the underlying parameter. Let for all $x\in (\rr^2)^N$, $\pi_H(x)$ the orthogonal projection on $H$. Writing 
 $$
 \gamma(x)=\gamma_N(x) = (S_{\ig 1,N\id}(x),\dots,S_{\ig 1,N\id}(x))\in (\rr^2)^N,
 $$
  we see that $\pi_H(x) = x-\gamma (x)$. Moreover for all $x\in (\rr^2)^N$, we set $\pi_{x^\perp} = I_d - x\otimes x / \|x\|^2$ the orthogonal projection on $x^\perp$. We are now ready to give the following. Recall Lemma \ref{lem:dimension} for the definition of $k_0$.
\begin{defi}[$KS$ and $SKS$-processes]\label{def:KS/KSK}
(a) For all $\mu\in \cP((E_{k_0})_\triangle)$, we say that a tuple 
$$
(\Omega,\cF,(\cF)_{t\ge 0},\PP, (B_t)_{t\ge 0},(X_t)_{t\ge 0})
$$
is a $KS_\mu(\theta,N)$-process or equivalently a $KS(\theta,N)$-process starting with initial distribution $\mu$ if
\begin{itemize}
\item $(\Omega,\cF,(\cF)_{t\ge 0},\PP)$ is a filtered probability space on which $(B_t)_{t\ge 0}$ is a $2N$-dimensional Brownian motion,
\item the process $(X_t)_{t\ge 0}$ takes its values in $(E_{k_0})_\triangle$, is continuous, $(\cF_t)_{t\ge 0}$-adapted and killed upon exiting $E_{k_0}$ where $k_0 = \lceil 2N/\theta\rceil$,
\item we have $\int_0^\infty \indiq_{X_t \notin E_2\cup\{\triangle\}} \dd t = 0$ a.s.,
\item $X_0$ is $\mu$-distributed,
\item for all $T>0$, there exists $f :\rr_+\to \rr_+$ which is bounded by above and by below by a positive constant such that $\int_0^T \|b_{\e f(R_{\ig 1,N\id}(X_s))}(X_s)\|\dd s <\infty$ a.s. and
\begin{align}\label{eq:EDSproba}
\sup_{t\in [0,T\wedge \tau_{k_0})}\Big\|X_t - X_0 - B_t -\int_0^t b_{\e f(R_{\ig 1,N\id}(X_s))}(X_s)\dd s\Big\|\to 0 \quad \mbox{ in probability as } \e \to 0,
\end{align}
where for all $x\in E_2$, all $\e>0$, $b_\e(x) = (b^1_\e(x),\cdots,b^N_\e(x))$ with for all $i\in \ig 1,N\id$,
\begin{align}\label{eq:bepsilon}
b^i_{\e} (x) = b^i(x)\prod_{K\in \cB^N_\theta : i\in K} \indiq_{R_K(x)>\e}
\end{align}
with $\cB^N_\theta$ defined in \eqref{eq:btheta} and for all $i\in \ig 1,N\id$,
\begin{align}\label{eq:drift}
b^i(x) = -\frac\theta N\sum_{j\ne i} \frac{x^i-x^j}{\|x^i-x^j\|^2}.
\end{align}
Finally, we set $b_0(x)= b(x) = (b^1(x),\dots b^N(x))$. 
\end{itemize}

\vip

(b) For all $\mu\in \cP((\sS_N)_\triangle)$, we say that a tuple 
$$
(\Omega,\cF,(\cF)_{t\ge 0},\PP, (B_t)_{t\ge 0},(U_t)_{t\ge 0})
$$
is a $SKS_\mu(\theta,N)$-process or equivalently a $SKS_\mu(\theta,N)$-process with initial distribution $\mu$ if
\begin{itemize}
\item $(\Omega,\cF,(\cF)_{t\ge 0},\PP)$ is a filtered probability space on which $(B_t)_{t\ge 0}$ is a $2N$-dimensional Brownian motion,
\item the process $(U_t)_{t\ge 0}$ takes its values in $(\sS_N\cap E_{k_0})_\triangle$, is continuous, $(\cF_t)_{t\ge 0}$-adapted and killed upon exiting $\sS_N\cap E_{k_0}$,
\item we have $\int_0^\infty \indiq_{U_t \notin E_2\cup\{\triangle\}} \dd t = 0$ a.s.,
\item $U_0$ is $\mu$-distributed,
\item for all $T>0$, we have $\int_0^T \|b_{\e}(U_s)\|\dd s <\infty$ a.s. and
\begin{align}\label{eq:EDSprobaU}
\sup_{t\in [0,T\wedge \tau_{k_0})}\Big\|U_t - U_0 - \int_0^t \pi_{U_s^\perp}\pi_{H_N}\dd B_s -\int_0^t\pi_{U_s^\perp}\pi_{H_N}b_{\e}(U_s)\dd s +\frac{2N-3}{2} \int_0^tU_s\dd s \Big\|\to 0,
\end{align}
in probability as $\e \to 0$.
\end{itemize}
\end{defi}
\begin{rk}\label{rk:dimensions}
(i) The requirement of the existence of a bounded $f_\eta$ in \eqref{eq:EDSproba} may seem artificial since the most natural choice for a principal value would be to take $f_\eta =1$. However for technical reasons due to the passage from the $KS$-process to the $SKS$-process (see Section \ref{sec:decompo}), one needs more flexibility in the definition of a $KS$-process. Observe that adding this in the definition improves slightly our uniqueness result. 

\vip

(ii) It is shown in \cite[Subsection 3.13]{ft} that the SDE of the particle system cannot have a solution in the classical sense because the drift term would explode. An intuitive way to think about it is the following : consider a subset of particles indexed by $K\in \cB^N_\theta$ colliding and since these particles are far away from the other, we do the approximation that they don't interact (this will be described rigorously in Proposition \ref{prop:girsanoveta}). Then, the empirical variance of particles indexed by $K$  would behave like squared Bessel processes of dimension $d_{\theta,N}(|K|)\in (0,1)$. Thus, defining the particle system should involve similar difficulties as when defining a Bessel process of dimension belonging to $(0,1)$ such as subtracting the local time at $0$ in a proper way in the SDE, see \cite[Section 3]{z}. This is exactly the role of the principal value we introduce. 

\vip

(iii) Observe that according to Lemma \ref{lem:dimension}, the classical Keller-Segel SDE \eqref{eq:EDS} is recovered on the time interval $[0,\tau_{k_1})$. In particular, if $d_{\theta,N}(N)>1$, then this gives a global in time solution of the Keller-Segel SDE in the classical sense. 

\vip

(iv) If for all $k\in \ig 2,N\id$ we have $d_{\theta,N}(k)>0$, then $\lim_{\eta \to 0}\tau^\eta_{k_0} = \infty$ a.s.

\vip

(v) If $x\in (E_{k_0})_\triangle$ and $\mu = \delta_x$, then we write $KS_x(\theta,N)$ and speak about a $KS(\theta,N)$-process starting at $x$. The same comment holds for $SKS(\theta,N)$-processes.
\end{rk}
We introduce what will be respectively the density with respect to Lebesgue on $(\rr^2)^N$ of the reversible measures of a $KS(\theta,N)$-process and the density with respect to the surface measure on $\sS_N$ of the reversible measure of a $SKS(\theta,N)$-process :
\begin{align}\label{eq:invmeasureX}
\mbox{ for all } x \in E_2, \quad \mu^N_\theta (x) = \prod_{1\le i\ne j\le N} \|x^i-x^j\|^{-\theta/N},
\end{align}
and 
\begin{align}\label{eq:invmeasureU}
\mbox{ for all } u \in E_2\cap \sS_N, \quad \beta^N_\theta (u) = \prod_{1\le i\ne j\le N} \|u^i-u^j\|^{-\theta/N}.
\end{align}
This is intuitively clear for a $KS(\theta,N)$-process since it is morally a Langevin process on $(\rr^2)^N$ with potential
$$
\frac\theta N\sum_{1\le i\ne j\le N} \log(\|x^i-x^j\|).
$$
We will abusively confuse the densities $\mu^N_\theta$, $\beta^N_\theta$ and their associated measure. We will need the following result.
\begin{lemma}\label{lem:radon}
Recall that $k_0 = \lceil 2N/\theta\rceil$. Then if $N>\theta$,

\vip

(i) the measure $\mu^N_\theta$ is Radon on $E_{k_0}$,

\vip

(ii) the measure $\beta^N_\theta$ is Radon on $\sS_N\cap E_{k_0}$, in particular, if $k_0\ge N$, then $\beta^N_\theta (\sS_N) < \infty$.

\vip

(iii) Let $k\ge 2$ such that $d_{\theta,N}(k)>1$, then for all $\cC$ compact set of $E_{k+1}$,
$$
\int_{\cC} \frac{1}{\|x^1-x^2\|}\mu^N_\theta(\dd x) <\infty \quad \mbox{ and } \quad \int_{\sS_N\cap \cC} \frac{1}{\|u^i-u^j\|} \beta^N_\theta(\dd u) <\infty.
$$
\end{lemma}
Points (i) and (ii) are already proven in \cite[Proposition A.1-Proposition A.3]{ft}. We add the proof of (iii) in the Appendix \ref{app:preuveradon}.

\vip

Before stating Theorem \ref{mainresult}, which is the main theorem of this paper, we need to introduce an eternal version of the $SKS$-process since it will be an important part of the construction of a $KS$-process.
\begin{defi}\label{def:eternalSKS}
Suppose that $k_0 \ge N$. We say that a tuple 
$$
(\Omega,\cF,(\cF)_{t\in\rr},\PP, (B_t)_{t\in \rr},(U_t)_{t\in \rr})
$$
is an eternal $SKS(\theta,N)$-process if for all $t_0 \in \rr$, the tuple 
$$
(\Omega,\cF,(\cF)_{t\ge t_0},\PP, (B_t-B_{t_0})_{t\ge t_0},(U_t)_{t\ge t_0})
$$
is a $SKS(\theta,N)$-process with $U_{t_0}$ following the law $Z_{N,\theta}^{-1}\beta^N_\theta \dd \sigma$ on $\sS_N$ where $\sigma$ is the surface measure on $\sS_N$ and $Z_{N,\theta}$ is a normalizing constant (see Lemma \ref{lem:radon}).
\end{defi}

We recall that if $E$ is a locally compact separable metrizable space and $m$ is a Radon measure on $E$, we say that a Markov process $\mM = (\Omega,\cF,(\cF_t)_{t\ge 0},(Z_t)_{t\ge 0},(\PP_x)_{x\in E_\triangle})$ is $m$-symmetric if, for all $f,g :E \to \rr_+$ we have
$$
\int_E \E_x[f(Z_t)\indiq_{Z_t\ne \triangle}]g(x) m(\dd x) = \int_E f(x)\E_x [g(Z_t)\indiq_{Z_t\ne \triangle}] m(\dd x),
$$
where $\E_x$ denotes the expectation with respect to $\PP_x$. Moreover if $\mM$ is a continuous $m$-symmetric strong Markov process, it is associated in the sense of \cite{f} to a Dirichlet space on $L^2(E,m)$, see \cite[Appendix A]{ft}. Moreover, we recall from \cite[Appendix A]{ft} that for any Dirichlet space $(\cE,\cF)$ on $L^2(E,m)$, a core is a set $\cC\subset C_c^\infty (E)\cap\cF$ such that $\cC$ is dense in $C_c^\infty (E)$ for the uniform convergence topology and dense in $\cF$ for the norm induced by the scalar product $\cE_1 := \cE(\cdot, \cdot) + (\cdot,\cdot)_{L^2(E,m)}$.

\vip

We denote by $\gK_N$ the set of the partitions of $\ig 1,N\id$ and for all $p\in \nn^*$, $\gK^p_N$ the set of partitions $\cK$ of $\ig 1,N\id$ such that for all $K\in \cK$, $|K|\le p$.  We introduce for all $\cK \in \gK_N$,
\begin{align*}
    \cC_\cK = \Big\{ \varphi : x&\in (\rr^2)^N\mapsto \psi((R_K(x),S_K(x))_{K\in \cK}) \\
    &\mbox{ with } \psi \mbox{ smooth and compactly supported such that } \varphi \in C^\infty_c (E_\cK) \Big\},
    \end{align*}
    where
    $$
    E_\cK = \Big\{ x\in (\rr^2)^N : \inf_{K\ne L\in \cK} \inf_{i\in K,j\in L} \|x^i-x^j\| >0 \Big\}.
    $$
We are now ready to state our main result.
\begin{thm}[$\mathrm{Main Theorem}(\theta,N)$]\label{mainresult}
Suppose that $N>3\theta$ and $\mathrm{(H)}$ holds.

\vip

(i) For all $\mu\in \cP((E_{k_0})_\triangle)$, there exists a unique in law $KS_\mu(\theta,N)$-process, we denote by $\PP^{X,N,\theta}_\mu$ its law (and $\PP^{X,N,\theta}_x$ if $\mu = \delta_x$ with $x\in (E_{k_0})_\triangle$). Such a process has a.s. a left-limit at $\tau_{k_0}$ for the usual topology on $(\rr^2)^N$. Moreover, one can build a filtered measure space $(\Omega,\cF,(\cF_t)_{t\ge 0})$ such that $\mM^{X,N,\theta} = (\Omega,\cF,(\cF_t)_{t\ge 0}, (X_t)_{t\ge 0}, (\PP^{X,N,\theta}_x)_{x\in (E_{k_0})_\triangle})$ strong Markov process such that for all $A\in \cB(C(\rr_+,(E_{k_0})_\triangle))$, $x\in (E_{k_0})_\triangle\mapsto \PP_x^{X,N,\theta}(A)$ is measurable. Moreover, it is Feller in the sense that for all $\varphi \in C_0(E_{k_0})$, 
\begin{gather*}
    x\in (E_{k_0})_\triangle \mapsto \E_x^{X,N,\theta}[\varphi (X_t)] \in C_0(E_{k_0})\quad \mbox {and} \quad \lim_{t\to 0} \sup_{x\in (E_{k_0})_\triangle} |\E_x^{X,N,\theta}[\varphi (X_t)] - \varphi (x)|= 0. 
\end{gather*}
Finally $\mM^{X,N,\theta}$ is $\mu^N_\theta$-reversible so that it is associated to a Dirichlet space $(\cE_X^N,\cF_X^N)$ on $L^2(E_{k_0},\mu^N_\theta)$ which is given by
\begin{align*}
\cF_X^N &= \overline{ C_c^\infty (E_{k_0})}^{\cE_{X,1}^N},\\
\mbox{ for all } f, g \in \cF_X^N,& \quad \cE_X^N(f,g) = \frac12 \int_{E_{k_0}} \nabla f \cdot \nabla g \dd \mu^N_\theta.
\end{align*}
Moreover, $\sum_{\cK \in \gK^{k_0-1}_N}\cC_\cK$ is a core for $(\cE^N_X,\cF^N_X)$.  

\vip

(ii) For all $\mu\in \cP((E_{k_0}\cap\sS_N)_\triangle)$, there exists a unique in law $SKS_\mu(\theta,N)$-process, we denote by $\PP^{U,N,\theta}_\mu$ its law (and $\PP^{U,N,\theta}_u$ if $\mu = \delta_u$ with $u\in (E_{k_0}\cap\sS_N)_\triangle$). Such a process has a.s. a left-limit at $\tau_{k_0}$ for the usual topology on $\sS_N$. Moreover, one can build a filtered measure space $(\Omega,\cF,(\cF_t)_{t\ge 0})$ such that $\mM^{U,N,\theta} = (\Omega,\cF,(\cF_t)_{t\ge 0}, (X_t)_{t\ge 0}, (\PP^{U,N,\theta}_u)_{u\in (E_{k_0}\cap\sS_N)_\triangle})$ is a strong Markov process such that for all $A\in \cB(C(\rr_+,(\sS\cap E_{k_0})_\triangle))$, $u\in (\sS\cap E_{k_0})_\triangle\mapsto \PP_u^{U,N,\theta}(A)$ is measurable. Moreover, it is Feller in the sense that for all $\varphi \in C_0(E_{k_0}\cap \sS)$, 
\begin{gather*}
    u\in (E_{k_0}\cap\sS)_\triangle \mapsto \E_u^{U,N,\theta}[\varphi (U_t)] \in C_0(E_{k_0}\cap \sS)\quad \mbox {and} \quad  \lim_{t\to 0}\sup_{u\in (E_{k_0}\cap \sS)_\triangle} |\E_u^{U,N,\theta}[\varphi (U_t)] - \varphi (u)| = 0. 
\end{gather*}
Finally, $\mM^{U,N,\theta}$ is $\beta^N_\theta$-reversible so it is associated to a Dirichlet space $(\cE_U^N,\cF_U^N)$ on $L^2(E_{k_0}\cap \sS_N,\beta^N_\theta)$ which is given by
\begin{align*}
\cF_U^N &= \overline{ C_c^\infty (E_{k_0}\cap \sS_N)}^{\cE_{U,1}^N},\\
\mbox{ for all } f, g \in \cF_U^N,& \quad \cE_U^N(f,g) = \frac12\int_{E_{k_0}\cap \sS_N} \nabla_{\sS_N} f \cdot \nabla_{\sS_N} g \dd \beta^N_\theta,
\end{align*}
where $\nabla_{\sS_N}f = \nabla [f\circ \Phi]$ in the sense of distributions with 
$$
\mbox{ for all } x\in E_N, \quad \Phi (x) = \pi_{H_N}(x)/\|\pi_{H_N}(x)\|.
$$
\end{thm}
\begin{rk}\label{rk:siKE2}
   (i) Observe that if $\cK = \{\{1\},\dots,\{N\}\}$ (which belongs to $\gK^1)$, then $\cC_\cK = C_c^\infty (E_2)$. Of course since $C^\infty_c(E_2)$ is not a sufficiently rich space of test functions to describe the dynamic because the collisions are excluded, we need to add the linear combination with other test functions that describe the behavior of our process around collisions. In the Dirichlet form theory language, we design a {\it core} of the Dirichlet form which is more adapted to the dynamic.

   \vip

   (ii) While it is plausible that $\sum_{\cK\in \gK_N^{(k_0-1)\wedge(N-1)}} \cC^\sS_\cK$ is a core for $(\cE^N_U,\cF^N_U)$, where for all $\cK\ne \{\ig 1,N\id\}$,
$$
\cC_{\cK}^\sS = \Big\{ \varphi : x\in \sS \mapsto \psi ((R_K(u),S_K(u))_{K\in \cK} \mbox{ with } \varphi \in C_c^\infty (E_\cK)\Big\}, 
$$
this property is not required for our proofs. We therefore omit its verification for the sake of brevity.
\end{rk}
In addition we can describe the dynamic of the $KS$ and $SKS$-processes in terms of some analogue of the Ito's formula. Let us first introduce the infinitesimal generators of a $KS(\theta,N)$-process and a $SKS(\theta,N)$-process : $\cL^X := \cL^{X,N,\theta}$ and $\cL^U := \cL^{U,N,\theta}$ are defined through
\begin{gather}
    \mbox{ for all } \varphi \in C^2((\rr^2)^N), \mbox{ all } x\in E_2, \quad \cL^X\varphi (x) = \frac{1}{2} \Delta \varphi (x) -\frac{\theta}{N} \sum_{1\le i\ne j\le N} \frac{x^i-x^j}{\|x^i-x^j\|^2} \cdot \nabla_{x_i} \varphi (x) \label{eq:infinitesimalgeneratorX}\\
    \mbox{ for all } \varphi \in C^2(\sS), \mbox{ all } u\in \sS\cap E_2, \quad \cL^U\varphi (u) = \frac{1}{2} \Delta (\varphi\circ \Phi) (u) -\frac{\theta}{N} \sum_{1\le i\ne j\le N} \frac{u^i-u^j}{\|u^i-u^j\|^2} \cdot \nabla_{u_i} (\varphi\circ \Phi)  \label{eq:infinitesimalgeneratorU}
\end{gather}
with $\Delta$ and $\nabla$ in \eqref{eq:infinitesimalgeneratorU} understood as the operators on $(\rr^2)^N$. We also introduce their smoothed version $\cL_\alpha^X$ and $\cL^U_\alpha$ for every $\alpha >0$ defined through,
\begin{align}
    \mbox{ for all } \varphi \in C^2((\rr^2)^N), \mbox{ all } x\in (\rr^2)^N, \quad \cL^X_\alpha\varphi (x) = \frac{1}{2} \Delta \varphi (x) -\frac{\theta}{N} \sum_{i=1}^N\sum_{j\ne i} \frac{x^i-x^j}{\|x^i-x^j\|^2+\alpha} \cdot \nabla_{x_i} \varphi (x), \label{eq:smoothedinfinitesimalgeneratorX}\\
     \mbox{ for all } \varphi \in C^2(\sS), \mbox{ all } u\in \sS, \quad \cL^U_\alpha\varphi (u) = \frac{1}{2} \Delta (\varphi \circ \Phi) (u) -\frac{\theta}{N} \sum_{i=1}^N \sum_{j\ne i}\frac{u^i-u^j}{\|u^i-u^j\|^2+\alpha} \cdot \nabla_{u_i} (\varphi\circ \Phi) .\label{eq:smoothedinfinitesimalgeneratorU}
\end{align}
\begin{rk}
    Observe that if $x\in E_N$ with $x = \psi(m,r,u)$ where $m\in \rr^2$, $r>0$ and $u\in \sS$, it is classical that $\Delta (\varphi \circ \Phi)(x) = r^{-1}\Delta_\sS\varphi (u)$ with $\Delta_\sS$ exactly the operator of Laplace-Beltrami on $\sS$. Moreover we will sometimes denote $\nabla_\sS\varphi (u) = r^{1/2}\nabla (\varphi \circ \Phi)(x)$. We thus see that $\cL^U$ is a natural analogue of $\cL^X$ on $\sS$. 
\end{rk}
We can now state the following.
\begin{prop}[$\mathrm{Ito}(\theta,N)$]\label{prop:itoformula}
Suppose that $N>3\theta$ and $\mathrm{(H)}$ holds.

\vip

(i) Let $\mu\in \cP((E_{k_0})_\triangle)$ and $(\Omega,\cF,(\cF)_{t\ge 0},\PP, (B_t)_{t\ge 0},(X_t)_{t\ge 0})$ be a $KS_\mu(\theta,N)$-process. Then for any $\varphi \in C^2((\rr^2)^N)$ such that for any $\eta >0$ and $M>0$,
   $$
   \sup_{\alpha \in (0,1]}\sup_{x\in B(0,M)} \|\cL^X_\alpha \varphi (x)\| \le C_M,
   $$
   for $C_M<\infty$ only depends on $M$, then for all $t\in [0,\tau_{k_0})$,
   \begin{align}
      \varphi (X_{t}) = \varphi (X_0) +\int_0^{t} \nabla \varphi (X_s)\dd B_s + \int_0^{t} \cL^X\varphi (X_s)\dd s.
   \end{align}

   \vip

   (ii) Let $\mu\in \cP(\sS)$ and $(\Omega,\cF,(\cF)_{t\ge 0},\PP, (B_t)_{t\ge 0},(U_t)_{t\ge 0})$ ,be a $SKS_\mu(\theta,N)$-process. Then for any $\varphi \in C^2(\sS)$ such that for any $M>0$,
   $$
   \sup_{\alpha \in (0,1]}\sup_{u\in B(0,M)} \|\cL^U_\alpha \varphi (u)\| \le C_M,
   $$
   for $C_M<\infty$ only depends on $M$, then for all $t\in [0,\tau_{k_0}(U))$,
   \begin{align}\label{eq:itoU}
      \varphi (U_t) = \varphi (U_0) +\int_0^{t} \nabla (\varphi\circ \Phi) (U_s)\dd B_s + \int_0^{t} \cL^U\varphi (U_s)\dd s.
   \end{align}
\end{prop}
\begin{rk}\label{rk:coordonnéespolaires}
    (i) Observe that for all $\varphi \in C^2(\sS)$, all $u\in \sS$, 
    \begin{align}\label{eq:gradspherique}
    \nabla (\varphi \circ \Phi)(u) = \pi_H(\pi_{u^\perp}(\nabla \varphi (u))).
    \end{align}
    Moreover, if $x\in H\setminus \{0\}$, setting $\rho=\|x\|$ and $u = x/\rho$, we get for all $f :H\to \rr$ regular enough
    $$
    \Delta f (\rho,u) = \frac{\Delta_\sS f (\rho,u)}{\rho^2} + \frac{2N-3}{\rho} \partial_\rho f (\rho,u) + \partial^2_\rho f(\rho,u).
    $$
    Thus, taking $f(\rho,u) = \varphi (\rho u)$, we get,
    $$
    \Delta \varphi (x) = \frac{\Delta_\sS\varphi (\rho u)}{\rho^2} + \frac{2N-3}{\rho} u\cdot \nabla \varphi (\rho u) + u\cdot H\varphi (\rho u) (u),
    $$
    where $H\varphi$ denotes the Hessian of $\varphi$. Thus, taking $\rho =1$,
    \begin{align}\label{eq:laplacienspherique}
    \Delta_\sS \varphi (u) = \Delta \varphi (u) - (2N-3)u\cdot \nabla \varphi (u)- u\cdot H\varphi (u) (u). 
    \end{align}
    The identities \eqref{eq:gradspherique} and \eqref{eq:laplacienspherique} show that \eqref{eq:itoU} is coherent with \eqref{eq:EDSprobaU}.

    \vip

    (ii) The condition in Proposition \ref{prop:itoformula} is not claimed to be sharp. It would be interesting to understand if a more general Ito's formula could be obtained using the compensation of the angle as we do while defining the drift. 
\end{rk}
A crucial result in our paper is the following decomposition theorem. This is the key to prove Theorem \ref{mainresult} and justify the relevance of the induction reasoning we will use for the construction of a $KS$-process with $N-1$ particles to build a $KS$-process with $N$ particles (see Section \ref{sec:structure} for the strategy of the proof), since morally a $SKS$-process with $N$ particles behaves like a $KS$-process with less particles thanks to three degrees of freedom lost in the projection of the manifold $\sS$. 
\begin{thm}[$\mathrm{Decomposition}(\theta,N)$]\label{thm:decompo} Suppose that $N>3\theta$ and $\mathrm{(H)}$ hols. Let $x\in E_{k_0}$ and 
    $$
    (\Omega,\cF,(\cF_t)_{t\ge 0},\PP,(B_t)_{t\ge 0}, (X^N_t)_{t\ge 0})
    $$
    a $KS_x(\theta,N)$-process. 
    
    \vip
    
    (i) In the case where $r := R_{\ig 1,N\id}(x)>0$, we can uniquely write $x = \psi (m,r,u):= \gamma (m) + \sqrt{r} u$ with $m\in \rr^2$, $\gamma (m) = (m,\dots,m) \in (\rr^2)^N$ and $u\in E_{k_0}\cap \sS_N$. Thus 
    \begin{itemize}
\item $(M_t)_{t\ge 0} := (S_{\ig 1,N\id}(X^N_{t}))_{t\ge 0}$ is a $2$-dimensional Brownian motion with diffusion coefficient $N^{-1/2}$ starting at $m$,
\item $(D_t)_{t\ge 0}:= (R_{\ig 1,N\id}(X^N_t))_{t\ge 0}$ is the squared Bessel process starting at $r$ with dimension $d_{\theta,N}(N)$ stopped when it reaches $0$ and driven by the standard Brownian motion $(W_t)_{t\ge 0}$ defined as
\begin{align*}
\mbox{ for all } t\ge 0, \quad W_t = \int_0^t \indiq_{R_{\ig 1,N\id}(X^N_s)>0} \Phi(X^N_s)\cdot \dd \pi_{H_N}(B)_s,
\end{align*}
\item Setting $(U_t)_{t\ge 0} = (\Phi(X^N_{\rho_t}))_{t\ge 0}$, $(\Omega,\cF,(\cF_{\rho_t})_{t\ge 0}, (\tilde B_t)_{t\ge 0}, (U_t)_{t\ge 0})$ is a $SKS_u(\theta,N)$-process, where for all $t\ge 0$, $A_t = \int_0^t D_s^{-1}\dd s$ (which equals $\infty$ if $t\ge \tau_N$), $\rho_t = \inf\{s\ge 0 : A_s \ge t\}$ (which is a stopping time) and
$$
\tilde B_t = \int_{0}^{\rho_t} D_s^{-1/2} \dd B_s. 
$$
\end{itemize}
Moreover the family $((M_t)_{t\ge 0}, (D_t)_{t\ge 0}, (U_t)_{t\ge 0})$ is mutually independent. 

\vip

(ii) Suppose that $k_0>N$. Since $\cZ = \{ t\ge \tau_N : D_t = 0\}$ is a.s. a closed set, we can write $(\tau_N,\infty)\cap \cZ^c = \cup_{i\ge 1}(a_i,b_i)$ with the intervals being disjoint and the application $i \mapsto (a_i,b_i)$ measurable. Thus for all $i\ge 1$, setting for all $t\in \rr$, $U^i_t = \Phi ( X^N_{\rho^i_t})$ with for all $t\ge 0$, $A^i_t = \int_{(a_i+b_i)/2} D_s^{-1}\dd s$ (which equals $-\infty$ if $t\le a_i$ and $\infty$ if $t\ge b_i$), $\rho^i_t = \inf \{ s\ge (a_i,b_i) : A^i_s\ge t\}$, we have that knowing $(D_t)_{t\ge 0}$, the the family $\{(U^i_t)_{t\in \rr}, i\ge 1\}$ is i.i.d., independent from $(M_t)_{t\ge 0}$ with $(U^1_t)_{t\ge 0}$ distributed as an eternal $SKS(\theta,N)$-process.
\end{thm}
In the last result, we implicitly used Lemma \ref{lem:horlogeexplosion}.

\vip

We will show some additional properties of the $KS$-process.  In order to explain that, we will show the following lemma.
\begin{lemma}\label{lem:tpsarretsimple}
Let $\eta>0$. There exists a family $(\sigma^\eta_k)_{k\ge 0}$ of applications from $C(\rr_+,(\rr^2)^N\cup\{\triangle\})$ to $\rr_+\cup \{ \infty\}$ and a family $(\cK^\eta_k)_{k\ge 0}$ of elements of $\gK_N$ such that for all continuous $(\rr^2)^N$-valued random process $(Z_t)_{t\ge 0}$ defined on a filtered probability space $(\Omega,\cF,(\cF_t)_{t\ge 0},\PP)$, $\lim_{k\to \infty} \sigma^\eta_k (Z) = \inf\{ t\ge 0 : Z_t = \triangle\}$, where $\inf \emptyset = \infty$ by convention, $\sigma^\eta_0(Z) = 0$ and for all $k\ge 0$,  
\begin{itemize}
\item $\sigma^\eta_k(Z)$ is a stopping time,
\item $\cK^\eta_k(Z)$ is $\cF_{\sigma^\eta_k}$-measurable,
\item a.s. $\sigma^\eta_{k+1}(Z)-\sigma^\eta_k(Z) >0$,
\item a.s., for all $K\ne L \in \cK^\eta_k(Z)$, all $i\in K$, all $j\in L$, we have either $\sigma^\eta_k(Z) = \inf \{ t\ge 0 : Z_t = \triangle\}$ or 
$$
\mbox{ for all }t\in [\sigma^\eta_k(Z), \sigma^\eta_{k+1}(Z)],\quad \|Z^i_t-Z^j_t\|^2 \ge c\eta^{1/(|K|\vee |L|+1)},
$$
with $c>0$ a positive constant.
\end{itemize}
\end{lemma}
We refer to Section \ref{sec:config} for an explicit construction of these sequences. This quantifies the fact that at each time, we can either say that particles are all close together or not and in this case we can regroup them as disjoint subgroups of particles which are considered close from each other and far away from the other subgroups. This fact is not trivial, in particular the relation ''being close as an $\eta$ away'' is not an equivalence relation, see Remark \ref{rk:equivalencerelation}.

\vip

We now define the stopping time which will play the role of an adapted $\eta$-approximation of $\tau_{k_0\wedge N}$.
\begin{defi}\label{def:sigmastar}
We set for all $\eta >0$, all continuous $E_{k_0}\cup\{\triangle\}$-valued process $(Z_t)_{t\ge 0}$, 
$$
\sigma^\eta_*(Z) = \inf \Big\{ \sigma^\eta_k(Z) : R_{\ig 1,N\id}(Z_{\sigma^\eta_k(Z)}) \le \eta^{1/N} \quad \mbox{ or } \quad \|Z_{\sigma^\eta_k(Z)}\|\ge \eta^{-1/N} \Big\}\wedge \inf\{ t\ge 0 : Z_t = \triangle\},
$$
where by convention $\inf \emptyset = \infty$.
\end{defi}
We  use Lemma \ref{lem:tpsarretsimple} to define the notion of $\eta$-$pre$-$KS$ process. 
\begin{defi}\label{def:etapreKS}
Let $\eta \in (0,\eta_0)$, $\mu \in \cP((E_{k_0})_\triangle)$ and $(\Omega,\cF,(\cF_t)_{t\ge 0},\PP)$ a filtered probability space on which there exists a $(\cF_t)_{t\ge 0}$-adapted $(\rr^2)^{N}$-valued standard Brownian motion $(B_t)_{t\ge 0}$. We say that $(\Omega,\cF,(\cF_t)_{t\ge 0},\PP,(B_t)_{t\ge 0},(X_t)_{t\ge 0})$ is a $\eta$-$pre$-$KS_\mu(\theta,N)$-process driven by $(B_t)_{t\ge 0}$ if $X_0$ is $\mu$-distributed, it is a continuous $(\cF_t)_{t\ge 0}$-adapted $(\rr^2)^{N}\cup\{\triangle\}$-valued process stopped at $\sigma^\eta_*$ and for all $k\ge 0$, all $K\in \cK^\eta_k$, knowing $\cF_{\sigma^\eta_k}$, on the event $\{ \sigma^\eta_k < \sigma^\eta_*\}$, $((X^i_{t\wedge \sigma^\eta_{k+1}})_{t\ge \sigma^\eta_k})_{i\in K}$ is a $KS(\theta |K|/N, |K|)$-process starting at $(X^i_{\sigma^\eta_k})_{i\in K}$, stopped at $\sigma^\eta_1$ driven by $((B^i_t-B^i_{\sigma^\eta_k})_{t\ge \sigma^\eta_k})_{i\in K}$.
We write its law $\tilde \PP^{\eta,N,\theta}_\mu$ (and write $\tilde \PP^{\eta,N,\theta}_\mu$ if $\mu = \delta_x$ with $x\in (E_{k_0})_\triangle$). 
\end{defi}

\begin{rk}\label{rk:unicitepreKS}
(i) Such a process has to be informally understood as follows: at the initial configuration, one can define clusters of particles which are ''$\eta$-close''. We let each cluster evolve as independent KS-processes with suitable parameters until the clusters are not the same, for example if one cluster split in several clusters or if two clusters fuse. Then the process is stopped at this time, new partition of subgroups of particles are identified and the same procedure continues until all the particles are close. The rescaling of the parameter $\theta$ to $\theta |K|/N$ for a cluster of particles indexed by $K\subset \ig 1,N\id$ follows from rewriting the interaction coefficient $\theta/N$ in \eqref{eq:EDS} as
$$
\frac{\theta}{N}=\frac{\theta |K|/N}{|K|}.
$$
Indeed, when considering this cluster in isolation, $|K|$ plays the role of the (fictitious) total number of particles. See \cite[Subsection 3.3]{ft} for a more detailed explanation.

\vip

 (ii) Observe that a $\eta$-$pre$-$KS(\theta,N)$-process is not rigorously a Markov process. For example, take $x\in (\rr^2)^N$ such that $R_{\ig 1,N\id}(x)> \eta^{1/N}$, then $\cK^\eta(x) \ne \{\ig 1,N\id\}$. One can show that there exists $t>0$ sufficiently small such that $\tilde\PP^{\eta,N+1,\theta}_x(A_t) >0$ where
$$
A_t = \Big\{ R_{\ig 1,N\id}(X_t) < \eta^{1/N} \quad \mbox{ and } \quad \inf \Big\{ R_{\ig 1,N\id}(X_s) , s\in [0,t]\Big\} > \frac{1}{2} \eta^{1/N}\Big\},
$$
with $(X_t)_{t\ge 0}$ the canonical process. In this case, $\cK^\eta(X_t) = \{\ig 1,N\id \}$, so that under $\tilde \PP^{\eta,N,\theta}_{X_t}$, $\sigma^\eta_* = 0$ a.s., but this is not the case for $(X_{s+t})_{s\ge 0}$.

\vip

(iii) However, by construction, this is still true that for all $k\ge 0$, $\tilde \PP^{\eta,N,\theta}_{X_{\sigma^\eta_k}}$ is the law of $(X_{t+\sigma^\eta_k})_{t\ge 0}$.

\vip

(iv) Observe that since strong existence or pathwise existence is not known for the $KS$-process, it is not clear whether two $KS$-processes driven by independent Brownian motion are measurable function of their respective driving Brownian motion. This implies that independence of such processes is not clear and this is the main obstacle to show uniqueness in law of $\eta$-$pre$-$KS$-processes. However, we will succeed showing such a result in Lemma \ref{lem:etapreKS} using a crucial more general result proved in Appendix \ref{sec:unicitefaible}.

\vip

(v) One could worry about the hypothesis $2N/\theta \notin \{ n - (n-4)/(2n-5), n\in \nn \setminus \{0,1,2,3\}\}$ to still hold when considering a $KS$-process with the different parameters $\tilde \theta :=\theta |K| /N$ and $\tilde N :=|K|$ in the proof because of the induction reasoning. Fortunately, this property is preserved through all the induction since $N/\theta = \tilde N/\tilde \theta$. 
\end{rk}

The following Proposition is crucial and motivates the induction argument that we will present in the next section and which structure this whole paper.
\begin{prop}\label{prop:girsanovpropre}
Suppose that $N>3\theta$ and $\mathrm{(H)}$ holds.

\vip

(i) Let $\mu\in \cP((E_{k_0})_\triangle)$ and $(\Omega, \cF,(\cF_t)_{t\ge 0}, \PP, (B_t)_{t\ge 0}, (X^N_t)_{t\ge 0})$ a $KS_\mu(\theta,N)$-process. For all $T>0$, all $\eta >0$, there exists a probability measure $\QQ^{T,\eta}_\mu$ such that under $\QQ^{T,\eta}_\mu$, $(X^N_{t\wedge \sigma^\eta_*})_{t\ge 0}$ is a $\eta$-$pre$-$KS_\mu(\theta,N)$-process and there exist two positive constants $c_{T,\eta}$, $C_{T,\eta}$ only depending on $T,\eta$ such that
$$
c_{T,\eta} \le \E \Big[ \Big( \frac{\dd \QQ^{T,\eta}_\mu}{\dd \PP} \Big)^2 \Big]\le C_{T,\eta},
$$

\vip

(ii) the applications $\sigma^\eta_k$ for all $k\ge 0$, $\sigma^\eta_*$ and $\tau_{k_0}$ are continuous $\cL_\mu((X^N_{t\wedge \sigma^\eta_* })_{t\ge 0})$-a.s. and for all $t>0$, 
$$
\PP(\sigma^\eta_*((X^N_{t\wedge \sigma^\eta_*})_{t\ge 0}) = t) = \PP(\tau_{k_0}((X^N_{t\wedge \sigma^\eta_*})_{t\ge 0})=t) = 0.
$$ 
\end{prop}
This explains why morally, until the first collision implying all the particles, a $KS(\theta,N)$-process can be built from some $KS(\theta,k)$-processes with $k<N$ far away from each other. The interactions between these subgroups of particles are bounded because of the distance so that one can add or suppress them thanks to the Girsanov theorem. 

\vip

An amusing consequence of the construction we will operate is the following:
\begin{prop}\label{prop:densité}
    For all $\theta >0$, $N\ge 2$ such that $\mathrm{(H)}$ holds and $k_0(\theta,N)>N$, all $m\in \rr^2$, all $t>0$, the law $\pi_t\#\PP^{X,N,\theta}_{\gamma (m)}$ has a density with respect to Lebesgue given by 
    \begin{align*}
       \rho^N_t(x^1,\dots,x^N) = \frac{e^{-\|x -\gamma (m)\|^2/(2t)}}{Z2^{\nu+1}\pi t^{2+\nu} \Gamma (\nu+1)}\mu^N_\theta(x),
    \end{align*}
    where $\nu = d_{\theta,N}(N)/2-1$, $Z$ is the normalization constant of $\beta^N_\theta$ and $\pi_t$ is the evaluation at time $t$. 
\end{prop}
The proof is done in Appendix \ref{app:densité}. In other words, this result explains that for any initial configuration where all the particles start at the same point, the angular part of the process is stationary and then becomes independent from the radius, making the law explicit to compute.

\section{Structure of the paper and idea of the construction}\label{sec:structure}
We will prove Theorem \ref{mainresult} with an induction argument on $N$, the number of particles. The base case is straightforward since for all $\theta >0$ a $KS(\theta,1)$ is a Brownian motion and the rest of the paper is devoted to prove the inductive step, trying to build a $KS(\theta,N)$-process for every $\theta >0$ knowing how to build a $KS(\theta',k)$ for all $\theta'>0$ and all $k\le N-1$. 

\vip

In Section \ref{sec:config} we prove a more precise version of Lemma \ref{lem:tpsarretsimple}.

\vip

Proposition \ref{prop:girsanovpropre} is going to be the center of attention in Section \ref{sec:killedEN}, where we define uniquely in law a $KS(\theta,N)$-process until the first time at which all the particles meet. As sketched in the previous Section, the idea is to first build a $\eta$-pre-$KS(\theta,N)$ for which existence and uniqueness will come from the induction hypothesis. Using a Girsanov transform (see Proposition \ref{prop:girsanoveta}) we will recover existence and uniqueness of the killed $KS(\theta,N)$-process when exiting $E_N$. 

\vip

Section \ref{sec:decompo} will be devoted to show that one can deduce from the previous section the existence of a unique in law $SKS(\theta,N)$-process and Theorem \ref{thm:decompo}-(i) at least up to the first time all the particles collide. 

\vip

In Section \ref{sec:constructionKS}, we build a candidate to be a $KS(\theta,N)$-process in the sense of Definition \ref{def:KS/KSK}, intuitively motivated by Theorem \ref{thm:decompo}-(ii) . In this construction we will need to use eternal $SKS(\theta,N)$-processes which will exist thanks to the results established in Section \ref{sec:ergodicity}. Indeed, in this section we will prove Proposition \ref{prop:ergodicity}, giving ergodicity of the $SKS(\theta,N)$-process. The existence of a unique in law eternal $SKS(\theta,N)$-process will follow thanks to the Kolmogorov's extension theorem.

 \vip
 
  We will conclude the induction argument by showing that the candidate process built in Section \ref{sec:constructionKS} is unique in law in Section \ref{sec:uniqueness}, is strong Markov and Feller in Section \ref{sec:strongMarkov}. Finally we will prove Proposition \ref{prop:itoformula} in Section \ref{sec:ito} and that our notion of solution is associated to the Dirichlet space described in Theorem \ref{mainresult} in Section \ref{sec:dirichlet}.
  
  \vip
  
  Appendix \ref{app:estimations} is devoted to prove several technical but important intermediate results for the inductive argument. In Section \ref{sec:unicitefaible} we show a crucial result about uniqueness in law which is interesting in itself, which we will use to prove independence of two Keller-Segel processes driven by independent Brownian motions. Appendix \ref{app:otherresult} is devoted to prove other auxiliary or classical results used in the paper. 

\vip

Before detailing what is the induction hypothesis, let us first define the following motivated by Remark \ref{rk:dimensions}.
\begin{defi}\label{def:collisionalapprox}
For every $N\ge 1$, every $\theta>0$, we call a family of collisional approximation every family of function
$$
\{\chi^\e_K, \e >0,K\in \cB^N_\theta\},
$$
where $\cB^N_\theta$ has been defined in Definition \ref{def:KS/KSK}, such that for every $\e >0$, every $K\in \cB^N_\theta$, $\chi^\e_K$ is an application from $(\rr^2)^N\times \rr^{d_K}$ to $[0,1]$ for some $d_K\in \nn^*$ and for all $(x,z) \in (\rr^2)^N\times \rr^{d_K}$ we have
\begin{itemize}
\item the existence of $\beta_K>\alpha_K >0$ such that $$
\indiq_{R_K(x)\ge \alpha_K\e}\ge \chi^\e_K(x,z) \ge\indiq_{R_K(x)\ge \beta_K\e},
$$
\item for all $K\in \cB^N_\theta$, $\chi^\e_K(x,z)$ non decreases pointwise to $\indiq_{R_K(x)>0}(x)$ as $\e \to 0$,
\item  for all $K\in \cB^N_\theta$, $\chi^\e_K(x,z) = \chi^\e_K(\cS_K(x),z)$, where $\cS_K :   (\rr^2)^N\to (\rr^2)^N$ such that for all $x\in (\rr^2)^N$, all $i\notin K$, $(\cS_K(x))^i = x^i$ and for all $i\in K$, $(\cS_K(x))^i = 2S_K(x) -x^i $.
\end{itemize}
Recall that according to Remark \ref{rk:dimensions}, for all $K\in \cB^N_\theta$, $|K| = k_1 = \lceil 2N/\theta \rceil -1$.
\end{defi}
The last definition is simply a slight extension of the applications of the type $\indiq_{R_K(x)\ge \alpha}$, in order to have a more flexible notion of principal value, see Subsection \ref{subsec:bessel}.

\vip

We are now ready to give the induction hypothesis. We set for all $N\ge 1$,

\vip

$\cP(N)$:'' for all $\theta >0$, all $k\in \ig 1,N\id$,

\vip

(a) $\mathrm{MainTheorem}(\theta,k)$, $\mathrm{Decomposition} (\theta,k)$ and $\mathrm{Ito}(\theta,k)$ hold,

\vip

(b) if $k_1(\theta,k)>k$ (see Remark \ref{rk:dimensions} for the definition of $k_1(\theta,k)$), then 
$$
\mbox{ for all } K>0, \quad \sup_{x\in B(0,K)} \E^{X,k,\theta}_x \Big[ \int_0^1 \|b(X_s)\|\dd s \Big] < \infty,
$$
and if $k_1(\theta,k) = k$, then,
$$
\mbox{ for all } K>0, \quad \sup_{x\in B(0,K)} \E^{X,k,\theta}_x \Big[ \int_0^{\tau_k\wedge 1} \|b(X_s)\|\dd s \Big] < \infty,
$$
where $\E^{X,k,\theta}_x$ is the expectation relative to $\PP^{X,k,\theta}_x$ and $(X_s)_{s\ge 0}$ is the canonical process,

\vip

(c) for all familly $((Z^K_t)_{t\ge 0})_{K\in \cB^k_\theta}$ of measurable processes independent from $(X^k_s)_{s\ge 0}$ such that for all $K\in \cB^k_\theta$, $(Z^K_t)_{t\ge 0}$ is $\rr^{d_K}$-valued (with $d_K$ as in Definition \ref{def:collisionalapprox}) and every family of collisional approximation $\{\chi^\e_K, \e >0,K\in \cB^k_\theta\}$, we have for all $i\in \ig 1,k\id$, all $T>0$, 
$$
\sup_{t\in [0,T\wedge \tau_{k_0})} \Big\| \int_0^t b^i(X^k_s)\prod_{K\in \cB^k_\theta : i\in K}\chi^\e_K(X^{k}_s,Z^K_s)\dd s - \int_0^t b^i_\e(X^k_s)\dd s\Big\| \to 0,
$$
in probability as $\e \to 0$,

\vip

(d) for all $x\in E_{k_0}$, all $\eta>0$, $r>0$, $\tau^\eta_{k_0}$ and $\tau_{k_0}$ are continuous for the one-point compactification topology $\PP^{X,k,\theta}_x$-a.e. and for all $s>0$,
\begin{gather*}
\PP^{X,k,\theta}_x(\tau^\eta_{k_0}((X^k_t)_{t\ge 0}) = s) = \PP^{X,k,\theta}_x(\tau_{k_0}((X^k_t)_{t\ge 0}) = s) = 0,\\
\mbox{ if } k_0<N, \quad\PP^{X,k,\theta}_x(\lim_{t\to \tau_{k_0}^-} R_{\ig 1,k\id}(X^k_t) = r) = 0,
\end{gather*}
where we recall $\mathrm{MainTheorem}(\theta,k)$ for the definition of $\PP^{X,k,\theta}_\mu$ (see Theorem \ref{mainresult}) and $(X^k_t)_{t\ge 0}$ is the canonical process.''

\vip

First if $\cP(N)$ is true for all $N\ge 1$, then Theorem \ref{mainresult}, Proposition \ref{prop:itoformula}, Theorem \ref{thm:decompo} and Proposition \ref{prop:girsanovpropre} are true. We see that $\cP(1)$ is clearly true if we set by convention that a $SKS(\theta,1)$-process is the process constantly equal to $0$, the assertions (b), (c) and (d) being empty. Let $N\ge 2$ and suppose that $\cP(N-1)$ is true. We now show that $\cP(N)$ holds, which is now the goal of this paper.

\section{Configurations}\label{sec:config} 
In this section we explain how until the first time where all the particles meet, it is always possible in a sense to separate particles in many subgroups as suggested in Lemma \ref{lem:tpsarretsimple}. We begin with a helpful identity about the empirical variances.
\begin{lemma}\label{lem:dispersionformula}
For all $m\in \ig 1,N\id$, all $\cK = (K_i)_{i\in \ig 1,m \id} \in\gK_N$, all $x\in (\rr^2)^N$, we have
$$
R_{\ig 1,N\id}(x) = \|x\|^2 - N \|S_{\ig 1,N\id}(x)\|^2 \quad \mbox{ and }\quad R_{\ig 1,N\id}(x) = \sum_{i=1}^m R_{K_i}(x) + R_\cK (x),
$$
where $R_\cK(x) = (2N)^{-1}\sum_{K\ne L\in \cK}|K||L| \|S_K(x)-S_L(x)\|^2$.
\end{lemma}

\begin{proof}
Fix $\cK = (K_i)_{i\in \ig 1,m\id}\in \gK_N$ and $x\in (\rr^2)^N$. Consider the probability measure
$$
\mu = \frac{1}{N} \sum_{i=1}^N\delta_{(x_i,\hat K(i))},
$$
where $\hat K (i)$ is the only element of $\cK$ containing $i$. Let $(X,Y)$ be a $\ig 1,N\id \times \cK$-valued random variable with law $\mu$. The first equality comes from the fact that $V_\mu[X] = \E_\mu[\|X\|^2] - \|\E_\mu[X]\|^2$. The second equality follows from the fact that $V_\mu(X) = \E_\mu [ V_\mu [X|Y]] + V_\mu [\E_\mu [X|Y]]$ with $\E_\mu$ and $V_\mu$ relative to $\mu$. Indeed, first notice that $V_\mu [X] = N^{-1}R_{\ig 1,N\id}(x)$. Then we conclude since we have
\begin{align*}
\E_\mu [X|Y] = S_Y(x) \quad \mbox{ and } \quad V_\mu [X|Y] = |\hat K(Y)|^{-1} R_Y(x),
\end{align*}
so that
\begin{align*}
V_\mu[\E_\mu[X|Y]] = N^{-1} R_\cK(x)\quad \mbox{ and } \quad \E_\mu [ V_\mu [X|Y]] = N^{-1} \sum_{i=1}^m R_{K_i}(x).
\end{align*}
\end{proof} 

We now show a slight variation of \cite[Lemma 13]{ft}.
\begin{lemma}\label{campingcar-unleashedversion}
There exists constants $\eta_0>0$, $c>0$ only depending on $N$ such that for all $\eta \in (0,\eta_0)$, all $\cK\in \gK_N$, all $x\in (\rr^2)^N$, if
\begin{gather*}
\mbox{for all } K\in \cK, \quad R_K(x) \le 2\eta^{1/|K|}\\  \mbox{ for all subfamily } (K_i)_{i\in \ig 1,p\id} \mbox{ of } \cK \mbox{ with } p\ge 2, \quad R_{K_1\cup\dots\cup K_p} \ge \frac12 \eta^{1/(\max_{1\le i\le p}|K_i|+1)},
\end{gather*}
then for all $K\ne L \in \cK$, all $i\in K$, $j\in L$, $\|x^i-x^j\|^2 \ge c\eta^{1/(|K|\vee |L|+1)}$.
\end{lemma}

\begin{proof}
Let $\eta_0>0$ to be chosen later. Fix $x\in (\rr^2)^N$ and $\cK\in \gK_N$ satisfying the assumptions. Take $K\ne L \in \cK$ and $i\in K$ and $j\in L$. Then, according to the assumptions and to Lemma \ref{lem:dispersionformula}, we get
\begin{align*}
\frac12\eta^{1/(|K|\vee |L|+1)} \le R_{K\cup L}(x) &= R_K(x) + R_L(x) + \frac{|K||L|}{|K|+|L|} \| S_K(x) - S_L(x)\|^2\\
&\le 2\eta^{1/|K|} + 2\eta^{1/|L|} + N \| S_K(x) - S_L(x)\|^2
\end{align*}
So that 
$$
\| S_K(x) - S_L(x)\|^2 \ge \frac{1}{N} (2^{-1}\eta^{1/(|K|\vee |L|+1)} - 2\eta^{1/|K|} - 2\eta^{1/|L|}).
$$
We deduce that for all $i\in K$, $j\in L$,
\begin{align*}
\|x^i-x^j\| &\ge \|S_K(x) -cS_L(x)\| - R_K(x) - R_L(x)\\
&\ge \frac{1}{N} (2^{-1}\eta^{1/(|K|\vee |L|+1)} - 2\eta^{1/|K|} - 2\eta^{1/|L|}) - 2\eta^{1/|K|} - 2\eta^{1/|L|}
\\
&\ge c\eta^{1/(|K|\vee |L|+1)},
\end{align*}
for all $\eta \in (0,\eta_0)$ with $\eta_0$ small enough only depending on $N$ and $c=1/(4N)$. 
\end{proof}
\bla
\begin{lemma}\label{lem:releq}
Fix $N\ge 1$. There exists $\eta_1>0$ only depending on $N$ such that if $0<\eta \le \eta_1$, then for all $x\in (\rr^2)^N$, all $K_1$, $K_2\in \cP^\eta_x$ such that $|K_1|\ge |K_2|$ where
$$
\cP^\eta_x = \Big\{ K\subset \ig 1,N\id : R_K(x) \le \eta^{1/|K|}\quad \mbox{ and } \quad \min_{i\notin K} R_{K\cup \{ i\}}(x) > \eta^{1/(|K|+1)} \Big\},
$$
we have either $K_1\cap K_2 = \emptyset$ or $K_2\subset K_1$.
\end{lemma}
\begin{proof}
Let $\eta>0$ to be chosen later, let $x\in (\rr^2)^N$ and $K_1$, $K_2 \in \cP^\eta_x$ such that $|K_1|\ge |K_2|.$ According to Lemma \ref{campingcar-unleashedversion} applied to $\cK_i = \{ K_1\} \cup \{ i\}$ for all $i\notin K_1$, we get that if $0<\eta\le \eta_1$, then for all $i\in K_1$, all $j\notin K_1$, $\|x^i-x^j\|^2 \ge c \eta^{1/(|K_1|+1)}  $. If $i\in K_1\cap K_2$, then taking any $j\in K_2$, we get 
$$
\|x^i-x^j\|^2 \le |K_2| R_{K_2}(x)\le N \eta^{1/|K_2|} \le N \eta^{1/|K_1|} ,
$$ 
according to the definition of $R_{K_2}$ and the fact that $K_2 \in \cP^\eta_x$. Then if $j\notin K_1$, this implies that
$$
c\eta^{1/(|K_1|+1)} \le N \eta^{1/|K_1|},
$$
i.e that
$$
cN^{-1} \le \eta^{1/[|K_1|(|K_1|+1)]} \le \eta^{1/[N(N+1)]},
$$
which does not hold if $\eta$ less than some $\tilde \eta_1$ only depending on $N$ and $c$ and small enough. We conclude choosing $\eta_1 = \eta_0 \wedge \tilde \eta_1$.
\end{proof}
For a fixed $N\ge 1$, we define for all $x\in (\rr^2)^N$, 
\begin{align}\label{eq:configeta}
\cK^\eta(x) := \{ K \in \cP^\eta_x \mbox{ maximal for the inclusion}\},
\end{align}
which is meaningful thanks to Lemma \ref{lem:releq}.
\begin{lemma}
For all $0<\eta<\eta_1$, all $x\in (\rr^2)^N$, $\cK^\eta(x)\in \gK_N$.
\end{lemma}
\begin{proof}
Let $0<\eta<\eta_1$. We first prove that $\cup_{K\in \cK^\eta(x)}K = \ig 1,N\id$. By symmetry, it is sufficient to show that there exists $K\in \cP^\eta_x$ that contains $1$, which we now do. If $\{1\} \in \cP^\eta_x$ then we are done. If not, there exists $i_1\in \ig 1,N\id$ such that $R_{\{1,i_1\}}(x) \le \eta^{1/2}$. If $\{1,i_1\} \in \cP^\eta_x$ it is done, and in the contrary there exists $i_2 \in \ig 1,N\id \setminus \{1,i_1\}$ such that $R_{\{1,i_1,i_2\}}(x) \le \eta^{1/3}$. If this procedure stop at a point we are done. If not then we end up with $R_{\ig 1,N\id}(x) \le \eta^{1/N}$, implying that $\ig 1,N\id \in \cP^\eta_x$.  

\vip

We now show that $\cK^\eta(x)$ is constituted of disjoint sets. Let $K_1 \ne K_2 \in \cK^\eta (x)\subset \cP^\eta_x$ with, for example, $|K_1|\ge |K_2|$. If $K_1\cap K_2\ne \emptyset$, then according to Lemma \ref{lem:releq}, $K_2\subset K_1$, with $K_1\ne K_2$, which contradicts the maximality of $K_2$. This implies that $K_1\cap K_2 = \emptyset$ and this ends the proof.
\end{proof}

For any continuous $E_{k_0}\cup\{\triangle\}$-valued process $(Z_t)_{t\ge 0}$, we define the families of stopping time $(\sigma^\eta_k(Z))_{k\ge 0}$ and of set $(\cK^\eta_k)_{k\ge 0}$ as follow : $\sigma^\eta_0(Z) = 0$, $\cK^\eta_0 = \cK^\eta(x)$ and for all $k\ge 0$,
\begin{align}\label{config}
\sigma^\eta_{k+1}(Z) = \inf \Big\{ &t\ge \sigma^\eta_k : Z_t = \triangle \quad \mbox{ or } \quad R_K(Z_t)\ge 2\eta^{1/|K|} \mbox{ for some } K\in \cK^\eta_k(Z) \\
&\mbox{ or } R_{K_1\cup\dots\cup K_p}(Z_t)\le \frac12 \eta^{1/(|K_1|+\dots+|K_p|)} \mbox{ for some } (K_i)_{i\in \ig 1,p\id} \subset \cK^\eta_k(Z)\Big\},\notag\\
\cK^\eta_k(Z) = \cK^\eta(&Z_{\sigma^\eta_k}) \mbox{ if } Z_{\sigma^\eta_k} \ne \triangle \mbox{ and } \cK^\eta_k(Z) = \triangle \mbox{ otherwise},\label{etaconfig}
\end{align}
with the convention $\inf \emptyset = \infty$. We say that for all $k\ge 0$, $\cK^\eta_k(Z)$ is the $k$-th $\eta$-configuration and $\sigma^\eta_k(Z)$ is the $k$-th time of change of $\eta$-configuration. Moreover we set
$$
\sigma^\eta_*  = \inf \Big\{ \sigma^\eta_k : R_{\ig 1,N\id}(Z_{\sigma^\eta_k}) \le \eta^{1/N} \quad \mbox{ or } \quad \|Z_{\sigma^\eta_k}\|\ge \eta^{-1/N} \Big\}\wedge \inf\{ t\ge 0 : Z_t = \triangle\},
$$
\begin{lemma}\label{lem:particulesloins}
For all $0<\eta<\eta_1$, all continuous $E_{k_0}\cup\{\triangle\}$-valued process $(Z_t)_{t\ge 0}$, for all $k\ge 0$ such that $\cK^\eta_k\ne \triangle$, all $K\in \cK^\eta_k$, all $t\in [\sigma^\eta_k, \sigma^\eta_{k+1}]$, all $K\in \cK^\eta_k$, we have 
$$
R_K(Z_t) \le 2\eta^{1/|K|} \quad \mbox{ and for all } i\in K, \mbox{ all } j\in L \mbox{ with } L\in \cK^\eta_k(Z)\setminus \{ K\}, \quad \|Z^i_t-Z^j_t\|^2 \ge c\eta^{1/(|K|\vee |L|+1)}.
$$
\end{lemma}
\begin{proof}
This is direct by applying Lemma \ref{campingcar-unleashedversion}.
\end{proof}
All in all, Lemma \ref{lem:tpsarretsimple} has been proven and precised.
\begin{rk}
Since $\eta \mapsto \tau^\eta_{k_0\wedge N}$ is a non increasing application, we can define
$$
\tau_{k_0\wedge N} := \lim_{\eta\to 0} \tau^\eta_{k_0\wedge N},
$$
which takes its value in $\rr_+\cup\{\infty\}$. Recalling Definition \ref{def:sigmastar}, it is worthy noticing that for all $\eta \in (0,\eta_1)$,
$$
 \tau^\eta_{k_0\wedge N}\le \sigma^\eta_* \le \tau_{k_0\wedge N}.
$$ 
This implies in particular that $\sigma^\eta_* \to \tau_{k_0\wedge N}$ as $\eta \to 0$. 
\end{rk}
\section{Construction of the $KS$-Process killed upon reaching $E_N$}\label{sec:killedEN}
In this section we prove the existence and uniqueness in law of $KS(\theta,N)$-process killed upon reaching $E_{k_0\wedge N}$. We show in addition that this is Feller and strong Markov which is associated with the desired Dirichlet form and that the empirical mean and the empirical dispersion (up to a multiplicative constant) are respectively, until the exiting time of $E_{k_0\wedge N}$, a Brownian motion and a squared Bessel process as stated in Theorem \ref{thm:decompo}-(i).

\vip

Let $\Lambda '$ be the set of increasing functions $\lambda$ mapping $\rr_+$ into $\rr_+$ (in particular, $\lambda (0) = 0$, $\lim_{t\to \infty} \lambda (t) = \infty$ and $\lambda$ is continuous). Let $\Lambda$ denote the set of Lipschitz-continuous functions $\lambda \in \Lambda '$ such that
$$
\gamma (\lambda) = \sup_{t>s \ge 0} \Big| \log\Big( \frac{\lambda (t) - \lambda (s)}{t-s}\Big) \Big|< \infty. 
$$
If $O$ is an open set of $(\rr^2)^N\times (\rr^2)^N$, we equip $\cD(\rr_+,O_\triangle)$ the space of c\`adl\`ag functions from $\rr_+$ to $O_\triangle$ with the classical metric $d$ defined for all $f$, $g \in \cD(\rr_+,(\rr^2)^N\times (\rr^2)^N)$ by, 
\begin{align}\label{distancebp}
d(f,g) = \inf_{\lambda \in \Lambda} \Big[ \gamma (\lambda) \vee \int_0^\infty e^{-u}\sup_{t\ge 0} d_{O_\triangle}(f(t\wedge u),g(\lambda(t) \wedge u)) \dd u \Big],
\end{align}
where $d_{O_\triangle}$ denotes the bounded metric inducing the topology of $O_\triangle$, the one-point compactification of $O$, see \cite[Proposition 7.1.15., page 187]{Cohn} (if $d_{O_\triangle}$ is not bounded, one can always consider $d_{O_\triangle}/(1+d_{O_\triangle})$ instead). 

\vip

The following lemma is a technical preliminary to Proposition \ref{prop:girsanoveta} but the proof, which is postponed to Appendix \ref{app:preuvepropetapreKS}, is not essential to understand the main arguments. 

\vip

\begin{lemma}\label{lem:etapreKS}
Let $\theta >0$, $\eta \in (0,\eta_1)$ and $\mu\in \cP((E_{k_0}\cap E^\eta_{N})_\triangle)$. Then the following holds.

\vip

(i) There exists a unique in law $\eta$-$pre$-$KS_\mu(\theta,N)$-process. We denote its law by $\tilde\PP^{\eta,N,\theta}_\mu$.

\vip

(ii) For all $k\ge 1$, the mappings $\sigma^\eta_k\indiq_{\sigma^\eta_{k-1} < \sigma^\eta_*}$,  $\indiq_{\sigma^\eta_{k-1}<\sigma^{\eta}_*}$ and $\indiq_{\sigma^\eta_{k-1}<\sigma^\eta_*}\cK^\eta_k$ are continuous for the one-point compactification topology $\tilde\PP^{\eta,N,\theta}_\mu$-a.e. Moreover for all $t>0$, 
$$
\tilde \PP^{\eta,N,\theta}_\mu(\sigma^\eta_k = t) = 0.
$$

\vip

(iii) The application $\sigma^\eta_*$ is continuous for the one-point compactification topology $\tilde\PP^{\eta,N,\theta}_\mu$-a.e. Moreover for all $t>0$, 
$$
\tilde \PP^{\eta,N,\theta}_\mu(\sigma^\eta_* = t) = 0.
$$
In addition, if $k_0\le N$, then $R_{\ig 1,N\id}(X_t)$ a.s. converges as $t\to \tau_{k_0}^-$ for the usual topology on $\rr$, and for all $r>0$,
$$
\tilde \PP_\mu^{\eta,N,\theta}\Big[ \lim_{t\to\tau_{k_0}^-} R_{\ig 1,N\id}(X_t) = r\Big] = 0.
$$

\vip

(iv) Let $x_\infty\in (E_{k_0}\cap E^\eta_N)_\triangle$ and let $(x_n)_{n\ge 0}$ a sequence of elements of $(E_{k_0}\cap E^\eta_N)_\triangle$ converging to $x_\infty$. Then $\tilde \PP^{\eta,N,\theta}_{x_n}\to \tilde \PP^{\eta,N,\theta}_{x_\infty}$ weakly as $n\to \infty$, with respect to the metric space $(C(\rr_+,(E_{k_0}\cap E^\eta_N)_\triangle,d)$.

\vip

(v) For all $\varphi \in C_0(E_{k_0}\cap E^\eta_N)$, we have
$$
\tilde P^{N,\eta}_t \varphi := \tilde \E^{\eta,N,\theta}[\varphi (X_t)\indiq_{t<\sigma^\eta_*}] \in C_0(E_{k_0}\cap E^\eta_N),
$$
where $(X_t)_{t\ge 0}$ denotes the canonical process and we recall that $C_0(E_{k_0}\cap E^\eta_N)$ is the space of continuous real-valued functions on $E_{k_0}\cap E^\eta_N$ that vanishes at the boundary of $E_{k_0}\cap E^\eta_N$. Furthermore, 
$$
\sup_{x\in E_{k_0}\cap E^\eta_N} |\tilde P^{N,\eta}_t\varphi (x) - \varphi (x)| \to 0 \quad \mbox{ as } t\to 0.
$$

\end{lemma}
We now prove a crucial lemma for the induction argument: it merely says that if several subgroups of particles are far away from each other, they behave like independent Keller-Segel processes. In a more precise way, we will prove that a $\eta$-$pre$-$KS(\theta,N)$-process and a $KS(\theta,N)$-process stopped at $\sigma^\eta_*$ have the same law up to a sufficiently regular change of probability. 

\begin{prop}\label{prop:girsanoveta}
 Let $N\ge 1$, $\theta>0$, $\mu\in \cP((\rr^2)^N)$. Then
 
 \vip
 
(i) there exists a unique in law $KS_\mu(\theta,N)$-process killed upon exiting $E_N$, we write its law $\PP^{X,N,\theta,*}_\mu$,

\vip

(ii) let $(\Omega,\cF,(\cF_t)_{t\ge 0},\PP,(B_t)_{t\ge 0},(X^N_t)_{t\ge 0})$ be a $KS_\mu(\theta,N)$-process killed upon exiting $E_N$, driven by a $2N$-dimensional Brownian motion $(B_t)_{t\ge 0}$. Then 

\vip

(a) for all $T>0$, all $\eta >0$, all family $((Z^K_t)_{t\ge 0})_{K\in \cB^N_\theta}$ of processes independent from $(X^N_t)_{t\ge 0}$ such that for all $K\in \cB^N_\theta$, $(Z^K_t)_{t\ge 0}$ is $\rr^{d_K}$-valued for some $d_K\in \nn^*$ and all family of collisional approximation $
\{\chi^\e_K, \e >0, K\in \cB^N_\theta\}$, we have for all $i\in \ig 1,N\id$,
\begin{align*}
\sup_{t\in [0,T\wedge \sigma^\eta_*\wedge \tau^\eta_{k_0}]} \Big\|X^{N,i}_t&-x^i - B^i_t-\int_0^t b^i(X^N_s)\prod_{K\in \cB^N_\theta : i\in K}\chi^\e_K(X^N_s, Z^K_s)\dd s \Big\| \to 0,
\end{align*}
in probability as $\e \to 0$,

\vip

(b)for all $\varphi \in C^\infty_c(E_N)$ such that
$$
\sup_{\alpha \in (0,1]}\sup_{x\in (\rr^2)^N} |\cL_\alpha^X \varphi (x) | < \infty,
$$
we have for all $t\in [0,\tau_{k_0\wedge N})$,
\begin{align*}
\varphi(X^N_{t})&=\varphi (X^N_0) + \int_0^t \nabla \varphi (X^N_s) \cdot\dd B_s + \int_0^t \cL^X \varphi (X^N_s)\dd s,
\end{align*}

\vip

(c) for all $T>0$, all $\eta>0$, there exists a probability $\dd\QQ^{T,\eta}_\mu$ such that under $\QQ^{T,\eta}_\mu$, $(X^N_{t\wedge \sigma^\eta_* })_{t\ge 0}$ is a $\eta$-$pre$-$KS_\mu(\theta,N)$-process driven by $(\tilde B_t)_{t\ge 0}$ where for all $t\in [0,\sigma^\eta_*]$,
\begin{align*}
\tilde B_t = B_t -\sum_{k\ge0} \int_{\sigma^\eta_k}^{\sigma^\eta_{k+1}} \nabla G_{\cK^\eta_k}(X^N_s)\dd s,
\end{align*}
with for all $x\in E_2$, all $\cK$ partition of $\ig 1,N\id$, 
\begin{align}\label{eq:GK}
\quad G _\cK = \frac{\theta}{2N} \sum_{K\ne L\in \cK} \sum_{i\in K,j\in L} \log(\|x^i-x^j\|).
\end{align}
Moreover, there exists $C_{N,\theta,\eta}\in (0,\infty)$ only depending on $N$, $\theta$ ,$\eta$  and a $\PP^{X,N,\theta,*}_\mu$-a.e. continuous application $F_{T,\eta} : C([0,T],\rr)\to \rr$ for the uniform convergence topology such that $\dd \QQ^{T,\eta}_\mu = F_{T,\eta} \dd \PP^{X,N,\theta,*}_\mu$ and
\begin{align}\label{eq:boundonF}
    \exp (-C_{N,\theta,\eta}T) \le F_{T,\eta} \le \exp( C_{N,\theta,\eta}T).
\end{align}

\end{prop}
\begin{proof}
Fix $\mu \in \cP((\rr^2)^N)$. According to the Kolmogorov's extension theorem and the fact that $\sigma^\eta_* \to \tau_{k_0\wedge N}$ as $\eta \to 0$, it is sufficient to prove that for all $T>0$, all $\eta >0$, there exists a unique in law $KS_\mu(\theta,N)$-process killed at $\sigma^\eta_*\wedge \tau^\eta_{k_0}$ satisfying the above properties, which we now do. Fix $T>0$ and $\eta >0$.

\vip

{\it Step 1.} We show that there exists a $KS_\mu(\theta,N)$-process killed at $\sigma^\eta_*$. According to Lemma \ref{lem:etapreKS}, there exists $(\Omega,\cF,(\cF_t)_{t\ge 0},\PP)$ a filtered probability space on which are defined a $\eta$-$pre$-$KS_\mu(\theta,N)$-process $(\tilde X^N_t)_{t\ge 0}$ and its driving Brownian motion $(B_t)_{t\ge 0}$.

\vip

{\it Step 1.1} We define a probability $\tilde \QQ^{T,\eta}$ under which 
$$
(\tilde B_t)_{t\in [0,T]} := ((B^{i,\ell}_t - \langle B^{i,\ell}, \tilde L \rangle _t)_{i\in \ig 1,N \id, \ell \in \{ 1,2 \}})_{t\in [0,T]}
$$
is a $2N$-dimensional Brownian motion where for all $t\ge 0$, 
\begin{align*}
\tilde L_t &= -\frac{\theta}{N}\sum_{k\ge 0} \int_{\sigma^\eta_k}^{\sigma^\eta_{k+1}} \indiq_{s\le t\wedge\sigma^\eta_*}\sum_{K,L\in \cK^\eta_k, K\ne L} \sum_{i\in K, j\in L} \frac{\tilde X^{N,i}_s-\tilde X^{N,j}_s}{\|\tilde X^{N,i}_s-\tilde X^{N,j}_s\|^2}\cdot(\dd B^i_s-\dd B^j_s)\\
&= -\sum_{k\ge 0} \int_{\sigma^\eta_k}^{\sigma^\eta_{k+1}} \indiq_{s \le t\wedge \sigma^\eta_*}\nabla G_{\cK^\eta_k}(\tilde X^N_s) \cdot \dd B_s.
\end{align*}
For all $k\ge 0$, all $t\in [\sigma^\eta_k, \sigma^\eta_{k+1}]\cap [0,\sigma^\eta_*]$, all $i\in K$, $j\in L$ with $K,L\in \cK^\eta_k$ we have $\|\tilde X^{N,i}_t - \tilde X^{N,j}_t\| \ge c \eta^{1/(|K|\vee |L|+1)}$ by definition of $(\sigma^\eta_k)_{k\ge 0}$ and Lemma \ref{campingcar-unleashedversion}. We deduce that $\langle L \rangle _T$ is bounded by $C_{N,\theta,\eta}T$ with $C_{N,\theta,\eta}$ some deterministic constant depending on $N$, $\theta$ and $\eta$. Since the Novikov condition is fulfilled, we can apply the Girsanov's theorem and obtain the result with $\tilde\QQ^{T,\eta} := \cE(L)_T\dd \PP$. Moreover, since $(\cE(L)_T)^2 = \cE(2L)_T \exp (\langle L\rangle_T)$, we get again with the Novikov criterion and the fact that $\langle L\rangle_T$ is bounded that there exist $\tilde c_{T,\eta}$ and $\tilde C_{T,\eta}$ positive finite constants such that
\begin{align*}
\exp( - C_{N,\theta,\eta}T) \le \E \Big[ \Big( \frac{\dd \tilde\QQ^{T,\eta}}{\dd \PP}\Big)^2\Big] \le \exp( C_{N,\theta,\eta}T).
\end{align*}
Observe that since the bound on $\langle L\rangle_T$ does not depend on $\mu$, it is the same for $C_{N,\theta,\eta}$.

\vip

{\it Step 1.2.} Let $\eta >0$ and $T>0$. We show that $(\tilde X^N_t)_{t\in [0,\sigma^\eta_*\wedge T]}$ under $\tilde \QQ^{T,\eta}$ is a $KS_\mu(\theta,N)$-process stopped at $\sigma^\eta_*\wedge T$ for which (ii)-(a) is satisfied. Fix a family of collisional approximation $\{ \chi^\e_K, \e >0, K\in \cB^N_\theta\}$ and family $((Z^K_t)_{t\ge 0})_{K\in \cB^N_\theta}$ of $\rr^{d_K}$-valued processes with $d_K\in \nn^*$ independent from all the other variables involved here. We set for all $t\in [0,\sigma^\eta_*]$, $R^{i,\e}_t = R^{1,i,\e}_t + R^{2,i,\e}_t$ where
\begin{gather}\label{eq:approsurcoll}
R^{1,i,\e}_t \!\!= \!\tilde X^{N,i}_{t\wedge\tau^\eta_{k_0}} \!\!-\tilde X^{N,i}_0\!\! - \!B^i_{t\wedge\tau^\eta_{k_0}}\!\! -\!\!\frac\theta N\sum_{k\ge 0}\int_{\sigma^\eta_k}^{\sigma^\eta_{k+1}} \!\!\indiq_{s\le t \wedge\tau^\eta_{k_0}} \!\!\sum_{j\in K^i_k} \frac{\tilde X^{N,i}_s - \tilde X^{N,j}_s}{\|\tilde X^{N,i}_s - \tilde X^{N,j}_s\|^2}\!\!\!\!\prod_{K\in \cB^N_\theta : K\subset K^i_k}\!\!\!\!\chi^\e_K (\tilde X^N_s,Z^K_s)\dd s,\\
R^{2,i,\e}_t = -\frac\theta N\sum_{k\ge 0}\int_{\sigma^\eta_k}^{\sigma^\eta_{k+1}} \indiq_{s\le t \wedge\tau^\eta_{k_0}} \sum_{j\notin K^i_k} \frac{\tilde X^{N,i}_s - \tilde X^{N,j}_s}{\|\tilde X^{N,i}_s - \tilde X^{N,j}_s\|^2}\Big[\prod_{K\in \cB^N_\theta : K\subset K^i_k}\chi^\e_K (\tilde X^N_s,Z^K_s)-1\Big]\dd s,
\end{gather}
where $K^i_k$ is the element of the partition $\cK^\eta_k$ containing $i$. Thanks to the definition of $(\sigma^\eta_k)_{k\ge 0}$, we get that for $\e$ less than a positive deterministic constant, for all $k\ge 0$, all $s\in [\sigma^\eta_k,\sigma^\eta_{k+1}]$, $\chi^\e_K(\tilde X^N_s, Z^K_s) = 1$ for every $K\in \cB^N_\theta$ not contained in an element of $\cK^\eta_k$. Recalling the definition of $(\tilde B_t)_{t\in [0,T]}$ in Step 1.1, this implies that for all $t\in [0,\sigma^\eta_*]$,
\begin{align*}
    R^{i,\e}_t = \tilde X^{N,i}_{t\wedge\tau^\eta_{k_0}} -\tilde X^{N,i}_0 - \tilde B^i_{t\wedge\tau^\eta_{k_0}} -\int_0^{t\wedge \tau^\eta_{k_0}} b^i(\tilde X^N_s) \prod_{K\in \cB^N_\theta : i\in K} \chi^\e_K(\tilde X^N_s, Z^K_s) \dd s.
\end{align*}

Fix $\nu >0$ and $\alpha >0$. According to Step 1.1 we get by the Cauchy-Schwarz inequality the existence of a positive and finite constant $\tilde C_{T,\eta}$ such that, 
\begin{align}\label{ttt}
\tilde \QQ^{T,\eta} [ \sup_{t\in [0,T]}\| R^{i,\e}_t\|>\nu] \le (\tilde C_{T,\eta})^{1/2} \Big( \PP [ \sup_{t\in [0,T]}\| R^{i,\e}_t\|>\nu]\Big)^{1/2}.
\end{align}
It is clear that $\PP [ \sup_{t\in [0,T]}\| R^{2,i,\e}_t\|>\nu] \to 0$ as $\e \to 0$. Indeed, it exists a constant $C\in (0,\infty)$ and $\beta \in (0,\infty)$ according to Definition \ref{def:collisionalapprox} such that
$$ 
|R^{2,i,\e}_t|\le C\sum_{k\ge 0}\int_{\sigma^\eta_k}^{\sigma^\eta_{k+1}} \indiq_{s\le t \wedge \sigma^\eta_*\wedge\tau^\eta_{k_0}}\!\! \prod_{K\in \cB^N_\theta : K\subset K^i_k}\!\!\indiq_{R_K(\tilde X^N_s) \le \beta\e}\dd s \le C\sum_{K\in \cB^N_\theta} \int_0^{t\wedge \sigma^\eta_*\wedge\tau^\eta_{k_0}} \indiq_{R_K(\tilde X^N_s)\le \beta \e } \dd s.
$$
Thus, by the monotone convergence theorem and the third point in Definition \ref{def:KS/KSK}, we have $R^{2,i,\e}_t\to 0$ as $\e \to 0$ a.s. Consequently, it only remains to show that $\PP [ \sup_{t\in [0,T]}\| R^{1,i,\e}_t\|>\nu] \to 0$ as $\e \to 0$ in order to conclude this step, which we now do.

\vip

We have for all $p\ge 0$,
\begin{align*}
\PP [ \sup_{t\in [0,T]}\| R^{1,i,\e}_t\|>\nu] \le \sum_{k=0}^p \PP [ \sup_{t\in [\sigma_k^\eta \wedge T,\sigma_{k+1}^\eta \wedge T]}\| R^{1,i,\e}_t\|>\nu] +  \PP [ \sigma^\eta_p < T],
\end{align*}
By continuity of a $KS(\tilde \theta,\tilde N)$-process for all $\tilde \theta>0$, all $\tilde N \le N$ on $[0,\tau^\eta_{k_0}]$, we get that there exists $p\ge 0$ such that $\tilde\PP^{\eta,N,\theta}_\mu[ \sigma^\eta_p\le T] \le \alpha$.
Observe that for all $K\subset \ig 1,N\id$, thanks to the expression of $k_1$ (see Remark \ref{rk:dimensions}) $k_1(\theta,N) = k_1( \theta |K|/N,|K|)$, so that the set
$$
\{ \chi^\e_L, \e >0, L\subset K : |L| = k_1(|K|\theta/N,|K|)\}
$$ 
is a family of collisional approximation (seen as application from $(\rr^2)^{|K|}\times (\rr^2)^{N-|K| + d_{|K|}}$ to $ \rr_+$). According $\cP(N-1)$-(c), we get that for all $k\in \ig 0,p\id$, if $\e$ is small enough then $\PP [ \sup_{t\in [\sigma^\eta_k\wedge T,\sigma^\eta_{k+1}\wedge T]}\| R^{1,i,\e}_t\|>\nu]\le \alpha/p$. This implies easily that
$$
\PP \Big[ \sup_{t\in [0,T]}\| R^{1,i,\e}_t\|>\nu \Big]\le 2\alpha.
$$

\vip

Thus, together with \eqref{ttt} we get that
\begin{align}\label{qqq}
\mbox{under }\tilde \QQ_{T,\eta}, \quad \sup_{t\in [0,T]}\| R^{i,\e}_t\| \to 0 \quad \mbox{ in probability as } \e \to 0.
\end{align} 

\vip

{\it Step 2.} We now prove (ii)-(b). Since uniqueness will be shown in Step 3, it is sufficient to prove it for the particular process we considered in Steps 1 and 2. Let $T>0$, $\eta >0$ and $\varphi \in C^\infty_c(E_N)$ such that
$$
\sup_{\alpha \in (0,1]}\sup_{x\in (\rr^2)^N} |\cL_\alpha ^X\varphi (x) | < \infty. 
$$
First, this implies that for every partition $\cK$ of $\ig 1,N\id$, we have for all $\delta >0$,
$$
\sup_{\alpha \in (0,1]} \sup_{x\in E^\delta_\cK} \Big|\frac12 \Delta \varphi (x) -\frac{\theta}{N} \sum_{K\in \cK} \sum_{i\in K} \nabla_{x_i} \varphi (x) \cdot \sum_{j\in K} \frac{x_i - x_j}{\|x_i-x_j\|^2 + \alpha} \Big| < \infty,
$$
since $\sup_{\alpha \in (0,1]}\sup_{x\in E^\delta_\cK} |\cL_\alpha ^X \varphi (x) | \le \sup_{\alpha \in (0,1]}\sup_{x\in (\rr^2)^N} |\cL_\alpha ^X\varphi (x) |$
 with
 $$
 E^\delta_\cK = \Big\{ x\in (\rr^2)^N : \min_{K\ne L\in \cK} \min_{i\in K,j\in L} \|x_i-x_j\|^2\ge \delta\Big\} \cap B(0,\delta^{-1}).
$$
It is thus possible for all $\cK\in \gK_N$ and $\delta >0$ to modify $\varphi$ outside of $E^\delta_\cK$ in a way that both $\varphi$ is smooth and $\cL^X\varphi$ is bounded on $(\rr^2)^N$.
As a consequence, one can use for all $k\ge 0$, $\mathrm{Ito}(\theta,N-1)$ from $\cP(N-1)$-(a) on the interval $[\sigma^\eta_k,\sigma^\eta_{k+1}]$, and the definition of $(\tilde X^N_t)_{t\in [0,\sigma^\eta_*\wedge T]}$ under $\PP$, we get by summing on $k$ that for all $t\in [0,\sigma^\eta_*\wedge \tau^\eta_{k_0}]$
\begin{align*}
\varphi (\tilde X^N_{t}) =& \varphi (\tilde X^N_0) + \int_0^{t} \nabla \varphi (\tilde X^N_s)\cdot \dd B_s + \frac12\int_0^{t} \Delta \varphi (\tilde X^N_s)\dd s \\
&- \frac{\theta}{N} \sum_{k\ge 0} \int_{\sigma^\eta_k}^{\sigma^\eta_{k+1}}  \indiq_{s\le t}\sum_{K\in \cK^\eta_k} \sum_{i\in K} \nabla_{x_i} \varphi ( \tilde X^N_s) \cdot \sum_{j\in K\setminus\{i\}}\frac{\tilde X^{N,i}_s - \tilde X^{N,j}_s}{\|\tilde X^{N,i}_s - \tilde X^{N,j}_s\|^2} \dd s.
\end{align*}
One could wonder about the possibility to apply simulatenously $\mathrm{Ito}(\theta,N-1)$ on $\varphi$ depending on all the clusters. Let $k\ge 0$ and write $\cK^\eta_k = (K_1,\dots, K_m)$. Take also $(t_i)_{i\in \ig 1,M\id}$ a discretization of $[\sigma^\eta_k,\sigma^\eta_{k+1}]$. We set for all $p\in \ig 0,m\id$, all $i\in \ig 1,M-1\id$,
$$
\varphi^{p,i}_k = \varphi ((\tilde X^{N,i}_{t_{i+1}})_{i\in K_1},\dots,(\tilde X^{N,i}_{t_{i+1}})_{i\in K_p},(\tilde X^{N,i}_{t_i})_{i\in K_{p+1}},\dots (\tilde X^{N,i}_{t_i})_{i\in K_m}).
$$
Observe that for all $i\in \ig 1,M-1\id$, $\varphi ( \tilde X^N_{t_{i+1}}) = \varphi^{m,i}_k$, $\varphi (\tilde X^N_{t_i}) = \varphi^{0,i}_k(t)$ and
$$
\varphi(\tilde X^N_{t_{i+1}}) - \varphi (\tilde X^N_{t_i}) = \sum_{p=0}^{m-1} \varphi_k^{p+1,i} - \varphi_k^{p,i}.
$$
It thus suffices to apply $\mathrm{Ito}(\theta,N-1)$ for each increment, use independence between the clusters and pass to the limit as the mesh of the discretization tends to $0$. This is classical thanks to the uniform bound on the drift. 

\vip

Finally, according to Step 1.1, we get that under $\tilde \QQ^{T,\eta}$, for all $t\in [0,\sigma^\eta_*\wedge\tau^\eta_{k_0}]$,
$$
\varphi (\tilde X^N_{t}) = \varphi(\tilde X^N_0) +\int_0^{t} \nabla \varphi (\tilde X^N_s)\cdot \dd \tilde B_s + \int_0^{t} \cL^X \varphi (\tilde X^N_s)\dd s.
$$
This implies the result since $\lim_{\eta\to 0} \sigma^\eta_*\wedge \tau^\eta_{k_0} = \tau_{k_0\wedge N}$ a.s.

\vip

{\it Step 3.} In this step, we show the uniqueness result and regularity of the application $F_{T,\eta}$ of the statement. Fix $(\Omega^1,\cF^1,(\cF^1_t)_{t\ge 0}, \PP^1,(B^1_t)_{t\ge 0},(X^1_t)_{t\ge 0})$ and $(\Omega^2,\cF^2,(\cF^2_t)_{t\ge 0}, \PP^2,(B^2_t)_{t\ge 0},(X^2_t)_{t\ge 0})$ two $KS_\mu(\theta,N)$-processes. We define for all $n\in \{1,2\}$, $\QQ^{n,T,\eta}_\mu = \cE(L^n_T)\dd \PP^n$ where for all $t\ge 0$,
$$
L^n_t = \sum_{k\ge 0} \int_{\sigma^\eta_k}^{\sigma^\eta_{k+1}} \indiq_{s\le t\wedge \sigma^\eta_*}\nabla G_{\cK^\eta_k}( X^n_s)\cdot\dd B^n_s.
$$
Checking the Novikov criterion as in Step 2.1, we get that for all $n\in \{1,2\}$, $(X^n_{t\wedge\sigma^\eta_*\wedge \tau^\eta_{k_0}})_{t \in [0, T]}$ under $\QQ^{n,T,\eta}_\mu$ is a $\eta$-$pre$-$KS_\mu(\theta,N)$ process. We now aim at using the analogue of Itô's formula on the exponential martingale of the change of probability in order to get rid of the stochastic integral and conclude in Step 3.2.

\vip

{\it Step 3.1}
 According to \cite[Lemma 12]{ft}, for all $\cK$ partition of $\ig 1,N\id$, all $\delta >0$, there exists $\Gamma^\delta_\cK \in C_c^\infty ((\rr^2)^N)$ such that 
$$
\mathrm{Supp}\Gamma^\delta_\cK \subset E^{\delta'}_\cK \quad \mbox{ for some } \delta'<\delta, \quad \Gamma^\delta_\cK =1 \mbox{ on } E^{\delta}_\cK \quad \mbox{ and } \quad \sup_{\alpha \in (0,1]} \sup_{x\in (\rr^2)^N} |\cL^X\Gamma^\delta_\cK(x)|<\infty. 
$$
Fix $\delta >0$ and $\cK\in \gK_N$. We show in this step that
$$
\sup_{\alpha \in (0,1]} \sup_{x\in (\rr^2)^N} |\cL_\alpha ^X (G_\cK \Gamma^\delta_\cK)(x)| <\infty. 
$$
Let $\alpha \in (0,1]$ and $x = (x_1,\dots, x_N)\in (\rr^2)^N$. It is first classical that
$$
\cL^X_\alpha (G_\cK \Gamma^\delta_\cK)(x) = \Gamma^\delta_\cK\cL^X_\alpha G_\cK (x) +G_\cK(x)\cL^X_\alpha  \Gamma^\delta_\cK(x) +  \nabla G_\cK (x) \cdot \nabla \Gamma^\delta_\cK (x).
$$
The two last terms are bounded uniformly in $\alpha $ and $x$ by definition of $\Gamma^\delta_\cK$, it is thus sufficient to show that $\sup_{\alpha \in (0,1]}\sup_{x\in E^\delta_\cK} |\cL^X_\alpha G_\cK (x) |<\infty.$ Since $\Delta G_\cK$ is clearly bounded on $E^\delta_\cK$ we only need to bound the drift part. We have for all $x\in E^\delta_\cK$,
\begin{align*}
S:=& \sum_{i=1}^N \nabla_{x_i} G_\cK(x) \cdot \sum_{j\ne i} \frac{x_i-x_j}{\|x_i-x_j\|^2 + \alpha}\\
&= \sum_{K\in \cK}\sum_{i\in K}\sum_{k\notin K}\sum_{j\ne i}\frac{x_i-x_k}{\|x_i-x_k\|^2}\cdot  \frac{x_i-x_j}{\|x_i-x_j\|^2 + \alpha}\\
&=  \sum_{K\in \cK}\sum_{i\in K}\sum_{k\notin K}\sum_{j\in K}\frac{x_i-x_k}{\|x_i-x_k\|^2}\cdot  \frac{x_i-x_j}{\|x_i-x_j\|^2 + \alpha} + r(x),
\end{align*}
where $r$ a bounded function of $x \in E^\delta_\cK$. By symmetry we get for all $K\in \cK$, $k\notin K$,
$$
\sum_{i\in K}\sum_{j\in K}\frac{x_i-x_k}{\|x_i-x_k\|^2}\cdot  \frac{x_i-x_j}{\|x_i-x_j\|^2 + \alpha} =\frac12 \sum_{i\in K}\sum_{j\in K}\Big(\frac{x_i-x_k}{\|x_i-x_k\|^2} - \frac{x_j-x_k}{\|x_j-x_k\|^2}\Big)\cdot  \frac{x_i-x_j}{\|x_i-x_j\|^2 + \alpha}.
$$
Writing in complex notations we get
$$
\Big\|\frac{x_i-x_k}{\|x_i-x_k\|^2} - \frac{x_j-x_k}{\|x_j-x_k\|^2}\Big\| = \Big\|\frac{1}{\overline {x_i-x_k}} - \frac{1}{\overline {x_j-x_k}}\Big\| =\Big\| \frac{\overline{x_i-x_j}}{(\overline {x_i-x_k})(\overline{x_j-x_k})} \Big\| \le C_\delta'\|x_i-x_j\|,
$$
with $C_\delta'$ a constant only depending on $\delta.$.  All in all we get
\begin{align*}
|S|&\le \sup_{x\in E^\delta_\cK}\|r(x)\| + C_\delta'\sum_{K\in \cK} \sum_{i\in K} \sum_{k\notin K} \sum_{j\in K} \frac{\|x_i-x_j\|^2}{\|x_i-x_j\|^2 + \alpha } \\
&\le \sup_{x\in E^\delta_\cK}\|r(x)\| + C_\delta'\sum_{K\in \cK} |K|^2 (N-|K|) <\infty.
\end{align*}

{\it Step 3.2}
We conclude here. According to Step 3.1, the Definition of a $\eta$-pre-$KS$-process, and $\mathrm{Ito}(\theta,N-1)$ applied iteratively, we get that
\begin{align*}
   \sum_{k\ge 0}(G_{\cK^\eta_k}(X^1_{T\wedge \sigma^\eta_*\wedge\sigma^\eta_k}) - G_{\cK^\eta_{k+1}}(X^1_{T\wedge \sigma^\eta_*\wedge\sigma^\eta_{k+1}}))=& -\sum_{k\ge 0}\int_{\sigma^\eta_k}^{\sigma^\eta_{k+1}} \indiq_{s\le T\wedge\sigma^\eta_*}\nabla G_{\cK^\eta_k} (X^1_s)\cdot\dd B^1_s\\
   &-\sum_{k\ge 0}\int_{\sigma^\eta_k}^{\sigma^\eta_{k+1}} \indiq_{s\le T\wedge\sigma^\eta_*} \cL^XG_{\cK^\eta_k}(X^1_s)\dd s. 
\end{align*}
This implies that
\begin{align*}
    \frac{\dd \QQ^{1,T,\eta}}{\dd \PP^1} =& \exp \Big( \sum_{k\ge 0}(G_{\cK^\eta_k}(X^1_{T\wedge \sigma^\eta_*\wedge\sigma^\eta_k}) - G_{\cK^\eta_{k+1}}(X^1_{T\wedge \sigma^\eta_*\wedge\sigma^\eta_{k+1}})) \\
    &+  \sum_{k\ge 0} \int_{\sigma^\eta_k}^{\sigma^\eta_{k+1}}\indiq_{s\le T\wedge \sigma^\eta_*}( \cL^X G_{\cK^\eta_k}(X^1_s)- \frac12\|\nabla G_{\cK^\eta_k}(X^1_s)\|^2\dd s\Big)\\
:=&F_{T,\eta}((X^1_{t\wedge\sigma^\eta_*})_{t\in [0,T]}),
\end{align*}
where $F_{T,\eta} : C([0,T],(E_{k_0}\cap E^\eta_N)_\triangle)\to \rr$ is measurable. The same computation gives us that $\frac{\dd \QQ^{2,T,\eta}}{\dd \PP^2} = F_{T,\eta}((X^2_{t\wedge\sigma^\eta_*})_{t\in [0,T]})$. Thus if $H : C([0,T],(E_{k_0}\cap E^\eta_N)_\triangle)\to \rr$ is measurable and bounded, we have according to Lemma \ref{lem:etapreKS}-(i)
\begin{align*}
\E_{\PP^1}[H((X^1_{t\wedge\sigma^\eta_*})_{t\in [0,T]}) ] &= \E_{\QQ^{1,T,\eta}} [ (F_{T,\eta}((X^1_{t\wedge\sigma^\eta_*})_{t\in [0,T]}))^{-1}H((X^1_{t\wedge \sigma^\eta_*})_{t\in [0,T]}) ]\\
&=\E_{\QQ^{2,T,\eta}} [ (F_{T,\eta}((X^2_{t\wedge\sigma^\eta_*})_{t\in [0,T]}))^{-1}H((X^2_{t\wedge \sigma^\eta_*})_{t\in [0,T]}) ]\\
&=\E_{\PP^2}[H((X^2_{t\wedge\sigma^\eta_*})_{t\in [0,T]}) ],
\end{align*}
which proves the uniqueness. Since we obviously have \eqref{eq:boundonF}, we end the proof by showing that $F_{T,\eta} : C([0,T],(E_{k_0}\cap E^\eta_N)_\triangle) \to \rr$ is $\PP^{X,N,\theta,*}_\mu$-a.e. continuous for the uniform convergence topology. According to Lemma \ref{lem:etapreKS}-(ii) and the fact that $\cL_{\PP^1}((\indiq_{t<\sigma^\eta_*}X^1_t + \triangle \indiq_{t\ge \sigma^\eta_*})_{t\ge 0})$ and $\tilde \PP^{\eta,N,\theta}_\mu$ are equivalent, there exists a borelian set $A$ of $C([0,T],(E_{k_0}\cap E^\eta_N)_\triangle)$ such that $\PP^1( (\indiq_{t<\sigma^\eta_*}X^1_t + \triangle \indiq_{t\ge \sigma^\eta_*})_{t\in [0,T]} \in A) = 1$ and for all $\omega \in A$, the applications $\sigma^\eta_k\indiq_{\sigma^\eta_{k-1} < \sigma^\eta_*}$, $\indiq_{\sigma^\eta_{k-1} < \sigma^\eta_*}$, $\indiq_{\sigma^\eta_{k-1}<\sigma^\eta_*} \cK^\eta_k$ and $\sigma^\eta_*$ for all $k\ge 1$ are continuous at $\omega$ and $\int_0^\infty \indiq_{\omega_s \notin E_2\cup\{\triangle\}}\dd s = 0$. Let $\omega \in A$ and $\omega_n \to \omega$ for the uniform convergence topology on $C([0,T],(E_{k_0}\cap E^\eta_N)_\triangle)$. Since for all $k\ge 1$, $\indiq_{\sigma^\eta_{k-1} (\omega^n) < \sigma^\eta_*(\omega^n)}$ and $\indiq_{\sigma^\eta_{k-1} (\omega^n) < \sigma^\eta_*(\omega^n)}\cK^\eta_k(\omega^n)$ converges respectively to $\indiq_{\sigma^\eta_{k-1} (\omega) < \sigma^\eta_*(\omega)}$ and $\indiq_{\sigma^\eta _{k-1}(\omega) < \sigma^\eta_*(\omega)}\cK^\eta_k(\omega)$ as $n\to \infty$ by definition of $A$, we get
thanks to the discrete structure of indicators and $(\cK_k^\eta)_{k\ge 0}$ that for all $k\ge 1$, there exists $n_k(\omega)$ such that for all $n\ge n_k(\omega)$, 
\begin{gather*}
\indiq_{\sigma^\eta_{k-1} (\omega^n) < \sigma^\eta_*(\omega^n)} = \indiq_{\sigma^\eta_{k-1} (\omega) < \sigma^\eta_*(\omega)},\\ \indiq_{\sigma^\eta_{k-1}(\omega^n)<\sigma^\eta_*(\omega^n)} \cK^\eta_{k}(\omega^n) = \indiq_{\sigma^\eta_{k-1}(\omega)<\sigma^\eta_*(\omega)} \cK^\eta_{k}(\omega).
\end{gather*}
From this, we get that for all $k\ge 1$ and all $n_k(\omega)$, on the event $\{\sigma^\eta_k(\omega)<\sigma^\eta_*(\omega)\}$,
$$
|G_{\cK^\eta_k(\omega^n)} ( \omega^n_{\sigma^\eta_k(\omega^n)}) - G_{\cK^\eta_k(\omega)} ( \omega_{\sigma^\eta_k(\omega)})| = |G_{\cK^\eta_k(\omega)} ( \omega^n_{\sigma^\eta_k(\omega^n)}) - G_{\cK^\eta_k(\omega)} ( \omega_{\sigma^\eta_k(\omega)})|\le C \| \omega^n_{\sigma^\eta_k(\omega^n)} -\omega_{\sigma^\eta_k(\omega)} \|,
$$
which tends to $0$ as $n\to \infty$ since $\omega^n \to \omega$ in the uniform convergence topology and $\sigma^\eta_k(\omega^n) \to \sigma^\eta_k(\omega)$ as $n\to \infty$. Finally for all $k\ge 1$, all $n\ge n_k(\omega)$,
\begin{align*}
   \Big| \int_{\sigma^\eta_k(\omega^n)} ^{\sigma^\eta_{k+1}(\omega^n)} &\indiq_{s\le T\wedge \sigma^\eta_*}(\omega^n)( \cL^XG_{\cK^\eta_k(\omega)}(\omega^n_s) - \frac12 \| \nabla G_{\cK^\eta_k(\omega)}(\omega^n_s)\|^2))\dd s  \\
   &-\int_{\sigma^\eta_k(\omega)} ^{\sigma^\eta_{k+1}(\omega)} \indiq_{s\le T\wedge \sigma^\eta_*}(\omega)( \cL^XG_{\cK^\eta_k(\omega)}(\omega) - \frac12 \| \nabla G_{\cK^\eta_k(\omega)}(\omega_s)\|^2))\dd s \Big| \\
   &\le C\indiq_{\sigma^\eta_k(\omega) <\sigma^\eta_*(\omega)}\Big( \int_{\sigma^\eta_k(\omega)}^{\sigma^\eta_{k+1}(\omega)} |\cL^X G_{\cK^\eta_k(\omega)}(\omega^n_s) -\cL^X G_{\cK^\eta_k(\omega)}(\omega_s)|\dd s 
   \\
   &+\|\omega^n-\omega\| + |\sigma^\eta_k(\omega^n) - \sigma^\eta_k(\omega)|
   + |\sigma^\eta_{k+1}(\omega^n) - \sigma^\eta_{k+1}(\omega)| + |\sigma^\eta_*(\omega^n) - \sigma^\eta_*(\omega)| \Big),
\end{align*}
where $C$ is a positve constant. Since $\int_0^\infty \indiq_{\omega_s\notin E_2\cup\{\triangle\}}\dd s = 0$, for a.e. $s\in [\sigma^\eta_k(\omega), \sigma^\eta_{k+1}(\omega)]$, $\cL^X G_{\cK_k^\eta(\omega)}$ is continuous at $\omega_s$. Thus, the dominated convergence theorem gives that
$$
\int_{\sigma^\eta_k(\omega)}^{\sigma^\eta_{k+1}(\omega)} |\cL^X G_{\cK^\eta_k(\omega)}(\omega^n_s) -\cL^X G_{\cK^\eta_k(\omega)}(\omega_s)|\dd s \to 0 \quad \mbox{ as } n\to \infty.
$$
We conclude since all the other terms tends to $0$ as $n\to \infty$ by definition of $\omega$. 
\end{proof}
\begin{rk}\label{rk:girsanovopen}
Define for all $\eta >0$, all $\cK\in \gK_N$ the sequences of stopping time $(\tilde T^{\cK,\eta}_k)_{k\ge 0}$ and $(T^{\cK,\eta}_k)_{k\ge 1}$ by induction by setting $\tilde T^{\cK,\eta}_0 = 0$ and for all $k\ge 0$,
\begin{align*}
    T^{\cK,\eta}_{k+1} = \inf \Big\{ t\ge \tilde T^{\cK,\eta}_k : \min_{K\in \cK} \min_{i\in K,j\notin K} \|X^i_t - X^j_t\| \ge 2\eta \Big\}\wedge \sigma^\eta_*,\\
    \tilde T^{\cK,\eta}_{k+1} = \inf \Big\{ t\ge T^{\cK,\eta}_{k+1} : \min_{K\in \cK} \min_{i\in K,j\notin K} \|X^i_t - X^j_t\| \le \eta \Big\}\wedge \sigma^\eta_*. 
\end{align*}
    An alternate version of Proposition \ref{prop:girsanoveta}-(ii)-(c) (which we don't prove for the sake of conciseness and since it uses exactly the same arguments) is the following : 
    for all $T>0$, all $\eta >0$, all $\cK\in \gK_N$, $\dd\tilde \QQ^{\cK,T,\eta}_\mu := \cE(\tilde L^{\cK,\eta}_T)\dd \PP^{X,N,\theta,*}_\mu$ is a probability with
\begin{align*}
\tilde L^{\cK,\eta}_T = -\sum_{k\ge 0}\int_{T^{\cK,\eta}_k}^{\tilde T^{\cK,\eta}_{k+1}} \indiq_{s\le T}\nabla G_\cK (X^N_s)\cdot\dd B_s.
\end{align*}
Under $\tilde\QQ^{\cK,T,\eta}_\mu$, for all $k\ge 0$, 
$$
\Big\{((X^{N,i}_{(t+T^{\cK,\eta}_k)\wedge \tilde T^{\cK,\eta}_k})_{i\in K})_{t\ge 0},\quad K\in \cK \Big\}
$$
has the same law as $\{ (X^K_{t\wedge \tau})_{t\ge0}, K\in \cK \}$ where 
$$
\tau = \inf \{ t\ge 0 : \min_{K\ne L\in \cK} \min_{i\in \ig 1,|K|\id, j\in \ig 1,|L|\id} \|X^{K,i}_t - X^{L,j}_t\|\le \eta \Big\},
$$
the family $\{ (X^K_t)_{t\ge 0}, K\in \cK\}$ is mutually independent conditionally on $\sigma(\{X^K_0, K\in \cK\})$ and for all $K\in \cK$, the process  $(X^K_t)_{t\ge 0}$ is a $KS_{\mu_K}(\theta|K|/N, |K|)$-process with $\mu_K = \pi_K\#\mu$ where $\pi_K : x\in (\rr^2)^N \mapsto (x^i)_{i\in K} \in (\rr^2)^{|K|}$. 

\vip

Moreover, there exists $\tilde C_{\cK,N,\theta,\eta}\in (0,\infty)$ only depending on $\cK$, $N$, $\theta$ ,$\eta$  and a $\PP^{X,N,\theta,*}_\mu$-a.e. continuous application $\tilde F_{\cK,T,\eta} : C([0,T],(E_{k_0}\cap E^\eta_N)_\triangle)\to \rr$ for the uniform convergence topology such that $\dd \tilde\QQ^{\cK,T,\eta}_\mu = \tilde F_{\cK,T,\eta} \dd \PP^{X,N,\theta,*}_\mu$ and
\begin{align}\label{eq:boundonF2}
    \exp (-\tilde C_{\cK,N,\theta,\eta}T) \le \tilde F_{\cK,T,\eta} \le \exp( \tilde C_{\cK,N,\theta,\eta}T).
\end{align}
\end{rk}

We now explain why the barycenter and the up-to-a-constant-empirical measure of a $KS(\theta,N)$-process respectively behave like a Brownian motion and a squared Bessel process as suggested in Theorem \ref{thm:decompo}-(i).
\begin{prop}\label{prop:barycentreetbessel}
Let $\mu \in \cP((E_{k_0})_\triangle)$ and $(X^N_t)_{t\ge 0}$ be a $KS_\mu(\theta,N)$-process killed at $\tau_{k_0\wedge N}$ defined on a filtered probability space $(\Omega,\cF,(\cF_t)_{t\ge 0},\PP)$ with driving $2N$-dimensional Brownian motion $(B_t)_{t\ge 0}$. Define 
\begin{itemize}
\item the $2$-dimensional Brownian motion $(M_t)_{t\ge 0} = (S_{\ig 1,N\id}(X^N_0) + S_{\ig 1,N\id}(B_{t}))_{t\ge 0}$ with diffusive coefficient $N^{-1/2}$ starting at $S_{\ig 1,N\id}(X^N_0)$,
\item the squared Bessel process $(D_t)_{t\ge 0}$ starting at $R_{\ig 1,N\id}(X^N_0)$ with dimension $d_{\theta,N}(N)$ stopped when it reaches $0$ and driven by the standard Brownian motion $(W_t)_{t\ge 0}$ defined as
\begin{align*}
\mbox{ for all } t\ge 0, \quad W_t = \int_0^{t\wedge\tau_N}\Phi_{\sS_N}(X^N_s)\cdot \dd \pi_H(B)_s + \indiq_{t\ge \tau_N}(\tilde B_t - \tilde B_{\tau_N}),
\end{align*}
with $(\tilde B_t)_{t\ge 0}$ a standard $1$-dimensional Brownian motion independent from $(X^N_t,B_t)_{t\ge 0}$. 
\end{itemize}
We have

\vip (i) for all $t\in [0,\tau_{k_0\wedge N})$, $S_{\ig 1,N\id}(X^N_t) = M_t$,

\vip

(ii) for all $t\in [0,\tau_{k_0\wedge N})$, $R_{\ig 1,N\id}(X^N_t) = D_t$,

\vip

(iii) $(M_t)_{t\ge 0}$ and $(D_t)_{t\ge 0}$ are independent.  
\end{prop}
\begin{proof}
{\it Step 1.} We show (i) and (ii). Since $\sigma^\eta_* \to \tau_{k_0\wedge N}$, as $\eta \to 0$, it is sufficient to show (i) and (ii) on $[0,\sigma^\eta_*]$ for all $\eta >0$, which we now do. Fix $\eta >0$. We aim at applying Proposition \ref{prop:girsanoveta}-(ii)-(b) with $\varphi_1(x) = S_{\ig 1,N\id}(x)$ and $\varphi_2(x) = R_{\ig 1,N\id}(x)$.

\vip

\vip

Using symmetries, we have for all $\alpha \in (0,1]$, all $x\in (\rr^2)^N$,
$$
\cL^X_\alpha \varphi_1(x) = 0,
$$
which concludes (i). Then, since for all $i\in \ig 1,N\id$, $\nabla_{x_i} R_{\ig 1,N\id}(x) =2(x_i - S_{\ig 1,N\id}(x))$, we get
$$
\cL^X_\alpha \varphi_2 (x) = \frac{1}{2} \Delta R_{\ig 1,N\id} (x) -\frac{2\theta}{N} \sum_{i=1}^N (x_i - S_{\ig 1,N\id}(x)) \cdot \sum_{j\ne i}\frac{x_i-x_j}{\|x_i-x_j\|^2+\alpha}. 
$$
Then, we deduce that from the expression of $\nabla R_{\ig 1,N\id}$ that $\Delta R_{\ig 1,N\id} (x) = 4(N-1)$. Moreover, by using symmetries we end up with 
\begin{align*}
\cL^X_\alpha \varphi_2(x) &= 2(N-1) - \frac{\theta}{N} \sum_{i=1}^N\sum_{j\ne i} (x_i - x_j) \cdot \frac{x_i-x_j}{\|x_i-x_j\|^2+\alpha},\\
&=2(N-1) - \frac{\theta}{N} \sum_{i=1}^N\sum_{j\ne i} \frac{\| x_i-x_j \|^2}{\|x_i-x_j\|^2+\alpha}
\end{align*}
which is bounded uniformly in $\alpha$ and $x$. This same computation also gives $\cL^X\varphi_2(x) = (N-1)(2-\theta)$. Thus we can apply Proposition \ref{prop:girsanoveta}-(ii)-(b) and we get for all $t\in [0,\sigma^\eta_*) $
\begin{align*}
    R_{\ig 1,N\id} (X^N_{t}) = R_{\ig 1,N\id} (X^N_0) + \int_0^{t}\nabla R_{\ig 1,N\id}
 (X^N_s) \cdot \dd B_s + (N-1)(2-\theta) t.
 \end{align*}
We thus get for all $t\in [0, \sigma^\eta_*)$
\begin{align*}
    R_{\ig 1,N\id} (X^N_{t}) = R_{\ig 1,N\id} (X^N_0) + 2\int_0^{t}\sqrt{R_{\ig 1,N\id}(X^N_s)}\dd W_s+ (N-1)(2-\theta) t,
 \end{align*}
 recalling the definition of $(W_t)_{t\ge 0}$ in the statement. This concludes the proof of (ii) and of the Step. 

\vip

{\it Step 3.} We show (iii), but this is easy using pathwise uniqueness for squared Bessel processes (see for example \cite[Theorem 3.5 page 390]{ry}, the fact that for all $t\ge 0$, $\langle W, S_{\ig 1,N\id} (B)\rangle _ t = 0$ and the Lévy characterization theorem.
\end{proof}

The following result is an illustration that the construction we follow to build a $KS$-process gives a good understanding of it since we get regularity of the semigroup. 
\begin{prop}\label{prop:feller}
Set $k_* = k_0\wedge N$ if $d_{\theta,N}(N)<2$ and $k_* = N+1$ if $d_{\theta,N}(N) \ge 2$. There exists a measure space $ (\Omega,\cF,(\cF_t)_{t\ge 0})$ such that 
$$
\mM^{X,N,\theta,*} = (\Omega,\cF,(\cF_t)_{t\ge 0}, (X_t)_{t\ge 0}, (\PP^{X,N,\theta,*}_x)_{x\in (E_{k_*})_\triangle})
$$
is a Feller and strong Markov process where for all $x \in  (E_{k_*})_\triangle$, $\PP^{X,N,\theta,*}_x$ is the law of a $KS_x(\theta,N)$-process killed when exiting $E_{k_*}$. Moreover, for all $A\in \cB(C(\rr_+,(E_{k_0})_\triangle))$, $x\in (E_{k_*})_\triangle\mapsto \PP^{X,N,\theta,*}_x(A)$ is measurable, $\tau_{k_0\wedge N}$ is continuous for the one-point compactification topology $\PP^{X,N,\theta,*}_\mu$-a.e. and for all $t>0$, $\PP^{X,N,\theta,*}_\mu (\tau_{k_0\wedge N} = t) = 0$. 
\end{prop}
The proof for the case where $d_{\theta,N}(N)\ge 2$ use the same kind of argument adapted to this context and is slightly easier. Thus we do not give the proof in this case and focus in the case where $d_{\theta,N}(N)<2$. 
\begin{proof}
Adapting the proof of \cite[Corollary 2]{yw} slightly to our setting, where the process is defined as a solution to \eqref{eq:EDSproba} rather than to an SDE, and using Proposition \ref{prop:girsanoveta}-(i), $\mM^{X,N,\theta,*}$ is a strong Markov process such that for all $A\in \cB(C(\rr_+,(E_{k_0})_\triangle))$, $x\in (E_{k_*})_\triangle\mapsto \PP^{X,N,\theta,*}_x$(A) is measurable, so that it is sufficient to show that it is a Feller process. We write for all $\varphi \in C_0(E_{k_0\wedge N})$, all $x\in (E_{k_0\wedge N})_{\triangle}$, all $t\ge 0$,
$$
P^{N}_t\varphi (x) = \E^{X,N,\theta,*}_x[\varphi(X_t) \indiq_{t<\tau_{k_0\wedge N}}], 
$$
with $\E^{X,N,\theta,*}_x$ relative to $\PP^{X,N,\theta,*}_x$. Our goal now is to prove
\begin{align}\label{eq:Fellertime}
\mbox{ for all }\varphi \in C_0(E_{k_0\wedge N}),\quad  \sup_{x\in (E_{k_0\wedge N})_\triangle} |P^N_t\varphi (x) - \varphi (x)| \to 0, \quad \mbox{ as } t\to 0.
\end{align}
\begin{align}\label{eq:Fellerspace}
\mbox{ for all }t>0, \quad P^N_t C_0(E_{k_0\wedge N}) \subset C_0(E_{k_0\wedge N})
\end{align}
{\it Step 1.} We show that for all $\mu \in \cP((E_{k_0\wedge N})_\triangle)$, $\tau_{k_0\wedge N}$ is continuous for the one-point compactification topology  $\PP^{X,N,\theta,*}_\mu$-a.s. and for all $t>0$, $\PP^{X,N,\theta,*}_\mu (\tau_{k_0\wedge N} = t) = 0$.  

\vip

According to Proposition \ref{prop:barycentreetbessel}, we can write $\tau_{k_0\wedge N} = \tau_1\wedge \tau_2$ where
\begin{gather*}
\tau_1 =\lim_{\eta\to 0} \tau^\eta_1 = \lim_{\eta\to 0} \inf\{ t\ge 0 : D_t = \eta\},\\
\tau_2 = \lim_{\eta \to 0}\tau^\eta_2 = \lim_{\eta\to0}\inf\{ t\ge 0 : \mbox{ it exists } K\subset \ig 1,N\id : |K| \in \ig k_0,N\id : R_K(X_t)=\eta\}.
\end{gather*}
By classical results on squared Bessel processes (see \cite{ry}), there exists $A_1 \in \cF$ with $\PP^{X,N,\theta,*}_\mu(A_1) = 1$ and for all $\omega \in A_1$, $\tau_1$ and $\tau^\eta_1$ for all $\eta \in \QQ^*_+$ are continuous at $\omega$.  According to Lemma \ref{lem:etapreKS}-(iii) and Proposition \ref{prop:girsanoveta}-(ii)-(c), for all $\eta \in \QQ_+^*$, there exists $A^\eta_2\in \cF$ such that $\PP^{X,N,\theta,*}_\mu(A^\eta_2) = 1$ and for all $\omega \in A^\eta_2$, $\sigma^\eta_*$ is continuous at $\omega$.

\vip

We set $A = (A_1\cap \{ \tau_{k_0\wedge N} = \tau_1\})\cup \bigcup_{\eta\in \QQ^*_+} (A^\eta_2 \cap \{ \tau_{k_0\wedge N} < \tau^\eta_1\}) \in \cF $. Then $\PP^{X,N,\theta,*}_\mu (A) = 1$ and for all $\omega \in A$, $\tau_{k_0\wedge N}$ is continuous at $\omega$. Indeed, if there exists $\eta \in \QQ^*_+$ such that $\omega \in \{\tau_{k_0\wedge N} < \tau^\eta_{1}\}$, then $\tau_{k_0\wedge N}(\omega) = \sigma^\eta_*(\omega)$, thus $\sigma^\eta_*(\omega) < \tau^\eta_1(\omega)$. By continuity of $\sigma^\eta_*$ and $\tau^\eta_1$ at $\omega$, we get that for all $n$ large enough, $\tau_{k_0\wedge N}(\omega^n) = \sigma^\eta_*(\omega^n)$ so that $\tau_{k_0\wedge N}(\omega^n) = \sigma^\eta_*(\omega^n) \to \sigma^\eta_*(\omega) = \tau_{k_0\wedge N}(\omega)$.

In the other scenario where $\tau_{k_0}(\omega) = \tau_1(\omega)$, we set for all $\delta >0$,
$$
H_\delta : x=(x_t)_{t\ge 0} \in C(\rr_+, (E_{k_0\cap E_N)_\triangle}) \mapsto \inf_{t\in [0,\tau^1(x)-\delta)} \min_{K\subset \ig 1,N\id : |K|=k_0} R_K(x_t),
$$ 
which is continuous at $\omega$ by definition of $A$. Thus $H_\delta(\omega^n) \to H_\delta (\omega)$ as $n\to \infty$. Since $H_\delta(\omega)>0$ for all $\delta >0$ and $\delta\mapsto H_\delta (\omega)$ is continuous, we deduce that for all $\delta>0$, there exists $\e>0$ such that $\liminf_{n\to \infty}H_\delta (\omega^n) \ge \e$, so that there exists $\e'$ such that $\liminf_{n\to \infty} \tau^{\e'}_{k_0}(\omega^n) \ge \tau_1(\omega) - \delta$. However, $\tau_{k_0}\ge \tau^\e_{k_0}$ so that $\liminf_{n\to \infty}\tau_{k_0}\ge \tau_1 - \delta$. Letting $\delta \to 0$ and recalling that $\tau_{k_0}\le \tau_N$, we get that $\tau_{k_0}(\omega^n)\to\tau_1(\omega)$. This together with $\tau_1(\omega^n)\to \tau_1(\omega)$ as $n\to \infty$ gives the result. 

\vip

We can conclude according to Lemma \ref{lem:etapreKS}-(iii), Proposition \ref{prop:girsanoveta}-(ii)-(c) and the fact that it is classical that $\PP^{X,N,\theta,*}_\mu (\tau_1 = t) = 0$ for all $t>0$, see \cite{ry}.

\vip

{\it Step 2.} We show that for all $t\ge 0$, all $\varphi \in C_0(E_{k_0\wedge N})$, $P^N_t \varphi \in C_0(E_{k_0\wedge N})$. The proof for the case where $d_{\theta,N}(N)\ge 2$ would use the same kind of argument adapted to the context, so that we don't do the proof in this case for the sake of conciseness. 

\vip

Fix $\varphi \in C_0(E_{k_0\wedge N})$, $t>0$ and a sequence $(x_n)_{n\ge 0}$ of elements of $(E_{k_0\wedge N})_\triangle$ converging to some $x_\infty \in (E_{k_0\wedge N})_\triangle$. We show that $P^N_t\varphi (x_n) \to P^N_t\varphi (x_\infty)$ as $n\to \infty$.

\vip

{\it Step 2.1.} We treat the case where there exists $\eta >0$ such that for all $n\ge 0$, $x_n\in E^\eta_N$ (observe that $x_\infty$ can be equal to $\triangle$). 

\vip

Fix $\psi \in C_c(E_{k_0\wedge N})$. There exists $\eta >0$ such that $\mathrm{Supp}\psi \subset E_{k_0}\cap E^\eta_N$. According to Proposition \ref{prop:girsanoveta}-(ii)-(c) with $T=t$, for every $\eta'$ small enough, there exists $F_{T,\eta'} : C([0,T],(E_{k_0}\cap E^{\eta'}_N)_\triangle)\to \rr$ bounded by above and below by positive constants and continuous for the uniform convergence topology such that for all $n\in \nn\cup\{\infty\}$, 
$$
P^{N,\eta'}_t\psi (x_n) := \tilde \E^{\eta,N,\theta}_{x_n}[\indiq_{t<\sigma^{\eta'}_*}\psi (X_t)(F_{T,\eta'}((X_t)_{t\ge 0}))^{-1}],
$$ 
where $(X_t)_{t\ge 0}$ is the canonical process. 
We thus conclude that $P^{N,\eta'}_t\psi(x_n) \to P^{N,\eta'}_t\psi (x)$ as $n\to \infty$ by Lemma \ref{lem:etapreKS}-(iv). Moreover, we have
\begin{align*}
    |P^{N}_t\psi(x_n) - P^{N}_t\psi(x)| \le |P^{N,\eta'}_t\psi(x_n) - P^{N,\eta'}_t\psi(x)| + 2\|\psi\|_\infty \sup_{z\in (E_{k_0}\cap E^\eta_N)_\triangle}\PP^{\eta,N,\theta}_z[\sigma^{\eta'}_*\le t < \tau_{k_0\wedge N}].
\end{align*}
However, using the first Dini's theorem, we get that $\sup_{z\in E_{k_0}\cap E^\eta_N}\PP^{\eta,N,\theta}_z[\sigma^{\eta'}_*\le t < \tau_{k_0\wedge N}] \to 0$ as $\eta'\to 0$ since for all $z\in (E_{k_0}\cap E^\eta_N)_\triangle$, $\PP^{\eta,N,\theta}_z[\sigma^{\eta'}_*\le t < \tau_{k_0\wedge N}] \to 0$ as $\eta' \to 0$, and for every $\eta>0$ fixed, the application $z \in (E_{k_0}\cap E^\eta_N)_\triangle \mapsto \PP^{\eta,N,\theta}_z[\sigma^{\eta'}_*\le t < \tau_{k_0\wedge N}]$ is continuous for the one point compactification topology. Thus we get that $P^N_t\psi (x_n) \to P^N_t\psi (x)$ as $n\to \infty$. We finally conclude that $P^N_t\varphi (x_n) \to P^N_t\varphi (x)$ as $n\to \infty$ using the density of $C_c(E_{k_0}\cap E^\eta_N)$ in $C_0(E_{k_0}\cap E^\eta_N)$ for the uniform convergence topology. 

\vip

{\it Step 2.2.} We conclude the Step 2 by treating the case where $x = \triangle$ and $R_{\ig 1,N\id}(x_n)\to 0$ as $n\to \infty$. Let us show that $P^N_t\varphi (x_n) \to 0$ as $n\to \infty$. 

\vip

According to Proposition \ref{prop:barycentreetbessel}, there exists $(D^n_t)_{t\ge 0}$ a sequence of squared Bessel processes driven by the same Brownian motion respectively starting at $R_{\ig 1,N\id}(x_n)$ with dimension $d_{\theta,N}(N)$ stopped when it reaches $0$ such that
$$
|P^N_t \varphi (x_n)|  \le \|\varphi\|_\infty \PP[t<\tau_n], \quad \mbox{ with } \quad \tau_n = \inf\{t\ge 0 : D^n_t = 0\}.
$$
Using the comparison theorem \cite[Theorem 2.7 page 394]{ry}, we get that 
$$
|P^N_t \varphi (x_n)|  \le \|\varphi\|_\infty \PP[t<\tilde \tau_n], \quad \mbox{ with } \quad \tilde\tau_n = \inf\{t\ge 0 : \tilde D^n_t = 0\},
$$
with $(\tilde D^n_t)_{t\ge 0}$ the family of squared Bessel process with dimension $1\vee d_{\theta,N}(N) \in [1,2)$ starting at $R_{\ig 1,N\id}(x_n)$ all driven by the same Brownian motion as the one of the family $((D^n_t)_{t\ge 0})_{n\ge 0}$. This is classical by scaling and \eqref{eq:densitetau} that $\PP_{x_n} [t<\tilde \tau_n] =\PP_{1} [t/x_n<\tilde \tau_n] \le C (x_n/t)^{1-(1\vee d_{\theta,N}(N))/2}\to 0$ as $n\to \infty$. This allows us to conclude thanks to the Markov inequality.

\vip

{\it Step 3.} We prove \eqref{eq:Fellertime}. Fix $\varphi \in C_0(E_{k_0\wedge N})$. 

\vip

{\it Step 3.1.} We first prove that for all $\eta>0$, all $\psi \in C_c(E_{k_0}\cap E^\eta_N)$,
\begin{align}
\sup_{x\in E_{k_0}\cap E^\eta_N} \Big|\E^{X,N,\theta,*}_x[\psi(X_t) \indiq_{t< \sigma^\eta_*}] - \psi (x)\Big| \to 0 \quad \mbox{ as } t\to 0.
\end{align}
According Proposition \ref{prop:girsanoveta}-(ii)-(c) applied with $T=t$, we have
\begin{align*}
\sup_{x\in E_{k_0}\cap E^\eta_N} \Big|\E^{X,N,\theta,*}_x[\psi(X_t) \indiq_{t< \sigma^\eta_*}] - \psi(x)\Big| \le& \sup_{x\in E_{k_0}\cap E^\eta_N} \Big|\tilde\E^{\eta,N,\theta}_x\Big[\psi(X_t) \indiq_{t< \sigma^\eta_*}\Big( F_{T,\eta}((X_t)_{t\ge 0})- 1 \Big) \Big]\Big|\\
&+\sup_{x\in E_{k_0}\cap E^\eta_N} \Big| \tilde \E^{\eta,N,\theta}[\psi(X_t) \indiq_{t< \sigma^\eta_*}] - \psi (x) \Big|.
\end{align*}
Since the second term of the RHS converges to $0$ as $t\to 0$ according to Lemma \ref{lem:etapreKS}-(iv), it remains to show that the first term of the RHS converges as well to $0$ as $t\to 0$. Using the Cauchy-Schwarz inequality, one gets for all $x\in E_{k_0} \cap E^\eta_N$,
\begin{align*}
\Big|\tilde\E^{\eta,N,\theta}_x\Big[\psi(X_t) \indiq_{t< \sigma^\eta_*}\Big( F_{T,\eta}((X_t)_{t\ge 0}) - 1 \Big) \Big]\Big|\le \|\psi\|_\infty \Big(\tilde\E^{\eta,N,\theta}_x\Big[\Big( F_{T,\eta}((X_t)_{t\ge 0}) - 1 \Big)^2 \Big]\Big)^{1/2}.
\end{align*}
Expanding the square, using that $F_{T,\eta}\dd \tilde\PP_x^{\eta,N,\theta}$ is a probability measure and \eqref{eq:boundonF}, we get that for all $x\in E_{k_0}\cap E^\eta_N$, 
$$
\tilde \E^{\eta,N,\theta}_x\Big[\Big( F_{T,\eta}((X_t)_{t\ge 0}) - 1 \Big)^2 \Big] = \tilde\E^{\eta,N,\theta}_x\Big[F_{T,\eta}((X_t)_{t\ge 0})^2 \Big] -1 \le \exp ( 2C_{N,\theta,\eta}t) -1,
$$
which conclude this step since the last bound is independent from $x$ and converges to $0$ as $t\to 0$.

\vip

{\it Step 3.2.} We conclude Step 2 here. Let $\psi \in C_c(E_{k_0\wedge N})$. There exists $\eta >0$ such that $\mathrm{Supp}\psi \subset E_{k_0}\cap E^\eta_N$. We write
\begin{align}\label{eq:majfeller}
\Delta_t :=\sup_{x \in E_{k_0\wedge N}}|\E^{X,N,\theta,*}_x[\psi(X_t) \indiq_{t< \tau_{k_0\wedge N}}] - \psi (x) |\le& E^1_t + E^2_t,
\end{align}
where
\begin{align*}
E^1_t = \sup_{x \in E_{k_0\wedge N}}|\E^{X,N,\theta,*}_x[\psi(X_t) \indiq_{t< \sigma^\eta_*}] - \psi (x) |,\\
E^2_t = \sup_{x \in E_{k_0\wedge N}}\E^{X,N,\theta,*}_x[\sup_{s\in [0,t]}|\psi(X_s) \indiq_{\sigma^\eta_*\le t}]|.
\end{align*}
Observe that for all $x\in E_{k_0} \cap (E^\eta_N)^c$, one has $\sigma^\eta_* = 0$ so that, recalling the definition of $\eta$,
$$
E^1_t \le \sup_{x \in E_{k_0}\cap E^\eta_N}|\E^{X,N,\theta,*}_x[\psi(X_t) \indiq_{t< \sigma^\eta_*}] - \psi (x) | \to 0 \quad \mbox{ as } t\to 0,
$$
according to Step 3.1. Finally, we conclude  that $\Delta_t \to 0$ as $t\to 0$ since $E^2_t \to 0$ as $t\to 0$ because of the first Dini's theorem. Indeed, $(E_{k_0\wedge N})_\triangle$ is a compact set for the one point compactification topology and setting the application
$$
f : (t,x)\in \rr_+\times(E_{k_0\wedge N})_\triangle \mapsto \E^{X,N,\theta,*}_x\Big[\sup_{s\in [0,t]}|\psi(X_s) \indiq_{\sigma^\eta_*\le t}|\Big],
$$
we have for all $t\in \rr_+$, all $x\in (E_{k_0\wedge N})_\triangle$, $f(t,\cdot)$ continuous according to Lemma \ref{lem:etapreKS}-(iii) and $f(\cdot,x)$ nondecreasing with $f(0,x) = 0$ because $\mathrm{Supp}\psi \subset E_{k_0}\cap E^\eta_N$. Finally we conclude \eqref{eq:Fellertime} by density of $C_c(E_{k_0\wedge N})$ in $C_0(E_{k_0\wedge N})$ for the uniform convergence topology. 
\end{proof}
We will describe our process in terms of Dirichlet forms in order to bridge our work with \cite{ft}. We will frequently rely on \cite[Appendix B.]{ft} in the sequel where all the results of the Dirichlet form theory we will use are summarized in a concise way. Before that we need the following lemma as a preparation. We set for all $K\subset \ig 1,N\id$, all $u\in \sS_{|K|}$ the application $\phi^K_{u} : (\rr^2)^N \to (\rr^2)^N$ that maps $x\in (\rr^2)^N$ to the vector $y\in (\rr^2)^N$ such that
$$
\mbox{ for all } i\notin K,\quad  y_i = x_i \quad \mbox{ and for all } i\in K, \quad y_i = S_K(x) + R_K(x) u_i. 
$$
Finally, recall Theorem \ref{mainresult} for the definition of $\cE^N_{X,1}$ and $\cC_\cK$ for all $\cK\in \gK_N$.
\begin{lemma}\label{lem:core}
We have 
$$
\overline{C_c^\infty(E_{k_0\wedge n})}^{\cE^N_{X,1}} = \overline{\sum_{\cK \in \gK^n_N}\cC_\cK}^{\cE^N_{X,1}}.
$$
\end{lemma}
\begin{proof}
{\it Step 1.} We set for all $(K_p)_{1\le p\le m}$ and $(L_q)_{1\le q\le m'}$ such that $\{ K_p, 1\le p\le m\} \cup \{ L_q, 1\le q\le m'\} \in \gK_N$, 
\begin{align*}
\cC_{(K_p),(L_q)} =&C_c^\infty (E_{k_0}) \cap C_c^\infty (E_{(K_p)\cup (L_q)})\\
&\cap \Big\{ \varphi : x\in (\rr^2)^N\mapsto \psi ((S_{K_p}(x), R_{K_p}(x))_{1\le p\le m}, (x^k)_{k\in \cup_{q=1}^{m'}L_q}) \} .
\end{align*}
We first prove that for all family $(K_p)_{1\le p\le m}$, $(L_q)_{q\in \ig 1,m'\id}$ such that $\{ K_p, p\in \ig 1,m\id\} \cup \{ L_q, q\in \ig 1,m'\id\} \in \gK_N$ we get that for all $L_0\in (L_q)$
$$
\cC_{(K_p)\cup\{L_0\}, (L_q)\setminus \{L_0\}} + \sum_{L\subset L_0} \cC_{(K_p),((L_q)\setminus \{L_{0}\})\cup \{ L,L_{0}\setminus L\} } \mbox{ is dense in } \cC_{(K_p),(L_q)} \mbox{ for the } \cE^N_{X,1} \mbox{ norm}.
$$
Fix $(K_p)_{1\le p\le m}$, $(L_q)_{1\le q\le m'}$ and $L_0$ as above and $\varphi \in \cC_{(K_p), (L_q)}$. Let $u\in \sS_{|L_0|}$. We set 
\begin{gather*}
    \varphi_{L_{0}} = \varphi \circ \phi_u^{L_{0}} \in \cC_{(K_p)\cup\{L_{0}\}, (L_q)\setminus \{L_{0}\}},\\
    \psi_\e = (\varphi - \varphi _{L_{0}})(1-\chi( R_{L_{0}}/\e)) \in \sum_{L\subset L_{0}} \cC_{(K_p), ((L_q)\setminus \{L_{0}\} \cup \{ L,L_{0}\setminus L\}},
\end{gather*}
where $\chi : \rr_+ \to [0,1]$ smooth such that $\chi(r) = 1$ for all $r\in [0,1]$ and $\chi (r) = 0$ if $r\ge 2$. Indeed, $\mathrm{Supp}(1-\chi(R_{L_0}/\e)) \subset \{ x\in (\rr^2)^N : R_{L_0}(x)>0\}$ so that it can be decomposed as $\sum_{L\subset L_0} \chi_L$ with for all $L\subset L_0$, $\mathrm{Supp}\chi_L \subset \{x\in (\rr^2)^N : \inf_{i\in L, j\in L_0\setminus L} \|x^i-x^j\| >0\}$. 

\vip

It is thus sufficient to show that $\varphi - \varphi_{L_0} - \psi _\e \to0$ as $\e \to 0$ for the $\cE^N_{X,1}$ norm. Since $\varphi - \varphi_{L_0} - \psi_\e = (\varphi -\varphi _{L_0})\chi ( R_{L_0}/\e)$ we have $\cE^N_{X,1}(\varphi - \varphi_{L_0} - \psi_\e,\varphi - \varphi_{L_0} - \psi_\e) \le 2E^\e_1 +2E^\e_2 + E^\e_3$ where
\begin{gather*}
    E^\e_1 = \int_{E_{k_0}} \|\nabla \varphi - \nabla \varphi_{L_0}\|^2 (\chi (R_{L_0}/\e))^2\mu^N_\theta ,\\
    E^\e_2 = \int_{E_{k_0}} ( \varphi - \varphi_{L_0})^2\frac{\| \nabla R_{L_0}\|^2}{\e^2} (\chi' (R_{L_0}/\e))^2\mu^N_\theta ,\\
    E^\e_3 =\int_{E_{k_0}}( \varphi - \varphi_{L_0})^2 (\chi (R_{L_0}/\e))^2\mu^N_\theta .
\end{gather*}
Since $\varphi$ and $\varphi_{L_0}$ are smooth and compactly supported in $E_{k_0\wedge n}$ and that $\mu^N_\theta$ is Radon on $E_{k_0}$ thanks to Lemma \ref{lem:radon}-(i), it is clear by the dominated convergence theorem that $E^\e_1 + E^\e_3\to 0$ as $\e \to 0$. Next we treat $E^\e_2$. We have for $\e>0$ small enough,
\begin{align*}
    E^\e_{2} \le \frac{C}{\e}\int_{E_{k_0}} \indiq_{R_{L_0}\le 2\e} ( \varphi - \varphi_{L_0})^2\mu^N_\theta,
\end{align*}
with $C$ a constant independent of $\e$ allowed to change from line to line, since $\|\nabla R_{L_0}\|^2 =4 R_{L_0}$. Moreover, if $x\in (\rr^2)^N$ is such that $R_{L_0}(x) \le 2\e$ then by regularity of $\varphi$,
$$
|\varphi (x) - \varphi_{L_0}(x)|^2 \le C\| D\varphi\|^2_\infty \indiq_{\mathrm{Supp}(\varphi - \varphi_{L_0})} ( R_{L_0}(x) + (R_{L_0}(x))^2) \le C\| D\varphi\|^2_\infty \indiq_{\mathrm{Supp}(\varphi - \varphi_{L_0}) } \e.
$$
We thus conclude that 
$$
E^\e_{2}\le C\|D\varphi\|^2_\infty\int _{\mathrm{Supp}\varphi} \indiq_{R_{L_0}\le 2\e}\mu^N_\theta \to 0 \quad \mbox{ as } \e \to 0,
$$
according to the dominated convergence theorem because $\mathrm{Supp}\varphi$ is compact in $E_{k_0}$ and $\mu^N_\theta$ is Radon on the same space. 

\vip

{\it Step 2} We conclude here, but the result is clear by iterating the previous step since $C_c^\infty (E_{k_0}) = \cC_{\emptyset,\{\ig 1,N\id\}}$ and for all $\cK \in \gK_N$, $\cC_\cK = \cC_{\cK,\emptyset}$. 
\end{proof}

\begin{prop}\label{prop:dirichlettué}
The process $\mM^{X,N,\theta,*}$ defined in Proposition \ref{prop:feller} is a $\mu^{N}_\theta |_{E_{k_0\wedge N}}$-symmetric strong Markov process so is associated to the Dirichlet space $(\cE^{N,*}_X,\cF^{N,*}_X)$ on $L^2(E_{k_0\wedge N},\mu^{N}_\theta)$ in the sense of \cite{f} which is given by
\begin{align*}
\cF^{N,*}_X = \overline{C^\infty_c(E_{k_0\wedge N})}^{\cE^{N,*}_{X,1}},\quad \mbox{ and for all } f,g \in \cF^{N,*}_X, \quad \cE^{N,*}_X(f,g) = \frac12\int_{E_{k_0\wedge N}}\nabla f \cdot \nabla g \dd \mu^{N}_\theta,
\end{align*}
where $\cE^{N}_{X,1}(\cdot, \cdot) = \cE^{N}_X + (\cdot, \cdot)_{L^2(\E_{k_0\wedge N}, \mu^{N}_\theta)}$, see Theorem \ref{mainresult}. Moreover, $\sum_{\cK\in \gK_N}\cC_\cK$ is a core for $(\cE^{N,*}_X,\cF^{N,*}_X)$.
\end{prop}

\begin{proof}
 According to \cite[Proposition 2]{ft}, there exists a $(E_{k_0})_\triangle$-valued diffusion 
 $$
 \xX = (\Omega^X,\cF^X,(X_t)_{t\ge 0},(\PP^X_x)_{x\in (E_{k_0})_\triangle})
 $$
 which is $\mu^N_\theta$-symmetric  with regular Dirichlet space $(\cE^N_X,\cF^N_X)$ defined in Theorem \ref{mainresult}. We set
   $$
   \tau_{k_0\wedge N} = \inf\{ t\ge 0 : X_t \notin E_{k_0\wedge N}\}.
   $$
   According to \cite[Lemma B.6]{ft} (which directly comes from Theorems 4.4.2 page 173 and 4.4.3-(i) page 174 in \cite[Section 4.4]{f}), we get that 
   $$
   (\Omega^X,\cF^X, (\indiq_{t<\tau_{k_0\wedge N}}X_t + \triangle \indiq_{t\ge \tau_{k_0\wedge N}})_{t\ge 0}, (\PP^X_x)_{x\in (E_{k_0\wedge N})_\triangle})
   $$
   is a $\mu^N_\theta|_{E_{k_0\wedge N}}$-symmetric $(E_{k_0\wedge N})_\triangle$-valued diffusion with regular Dirichlet space $(\cE^{N,*}_X,\cF^{N,*}_X)$.
   \vip

   The goal of the proof is to show that there exists a properly exceptional set $\cN$ such that for all $x\in E_{k_0\wedge N}\setminus \cN$, the law of $(X^N_t)_{t\ge 0}$ under $\PP^{X,N,*}_x$ equals the law of $(X_t)_{t\ge 0}$ under $\PP^{X,*}_x$ with obvious notations.

\vip

We set for all $\eta>0$,
\begin{align}\label{eq:decoupetransition}
E^\eta = \bigcup_{\cK\in \gK_N \setminus \{\ig 1,N\id \}}E^\eta_\cK,
\end{align}
with 
\begin{align*}
E^\eta_\cK =& \{ x \in (\rr^2)^N : \mbox{ for all } K\in \cK, R_K(x)\le 2\eta^{1/|K|} \\
&\mbox{ and } \mbox{ for all } (K_i)_{i\in \ig 1,m\id} \subset \cK \mbox{ with } m\ge 2, R_{K_1\cup \dots \cup K_p}(x) \ge \frac12 \eta^{1/|K_1\cup\dots\cup K_p|}\}.
\end{align*}

\vip

{\it Step 1.}  We show that for all $\cK \in \gK_N \setminus \{\ig 1,N\id \}$, all $\eta >0$, the process $(\indiq_{t<\rho^\eta_\cK}X^{N}_t + \triangle \indiq_{t\ge \rho^\eta_\cK})_{t\ge 0}$ where 
$$
\rho^\eta_\cK := \inf \{ t\ge 0 : X^N_t \notin E^\eta_\cK\},
$$
is $\mu^{N}_\theta |_{E^\eta_\cK}$-symmetric and is associated with the Dirichlet space $(\cE^\eta_\cK,\cF^\eta_\cK)$ given by
$$
\cF^\eta_\cK = \overline{C^\infty_c(E^\eta_\cK)}^{(\cE^\eta_\cK)_1} \quad \mbox{ and for all } \varphi,\psi \in \cF^\eta_\cK, \quad \cE^\eta_\cK(\varphi, \psi) = \frac12 \int_{E^\eta_\cK}\nabla \varphi \cdot \nabla \psi \dd \mu^{N}_\theta.
$$

Let us fix $\eta >0$, $\cK\in \gK_N \setminus \{\ig 1,N\id \}$ and $(\Omega,\cF,(\cF_t)_{t\ge 0},\PP,(B_t)_{t\ge 0},(\tilde X^N_t)_{t\ge 0})$ a $\eta$-pre-$KS(\theta,N)$-process. According to $\cP(N-1)$-(a), \cite[Lemma B.5.]{ft} and \cite[Lemma B.6.]{ft}, we get that,
$$
(\indiq_{t<\rho^\eta_\cK}\tilde X^N_t+\triangle\indiq_{t\ge \rho^\eta_\cK})_{t\ge 0}
$$
is a $\tilde\mu^{N}_{\theta,\cK}|_{E^\eta_\cK}$-symmetric diffusion where for all $A\in \cB((\rr^2)^N)$,
$$
\tilde\mu^{N}_{\theta,\cK}(A) = \int_{A\cap E^\eta_\cK}   \prod_{K \in \cK} \prod_{i\ne j\in K}\|x^i-x^j\|^{-\theta/N} \dd x,
$$
and it is associated with the Dirichlet space $(\tilde\cE^\eta_\cK,\tilde\cF^\eta_\cK)$ given by 
$$
\tilde\cF^\eta_\cK = \overline{C^\infty_c(E^\eta_\cK)}^{\tilde\cE^\eta_{\cK,1} }\quad \mbox{ and for all } \varphi,\psi \in \cF, \quad \tilde\cE^\eta_\cK(\varphi, \psi) = \frac12 \int_{E^\eta_\cK}\nabla \varphi \cdot \nabla \psi \dd \tilde\mu^{N}_{\theta,\cK},
$$
with $\tilde\cE^\eta_{\cK,1}(\cdot,\cdot) = \cE(\cdot,\cdot) + (\cdot,\cdot)_{L^2(E^\eta_\cK,\tilde \mu^{N}_{\theta,\cK})}$.

\vip

We conclude this step according to Proposition \ref{prop:girsanoveta}-(ii)-(c), \cite[Lemma B.7.]{ft} and the fact that 
$$
\mu^{N}_\theta |_{E^\eta_\cK} = \tilde \mu^{N}_{\theta,\cK} \exp (- 2G_\cK) |_{E^\eta_\cK},
$$
where we recall \eqref{eq:GK} for the definition of $G_\cK$.

\vip

{\it Step 2.} We show in this step that for all $\eta >0$, there exists an exceptional set $\cN^\eta$ (relatively to the Dirichlet space $(\cE^{N,*},\cF^{N,*})$, see \cite[Appendix B.]{ft}) such that for all $x\in E^\eta\setminus \cN^\eta$, all $\cK \in \gK_N\setminus \{\ig 1,N\id\}$ such that $x\in E^\eta_\cK$, the law of $(X^N_{t\wedge \rho^\eta_{\cK}})_{t\ge 0}$ under $\PP^{X,N,*}_x$ equals the law of $(X_{t\wedge \rho^\eta_{\cK}})_{t\ge 0}$ under $\PP^{X,*}_x$. 

\vip

Using \cite[Theorem 4.2.8 p.167]{f} and Step 1, we get the existence for every $\cK\in \gK_N \setminus \{\ig 1,N\id \}$ of a properly exceptional set $\cN^\eta_\cK$ (relatively to the Dirichlet space $(\cE^\eta_\cK,\cF^\eta_\cK)$) such that the for all $x\in E^\eta_\cK\setminus \cN^\eta_\cK$, the above equality in law holds. It is thus sufficient to show that $\cup_{\cK \in \gK_N \setminus \{\ig 1,N\id \}}\cN^\eta_\cK$ is exceptional relatively to $(\cE^{N,*}_X, \cF^{N,*}_X)$. However it is clear that for any $\cK\in \gK_N \setminus \{\ig 1,N\id \}$, $\cN^\eta_K$ is exceptional relatively to $(\cE^{N,*}_X, \cF^{N,*}_X)$ according to \cite[Theorem 4.4.3-(ii) p.174]{f}. The fact that the finite union stays a properly exceptional set follows from the definition of an exceptional set see Appendix \cite[Appendix B.]{ft}. 

\vip

{\it Step 3.} We prove that for all $\eta >0$, all $x\in E^{\eta}\setminus  \cN^{\eta}$, the law of $(X^N_{t\wedge \sigma^{\eta}_*(X^N)})_{t\ge 0}$ under $\PP^{X,N,*}_x$ equals the law of $(X_{t\wedge \sigma^{\eta}_*(X)})_{t\ge 0}$ under $\PP^{X}_x$. It is sufficient to show by induction that $\cP(k)$ :''for all $x\in E^{\eta}\setminus \cN^{\eta}$ the law of $(X^N_{t\wedge \sigma^{\eta}_k(X^N)})_{t\ge 0}$ under $\PP^{X,N,*}_x$ equals the law of $(X_{t\wedge \sigma^{\eta}_k(X)})_{t\ge 0}$ under $\PP^{X}_x$.'' Of course $\cP(0)$ is true and we suppose that $\cP(k)$ holds for some $k\ge0$. We now prove $\cP(k+1)$. 

\vip

Observe that $x\in E^\eta_{\cK^\eta(x)}$, (recalling \eqref{eq:configeta} for the definition of $\cK^\eta (x)$). Thus, the result easily follows from the definition of a properly exceptional set, the definition of $\cN^{\eta}$ according to Step 2, the fact that for all $x\in E^{\eta}$, $\PP^{X,N,*}_x[\sigma^\eta_1 = \rho^{\eta}_{\cK^\eta(x)}]=1$ and the strong Markov property.

\vip
{\it Step 4.} We conclude the proof here reasoning exactly as in Step 2, by setting $\cN = \cup_{\eta \in \QQ_+^*}\cN^{\eta}$ and recalling that for all $x\in C(\rr_+, (\rr^2)^N_\triangle)$, $\sigma^\eta_*(x)\to \tau_{k_0\wedge N}(x)$ as $\eta \to 0$. Observe that it is clear that $\sum_{\cK\in \gK_N}\cC_\cK$ is a core for the Dirichlet space thanks to Lemma \ref{lem:core}.
\end{proof}
\begin{rk}
    One could wonder why we proved the Itô's formula for $KS$-processes in Proposition \ref{prop:girsanoveta}-(ii)-(b) by hand while applying \cite[Lemma B.2 page 59]{ft} gives directly a notion of Itô's formula using the Dirichlet forms theory point of view. Firstly, the point of view in \cite{f} would only hold for quasi every point while we prove it for any starting point. Secondly, even if the process had a Feller property that would allow to replace the ''quasi everywhere properties'' by ''everywhere'' properties (see for example \cite[Condition (4.2.9) page 165]{f}), this would lead to a circular reasoning since we need the use of the It\^o's formula for every starting point in order to get our notion of Feller property. 
\end{rk}
The following Lemma is adapted from \cite[Proposition 16, Lemma 18]{ft}. The proof is exactly the same, but in our case, since our process is well defined for every starting point, we get also the result for all starting point. 
\begin{lemma}\label{lem:limitusualtopo}
    Suppose that $k_0\le N$. For all $x \in E_{k_0}$, all $K\subset \ig 1,N\id$, we have that $\PP_x^{X,N,\theta}$- a.s, $X_{\tau_{k_0}^-} := \lim_{t\to \tau_{k_0}^-}X_t$ exists for the usual topology of $(\rr^2)^N$.
\end{lemma}
Observe that in the case where $k_0<N$, then $\PP_x^{X,N,\theta} = \PP_x^{X,N,\theta,*}$ for all $\mu \in \cP((E_{k_0})_\triangle$, since $\tau_{k_0}\le \tau_N$.

\vip

The following Lemma allows us to conclude the proof of $\cP(N)$-(d).
\begin{lemma}\label{lem:pN-1b}
    Suppose that $k_0\le N$. Then, for all $x \in E_{k_0}$, for all $r> 0$, 
    $$
    \PP_x^{X,N,\theta} (\lim_{t\to \tau_{k_0}^-} R_{\ig 1,N\id}(X_s) = r) = 0.
    $$
\end{lemma}
Notice that according to Proposition \ref{prop:barycentreetbessel}, the limit in Lemma \ref{lem:pN-1b} is well defined by continuity of squared Bessel processes. 
\begin{proof}
The following proof is very close to the one of Step 2.1 of Lemma \ref{lem:etapreKS}. Fix $x\in E_{k_0}$. Let $r>0$. We have
    \begin{align*}
        P:=  \PP_x^{X,N,\theta} (\lim_{t\to \tau_{k_0}^-} R_{\ig 1,N\id}(X_s) = r) \le \sum_{K\subset \ig 1,N\id : |K|\ge k_0} P_K,
    \end{align*}
    where 
    $$
    P_K = \PP_x^{X,N,\theta} \Big[ R_{\ig 1,N\id}(X_{\tau_{k_0}^-}) = r, \quad R_K(X_{\tau_{k_0}^-}) = 0\quad  \mbox{ and } \quad\min_{i\notin K}R_{K\cup\{i\}}(X_{\tau_{k_0}^-})>0 \Big].
    $$
    Of course, $P_{\ig 1,N\id} = 0$ since $r>0$. It only remains to show that for all $K\subset \ig 1,N\id$ with $|K|\in \ig k_0,N-1\id$, $P_K = 0$. Fix such $K$. We have $P_K\le \sum_{\alpha \in \QQ_+^*}\sum_{\eta \in \QQ_+^*}\sum_{T\in \nn} P_K^{\alpha,\eta,T}$ where
    $$
    P_K^{\alpha, \eta,T} = \PP^{X,N,\theta}_x \Big[ R_{\ig 1,N\id}(X_{\tau^\alpha_K}) = r,\quad \mbox{ for all } t\in [\alpha,\tau^\alpha_K), \min_{i\in K,j\notin K} \|X^i_t - X^j_t\|\ge \eta,\quad  \tau^\alpha_K = \tau_{k_0}\le T\Big],
    $$
    where $\tau^\alpha_K = \inf\{t\ge \alpha : R_K(X_t) = 0\}$, with the convention $\inf \emptyset = \infty$.

    \vip

    Fix $\alpha,\eta \in \QQ^*_+$ and $T\in \nn$. According to Remark \ref{rk:girsanovopen} and Lemma \ref{lem:dispersionformula}, there exists a positive constant $C$, a probability $\QQ$ and a family $\{ (R^K_t)_{t\ge 0},(S^K_t)_{t\ge 0},(R^{K^c}_t)_{t\ge 0},(S^{K^c}_t)_{t\ge 0}\}$ of processes mutually independent knowing $\cF_ \alpha$ such that 
    \begin{itemize}
        \item $(R^K_t)_{t\ge 0}$ and $(R^{K^c}_t)_{t\ge 0}$ are squared Bessel processes with respective dimension $d_{\theta,N}(|K|)$ and $d_{\theta,N}(|K^c|)$, with $R^{K^c}_0>0$,
        \item $(S^K_t)_{t\ge 0}$ and $(S^{K^c}_t)_{t\ge 0}$ are $2$-dimensional Brownian motion with respective diffusion coefficient $|K|^{-1/2}$ and $|K^c|^{-1/2}$,
        \item $P^{\alpha,\eta,T}_x \le C\QQ \Big[  R^{K^c}_{\tilde \tau^\alpha_N} + \frac{|K||K^c|}{N}\|S^K_{\tilde \tau^\alpha_N} - S^{K^c}_{\tilde \tau^\alpha_N}\|^2 = r\Big]$,
    \end{itemize}
    where $\tilde\tau^\alpha_K = \inf\{t\ge \alpha : R^K_t = 0\}$. However, since for example $(R^{K^c}_t)_{t\ge 0}$ and $\tilde \tau^\alpha_K$ are independent and for all $t> 0$, $R^{K^c}_t$ has a density, we get that $P^{\alpha,\eta,T}_K = 0$, which concludes.    
\end{proof}

\section{Decomposition result and construction of $SKS$-process}\label{sec:decompo}
This section is devoted to show $\mathrm{MainTheorem}(\theta,N)$-(ii) and $\mathrm{Decomposition}(\theta,N)$-(i). In the case where $k_0\ge N+1$, the results proved in this Section are important in their own right but also, together with Section \ref{sec:ergodicity}, will allow us to build a candidate to extend the already built $KS(\theta,N)$-process killed when exiting $E_N$ after its life-time.
\subsection{Existence of $SKS$-process and Decomposition before $N$-collision}\label{subsec:existenceSKS}
This subsection aims at building a $SKS(\theta,N)$-process for which the barycenter-radius-angle decomposition is fulfilled before the first $N$-collision (see $\mathrm{Decomposition}(\theta,N)$-(i)).
\begin{proof}[Proof of existence in $\mathrm{Main Theorem}(\theta,N)$-(ii)]
Let $N\ge 2$, $\theta>0$, $x\in E_{k_0\wedge N}$. 

\vip

We show for all $u\in \sS\cap E_{k_0}$ the existence of a $SKS_u(\theta,N)$-process
$$
 (\Omega,\cF,(\tilde\cF_t)_{t\ge 0}, \PP, (\tilde B_t)_{t\ge 0},(U^N_t)_{t\ge 0})
$$
and $\mathrm{Decomposition}(\theta,N)$-(i) at the same time. In fact, concerning existence, we show the stronger result which is that for all $T>0$, all $\rr_+^*$-valued process $(Z_t)_{t\ge 0}$ independent from $(U^N_t)_{t\ge 0}$, we have $\sup_{t\in [0,T\wedge \tau_{k_0}(U^N))} \| U^{N,\e,Z}_t-U^N_t\| \to 0$ in probability as $\e \to 0$, where
\begin{gather*}
U^{N,\e,Z}_t = u + \int_0^t\pi_{(U^N_s)^\perp}\pi_H \dd \tilde B_s+\int_0^t b_{\e Z_s} (U^N_s)\dd s  + \Big(\frac{\theta (N-1)}{2}-\frac{2N-3}{2}\Big)\int_0^t U^N_s\dd s,
\end{gather*}
see the paragraph above Definition \ref{def:KS/KSK} for the definitions of $\pi_H$ and $\pi_{(U^N_s)^\perp}$, \eqref{eq:bepsilon} for the definition of $b_{\e'}$ for all $\e'>0$ and recall \eqref{eq:tempsarret} for the definition of $\tau_{k_0}(U^N)$. Fix $u\in \sS_{k_0}$ and a $\rr_+^*$-valued process $(Z_t)_{t\ge 0}$ independent from  $(U^N_t)_{t\ge 0}$ and upperbounded by $c^{-1}$ and bounded from below by $c$ where $c$ is a deterministic positive constant. Let $(\Omega,\cF,(\cF_t)_{t\ge 0}, (B_t)_{t\ge 0}, (X^N_t)_{t\ge 0})$ a $KS_u(\theta,N)$-process, which is meaningful since $u\in \sS_{k_0} \subset E_{k_0\wedge N}$. We set for all $t\in [0, \tau_{k_0\wedge N})$, $\rho_t = A^{-1}_t$ with $A_t = \int_0^t(R_{\ig 1,N\id}(X^N_s))^{-1}\dd s$ and 
$$
U^N_t =\Phi (X^N_{\rho_t}),
$$
where $\Phi : x\in E^\eta_{k_0\wedge N}\mapsto \pi_H(x)/\|\pi_H(x)\| = (x-\gamma (x))/\|x-\gamma (x)\|$. We set for all $t\in [0,\tau_{k_0\wedge N})$, 
 $$
 \hat X^{N,\e}_t = x + B_t + \int_0^t b_{\e R_{\ig 1,N\id}(X^N_s)Z_s}(X^N_s) \dd s.
$$
We recall that according to Lemma \ref{lem:horlogeexplosion}, $\rho : [0, A_{\tau_{k_0}}) \to [0,\tau_{k_0\wedge N})$ is a.s. a bijection.
 
\vip
 
{\it Step 1} Fix $\e >0$. We set for all $t\ge 0$, $\tilde \cF_t = \cF_{\rho_t}$, which is well defined since $\rho_t$ is a stopping time. We show that there exists $(\tilde B_t)_{t\ge 0}$ a $(\tilde \cF_t)_{t\ge 0}$-adapted $2N$-dimensional Brownian motion such that for all $t\in [0,A_{\tau_{k_0}})$,
\begin{align}\label{eq:ff}
U^{N,\e}_t := \Big( \frac{R_{\ig 1,N\id}(\hat X^{N,\e}_{\rho_t})}{R_{\ig 1,N\id}(X^N_{\rho_t})} \Big)^{1/2}\Phi (\hat X^{N,\e}_{\rho_t}) =& u + \tilde B_t - \gamma (\tilde B_t) - \int_0^t U^{N,\e}_s U^N_s \cdot \dd \tilde B_s + \int_0^t b_{\e Z_s}(U^N_s)\dd s \notag\\
&+\Big( \frac{\theta (N-1)}{2} - \frac{2(N-1)}{2} + \frac32 \Big) \int_0^t U^{N,\e}_s\dd s - \int_0^t U^N_s\dd s.
\end{align}
 
\vip

Using the Itô's formula, we get for all $i\in \ig 1,N\id$, all $e\in \{ 1,2\}$, all $t\in [0,\tau_{k_0\wedge N})$,
\begin{align}\label{ccccc}
\Big( \frac{R_{\ig 1,N\id}(\hat X^{N,\e}_t)}{R_{\ig 1,N\id}(X^N_t)} \Big)^{1/2}\Phi ( \hat X^{N,\e}_t)^{i,e} =& u^{i,e} +I^1_t + I^2_t + I^3_t + I^4_t 
\end{align}
where $u^i = (u^{i,1},u^{i,2})$, $\Phi(\cdot)^i = (\Phi(\cdot)^{i,1},\Phi(\cdot)^{i,2})$ and 
\begin{gather*}
I^1_t = \int_0^t \frac{\dd (\hat X^{N,\e,i,e}_s - S_{\ig 1,N\id}(\hat X^{N,\e}_s)^e)}{(R_{\ig 1,N\id}(X^{N}_s))^{1/2}}\\
I^2_t = - \frac{1}{2} \int_0^t \frac{\hat X^{N,\e,i,e}_s-S_{\ig 1,N\id}(\hat X^{N,\e}_s)^e}{(R_{\ig 1,N\id}(X^{N}_s))^{3/2}}\dd R_{\ig 1,N\id}(X^{N}_s)\\
I^3_t = -\frac{1}{2} \int_0^t \frac{\dd \langle R_{\ig 1,N\id}(X^{N}_\cdot), X^{N,\e,i,e}_\cdot - S_{\ig 1,N\id}(X^{N,\e}_\cdot)^e \rangle _s}{R_{\ig 1,N\id}(X^{N}_s)^{3/2}}\\
I^4_t = \frac{3}{8} \int_0^t \frac{\hat X^{N,\e,i,e}_s - S_{\ig 1,N\id}(\hat X^{N,\e})^e}{(R_{\ig 1,N\id}(X^{N}_s))^{5/2}} \dd \langle R_{\ig 1,N\id}(X^N_\cdot)\rangle_s,
\end{gather*}
with
\begin{gather*}
X^{N,\e,i}_\cdot = ( X^{N,\e,i,1}_\cdot,X^{N,\e,i,2}_\cdot),\\
\hat X^{N,\e,i}_\cdot = (\hat X^{N,\e,i,1}_\cdot,\hat X^{N,\e,i,2}_\cdot),\\
S_{\ig 1,N\id}(X^{N,\e}_\cdot) = (S_{\ig 1,N\id}(X^{N,\e}_\cdot)^1,S_{\ig 1,N\id}(X^{N,\e}_\cdot)^2).
\end{gather*}

\vip

First, for all $t\in [0,\tau_{k_0\wedge N})$, we get according to Proposition \ref{prop:barycentreetbessel},
\begin{align*}
\langle R_{\ig 1,N\id}(X^{N}_\cdot), \hat X^{N,\e,i,e}_\cdot - S_{\ig 1,N\id}(\hat X^{N,\e}_\cdot)\rangle_t &\!=\! 2\Big\langle \int_0^\cdot (R_{\ig 1,N\id}(X^N_s))^{1/2}\Phi (X^N_s)\cdot \dd B_s ,  B^{i,e}_\cdot \!- \!S_{\ig 1,N\id}(B_\cdot)^e\Big\rangle_t\\
&\!=\!2\int_0^t (R_{\ig 1,N\id}(X^N_s))^{1/2}\Phi(X^N_s)^{i,e} \dd s,
\end{align*}
because $X^N_s \in E_2$ for a.e. $s \in [0,\tau_{k_0\wedge N})$ according to by definition and for all $x\in E_2$, $\Phi (x) \in \sS$ so that $\sum_{i=1}^N \Phi (x)^i = 0$. We deduce that 
\begin{align}
I^2_t +I^4_t =& - \int_0^t \frac{(\hat X^{N,\e,i,e}_s - S_{\ig 1,N\id}(\hat X^{N,\e})^e)\Phi (X^N_s) \cdot \dd B_s}{R_{\ig 1,N\id}(X^{N}_s)} \label{aaaaa1}\\
&+\Big(\frac{\theta(N-1)}{2}-\frac{2(N-1)}{2}+\frac32\Big)\int_0^t \frac{(\hat X^{N,\e,i,e}_s - S_{\ig 1,N\id}(\hat X^{N,\e})^e)}{(R_{\ig 1,N\id}(X^{N}_s))^{3/2}} \dd s, \notag \\
I^3_t =& -\int_0^t \frac{\Phi (X^{N}_s)^{i,e}}{R_{\ig 1,N\id}(X^{N}_s)}\dd s\label{aaaaa2},
\end{align}
Then, recalling \eqref{eq:bepsilon}, we have for all $x\in E_2$, all $\lambda >0$, $b_{\lambda R_{\ig 1,N\id}(x)}(x) = (R_{\ig 1,N\id}(x))^{-1/2} b_{\lambda}(\Phi (x))$. We deduce that for all $t\in [0,\tau_{k_0\wedge N})$,
\begin{align}\label{bbbbb}
I^1_t = \int_0^t \frac{\dd (B^{i,e}_s - S_{\ig 1,N\id}(B_s)^e)}{(R_{\ig 1,N\id}(X^{N}_s))^{1/2}} +\int_0^t \frac{b_{\e Z_s}(\Phi(X^N_s))^{i,e}}{R_{\ig 1,N\id}(X^N_s)}\dd s,
\end{align}
where $B^i_\cdot = (B^{i,1}_\cdot, B^{i,2}_\cdot)$. Gathering \eqref{ccccc}, \eqref{aaaaa1}, \eqref{aaaaa2} and \eqref{bbbbb}, specifying at $\rho_t$, we get for all $t\in [0,\tau_{k_0\wedge N})$,
\begin{align}\label{dd}
U^{N,\e}_t =& u + \int_0^{\rho_t} \frac{\dd (B_s - \gamma (B_s))}{(R_{\ig 1,N\id}(X^{N}_s))^{1/2}} -  \int_0^{\rho_t} \frac{U^{N,\e}_{A_s} U^N_{A_s} \cdot \dd B_s}{(R_{\ig 1,N\id}(X^N_s))^{1/2}}-\int_0^t U^{N}_s\dd s \notag\\
&+\int_0^t b_{\e Z_s}(U^N_s)\dd s+ \Big(\frac{\theta(N-1)}{2}-\frac{2(N-1)}{2}+\frac32\Big)\int_0^t U^{N,\e}_s\dd s,
\end{align}
where we have performed the change of variable $s' = \rho_s$ in the non-stochastic integrals and then written $s$ instead of $s'$. However, $
(\tilde B_t)_{t\ge 0} = (\int_0^{\rho_t} (R_{\ig 1,N\id}(X^N_s))^{-1/2}\dd B_s)_{t\ge 0}$ is a $2N$-dimensional Brownian motion. Indeed, $(\tilde B_t)_{t\ge 0}$ is a martingale and since $\rho_t$ is a stopping time, we get for all $t\in [0,A_{\tau_{k_0}})$,
\begin{align*}
\langle \tilde B \rangle_t = \Big\langle \int_0^\cdot \frac{\dd B_s}{(R_{\ig 1,N\id}(X^N_s))^{1/2}} \Big\rangle_{\rho_t} = A_{\rho_t} = t.
\end{align*}
Inserting in \eqref{dd}, we get \eqref{eq:ff}.

\vip

{\it Step 2} We show that for all $T\ge 0$, $\sup_{t\in [0,T\wedge \tau_{k_0\wedge N}(U^N))} \| U^{N,\e}_t - U^N_t\| \to 0$ in probability as $\e \to 0$. Indeed, let $\delta >0$ and $\nu >0$. According to the definitions of $U^{N,\e}_t$ and $U^N_t$ and the fact that $x\in (\rr^2)^N\mapsto x - \gamma (x)$ is Lipschitz, we get the existence of a constant $c_N$ only depending on $N$ such that
\begin{align}\label{kk}
P^{1,\nu}_{\e,T} :=& \PP \Big[ \sup_{t\in [0,T\wedge A_{\tau_{k_0 \wedge N}})} \| U^{N,\e}_t - U^N_t\|\ge \frac\nu 3\Big]\notag\\
\le& \PP \Big[ \sup_{t\in [0,\rho_T\wedge \tau_{k_0\wedge N})} \| \hat X^{N,\e}_t - X^N_t\|\ge c_N\inf_{s\in [0,\rho_T\wedge\tau_{k_0\wedge N)}}(R_{\ig 1,N\id}(X^N_s))^{1/2}\frac\nu 3\Big]\notag\\
\le& \PP \Big[ \sup_{t\in [0,M\wedge \tau_{k_0\wedge N})} \| \hat X^{N,\e}_t - X^N_t\|\ge c_NM^{-1/2}\frac\nu 3, \sup_{t\in [0,M\wedge \tau_{k_0\wedge N})} R_{\ig 1,N\id}(X^N_t)\le M\Big] \\
&+ \PP [ \rho_T\wedge \tau_{k_0\wedge N} \ge M]\notag + \PP \Big[\inf_{s\in [0,\rho_T\wedge \tau_{k_0\wedge N})}R_{\ig 1,N\id}(X^N_s) < M^{-1}\Big]\notag \\
&+ \PP\Big[\sup_{s\in [0,\rho_T\wedge \tau_{k_0\wedge N})}R_{\ig 1,N\id}(X^N_s) > M\Big]\notag,
\end{align}
where $M>0$ is chosen such that
$$
\PP[\rho_T\wedge \tau_{k_0\wedge N}\ge M]+\PP \Big[\inf_{s\in [0,\rho_T\wedge\tau_{k_0\wedge N})}R_{\ig 1,N\id}(X^N_s) < M^{-1}\Big]+ \PP\Big[\sup_{s\in [0,\rho_T\wedge \tau_{k_0\wedge N})}R_{\ig 1,N\id}(X^N_s) > M\Big] \le \delta,
$$ which is possible since according to Proposition \ref{prop:barycentreetbessel} and Lemma \ref{lem:horlogeexplosion}, 
\begin{gather*}
\inf_{s\in [0,\rho_T\wedge \tau_{k_0\wedge N})}R_{\ig 1,N\id}(X^N_s) >0 \mbox{ a.s.},\\ \rho_T\wedge \tau_{k_0\wedge N} +\sup_{s\in [0,\rho_T\wedge \tau_{k_0\wedge N})}R_{\ig 1,N\id}(X^N_s) <\infty \mbox{ a.s.}
\end{gather*}
Let for all $K\in \cB^N_\theta$,
\begin{align*}
    \chi^\e_K : (x,z)\in  (\rr^2)^{N} \times \rr \mapsto \indiq_{R_K(x) \ge \e f(x,z)}\\
    f (x,z) =(c^{-1}M)\wedge(R_{\ig 1,N\id}(x)z)\vee (cM^{-1}).
\end{align*}
The family $\{\chi^\e_K, K\in \cB^N_\theta\}$ is a family of collisional approximations. The only subtle point of the definition of a collisional family to check is the last one, which holds since $R_{\ig 1,N\id}$ is invariant under $\cS_K$ for all $K\in \cB^N_\theta$ since $R_{\ig 1,N\id} = R_K + R_{K^c} + |K||K^c|\|S_K-S_{K^c}\|^2/N$ according to Lemma \ref{lem:dispersionformula}. We can thus apply Proposition \ref{prop:girsanoveta}-(ii)-(a), concluding the step. 

\vip

{\it Step 3.} We set for all $\nu>0$, all $T>0$, 
$$
P^\nu_{\e,T} = \PP \Big[ \sup_{t\in [0,T\wedge A_{\tau_{k_0\wedge N}})}\| U^N_t - U^{N,\e,Z}_t\|\ge \nu\Big],
$$
and show that $P^\nu_{\e,T} \to 0$ as $\e \to 0$. According to Step 1 we get 
\begin{align}\label{ii}
P^\nu_{\e,T}\le P^{1,\nu}_{\e,T} + P^{2,\nu}_{\e,T} + P^{3,\nu}_{\e,T},
\end{align}
where we recall the definition of $P^{1,\nu}_{\e,T}$ in Step 2 above and
\begin{gather*}
P^{2,\nu}_{\e,T} = \PP \Big[\sup_{t\in [0,T\wedge A_{\tau_{k_0\wedge N}})} \Big\| \int_0^t (U^{N,\e}_s-U^N_s) U^N_s\cdot \dd \tilde B_s \Big\| > \frac\nu 3\Big],\\
P^{3,\nu}_{\e,T} = \PP \Big[\sup_{t\in [0,T\wedge A_{\tau_{k_0\wedge N} })} \Big\| \int_0^t (U^{N,\e}_s-U^N_s)\dd s \Big\| > c_{\theta,N}\frac\nu 3\Big], 
\end{gather*}
where $c_{\theta,N} = 2|\theta (N-1) - (2N-5)|^{-1}$. We suppose that $c_{\theta,N}$ is not equal to infinity, in the other case we simply set $P^{3,\nu}_{\e,T} = 0$. We easily have that 
\begin{align}\label{jj}
P^{3,\nu}_{\e,T} \le P^{1,c_{\theta,N}\nu T^{-1}}_{\e,T}.
\end{align}
Let $\delta >0$. We set for all $i\in \ig 1,N\id$, all $e\in \{ 1,2\}$, 
$$
(M^{i,e}_t)_{t\ge 0} :=(\int_0^{t\wedge A_{\tau^\eta_{k_0\wedge N}}} (U^{N,\e,i,e}_s-U^{N,i,e}_s)U^N_s\cdot \dd \tilde B_s)_{t\ge 0},
$$ 
where $U^{N,i}_\cdot = (U^{N,i,1}_\cdot,U^{N,i,2}_\cdot)$ and $U^{N,\e,i}_\cdot = (U^{N,\e,i,1}_\cdot,U^{N,\e,i,2}_\cdot)$. Since for all $i\in \ig 1,N\id$, all $e\in \{ 1,2\}$ $(M^{i,e}_t)_{t\ge 0}$ is a local martingale we get thanks to the Dubins-Schwarz theorem, even if it means to enlarge the probability space, the existence of $(W^{i,e}_t)_{t\ge 0}$ a $1$-dimensional Brownian motion starting at $0$ such that for all $t\ge 0$,
$$
\int_0^{t\wedge A_{\tau^\eta_{k_0\wedge N}}} (U^{N,\e,i,e}_s-U^{N,i,e}_s)U^N_s\cdot \dd \tilde B_s = W^{i,e}_{\langle M^{i,e}\rangle_t}.
$$
Because of the equivalences of the norms in finite dimension, there exists $C_N<\infty$ such that for all $t\in [0,A_{\tau_{k_0\wedge N}})$,
$$
\Big\| \int_0^{t}(U^{N,\e}_s - U^N_s)U^N_s\cdot \dd\tilde B_s\Big\| \le C_N\sum_{i=1}^N\sum_{i\in \{ 1,2\}} \Big| \int_0^{t}(U^{N,\e,i,e}_s - U^{N,i,e}_s)U^N_s\cdot \dd\tilde B_s\Big|,
$$
so that thanks to the triangular inequality,
\begin{align}\label{hh}
P^{2,\nu}_{\e,T} &\le \sum_{i=1}^N \sum_{e\in \{ 1,2\}} \PP \Big[\sup_{t\in [0,T\wedge A_{\tau_{k_0\wedge N})}} | W^{i,e}_{\langle M^{i,e}\rangle_t} | > \frac\nu {6NC_N}\Big]\\
&= \sum_{i=1}^N \sum_{e\in \{ 1,2\}}\PP \Big[\sup_{t\in [0,\langle M^{i,e}\rangle _{T\wedge A_{\tau_{k_0\wedge N} }})} | W^{i,e}_t | > \frac\nu {6NC_N}\Big] \notag \\
&\le \sum_{i=1}^N \sum_{e\in \{ 1,2\}} \PP \Big[\sup_{t\in [0,\alpha]} | W^{i,e}_t | > \frac\nu {6NC_N}\Big] + \PP\Big[\langle M^{i,e}\rangle _{T\wedge A_{\tau_{k_0\wedge N} }}\ge \alpha\Big], \notag
\end{align}
where $\alpha>0$ is small enough such that $\sum_{i=1}^N \sum_{e\in \{ 1,2\}}\PP \Big[\sup_{t\in [0,\alpha]} | W^{i,e}_t | > \frac\nu {6NC_N}\Big]\le \delta$. However, we have for all $i\in \ig 1,N\id$, all $e\in \{ 1,2\}$,
\begin{align}\label{gg}
\PP\Big[\langle M^{i,e}\rangle _{T\wedge A_{\tau_{k_0\wedge N}}}\ge \alpha\Big] = \PP \Big[ \int_0^{T\wedge A_{\tau_{k_0\wedge N}}} (U^{N,\e,i,e}_s - U^{N,i,e}_s)^2\dd s \ge \alpha \Big]\le P^{1,(\alpha (3T)^{-1})^{1/2}}_{\e,T}.
\end{align}
So that, combining \eqref{ii}, \eqref{jj}, \eqref{hh}, \eqref{gg}, the definition of $\alpha$ and Step 2, we get that
$$
P^\nu_{\e,T} \le P^{1,\nu}_{\e,T} +P^{1,c_{\theta,N}\nu T^{-1}}_{\e,T} + P^{1,(\alpha T^{-1})^{1/2}}_{\e,T} + \delta \le 2\delta,
$$
if $\e$ is sufficiently small, giving the result since, according to Lemma \ref{lem:horlogeexplosion}, $A_{\tau_{k_0\wedge N}(X^N)} = \tau_{k_0}(U^N) $. 
\end{proof}

\subsection{Uniqueness in law of $SKS(\theta,N)$-process and Strong Markov property}\label{subsec:uniquenessSKS}
We now turn our attention to the characterization result concerning $SKS$-processes. Two by-product of this Subsection together with Subsection \ref{subsec:existenceSKS} are $\mathrm{Decomposition}(\theta,N)$-(i) and the strong Markov property.
\begin{proof}[Proof of the uniqueness part in Theorem \ref{mainresult}-(ii)]
In Step 1, with a given $SKS(\theta,N)$-process, we build a $KS(\theta,N)$-process killed when exiting $E_N$ such that the decomposition result holds before $\tau_{k_0\wedge N}$. In Step 2 we show it implies uniqueness and $\mathrm{Decomposition}(\theta,N)$-(i).

\vip

{\it Step 1.} Let $x\in E_{k_0\wedge N}$. We set $m = S_{\ig 1,N\id}(x)$, $r = R_{\ig 1,N\id}(x)>0$ and $u = \Phi (x)$. Let $(\Omega,\cF,(\cF_t)_{t\ge 0},( B_t)_{t\ge 0},(U_t)_{t\ge 0})$ a $SKS_u(\theta,N)$-process with lifetime $\tau_{k_0}(U)$. Even if it means to enlarge the probability space, one can build $(M_t)_{t\ge 0}$ a $2$-dimensional Brownian motion starting at $m$ with diffusion coefficient $N^{-1/2}$ and $(D_t)_{t\ge 0}$ a squared Bessel with dimension $d_{\theta,N}(N)$ starting at $r$ and take these processes such that $((M_t)_{t\ge 0}, (D_t)_{t\ge 0}, (U_t)_{t\ge 0}$ is a family of mutually independent random variables. In this step we show that $(X^N_t)_{t\ge 0} := (M_t + D_t^{1/2} U_{A_t})_{t\in [0,\tau_{k_0\wedge N}(X^N))}$ is a $KS_x(\theta,N)$-process killed when exiting $E_N$. 

\vip

It is classical that $((D_{\rho_t})^{1/2})_{t\ge 0}$ satisfies the following SDE on $\rr_+$:
\begin{gather*}
\mbox{ for all } t\in [0, \tau_0), \quad (D_{\rho_t})^{1/2} = r^{1/2} + \int_0^t (D_{\rho_s})^{1/2} \dd W^D_s + \frac{d_{\theta,N}(N)-1}{2} \int_0^t (D_{\rho_s})^{1/2}\dd s,\\
\tau_0 = \inf\{t\ge 0 : D_t = 0\},
\end{gather*}
where $(W^D_t)_{t\ge 0}$ is a $1$-dimensional Brownian motion. We now set for all $t\ge 0$, all $\e >0$, $Y_t = (D_{\rho_t})^{1/2} U^{N}_t$ and $Y^\e_t = (D_{\rho_t})^{1/2} U^{N,\e}_t$ with 
$$
U^{N,\e}_t = u+ \int_0^t\pi_{(U^N_s)^\perp}\pi_H\dd B_s+\int_0^t b_\e (U^N_s)\dd s  + \Big(\frac{\theta (N-1)}{2}-\frac{2N-3}{2}\Big)\int_0^t U^N_s\dd s.
$$
Using the Itô's formula and the fact that $(D_{\rho_t})_{t\ge 0}$ and $( B_t)_{t\ge 0}$ are independent, we get for all $t\in [0, A_{\tau_0}\wedge \tau_{k_0}(U))$,
\begin{align*}
Y^\e_t =& r^{1/2}u + \int_0^t Y_s \dd W^D_s + \frac{d_{\theta,N}(N)-1}{2}\int_0^t Y_s\dd s + \int_0^t (D_{\rho_s})^{1/2} \pi_{(U^N_s)^\perp} \pi_H (\dd B_s )\\
& + \int_0^t (D_{\rho_s})^{1/2}b_\e(U^N_s)\dd s +\Big(\frac{\theta (N-1)}{2} -\frac{2N-3}{2} \Big)\int_0^tY_s\dd s,
\end{align*}
Since $\theta (N-1) - (2N-3) = -d_{\theta,N}(N)+1$ and for all $\ell >0$, $v \in \sS$ we have $b_\e(\ell v) = \ell ^{-1} b_{\e\ell^{-2}}(v)$, we get that for all $t\in [0, A_{\tau_0}\wedge \tau_{k_0}(U))$
\begin{align*}
Y^\e_t =& r^{1/2}u + \int_0^t Y_s \dd W^D_s +  \int_0^t \|Y_s\|\pi_{Y_s^\perp} \pi_H ( \dd B_s )\\
& + \int_0^t \|Y_s\|^2 b_{\e \|Y_s\|^2}(Y_s)\dd s.
\end{align*}
Observing that $\tau_{k_0}(X) = \rho_{\tau_{k_0}(U)}$, we see that for all $t\in [0, \tau_{k_0\wedge N}(X))$, $A_t = \int_0^t \|Y_s\|^{-2}\dd s$, and setting $X^{N,\e}_t = \gamma (M_t) + Y^\e_{A_t}$ we get
\begin{align*}
X^{N,\e}_t =& r^{1/2}u + \gamma (M_t) +  \int_0^{A_t} Y_s \dd W^D_s +  \int_0^{A_t} \|Y_s\| \pi_{Y_s^\perp} \pi_H ( \dd B_s )\\
&+\int_0^t b_{\e R_{\ig 1,N\id (X^N_s)}}(X^N_s)\dd s,
\end{align*}
where we did the change of variable $s'=A_s$ and write $s$ instead of $s'$ in the last integral. Moreover we have for all $T>0$, all $\eta >0$, $\sup_{t\in [0,T\wedge \tau^\eta_{k_0\wedge N})} \|X^N_t - X^{N,\e}_t\| \to 0$ in probability as $\e \to 0$. Indeed, let $T>0$, $\eta >0$ and $\nu >0$. Fixing $\delta >0$ we have
\begin{align*}
    \PP \Big[\sup_{t\in [0,T\wedge \tau^\eta_{k_0\wedge N}]} \|X^N_t - X^{N,\e}_t\|\ge \nu \Big] &\le  \PP \Big[\sup_{t\in [0,A_{T\wedge \tau^\eta_{k_0\wedge N}}]} \|U^N_t - U^{N,\e}_t\|\ge \sqrt2\eta^{1/(2N)}\nu\Big]\\
    &\le \PP \Big[\sup_{t\in [0,M]} \|U^N_t - U^{N,\e}_t\|\ge \sqrt{2}\eta^{1/(2N)}\nu\Big] + \delta
\end{align*}
with $M$ chosen such that $\PP [ A_{T\wedge \tau^\eta_{k_0\wedge N}} \ge M] \le \delta$, which is possible since $ A_{T\wedge \tau^\eta_{k_0\wedge N}}$ is a.s. finite. Thanks to Definition \ref{def:KS/KSK}-(b), we get that for all $\eta >0$, $(X^N_{t\wedge \tau^\eta_{k_0\wedge N}})_{t\ge 0}$ is a $KS_x(\theta,N)$-process stopped when exiting $E^\eta_{k_0\wedge N}$. Since $\tau^\eta_{k_0\wedge N}\to \tau_{k_0\wedge N}$ as $\eta \to 0$, we get the result. 

\vip

In order to conclude the step, it is then sufficient  to prove that there exists a $2N$-dimensional Brownian motion starting at $0$ such that for all $t\in [0, \tau_{k_0\wedge N})$, $\tilde M_t :=  \gamma (M_t - m) +  \int_0^{A_t} Y_s \dd W^D_s +  \int_0^{A_t} \|Y_s\| \pi_{Y_s^\perp} \pi_H ( \dd\tilde B_s ) = B_t$. However it is clearly a martingale and for $(i,e), (j,f) \in \ig 1,N\id \times\{ 1,2\}$, we get for all $t\in [0,\tau_{k_0\wedge N})$,
\begin{align}\label{eq:fgh}
\langle \tilde M^{i,e}, \tilde M^{j,f} \rangle_t =&\delta_{e=f}\frac{t}{N} + \Big\langle   \int_0^{\cdot} Y^{i,e}_s \dd W^D_s, \int_0^\cdot Y^{j,f}_s \dd W^D_s \Big\rangle_{A_t}\notag \\
&+\Big\langle  \int_0^{\cdot} \|Y_s\|\dd (\pi_{Y_s^\perp} \pi_H (  B_s)^{i,e}, \int_0^{\cdot} \|Y_s\|\dd (\pi_{Y_s^\perp} \pi_H (  B_s)^{j,f} )\Big\rangle _{A_t}\notag\\
=& \delta_{e=f}\frac{t}{N}+\int_0^{A_t} Y^{i,e}_s Y^{j,f}_s \dd s + \int_0^{A_t} \|Y_s\|^2 \dd \langle \pi _{Y_\cdot^\perp}\pi_H( B)^{i,e},\pi _{Y_\cdot^\perp}\pi_H( B)^{j,f} \rangle_{s},
\end{align}
where $\tilde M^i_\cdot = (\tilde M^{i,1}_\cdot, \tilde M^{i,2}_\cdot)$, $Y^i_\cdot = (Y^{i,1}_\cdot,Y^{i,2}_\cdot)$ and $B^i_\cdot = (B^{i,1}_\cdot, B^{i,2}_\cdot)$. Moreover, we have for all $t\in [0,\tau_{k_0\wedge N})$, since $\sum_{k=1}^N Y_\cdot^k = 0$, $\langle \pi _{Y_\cdot^\perp}\pi_H( B)^{i,e},\pi _{Y_\cdot^\perp}\pi_H( B)^{j,f} \rangle_t$ equals
\begin{align*}
    \Big\langle  B^{i,e} &- S_{\ig 1,N\id}( B)^e - \sum_{k=1}^N \frac{Y^k_\cdot}{\| Y_\cdot\|^2} \cdot B^k Y^{i,e}_\cdot,   B^{j,f} - S_{\ig 1,N\id}( B)^f - \sum_{k=1}^N\frac{ Y^k_\cdot}{\| Y_\cdot\|^2} \cdot B^k Y^{j,f}_\cdot\Big\rangle_t\\
    &= \delta_{e=f}\Big(\delta_{i=j} -\frac1N \Big) t - 2\int_0^t \frac{Y^{i,e}_s Y^{j,f}_s}{\| Y_s\|^2}\dd s + \int_0^t \frac{Y^{i,e}_s Y^{j,f}_s}{\|Y_s\|^2} \dd s. 
\end{align*}
The last equality together with \eqref{eq:fgh} gives
$$
\langle \tilde M^{i,e}, \tilde M^{j,f} \rangle_t  =\delta_{(i,e) = (j,f)}t,
$$
which concludes the proof by the Lévy characterization theorem.

\vip

{\it Step 2.} We conclude the proofs of $\mathrm{Decomposition}(\theta,N)$-(i) and the uniqueness part of point (ii) of $\mathrm{MainTheorem}(\theta,N)$. Fix $u\in \sS$, and $(\Omega,\cF,(\cF_t)_{t\ge 0}, \PP, (B_t)_{t\ge 0}, (U^N_t)_{t\ge 0})$ a $SKS_u(\theta,N)$-process. Even if it means to enlarge the probability space, we can define $(M_t)_{t\ge 0}$ a $2$-dimensional Brownian motion with diffusion coefficient $N^{-1/2}$ starting at $0$, $(D_t)_{t\ge 0}$ a squared Bessel process with dimension $d_{\theta,N}(N)$ starting at $1$ such that the family $\{(M_t)_{t\ge 0}, (D_t)_{t\ge 0}, (U^N_t)_{t\ge 0} \}$ is mutually independent. According to Step 2 and Proposition \ref{prop:girsanoveta} we have that the law of $(X^N_t)_{t\ge 0}$ under $\PP$ is $\PP^{X,N,*}_u$ with $u$ viewed as an element of $(\rr^2)^N$. This proves also $\mathrm{Decomposition}(\theta,N)$-(i).

\vip

Since $(U^N_t)_{t\ge 0}$ is a measurable transformation of $(X^N_t)_{t\ge 0}$, we get that the law of $(U^N_t)_{t\ge 0}$ is characterized. Moreover, the strong Markov property also holds with the same reasoning as in the very first lines of the proof of Proposition \ref{prop:feller}
\end{proof}

\subsection{Ito's formula for $SKS(\theta,N)$-process}\label{subsec:ItoSKS}
In this section, we prove a Ito's like formula for $SKS$-processes even though it is not necessarily a semi martingale.
\begin{proof}[Proof of $\mathrm{Ito}(\theta,N)$-(ii)]
Let $\mu \in \cP(\sS)$ and 
$$
(\Omega,\cF,(\cF_t)_{t\ge 0}, \PP, (B_t)_{t\ge 0}, (U^N_t)_{t\ge 0})
$$
a $SKS_\mu(\theta,N)$-process. Let $\varphi \in C^2(\sS)$ such that
$$
\sup_{\alpha \in (0,1]} \sup_{u \in \sS} |\cL^U_\alpha \varphi (u)| <\infty. 
$$
We fix $(M_t)_{t\ge 0}$, $(D_t)_{t\ge 0}$ and $(X^N_t)_{t\ge 0}$ as in Subsection \ref{subsec:uniquenessSKS} (in particular, $U^N_t = \Phi (X^N_{\rho_t})$ for all $t\in [0,\tau_{k_0}(U))$). Recalling \eqref{eq:smoothedinfinitesimalgeneratorU} and Remark \ref{rk:coordonnéespolaires}, we get that $\varphi \circ \Phi \in C^2(E_N)$ and
$$
\sup_{\alpha \in (0,1]} \sup_{x\in (\rr^2)^N} |\cL^X_\alpha (\varphi\circ \Phi ) (x)|<\infty,
$$
so that we can apply Proposition \ref{prop:girsanoveta}-(ii)-(b) to get for all $\eta >0$, all $t\in [0,\sigma^\eta_*)$,
\begin{align}\label{eq:itoXtoU}
\varphi\circ \Phi(X^N_{t})&=\varphi (u) + \int_0^{t} \nabla (\varphi\circ \Phi) (X^N_s) \cdot\dd B_s +  \int_0^{t} \cL^X (\varphi\circ \Phi) (X^N_s)\dd s. 
\end{align}
Using again Remark \ref{rk:coordonnéespolaires}, we get for all $x\in E_N$ (which we can write $\psi (m,r,u)$ with $m\in \rr^2$, $r>0$ and $u\in \sS$),
\begin{align*}
    \cL^X(\varphi \circ \Phi)(x) &= \frac{1}{2}\Delta (\varphi \circ \Phi)(x) - \frac{\theta}{N} \sum_{i=1}^N \sum_{j\ne i} \frac{x^i-x^j}{\| x^i-x^j\|^2} \cdot \nabla_{x^i} (\varphi \circ \Phi)\\
    &=\frac{1}{2r}\Delta_\sS (\varphi \circ \Phi)(u) - \frac{\theta}{Nr} \sum_{i=1}^N \sum_{j\ne i} \frac{u^i-u^j}{\| u^i-u^j\|^2} \cdot \nabla_\sS (\varphi \circ \Phi)((0,\dots, u_i,\dots,0)\\
    &= \frac{\cL^U(\varphi \circ \Phi)(u)}{r}.
\end{align*}
We can now rewrite \eqref{eq:itoXtoU} as
\begin{align*}
\varphi(U^N_{t})&=\varphi (u) + \int_0^{\rho_{t}} \frac{\nabla (\varphi\circ \Phi) (U^N_s)}{\sqrt{D_s}} \cdot\dd B_s + \int_0^{\rho_{t}} \frac{\cL^U (\varphi\circ \Phi) (U^N_s)}{D_s}\dd s. 
\end{align*}
Thus, performing the change of variable $ s = \rho_{s'}$ and writing $s'=s$ in the last integral, and setting for all $t\ge 0$, $\tilde B_t = \int_0^{\rho_t}\frac{\dd B_s}{\sqrt{D_s}}$, we get that for all $\eta>0$, all $t\in [0,\sigma^\eta_*)$,
\begin{align*}
\varphi(U^N_{t})&=\varphi (u) + \int_0^{t} \nabla (\varphi\circ \Phi) (U^N_s)\cdot\dd \tilde B_s + \int_0^{t} \cL^U (\varphi\circ \Phi) (U^N_s)\dd s,
\end{align*}
which gives the result since $(\tilde B_t)_{t\ge 0}$ is a $2N$-dimensional Brownian motion and $\rho_{\sigma^\eta_*\wedge \tau^\eta_{k_0}} \to \tau_{k_0}(U^N)$ a.s. as $\eta \to 0$. 
\end{proof}

\subsection{Feller Property and Dirichlet form of $SKS(\theta,N)$-process}
\begin{proof}

{\it Step 1.}
We prove that $\mM^{U,N,\theta} = (\Omega,\cF,(\cF_t)_{t\ge 0} , (U^N_t)_{t\ge 0}, (\PP^{U,N,\theta}_u)_{u\in (E_{k_0}\cap \sS_N)_\triangle})$ is a Feller process i.e. we need to show that for any $\varphi \in C_0(E_{k_0}\cap \sS_N)$, we have
 \begin{gather}
     P^{U,N}_t \varphi : u\in (E_{k_0}\cap \sS_N)_\triangle \mapsto \E^{U,N,\theta}_u[\varphi (U^N_t) \indiq_{t< \tau_{k_0}(U^N)}] \in C_0(E_{k_0}\cap \sS_N),\label{eq:FellerU1}\\
     \sup_{u\in (E_{k_0}\cap \sS_N)_\triangle} | P^{U,N}_t \varphi (u) - \varphi (u)| \to 0 \quad \mbox{ as } t\to 0. \label{eq:FellerU2}
 \end{gather}
We first show \eqref{eq:FellerU1}. Fix $t>0$ and let $(u_n)_{n\ge 0}$ a sequence of elements of $(E_{k_0}\cap \sS_N)_\triangle$ converging to some $u_\infty \in (E_{k_0}\cap \sS_N)_\triangle$. 

\vip

Thanks to Subsection \ref{subsec:uniquenessSKS}, even if it means to enlarge the probability space, one can build a sequence $((X^{N,n}_t)_{t\ge 0})_{n\ge 0}$ and $(X^{N,\infty}_t)_{t\ge 0}$ with respectively law $\PP^{X,N,\theta}_{u_n}$ and $\PP^{X,N,\theta}_{u_\infty}$ under $\PP$ such that for all $t\ge 0$, $U^{N,n}_t = \Phi (X^{N,n}_{\rho_t})$, $U^{N,\infty}_t = \Phi (X^{N,\infty}_{\rho_t})$ (with the usual notations). In particular $(U^{N,n}_t)_{t\ge 0}$ and $(U^{N,\infty}_t)_{t\ge 0}$ respectively starting from $u_n$ and $u_\infty$. Thus
$$
P^{U,N}_t \varphi(u_n) = \E_{u_n} [\varphi \circ \Phi (X^N_{\rho_t})] \quad \mbox{ and } \quad P^{U,N}_t \varphi(u_\infty) = \E_{u_\infty} [\varphi \circ \Phi (X^N_{\rho_t})].
$$
We now write for $\eta >0$,
\begin{align*}
    \limsup_{n\to \infty}|P^{U,N}_t \varphi(u_n) - P^{U,N}_t\varphi (u_\infty)| \le& \limsup_{n\to \infty}|\E [\varphi \circ \Phi (X^{N,n}_{\rho_t \wedge \tau^\eta_N})] - \E [\varphi \circ \Phi (X^{N,\infty}_{\rho_t \wedge \tau^\eta_N})]|\\
    &+ \|\varphi\|_\infty \Big( \limsup_{n\to \infty} \PP^{X,N,*}_{u_n}[\tau^\eta_N\le \rho_t] + \PP^{X,N,*}_{u_\infty}[\tau^\eta_N\le \rho_t]\Big). 
\end{align*}
Fix $\e >0$ and $\eta >0$ such that $\PP^{X,N,*}_{u_\infty}[\tau^\eta_N\le \rho_t] \le \e$. This is licit since $\rho_t < \tau_N$ a.s. Observe that $\tau^\eta_N$ and $\rho_t$ are continuous $\PP_{u_\infty}^{X,N,\theta,*}$-a.e. according to Proposition \ref{prop:barycentreetbessel}. This together with the fact that $(X^{N,n}_t)_{t\ge 0} \to (X^{N}_t)_{t\ge 0} $ in law as $n\to \infty$ according to \cite[Theorem 4.2.5 page 167]{ek} and setting $A = \{\tau^\eta_N\le \rho_t\}$,  we get
\begin{align*}
\limsup_{n\to \infty}\PP^{X,N,*}_{u_n}[\tau^\eta_N\le \rho_t]   &= \limsup_{n\to \infty}\E[\indiq_{(X^n_t)_{t\ge 0}\in A}]\\
&\le \lim_{n\to \infty}\E[\indiq_{\cup_{k\ge n} \{ (X^k_t)_{t\ge 0}\in A}]\\
&=\E[\indiq_{\cap_{n\ge 0}\cup_{k\ge n} \{ (X^k_t)_{t\ge 0} \in A\}}]\\
&\le\PP^{X,N,*}_{u_\infty}[\tau^\eta_N\le \rho_t] \le \e,
\end{align*}
where we used the monotone convergence theorem between the second and the third line. Finally we conclude by continuity $\PP_{u_\infty}^{X,N,*}$-a.e. for every fixed $\eta >0$ of the application
$$
(x_t)_{t\ge 0} \in C(\rr_+,(E_{k_0})_\triangle) \mapsto \varphi\circ \Phi (x_{\rho_t\wedge \tau^\eta(x)}).
$$

\vip

Let us now treat \eqref{eq:FellerU2}. Recall Remark \ref{rk:siKE2}-(ii) for the definition of $\cC_\cK^\sS$ for all $\cK \setminus \{\ig 1,N\id\}$. Fix $\cK \in \gK_N \setminus \{\ig 1,N\id \}$ and  $\varphi \in C_\cK^\sS$. We have $\sup_{\alpha \in (0,1]} \sup_{u\in \sS}|\cL^U_\alpha \varphi (u)| <\infty$. Indeed, this comes from the fact that according to similar computations as in Proposition \ref{prop:barycentreetbessel}, for all $K\in \cK$, 
$$
\sup_{\alpha \in (0,1]} \sup_{u\in \sS\cap E_\cK}(|\cL^U_\alpha S_K (u)| + |\cL^U_\alpha R_K(u)|) <\infty.
$$
We can thus apply Subsection \ref{subsec:ItoSKS} and apply the expectation in order to get for all $u\in \sS$,
\begin{align*}
    |P^{U,N}_t\varphi (u) - \varphi (u)| \le t \sup_{\alpha \in (0,1]}\sup_{u\in (E_{k_0}\cap\sS)_\triangle}|\cL^U_\alpha \varphi (u)| \to 0 \quad \mbox{ as } t\to 0. 
\end{align*}
We conclude for every $\varphi \in C_0(E_{k_0}\cap \sS)$ since $\sum_{\cK \in \gK_N \setminus \{\ig 1,N\id \}} \cC_\cK^\sS$ is dense in $C_0(E_{k_0}\cap \sS)$ for the uniform convergence topology. 

\vip

{\it Step 2.} We conclude the proof of $\mathrm{MainTheorem}(\theta,N)$-(ii) by proving that for all $\mu \in \cP( (E_{k_0}\cap \sS_N)_\triangle)$, the unique $SKS_\mu(\theta,N)$-process we have built is $\beta^N_\theta$-reversible and associated to the Dirichlet space $(\cE^N_U,\cF^N_U)$ defined in Theorem \ref{mainresult}. However it is a straight application of Proposition \ref{prop:dirichlettué} and \cite[Proposition 10]{ft}.
\end{proof}

We end this section by giving the following. 
\begin{prop}\label{prop:existuniqSKSeternel}
    Suppose that $k_0\ge N$. There exists an eternal $SKS(\theta,N)$-process as in Definition \ref{def:eternalSKS}. Moreover its law is unique.
\end{prop}
\begin{proof}
According to $\mathrm{MainTheorem}(\theta,N)$-(ii), Lemma \ref{lem:radon} and the fact that $k_0\ge N$, we have that $\beta^N_\theta$ normalized with a suitable constant $Z$ is an invariant probability measure of $SKS(\theta,N)$-process. We can thus build for all $n\in \nn$ a $SKS(\theta,N)$-process defined on $[-n,\infty)$ starting from $Z^{-1}\beta^N_\theta$. Thus the proof of the existence is classical using the Kolmogorov extension theorem. The uniqueness easily follows from Definition \ref{def:eternalSKS} and the uniqueness part of $\mathrm{MainTheorem}(\theta,N)$-(ii) proved in Section \ref{sec:decompo}.
\end{proof}

\section{Construction of a candidate}\label{sec:constructionKS}
We are finally in position to build a candidate for a $KS(\theta,N)$-process. In this Section one can additionally suppose that $N\le k_1$. Indeed, if $k_0\le N$, then one can not hope to extend the process after $\tau_N$. Thus, in this section we have $\cB^N_\theta \subset \{\ig 1,N\id\}$.

\vip

Finally, it is interesting to note that the cases $\cB^N_\theta=\emptyset$ or $\cB^N_\theta=\{\ig 1,N\id\}$ are also different in nature. In fact the need for excursion theory is mainly used for the case $\cB^N_\theta = \{\ig 1,N\id\}$ (which only concerns the case $k_1=N$, see Lemma \ref{lem:dimension}) for which the well posedness of the drift involved similar difficulties as for the well definition of a Bessel process with dimension lying in $(0,1)$. The case $\cB^N_\theta = \emptyset$ could be handled with a more classical approach since the construction in this case is more similar to the one of a Bessel process with dimension strictly greater than $1$. However, for the sake of consistency throughout the paper we use the same technique for both cases.  

\vip

As just explained, we will need the excursion theory for Bessel processes. We introduce a few notations and summarize what we will use. 

\vip

Let $(Z_t)_{t\ge 0}$ be a Bessel process with dimension $\delta \in (0,2)$. According to \cite{bertoin(1990)} for example, there exists $\{L_t^a, a\in \rr_+, t\ge 0\}$ a jointly continuous version of its local times, i.e. for every $\varphi$ bounded and measurable, we have a.s. 
\begin{align*}
\int_0^t \varphi (Z_s)\dd s = \int_{\rr_+} \varphi (a) L_t^a a^{\delta-1}\dd a.
\end{align*}
We set for all $t\ge 0$, $\tau^0_t = \inf\{ s\ge 0 : L^0_s > t\}$ the right-continuous inverse of the local time at $0$. Moreover we set the excursion process $(\gge^t)_{t\ge 0}$ where for all $t\ge 0$,
\begin{align*}
    \gge^t = (Z_{\tau^0_{t-} + s})_{s\in [0, \tau^0_t-\tau^0_{t-}]}.
\end{align*}
According to \cite[Theorem 3.1]{bertoin(1990)}, $(\gge^t)_{t\ge 0}$ is a Poisson point process with intensity measure $\dd t \otimes \gn$, where $\gn$ is a measure on 
\begin{align*}
\cE = \Big\{ f \in C(\rr_+, \rr_+) : f(0) = 0 \mbox{ and there exists } s_0>0 \mbox{ such that for all } s\ge s_0, f(s) = 0\Big\},
\end{align*}
and is called the Bessel excursion measure. For all $f\in \cE$, we call $\ell (f) = \sup \{ s> 0 : f(s) \ne 0\} \in (0,\infty)$ the length of $f$. Moreover, according to \cite[Theorem 1]{Pitman-Yor1996} and the paragraph that follows, one can write for all $A\in \cB(\cE)$,
\begin{align}\label{eq:mesureexcursionformule}
\gn(A) = c_\delta\int_0^\infty \frac{\dd \ell}{\ell^{2-\delta/2}} \int_{\cE_1} \indiq_{(\sqrt{\ell} \gge^{s/\ell})_{s\ge 0} \in A} \gn_1(\dd e),
\end{align}
where $c_\delta$ is a positive constant and $\cE_1 = \{ f\in \cE : \ell (f) = 1\}$ and $\gn_1$ is called the normalized Bessel excursion measure. Moreover $\gn_1$ is exactly the law of an excursion of the Bessel process with length $1$, i.e. the process $(\gge_t)_{t\in [0,1]}$ satisfying the SDE
\begin{align*}
  \mbox{ for all } t\in [0,1], \quad  \gge_t = B_t + \int_0^t\frac{3-\delta}{2\gge_s}\dd s - \int_0^t\frac{\gge_s}{1-s}\dd s,
\end{align*}
see \cite[page 37]{G.F.LawlerNotes}.

\vip

We set $\gM = \sum_{t \in J} \delta_{(t,\gge^t)}$ where $J = \{ t\ge 0 : \tau^0_t >\tau^0_{t-}\}$. It holds that $\tau^0_t = \int_0^t \int_{\cE} \ell (\gge) \gM (\dd s,\dd \gge)$. For any $\phi : \rr \to \rr_+$, any $t\ge 0$, we have
\begin{align}
\int_0^{\tau_t} \phi (Z_u)\dd u = \sum_{s\in J\cap [0,t]} \int_{\tau_{s^-}}^{\tau_s} \phi (Z_u)\dd s = \int_0^t \int_\cE \Big[ \int_0^{\ell (\gge)} \phi (\gge^u)\dd u \Big] \gM(\dd s,\dd \gge).
\end{align}

For all $x = \gamma (m) + r^{1/2}u \in (\rr^2)^N$ with $m\in \rr^2$, $r\ge 0$ and $u\in \sS$, we define $(M_t)_{t\ge 0}$ a $2$-dimensional Brownian motion with diffusive coefficient $N^{-1/2}$ starting at $m$, $(D_t)_{t\ge 0}$ a squared Bessel process with dimension $d_{\theta,N}(N)$ starting at $r$. We set for all $t\ge 0$, $Z_t = D_t^{1/2}$ which is a Bessel process with dimension $d_{\theta,N}(N)>0$ (because $k_0>N$) for which we re-use the above notations. Finally, let $(U^{N}_t)_{t\ge 0}$ be a $SKS_u(\theta,N)$-process and $(U^{N,i}_t)_{t\in \rr})_{i\ge 0}$ a countable collection of i.i.d eternal $SKS(\theta,N)$-processes. Let $\cZ = \{ t\ge 0 : D_t = 0\}$. Since $[\tau_N,\infty)\setminus\cZ$ is an open set of $\rr$, we can write it $\cup_{i\ge 0}(a_i,b_i)$ with the interval being disjoint. Moreover, we can always choose an indexation such that $i\in \nn \mapsto (a_i,b_i)$ is measurable (using for example the inverse of the local time of $(R_{\ig 1,N\id}(X_t))_{t\ge 0}$ to describe the sequence $(a_i,b_i)_{i\ge 0}$). We set for all $t\ge 0$, 
$$
X^N_t = \gamma (M_t) + (D_t)^{1/2}\Big( U^{N}_{A_t} \indiq_{t\in [0,\tau^0_0)} + \sum_{i\ge 0} U^{N,i}_{A^i_t}\indiq_{t\in (a_i,b_i)}\Big),
$$
where for all $t\in [0,\tau_N)$, $A_t = \int_0^tD_s^{-1}\dd s$ and for all $i\ge 0$, all $t\in (a_i,b_i)$, $A_t^i = \int_{(a_i+b_i)/2}^t D_u^{-1}\dd u$. Let us first observe that by construction, definition of a $SKS(\theta,N)$-process and properties of Bessel processes (see \cite[Proposition 1.5 page 442]{ry}),
\begin{align}\label{eq:pasdansE2}
    \int_0^\infty \indiq_{X^N_s \notin E_2\cup \{ \triangle \}} \dd s = 0 \quad \mbox{ a.s.}
\end{align}

\vip

Let us observe that in the case where $r=0$ and $d_{\theta,N}(N)\ge 2$ then there is only one excursion of infinite length, making it a particular case but easier to handle. In this case we set $A_t = \int_1^t D_s^{-1}\dd s$ for all $t>0$ and the process $(X^N_t)_{t\ge 0}$ analogously.

\vip

We begin with a result on the well definition of the candidate for being the drift of $(X^N_t)_{t\ge 0}$.
\begin{lemma}\label{lem:drifttendzero}
There exists a continuous process $(\gb_t)_{t\ge 0}$ such that for all $t>0$, all collisional approximation $\chi^\e_N$ (which can be equal to $1$ if $\cB^N_\theta = \emptyset$, i.e. $d_{\theta,N}(N)>1$, see Lemma \ref{lem:dimension}),
$$
\sup_{t\in [0,T]}\Big\|\int_0^t b(X^N_s) \chi^\e_N(X^N_s,Z_s) \dd s - \gb_t \Big\| \to 0
$$
in probability as $\e \to 0$, where $(Z_u)_{u\ge 0}$ is a $\rr^d$-valued process independent from all other random variables defined previously. In particular,
$$
\sup_{t\in [0,T]}\Big\|\int_0^t b(X^N_s)\indiq_{R_{\ig 1,N\id}(X^N_s)\ge \e} \dd s - \gb_t \Big\| \to 0
$$
\end{lemma}
\begin{proof}
{\it Step 1.} In this step we build the process $(\gb_t)_{t\ge 0}$.

\vip

{\it Step 1.1} We first build the drift between the first and the last instant of collision before time $t$ implying the $N$ particles using the excursion theory. We set $\cU = C(\rr,\sS_N)$ and $\Lambda$ the law of the eternal $SKS(\theta,N)$-process (see Proposition \ref{prop:existuniqSKSeternel}). Observe that since the measure $\beta^N_\theta$ is invariant under the transformation $u\in \sS_N\mapsto -u\in \sS_N$ and by definition and uniqueness in law of an eternal $SKS(\theta,N)$-process, we get that 
\begin{align}\label{eq:symduU}
    \mbox{ for all } \varphi : \cU\to \rr \mbox{ measurable bounded, } \quad \int_\cU \varphi (\gu) \Lambda (\dd \gu) =\int_\cU \varphi (-\gu) \Lambda (\dd \gu).
\end{align}
We set for all $\gge \in \cE$, all $\gu\in \cU$,
\begin{gather*}
    F(\gge,\gu) = \int_0^{\ell (\gge )} \gge^{-1}_s b(\gu_{A_s})\dd s \quad \mbox{ and } \quad \tilde F(\gge,\gu) = \int_0^{\ell (\gge )} \gge^{-1}_s \|b(\gu_{A_s})\|\dd s.
\end{gather*}
Observe that according to Lemma \ref{lem:estimedun}, we have for all $\e >0$, 
$$
\int_\cE\int_\cU\indiq_{\e \le \tilde F(\gge,\gu)\le 1} \|F(\gge,\gu)\|\gn(\dd \gge)\Lambda (\dd \gu) <\infty.
$$
Moreover, since for all $\gu\in \cU$, $\indiq_{\e \le \tilde F(\gge,-\gu)\le 1} F(\gge,-\gu) = - \indiq_{\e \le \tilde F(\gge,\gu)\le 1} F(\gge,\gu)$, the property \eqref{eq:symduU} gives us that
\begin{align*}
    \int_\cE\int_\cU \indiq_{\e \le \tilde  F(\gge,\gu)\le 1} F(\gge,\gu)\gn_1(\dd \gge ) \Lambda (\gu) = 0.
\end{align*}
Using now Lemma \ref{lem:estimedun} and the fact that $\| F\| \le \tilde F$, one can conclude that there exists a càdlàg process $(H_t)_{t\ge 0}$ (representing the small jumps) such that
$$
\sup_{t\in [0,T]} \Big\| \int_0^t\indiq_{\e\le \tilde F(\gge,\gu)\le 1}F(\gge,\gu)\gM(\dd s,\dd \gge,\dd \gu) - H_t \Big\| \to 0 \quad \mbox{ in } L^2,
$$
see \cite[Chapter II.3.]{iw}. Moreover, setting 
$$
E_T = \Big\{ s\in J\cap [0,T] : \tilde F(\gge ,\gu)\ge 1\Big\}, 
$$
Lemma \ref{lem:estimedun} tells us that $|E_T|$ is a.s. finite. We then set 
$$ \gb^e_t=H_t+\int_0^\infty\!\!\int_E\!\int_U\indiq_{\{\tau^0_s\le t\}}\indiq_{\{\tilde F(e,u)\ge1\}}F(e,u)\,\cM(\dd s,\dd e,\dd u),$$
where $H_t$ is the limit in $L^2$ of $\int_0^\infty\!\int_E\!\int_U\indiq_{\{\tau^0_s\le t\}}\indiq_{\{\e\le\tilde F(e,u)\le1\}}F(e,u)\,\cM(\dd s,\dd e,\dd u)$, i.e. the compensated sum over the excursions completed before time $t$. 

\vip

{\it Step 1.2} We conclude Step 1 here. Set for all $t\ge 0$,
$$
\gb_t = \int_0^{\tau^0_0\wedge t} b(X^N_s)\dd s + \gb^e_t + \int_{g_t}^tb(X^N_s)\dd s,
$$
where $g_t := \sup \{ \tau^0_s : \tau^0_s\le t\}$ with the convention that $g_t = 0$ if $\{ \tau^0_s : \tau^0_s\le t\} = \emptyset$. We also recall Step 1.1 for the definition of $(\gb^e_t)_{t\ge 0}$. To this aim it is sufficient to show that  
$$
\int_0^{\tau^0_0\wedge t} \| b(X^N_s)\|\dd s + \int_{g_t}^t \| b(X^N_s)\| \dd s < \infty \mbox{ a.s.}
$$
First, one can suppose $r>0$ because the other scenario can be proved in an easier way since the first integral would be equal to $0$.

\vip

We refer to Step 1 of Lemma \ref{lem:estimeduU} and its proof to get
\begin{gather*}
\int_0^{\tau^0_0}\| b(X^N_s)\|\dd s = \sqrt{\tau^0_0} \int_0^\infty Z^r_v\|b(U^N_v)\|\dd v,
\end{gather*}
with the same notation as in  Lemma \ref{lem:estimexcursionbessel} and $\mathrm{Decomposition}(\theta,N)$. Thus, taking the expectation, using the Fubini's theorem, independence between $(Z^r_v)_{v\ge 0}$ and $(U^N_v)_{v\ge 0}$ and Lemma \ref{lem:estimexcursionbessel}, we get the existence of positive and finite constants $C$ and $c$ such that
\begin{align*}
\E \Big[(\tau^0_0)^{-1/2}\int_0^{\tau^0_0} \| b(X^N_s)\|\dd s\Big]&= \int_0^\infty \E[ Z_v]\E[ \| b (U^{N}_v)\|] \dd v\\
&\le C \E \Big[ \int_0^\infty e^{-c\sqrt v} \|b(U^{N}_v)\|\dd v \Big]\\
&\le C\sum_{n=0}^\infty e^{-c\sqrt{n}} \E \Big[ \int_n^{n+1} \| b(U^{N}_v)\|\dd v\Big].
\end{align*}
According to the Markov property of $(U^N_t)_{t\ge 0}$ (see Subsection \ref{subsec:uniquenessSKS}), we get that
\begin{align*}
    \E \Big[(\tau^0_0)^{-1/2}\int_0^{\tau^0_0} \| b(X^N_s)\|\dd s\Big] \le C \sum_{n=0}^\infty e^{-c\sqrt n} \sup_{u\in \sS^{N}} \E_u \Big[ \int_0^1 \| b(U^{N}_v)\|\dd v\Big] 
\end{align*}
which is finite according to Lemma \ref{lem:estimeduU}. This concludes for the first term since in addition, $\sqrt{\tau_0^0}<\infty$ a.s. We now treat the second term. Observe that a.s. there exists $i_0\ge 0$ such that we have $g_t = a_{i_0}$. We can thus write by definition of $(X^N_t)_{t\ge0}$,
\begin{align*}
    \int_{g_t}^t \| b(X^N_s)\|\dd s = \int_{g_t}^{t} D_s^{-1/2} \| b(U^{N,i_0}_{A^{i_0}_s})\|\dd s.
\end{align*}
Moreover, by construction, $(U^{N,i_0}_s)_{s\ge 0}$ knowing $(D_s)_{s\ge 0}$ is distributed as an eternal $SKS(\theta,N)$-process, so by definition of such a process,
\begin{align*}
\E \Big[\int_{g_t}^t \| b(X^N_s)\|\dd s |\sigma((D_t)_{t\ge 0})\Big] &= \int_{g_t}^t D_s^{-1/2}\dd s \int_{\sS^{N}} \|b\| \dd \beta^N_\theta <\infty \mbox{ a.s.}
\end{align*}
according to Lemma \ref{lem:radon}-(iii). This thus implies the result.

\vip

{\it Step 2.} We prove the result in the easier case where $\cB^N_\theta = \emptyset$. One can suppose that $d_{\theta,N}(N)<2$ because otherwise, Step 1 is sufficient to conclude since the squared Bessel process never hits $0$ after time $0$. According to Lemma \ref{lem:estimedun} and \cite[Chapter II.3.]{iw}, we get that 
$$
\E\Big[\int_0^{g_t} \int\indiq_{\tilde F (\gge,\gu)\le 1} \tilde F(\gge,\gu) \gM(\dd s,\dd \gge, \dd \gu) \Big] <\infty.
$$
Moreover, Lemma \ref{lem:estimedun} still gives that $|E_T|<\infty$ a.s. Thus, using the estimates in Step 1.2, we get that for all $t\ge 0$,
$$
\E\Big[ \int_0^T \|b(X^N_s)\|\dd s\Big] <\infty.
$$
Thus, the dominated convergence theorem gives the result. The continuity of $(\gb_t)_{t\ge 0}$ follows directly from the integrability. 

\vip

{\it Step 3.} We prove the result in the more delicate case where $\cB^N_\theta = \{\ig 1,N\id \}$. Here, the drift is not absolutely integrable so we need to be more careful, using the compensations in the drift thanks to the angle being stationary. 

\vip

{\it Step 3.1.} Fix $T>0$, $\chi^\e_N$ a collisional approximations and $(Z_t)_{t\ge 0}$ a process independent from $(X^N_t)_{t\ge 0}$. We set for all $\e >0$,
\begin{align*}
    F_\e(s,\gge,\gu) = \int_0^{\ell (\gge)} \gge ^{-1}_ub(\gu_{A_u}) \chi^\e_{\ig 1,N\id}(\gge^{-1}_u\gu_{A_u},  Z_{s+u})\dd u. 
\end{align*}
This step is devoted to show that
\begin{align}\label{eq:smalljumps}
  \sup_{t\in [0,T]} \Big\| \int_0^t \int F_\e(\tau^0_{s-},\gge,\gu) \gM(\dd s, \dd \gge,\dd \gu) - \gb^e_{t} \Big\| \to 0 \quad \mbox{ in probability as } \e \to 0.  
\end{align}
For all $s\ge 0$, we have by definition of the family $\chi^\e_N$ (see Definition \ref{def:collisionalapprox}) and \eqref{eq:symduU},
\begin{align*}
   \int_\cE\int_\cU \int_{\tilde F(\gge,\gu)\le 1} F_\e (s,\gge,\gu)\gn (\dd \gge)\Lambda (\dd \gu) = - \int_\cE\int_\cU \int_{\tilde F(\gge,\gu)\le 1} F_\e (s,\gge,\gu)\gn (\dd \gge)\Lambda (\dd \gu) = 0.
\end{align*}
Thus, according to \cite[Chapter II.3.]{iw}, one can define for all $t\ge 0$, 
\begin{align*}
   N^\e_t := \int_0^t\int \indiq_{\tilde F(\gge,\gu)\le 1}( F_\e(\tau^0_{s-},\gge,\gu) - F(\gge,\gu)) \gM(\dd s, \dd \gge,\dd \gu). 
\end{align*}
In particular, $(N^\e_t)_{t\ge 0}$ is a $L^2$-martingale. Using the Doob inequality, we get
\begin{align*}
   \E\Big[\sup_{t\in [0,T]} \|N^\e_t\|^2\Big] \le 4 \E[\|N^\e_T\|^2] = 4\E\Big[\int_0^t \int_\cE\int_\cU \indiq_{\tilde F(\gge,\gu)\le 1}\|F_\e(\tau^0_{s-},\gge,\gu)-F(\gge,\gu)\|^2\dd s\gn_1(\dd \gge)\Lambda(\gu)\Big].
\end{align*}
Since $\sup_{s\ge 0}\|F_\e(s,\cdot)\| + \|F(\cdot)\| \le 2 \tilde F(\cdot)$, Lemma \ref{lem:estimedun} and the dominated convergence theorem gives
\begin{align}\label{eq:cvL2petitssauts}
  \sup_{t\in [0,T]} \Big\| \int_0^t\int \indiq_{\tilde F(\gge,\gu)\le 1} (F_\e(\tau^0_{s-},\gge,\gu)-F(\gge,\gu)) \gM(\dd s, \dd \gge,\dd \gu)  \Big\| \to 0 \quad \mbox{ in } L^2 \mbox{ as } \e \to 0.  
\end{align}

\vip

We now recall from Step 1 that $|E_T|<\infty$ a.s. which immediately gives
\begin{align*}
    \sup_{t\in [0,T]} \Big\| \int_0^t \indiq_{\tilde F(\gge,\gu)\ge 1} (F_\e(\tau^0_{s-},\gge,\gu) - F(\gge,\gu))\gM(\dd s,\dd \gge, \dd \gu) \Big\| \to 0 \quad \mbox{ a.s.},
\end{align*}
which together with \eqref{eq:cvL2petitssauts} concludes the step.

\vip

{\it Step 3.2.} We conclude the proof here. We have
$$
 \sup_{t\in [0,T]}\Big\|\int_0^t b(X^N_s)\chi^\e_{\ig 1,N\id}(X^N_s,Z_s) \dd s -\gb_t \Big\| \le I^1_\e + I^2_\e + I^3_\e,
$$
where 
\begin{gather*}
I^1_\e = \sup_{t\in [0,T]} \Big\| \int_0^{\tau^0_0\wedge t} b(X^N_s) (\chi^\e_{\ig 1,N\id}(X^N_s,Z_s)-1) \dd s \Big\|\\
I^2_\e = \sup_{t\in [0, L^0_T]} \Big\|\int_0^t\int_\cE\int_\cU(F_\e(\tau^0_{s-},\gge,\gu) - F(\gge,\gu) )\gM(\dd s,\dd \gge,\dd \gu) \Big\|\\
I^3_\e = \sup_{t\in [0,T]} \Big\|  \int_{g_t}^{t} b(X^N_s)  (\chi^\e_{\ig 1,N\id}(X^N_s,Z_s)-1)\dd s\Big\|,
\end{gather*}
According to Step 1.2, the fact that for all $\e>0$, all $K\in \cB^N_\theta$, $\|\chi^\e_K\|_\infty \le 1$ and the dominated convergence theorem, we have that $I^1_\e + I^3_\e \to 0$ a.s. as $\to 0$. Finally, fix $\nu,\eta >0$. We have
\begin{align*}
    \PP [ I^2_\e \ge \nu] \le \PP[L^0_T\ge L] + \PP\Big[  \sup_{t\in [0, L]} \Big\|\int_0^t\int_\cE\int_\cU(F_\e(\tau^0_{s-},\gge,\gu) - F(\gge,\gu) )\gM(\dd s,\dd \gge,\dd \gu) \Big\| \ge \nu\Big].
\end{align*}
Since $L^0_T < \infty$ a.s., there exists $L>0$ big enough such that $\PP[L^0_T\ge L]\le \eta$. Step 3.1 allows us to conclude that for $\e$ small enough, $\PP[I^2_\e \ge \nu] \le 2\eta$, which gives the convergence in probability. Moreover, since for all $\e >0$, the process 
$$
\Big( \int_0^t b(X^N_s)\chi^\e_{\ig 1,N\id}(X^N_s,Z_s) \dd s \Big)_{t\ge 0}
$$
is a.s. continuous, the uniform convergence implies that $(\gb_t)_{t\ge 0}$ is also a.s. continuous.
\end{proof}

\begin{lemma}\label{lem:EDSpv}
There exists a $2N$-dimensional Brownian motion $(B_t)_{t\ge 0}$ such that, 
$$
\mbox{ for all } t\ge 0, \quad X^N_t = x +B_t +\gb_t,
$$
where $(\gb_t)_{t\ge 0}$ is the process defined in Lemma \ref{lem:drifttendzero}. 
\end{lemma} 
\begin{proof}
We set for all $t\ge 0$, $ N_t = X^N_t - x - \gb_t$ and we show that $(N_t)_{t\ge0}$ is a $2N$-dimensional Brownian motion.

\vip

{\it Step 1.} We first show that for all $x\in E_{k_0}\cap E_N$, all $s\ge 0$, all $\e >0$, $(N_{(u +s)\wedge \tau^{s}_0}-N_{s})_{t\ge 0}$ is a continuous martingale with bracket equal to $u\wedge (\tau^s_0-s)$, where
$$
\tau^s_0 = \inf \{ u \ge s : D_u = 0\}. 
$$

\vip

According to the strong Markov property for Brownian motion, squared Bessel processes and $SKS(\theta,N)$-processes, we have that conditionally to $\cF_s$, $(S_{\ig 1,N\id}(X^N_{(u+s)\wedge \tau^s_0}))_{t\ge 0}$ is a $2$-dimensional Brownian motion stopped at $\tau^s_0$ starting at $S_{\ig 1,N\id}(X^N_s)$, $(R_{\ig 1,N\id}(X^N_{(u+s)\wedge \tau^s_0}))_{t\ge 0}$ is a squared Bessel process stopped at $\tau^s_0$ with dimension $d_{\theta,N}(N)$ starting at $R_{\ig 1,N\id}(X^N_s)$ and $(\Phi (X^N_{\rho^{s}_{(u+A^{s}_{s})}\wedge \tau^s_0}))_{u\ge 0}$ is a $SKS(\theta,N)$-process stopped at $A^{s}_{\tau^s_0}$ starting from $\Phi (X^N_s)$, where $\rho^{s}$ and $A^{s}$ are defined as usually, such that $A^s_s = 0$. 

\vip

However, according to $\mathrm{Decomposition}(\theta,N)$-(i) proved in Subsection \ref{subsec:uniquenessSKS}, $(X^N_{(u+s)\wedge \tau^s_0} )_{u\ge 0}$ is a $KS(\theta,N)$-process starting from $X^N_s$ and stopped when exiting $E_N$. Thus using Proposition \ref{prop:girsanoveta}-(i) and the fact that 
$$
g : (x_t)_{t\ge 0} \in C(\rr_+,((\rr^2)^N)_\triangle) \mapsto \Big(x_t - x_0 - \limsup_{\e\to 0} \int_0^t b_\e(x_s) \dd s\Big)_{t\ge 0},
$$
is measurable, we conclude the Step.

\vip

{\it Step 2.} We show that $(N_t)_{t\ge 0}$ is a continuous local martingale. Observe first that continuity is direct according to the definition of $(X^N_t)_{t\ge 0}$ and by Lemma \ref{lem:drifttendzero}. Then, fix $n\ge 0$ and $\tau_n = \inf \{t\ge 0 : N_t \ge n\}$ and $0\le s<t$, we want to prove that $\E[N_{t\wedge \tau_n} - N_{s\wedge \tau_n}|\cF_s] = 0$. 

\vip

According to Step 1, we get that $\E[ N_{\tau^s_0\wedge t\wedge\tau_n} - N_{s\wedge \tau_n} |\cF_s] = 0$, so that it remains to show that $\E[\indiq_{t\ge \tau^s_0}(N_{t\wedge \tau_n} - N_{\tau^s_0\wedge\tau_n})|\cF_s] = 0$. We can write $N_t = \tilde N_t + \gamma(M_t)$ for all $t\ge 0$. By classical properties of Brownian motion, it is clear that $\E[\indiq_{t\ge \tau^s_0}(\gamma(M_{t\wedge \tau_n}) - \gamma (M_{\tau^s_0\wedge\tau_n}))|\cF_s] = 0$. Moreover, by \eqref{eq:symduU}, by oddness of $(\tilde N_{t\wedge \tau_n} )_{t\ge 0}$ with respect to the eternal $SKS$-processes involved in its construction and by definition of $\tau^s_0$, we get $\E[\indiq_{t\ge \tau^s_0}(\tilde N_{t\wedge \tau_n} - \tilde N_{\tau^s_0\wedge \tau_n})|\cF_s] = -\E[\indiq_{t\ge \tau^s_0}(\tilde N_{t\wedge \tau_n} - \tilde N_{\tau^s_0\wedge \tau_n})|\cF_s]$, which concludes. Of course the same technique gives the result for all starting point $x \in E_{k_0}$ such that $R_{\ig 1,N\id}(x) = 0$. 

\vip

{\it Step 3.} We conclude. According to Steps 1 and 2, $(N_t)_{t\ge 0}$ is a continuous local martingale with for all fixed $i,j\in \ig 1,N\id$, $e_1,e_2\in \{1,2\}$, all $s\ge 0$, $\langle N^{i,e_1}, N^{j,e_2}\rangle_{(u+s)\wedge \tau^s_0} - \langle N^{i,e_1}, N^{j,e_2}\rangle_{s} = \delta_{(i,e_1) = (j,e_2)}u\wedge \tau^s_0$. Thus, we have for all $t\ge 0$,
\begin{align*}
    \langle N^{i,e_1},N^{j,e_2}\rangle_t &= \int_0^t \indiq_{s\notin \cZ}\dd \langle N^{i,e_1},N^{j,e_2}\rangle_s + \int_0^t \indiq_{s\in \cZ}\dd \langle N^{i,e_1},N^{j,e_2}\rangle_s\\
    &=\delta_{(i,e_1)=(j,e_2)} t + S^{i,j,e_1,e_2}_t,
\end{align*}
where $(S^{i,j,e_1,e_2}_t)_{t\ge 0} := ( \int_0^t \indiq_{s\in \cZ}\dd \langle N^{i,e_1},N^{j,e_2}\rangle_s)_{t\ge 0}$ is a continuous function. Thanks to the Lévy characterization, it suffices to show that for all $t\ge 0$, $S_t^{i,j,e_1,e_2} = 0$ a.s. Moreover, by Kunita-Watanabe, it is sufficient to show that for all $t\ge 0$, $S_t^{i,i,e_1,e_1} = 0$ a.s., which we now do. First, $(S^{i,i,e_1,e_1}_t)_{t\ge 0}$ is nondecreasing, only increasing on $\cZ$ which is a.s. a negligable set since $\int_0^\infty \indiq_{X^N_s \notin E_2\cup \{\triangle\}} \dd s = 0$ a.s. 

\vip

Since $(S^{i,i,e_1,e_1}_t)_{t\ge 0}$ only increases on $\cZ$, we get that $(S^{i,i,e_1,e_1}_{\tau^0_t})_{t\ge 0}$ is a continuous subordinator. Indeed, by construction of $(X^N_t)_{t\ge 0}$, it is clear that the increments are independent with stationary law. Thus, there exists $\kappa \ge0$ such that for all $t\ge 0$, $S^{i,i,e_1,e_1}_{\tau^0_t} = \kappa t$. Next, we see that by construction of $(X^N_t)_{t\ge 0}$ again, scaling properties of Brownian motion and squared Bessel processes, we have for all $\lambda >0$, $(\lambda^{-1/2}(X^N_{\lambda t + \tau^0_0} - X^N_{\tau^0_0}))_{t\ge 0}$ equals $(X^N_{ t + \tau^0_0} - X^N_{\tau^0_0})_{t\ge 0}$ in law. We thus deduce that for all $\lambda >0$, $(\lambda^{-1/2}(N_{\lambda t + \tau^0_0} - N_{\tau^0_0}))_{t\ge 0} $ equals $(N_{t + \tau^0_0} - N_{\tau^0_0})_{t\ge 0}$ in law, implying that $(\lambda^{-1} S^{i,i,e_1,e_1}_{\lambda t + \tau^0_0})_{t\ge 0}$ equals $(S^{i,i,e_1,e_1}_{ t + \tau^0_0})_{t\ge 0}$ in law (we recall that $S^{i,i,e_1,e_1}_{\tau^0_0}=0$). Moreover it is classical that for all $\lambda >0$, $(\lambda^{-2/(2-d_{\theta,N}(N))}\tau^0_{\lambda t})_{t\ge 0}$ has the same law as $(\tau^0_t)_{t\ge 0}$, see \cite{PitmanYor92} or \cite{MolchanovOstrovskii69}. All in all, we deduce that $(S^{i,i,e_1,e_1}_t)_{t\ge 0}$ is $2/(2-d_{\theta,N}(N))$-stable with $ 2/(2-d_{\theta,N}(N))>1$, so that $\kappa =0$, concluding the proof.

\end{proof}

\section{Ergodicity of the $SKS(\theta,N)$-process}\label{sec:ergodicity}
This section is devoted to show the following
\begin{prop}\label{prop:ergodicity}
    There exists $C,c>0$ positive and finite constants such that 
    $$
   \mbox{ for all }t\ge 0, \quad  \sup_{u\in \sS_N} \|\cL_u(U_t)- Z^{-1}\beta^N_\theta\|_{TV} \le Ce^{-ct},
    $$
    where for all $u\in \sS_N$, all $t\ge 0$, $\cL_u(U_t)$ is the law of a $SKS_u(\theta,N)$-process at time $t$ starting from $u$ and $Z$ is the normalizing constant such that $Z^{-1}\beta^N_\theta$ is a probability measure. 
\end{prop}
\begin{proof}
    {\it Step 1.} We prove that there exists $\e_N,M_N>0$ such that
    $$
   c:= \!\!\inf_{x\notin E^{\e_N}_2 : R_{\ig 1,N\id}(x) = 1} \PP^{X,N,\theta,*}_x\Big[X_1 \in E^1_2,\quad \mbox{ and for all } t\in [0,1], \quad M_N^{-1}\le R_{\ig 1,N\id}(X_t)\le M_N\Big] >0,
    $$
    recall Proposition \ref{prop:girsanoveta} for the definition of $\PP^{X,N,\theta,*}_x$. This step is the more technical of this proof, it explains why there is a uniform over some initial configuration positive probability to end up with a configuration where all particles are well separated from the other. Of course the issue is to deal with clusters of very close particles at time $0$ and understand how they can split in time $1$. We present informally our strategy : we will prove the result through an induction by showing that one can split all the clusters of particles with uniformly positive probability as follows: divide $[0,1]$ in time interval of size $1/N$. Roughly, during each of these intervals, 
    \begin{itemize}
        \item there is a positive probability that during the first $1/(2N)$ time, one cluster develops enough internal dispersion to break into several sub clusters,
        \item during the remaining $1/(2N)$ time, there is a positive probability for one of these cluster to go far away from all the other clusters. After $N$ steps, every particle will be isolated,
        \item  of course one needs in addition to keep track of the total dispersion to keep some control on the distance between particle when we will renormalize by the total dispersion to go back on the sphere.
    \end{itemize}
    Even though the idea is conceptually simple, it is technical because of all the (elementary but tedious) geometrical questions it raises.

    \vip
    
One can find a sequence $(a_n)_{n\in \ig 0,2N\id}$ such that $a_0 = 1$ and for all $n\in \ig 1,2N\id$, $ a_{n-1}+2 \le a_n/2$. The sequence $(a_n)_{n\in \ig 0,2N\id}$ will represent the size of the displacements of the clusters we decide to isolate at each step. It is chosen with a fast enough growth in order to compensate the maximal displacement performed at each induction step while preserving a uniform separation between already isolated clusters (see Figure \ref{fig:separations}). Define by induction the sequence $(C_n)_{n\in \ig 0,2N\id}$ by setting $C_0 = 1$ and for all $n\in \ig 0,2N-1\id$, 
$$
C_{n+1} =\max (C_n, 2 + 2C_n +4N(a_{2N}+1)^2 ). 
$$
The sequence $(C_n)_{n\in \ig 0,2N\id}$ represents the control we will have at each step of the total dispersion of the particles. Let us set for all $k\in \ig 0,2N\id$, all $\cK_1,\cK_2 \in \cP(\ig 1,N\id)$ such that $\cK_1\cup \cK_2\in \gK_N$,
\begin{align*}
G^k_{\cK_1,\cK_2} = \Big\{ x\in (\rr^2)^N :  \cK_1\cup\cK_2 \in \cP^{2^{-k}\e_N}_x, \min_{K\ne L, K\in \cK_1,L\in \cK_1\cup \cK_2}\| S_K(x) -S_L(x)\| \ge a_{2N-k} \\
\mbox{ and }\quad C_k^{-1}\le R_{\ig 1,N\id}(x) \le C_k\Big\},
\end{align*}
where we recall Lemma \ref{lem:releq} for the definition of $\cP^{2^{-k}\e_N}_x$ and $\e_N>0$ is to be chosen later. This set is informally the set of configurations such that the particles spatially regroup at the scale $2^{-k}\e_N$ according to the partition $\cK_1\cup\cK_2$ with the additional informations that the total dispersion stays under control and particles labelled by the subsets described by $\cK_1$ are very far away from all the particles which are not in their cluster. The asymmetry in the definition of $G^k_{\cK_1,\cK_2}$	is intentional: clusters in $\cK_1$ are already separated, whereas clusters in $\cK_2$ are allowed to remain close and will be progressively separated during the induction. 

\vip

We will prove by induction that
\begin{align*}
\cP'(n) : ''\inf_{x\notin E^{\e_N}_2 : R_{\ig 1,N\id}(x) = 1} \PP^{X,N,\theta,*}_x\Big[ X_{n/N} \in \bigcup_{\cK_1,\cK_2 : \cK_1\cup\cK_2\in \gK_N, |\cK_1| \ge n}G^{2n}_{\cK_1,\cK_2} \\
\mbox{ and for all } t\in [0,n/N] \quad C_{2n}^{-1}\le R_{\ig 1,N\id}(X_t)\le C_{2n} \Big]>0,''
\end{align*}
is true for all $n\in \ig 0,N\id$. At each induction step, either all clusters are already separated, or one cluster is extracted from $\cK_2$ and added to $\cK_1$ or one cluster of $\cK_1$ is split and one of them is added to $\cK_1$. Observe that the lower bound on the total dispersion is almost straightforward to get with $1/2$ instead of $C_{2n}^{-1}$ in the following proof. We thus do not detail this point for the sake of conciseness. Proving such induction will indeed conclude since $\cP'(N)$ implies the result with $M_N = C_{2N}$. 

\vip

First, $\cP'(0)$ is obviously true since $x\in G^0_{\emptyset, \cK^{\e_N}(x)}$. Now let $n\in \ig 0,N-1\id$ such that $\cP'(n)$ holds. We prove $\cP'(n+1)$. 

\vip

Let $x\notin E^{\e_N}_2$ such that $R_{\ig 1,N\id}(x) = 1$. By $\cP'(n)$, there exists $c_N^1>0$ (independent from $x$), $\cK_1,\cK_2$ such that $\cK_1\cup \cK_2 \in \gK_N$, $|\cK_1|\ge n$ and 
$$
\PP^{X,N,\theta,*}_x\Big[ X_{n/N} \in G^{2n}_{\cK_1,\cK_2} \quad \mbox{ and for all } t\in [0,n/N] \quad C_{2n}^{-1}\le R_{\ig 1,N\id}(X_t)\le C_{2n}\Big]\ge c_N^1. 
$$
Thus, using the strong Markov property (see Proposition \ref{prop:feller}), we get
\begin{align*}
\PP^{X,N,\theta,*}_x\Big[ X_{(n+1)/N}& \in \bigcup_{\cK_1,\cK_2 : \cK_1\cup\cK_2\in \gK_N, |\cK_1| \ge n+1}G^{2(n+1)}_{\cK_1,\cK_2}\\
&\mbox{ and for all } t\in [0,n/N] \quad C_{2(n+1)}^{-1}\le R_{\ig 1,N\id}(X_t)\le C_{2(n+1)}\Big] \\
    &\ge c_N^1\inf_{y\in G^{2n}_{\cK_1,\cK_2}}\PP^{X,N,\theta,*}_y \Big[ X_{1/N} \in \bigcup_{\cK_1',\cK_2' : \cK_1'\cup\cK_2'\in \gK_N, |\cK_1'| \ge n+1}G^{2(n+1)}_{\cK_1',\cK_2'} \\
    &\mbox{ and for all } t\in [0,n/N]\quad C_{2(n+1)}^{-1}\le R_{\ig 1,N\id}(X_t)\le C_{2(n+1)}\Big].
\end{align*}
Suppose first that $|\cK_1| = N$ (i.e. $\cK_1 = \ig 1,N\id$). This is the case where we have already succeeded to separate all the particles, we then only need to show that there is a positive probability that particles do not move significantly during the remaining time. We show that 
$$
\inf_{y\in G^{2n}_{\ig 1,N\id,\emptyset}}\!\!\PP^{X,N,\theta,*}_y \Big[ X_{1/N} \in G^{2(n+1)}_{\ig 1,N\id,\emptyset} \quad \!\!\mbox{ and for all } t\in [0,1/N] \quad \!\!C_{2(n+1)}^{-1}\le R_{\ig 1,N\id}(X_t)\le C_{2(n+1)}\Big]\!\!>0.
$$
However we have for all $y \in G^{2n}_{\ig 1,N\id,\emptyset}$,
\begin{align*}
    \PP^{X,N,\theta,*}_y \Big[ X_{1/N} \in G^{2(n+1)}_{\ig 1,N\id,\emptyset}\Big] \ge \PP^{X,N,\theta,*}_y \Big[ \max_{i\in \ig 1,N\id} \sup_{t\in [0,1/N]}\|X^i_t-y^i\|\le  \alpha_N\Big],
\end{align*}
with $\alpha_N$ chosen small enough and only depending on $N$. Applying Proposition \ref{prop:girsanoveta}-(ii)-(c) with $\cK = \ig 1,N\id$, $\eta = 1$ and $T = 1/N$, we get that
\begin{align*}
    \PP^{X,N,\theta,*}_y \Big[ X_{1/N} \in G^{2(n+1)}_{\ig 1,N\id,\emptyset}\Big] &\ge \exp(-C_{N,\theta,\eta}1/N) \QQ_y^{T,\eta} \Big[ \max_{i\in \ig 1,N\id} \sup_{t\in [0,1/N]}\|X^i_t-y^i\|\le  \alpha_N\Big]\\
    &=\exp(-C_{N,\theta,\eta}1/N) \tilde\PP \Big[ \max_{i\in \ig 1,N\id} \sup_{t\in [0,1/N]}\|B^i_t\|\le  \alpha_N\Big],
\end{align*}
where $((B^i_t)_{t\ge 0})_{i\in \ig 1,N\id}$ is a $2N$-dimensional Brownian motion starting from $0$ built on some probability space $(\tilde \Omega,\tilde \cF,\tilde\PP)$, which concludes this first case. 

\vip

We now suppose that $|\cK_1|<N$. We conclude this case in several steps. 

\vip

{\it Step 1.1} Suppose that there exists $K_0 \in \cK_1\cup\cK_2$ such that $|K_0|\ge 2$. In this case there is still a cluster of particle that one can split in several clusters. This step shows why there is a positive probability to see this happening. We show that $\inf_{y\in G^{2n}_{\cK_1,\cK_2}} P_y >0$ where 
\begin{align*}
    P_y =& \PP_y^{X,N,\theta,*}\Big[ \bigcup_{\cL\in \gK^{K_0}}\{ X_{1/(2N)} \in  G^{2n+1}_{\cK_1\setminus \{K_0\}, (\cK_2\setminus \{K_0\})\cup\cL} \\
    &\quad \mbox{ and for all } t\in [0,1/(2N)] \quad C_{2n+1}^{-1}\le R_{\ig 1,N\id}(X_t)\le C_{2n+1}\} \Big],
\end{align*}
with $\gK^{K_0}$ the set of partitions $\cL$ of $K_0$ such that $\cL \ne \{K_0\}$. Fix $y\in G^{2n}_{\cK_1,\cK_2}$. We have 
\begin{align*}
    P_y \ge \PP_y^{X,N,\theta,*}\Big[\bigcap_{K\in \cK_1\cup\cK_2} B_K\cap C_K\Big],
\end{align*}
where for all $K\in \cK_1\cup \cK_2$,
\begin{gather*}
 B_K = \Big\{  \sup_{t\in [0,1/(2N)]}\|S_K(X_t) - S_K(y)\| \le \alpha_N\Big\},
\end{gather*}
for all $K\in \cK_1\cup\cK_2\setminus \{K_0\}$
\begin{align*}
C_K =  \Big\{ \sup_{t\in [0,1/(2N)]} R_K(X_t)\le (2^{-2n}\e_N)^{1/|K|}+\alpha_N , R_K(X_{1/(2N)})\le (2^{-2n-1}\e_N)^{1/|K|} \Big\},
\end{align*}
and
\begin{align*}
 C_{K_0} = \Big\{ \sup_{t\in [0,1/(2N)]} R_{K_0}(X_t) \le (2^{-2n}\e_N)^{1/|K_0|} + \alpha_N,\quad  R_{K_0}(X_{1/2N}) \ge (2^{-2n}\e_N)^{1/|K_0|} +\alpha_N/2\Big\},
\end{align*}
where $\alpha_N>0$ is chosen small enough and only depending on $N$. Indeed, first, on the event $\bigcap_{K\in \cK_1\cup\cK_2} B_K\cap C_K$, if $\alpha_N$ and $\e_N$ are sufficiently small, there exists $\cL\in \gK^{K_0}$ such that $(\cK_1\cup\cK_2\setminus \{K_0\} )\cup\cL \in \cP_{X_{1/(2N)}}^{2^{-2n-1}\e_N}$. This is true since for all $K\in \cK_1\cup\cK_2$, all $i\in K$, the triangular inequality gives
\begin{align*}
    \|X^i_{1/(2N)} -y^i\| &\le \|X^i_{1/(2N)}-S_K(X_{1/(2N)})\| + \|S_K(X_{1/(2N)})- S_K(y)\| + \|S_K(y) - y^i\|\\
    &\le \sqrt{R_K(X_{1/(2N)})} + \alpha_N + \sqrt{R_K(y)}  \\
    &\le 2 ( (2^{-2n}\e_N)^{1/|K|} + \alpha_N)^{1/2} + \alpha_N,
\end{align*}
so that there exists a constant $C'_N$ only depending on $N$ such that for all $K\ne L\in \cK_1\cup\cK_2$, all $i\in K$, all $j\in L$,
$$
R_{K \cup \{j\}}(X_{1/(2N)}) \ge R_{K \cup \{j\}}(y) - C'_N \Big( 8 ( (2^{-2n}\e_N)^{1/|K|} + \alpha_N)^{1/2} + \alpha_N^2\Big) \ge (2^{-2n-1 }\e_N)^{1/(|K|+1)},
$$
providing $\alpha_N$ and $\e_N$ are small enough (we used that for all $x,y \in \rr$, $(x+y)^2\le 2(x^2+y^2)$). Moreover, according to  Lemma \ref{lem:dispersionformula}, we have for all $t\in [0,1/(2N)]$,
\begin{align*}
R_{\ig 1,N\id} (X_t) &= \sum_{K\in \cK_1\cup \cK_2} R_K(X_t) + \frac{1}{2N} \sum_{K\ne L \in \cK_1\cup\cK_2} |K| |L| \| S_K(X_t) - S_L(X_t)\|^2\\
&\le |\cK_1\cup\cK_2| ((2^{-2n}\e_N)^{1/N}+\alpha_N) + \frac{1}{N} \sum_{K\ne L \in \cK_1\cup\cK_2} |K||L| \|S_K(x) - S_L(x)\|^2 + 4N\alpha_N^2 \\
&\le N((2^{-2n}\e_N)^{1/N}+\alpha_N) + 2R_{\ig 1,N\id}(x) + 4N\alpha_N^2 \\
&\le 2 + 2C_k \le C_{k+1},
\end{align*}
by definition of $C_{k+1}$, if $\alpha_N$ and $\e_N$ are chosen only depending on $N$ and small enough. The lower bound on $R_{\ig 1,N\id}(X_t)$ follows easily after one iteration of the induction argument thanks to the cluster in $\cK_1$ which we do not detail since it is the same kind of argument as for the upperbound. Finally, it is clear that for all $K'\in \cK_1$, all $K\in \cK_1\cup \cK_2$ with $K'\ne K$, we have
$$
\inf_{t\in[0,1/(2N)]}\|S_{K'}(X_t) - S_K(X_t)\| \ge \|S_{K'}(z) - S_K(z)\| - 2\alpha_N\ge a_{2N-n} - 2\alpha_N \ge a_{2N-n-1}
$$
if $\alpha_N$ is chosen small enough since $a_{2N-n-1}+2\le a_{2N-n}/2$.

\vip

We now aim at using Remark \ref{rk:girsanovopen}. As a preparation, observe that on the event $\cap_{K\in \cK_1} B_K\cap C_K$, using the triangular inequality and Lemma \ref{campingcar-unleashedversion}, we get for all $K\ne K'\in \cK_1\cup\cK_2$, all $i \in K$, $j\in K'$, all $t\in [0,1/(2N)]$,
\begin{align}\label{eq:prepagirsanov}
\|X^i_t - X^j_t\| \ge& \|z^i-z^j\| - \|z^i-S_K(z)\| - \| S_K(z) - S_K(X_t)\|\\
&- \sqrt{R_K(X_t)} - \|z^j - S_{K'}(z)\| - \|S_{K'}(z)-S_{K'}(X_t)\| - \sqrt{R_{K'}(X_t)}\notag\\
\ge& \|z^i-z^j\|- 4((2^{-2n}\e_N)^{1/(|K|\vee |K'|)}+\alpha_N)^{1/2} -2\alpha_N\notag\\
\ge& \Big(c(2^{-2n}\e_N)^{1/(|K|\vee |K'|+1)} \Big)^{1/2}- 4((2^{-2n}\e_N)^{1/(|K|\vee |K'|)}+\alpha_N)^{1/2}  -2\alpha_N\ge c^2_N\notag,
\end{align}
with $c^2_N$ positive and independent from $K,K'$ and $k$ if $\e_N$ and $\alpha_N$ are chosen small enough and depending only on $N$. We can thus apply Remark \ref{rk:girsanovopen} with $\cK = \cK_1\cup\cK_2$, $\eta = c^2_N/2$ and $T = 1/(2N)$ and we get by independence,
\begin{align*}
    P_y &\ge \exp(-C_{N,\theta,\eta} /(2N)) \tilde\QQ_y^{\cK,T, \eta}[\cap_{K\in \cK_1\cup\cK_2} B_K\cap  C_K]\\
    &= \exp(-C_{N,\theta,\eta} /(2N)) \prod_{K\in \cK_1\cup\cK_2}\tilde\QQ_y^{\cK,T, \eta}[B_K] \tilde\QQ_y^{\cK,T, \eta}[C_K]\\
    &= \exp(-C_{N,\theta,\eta} /(2N)) \prod_{K\in \cK_1\cup\cK_2}\tilde \PP[\tilde B_K] \tilde\PP[\tilde C_K],
\end{align*}
where
\begin{gather*}
       \mbox{ for all } K\in \cK_1, \quad  \tilde B_K = \Big\{  \sup_{t\in [0,1/(2N)]}\|S^K_t - S_K(y)\| \le \alpha_N\Big\},\\
        \mbox{ for all } K\in \cK_1\setminus \{K_0\}, \quad  \tilde C_K= \Big\{ \sup_{t\in [0,1/(2N)]} R^K_t\le (2^{-2n}\e_N)^{1/|K|} +\alpha_N\Big\},\\
    \tilde C_{K_0} = \Big\{ \sup_{t\in [0,1/(2N)]} R^{K_0}_t \le (2^{-2n} \e_N)^{1/|K_0|}+\alpha_N,\quad  R^{K_0}_{1/(2N)} \ge (2^{-2n}\e_N)^{1/|K_0|} +\alpha_N/2\Big\},
\end{gather*}
where for all $K\in \cK_1\cup \cK_2$, $(S^K_t)_{t\ge 0}$ is a $2$-dimensional Brownian motion with diffusivity coefficient $|K|^{-1/2}$ starting from $S_K(y)$, $(R^K_t)_{t\ge 0}$ is a squared Bessel process with dimension $d_{\theta,N}(|K|)$ starting from $R_K(y)$ built on some probability space $(\tilde\Omega,\tilde \cF,\tilde \PP)$. Moreover, the family 
$$
\{ (S^K_t)_{t\ge 0}, (R^K_t)_{t\ge 0}, K\in \cK_1\cup\cK_2\}
$$
is mutually independent. We easily conclude this step thanks to well known properties of Brownian motion and Bessel processes, see \cite{ry}.

\vip

{\it Step 1.2}
 Let $\cK_3,\cK_4$ such that $\cK_3\cup\cK_4 \in \gK_N$ and $\cK_4\ne \emptyset$, $k\in \{2n,2n+1\}$ and $t_N >0$. Fix $z \in G^{k}_{\cK_3,\cK_4}$ and $L\in \cK_4$ such that $S_L(z)$ is an extremal point of the convex hull of $\cC := \{ S_K(z), K\in \cK_4 \}$. We conclude step 1 by showing that there exists $c^3_N>0$ independent of $z$, $\cK_3$, $\cK_4$, $L$, $k$ and $t_N$ such that
 \begin{align}\label{eq:finalduo}
\PP^{X,N,\theta,*}_z\Big[ X_{t_N} \in G^{2(n+1)}_{\cK_1 \cup\{L\},\cK_3 \setminus\{L\}} \quad \mbox{ and for all } t\in [0,t_N] \quad C_{2(n+1)}^{-1}\le R_{\ig 1,N\id}(X_t)\le C_{2(n+1)}\Big]\ge c^3_N. 
 \end{align}
This will indeed conclude since there are two cases : first for all $K\in \cK_1\cup\cK_2$, $|K|=1$, with $|\cK_1|<N$, which implies $\cK_2$ is non empty, and thus \eqref{eq:finalduo} with $t_N = 1/N$, $k=2n$ and $(\cK_3,\cK_4) = (\cK_1,\cK_2)$ concludes. In the second case, Step 1.1 combined with \eqref{eq:finalduo} with $t_N = 1/(2N)$, $k=2n+1$ and $(\cK_3,\cK_4) = (\cK_1\setminus \{K_0\}, (\cK_2\setminus \{K_0\})\cup \cL)$ with $\cL$ given by Step 1.1 concludes the proof. It thus only remains to show \eqref{eq:finalduo}. 

\vip

Fix $z\in G^k_{\cK_3,\cK_4}$ and $L$ as above. We denote $v \in \rr^2$ some vector such that $v-S_L(z)$ is unitary and for all $c\in \cC$, $\langle v-S_L(z), c-S_L(z)\rangle\le 0$ (which exists since $\cC$ is convex). We will first show that,
\begin{align*}
    \PP^{X,N,\theta,*}_z[X_{t_N} \in G^{k+1}_{\cK_1\cup\{L\},\cK_2\setminus \{L\}}] \ge \PP^{X,N,\theta,*}_z[ \cap_{K\in \cK_1\cup\cK_2} B_K\cap C_K  ],
\end{align*}
where for all $K\in \cK_1\cup\cK_2\setminus \{L\}$,
\begin{gather*}
     B_K = \Big\{ \sup_{t\in [0,t_N]} \|S_K(X_t)-S_K(z)\|\le \alpha_N \Big\} \\
     C_K = \Big\{ \sup_{t\in [0,t_N]}R_K(X_t)\le (2^{-k}\e_N)^{1/|K|}+\alpha_N\quad \mbox{ and } \quad  R_K(X_{t_N})\le (2^{-2(n+1)}\e_N)^{1/|K|} \Big\},
\end{gather*}
and
\begin{gather*}
     B_L = \Big\{ \sup_{t\in [0,t_N]} \|S_L(X_t)-S_L(z) - t_N^{-1}t(a_{2N-k-1}+1)(v-S_L(z))\|\le \alpha_N \Big\} \\
     C_L = \Big\{ \sup_{t\in [0,t_N]} R_L(X_t)\le (2^{-k}\e_N)^{1/|L|}+\alpha_N \quad \mbox{ and } \quad R_L(X_{t_N})\le (2^{-2(n+1)}\e_N)^{1/|L|} \Big\},
\end{gather*}
where $\alpha_N>0$ is chosen small enough. 

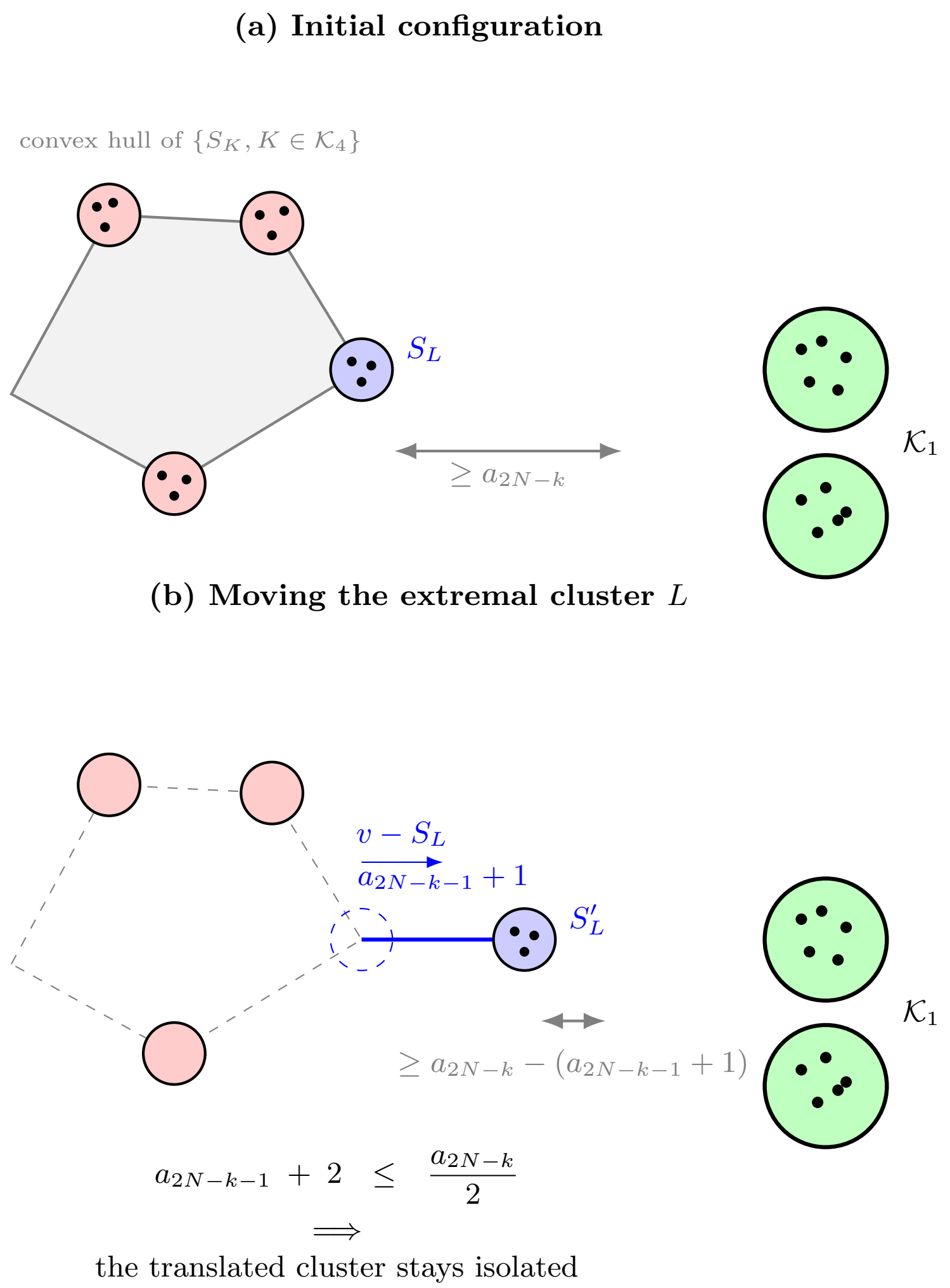
\begin{figure}[ht]
\centering
\begin{tikzpicture}[
    scale=0.85,
    >=Latex,
    every node/.style={font=\small}
]


\node at (5,4.5)
{\bf (a) Initial configuration};

\fill[gray!10]
(0,0)--(1.2,2.2)--(3.2,2.1)--(4.3,.3)--(2,-1.1)--cycle;

\draw[gray,thick]
(0,0)--(1.2,2.2)--(3.2,2.1)--(4.3,.3)--(2,-1.1)--cycle;

\draw[thick,fill=red!20]
(1.2,2.2) circle (.38);
\fill (1.05,2.3) circle (1.7pt);
\fill (1.25,2.35) circle (1.7pt);
\fill (1.15,2.05) circle (1.7pt);

\draw[thick,fill=red!20]
(3.2,2.1) circle (.38);
\fill (3.05,2.2) circle (1.7pt);
\fill (3.35,2.25) circle (1.7pt);
\fill (3.2,1.95) circle (1.7pt);

\draw[thick,fill=red!20]
(2,-1.1) circle (.38);
\fill (1.85,-1) circle (1.7pt);
\fill (2.15,-1.05) circle (1.7pt);
\fill (2,-1.25) circle (1.7pt);

\draw[thick,fill=blue!20]
(4.3,.3) circle (.38);
\fill (4.18,.4) circle (1.7pt);
\fill (4.42,.35) circle (1.7pt);
\fill (4.3,.15) circle (1.7pt);

\node[right,blue] at (4.7,.55) {$S_L$};

\draw[very thick,fill=green!25]
(10,.3) circle (.75);
\fill (9.7,.55) circle (2pt);
\fill (9.95,.65) circle (2pt);
\fill (10.25,.45) circle (2pt);
\fill (9.8,.15) circle (2pt);
\fill (10.15,.05) circle (2pt);

\draw[very thick,fill=green!25]
(10,-1.5) circle (.75);
\fill (9.7,-1.3) circle (2pt);
\fill (10,-1.15) circle (2pt);
\fill (10.25,-1.45) circle (2pt);
\fill (9.9,-1.7) circle (2pt);
\fill (10.15,-1.55) circle (2pt);

\node[right] at (10.8,-.6) {$\mathcal K_1$};

\draw[<->,gray,thick]
(4.7,-.7)--(7.5,-.7);

\node[gray] at (6.1,-1)
{$\geq a_{2N-k}$};

\node[gray] at (2.2,3.1)
{\scriptsize convex hull of $\{S_K,K\in\mathcal K_4\}$};


\begin{scope}[yshift=-7cm]

\node at (5,4.5)
{\bf (b) Moving the extremal cluster $L$};

\draw[dashed,gray]
(0,0)--(1.2,2.2)--(3.2,2.1)--(4.3,.3)--(2,-1.1)--cycle;

\draw[thick,fill=red!20]
(1.2,2.2) circle (.38);
\draw[thick,fill=red!20]
(3.2,2.1) circle (.38);
\draw[thick,fill=red!20]
(2,-1.1) circle (.38);

\draw[dashed,blue]
(4.3,.3) circle (.38);

\draw[very thick,blue,->]
(4.3,.3)--(6.3,.3);

\node[blue] at (5.3,1.05)
{$a_{2N-k-1}+1$};

\draw[thick,fill=blue!20]
(6.3,.3) circle (.38);
\fill (6.18,.4) circle (1.7pt);
\fill (6.42,.35) circle (1.7pt);
\fill (6.3,.15) circle (1.7pt);

\node[right,blue] at (6.7,.55) {$S_L'$};

\draw[very thick,fill=green!25]
(10,.3) circle (.75);
\fill (9.7,.55) circle (2pt);
\fill (9.95,.65) circle (2pt);
\fill (10.25,.45) circle (2pt);
\fill (9.8,.15) circle (2pt);
\fill (10.15,.05) circle (2pt);

\draw[very thick,fill=green!25]
(10,-1.5) circle (.75);
\fill (9.7,-1.3) circle (2pt);
\fill (10,-1.15) circle (2pt);
\fill (10.25,-1.45) circle (2pt);
\fill (9.9,-1.7) circle (2pt);
\fill (10.15,-1.55) circle (2pt);

\node[right] at (10.8,-.6) {$\mathcal K_1$};

\draw[<->,gray,thick]
(6.5,-.7)--(7.3,-.7);

\node[gray,align=center,anchor=north]
at (6.9,-.9)
{$\displaystyle
\geq a_{2N-k}-(a_{2N-k-1}+1)
$};

\node[align=center,text width=7cm]
at (4,-3.1)
{
$\displaystyle
a_{2N-k-1}+2\leq \frac{a_{2N-k}}2
$

\vspace{2mm}

$\Longrightarrow$

the translated cluster stays isolated
};

\draw[->,blue]
(4.3,1.25)--(5.3,1.25);

\node[blue] at (4.8,1.6) {$v-S_L$};

\end{scope}

\end{tikzpicture}
\caption{
Geometric illustration of the displacement argument.
The cluster $L$ is an extremal point of the convex hull of the remaining
clusters. It is translated in the outward normal direction by a distance
$a_{2N-k-1}+1$. Since the initial distance from the already separated
clusters $\mathcal K_1$ is at least $a_{2N-k}$, the inequality
$a_{2N-k-1}+2\le a_{2N-k}/2$ guarantees that the displacement does not
destroy the separation.
}
\label{fig:separations}
\end{figure}

Indeed, on the event $\cap_{K\in \cK_1\cup\cK_2} B_K\cap C_K$, we have for all $K\in \cK_2\setminus \{L\}$,,
\begin{align*}
    \|S_L(X_{t_N})-S_K(X_{t_N})\| \ge &\| S_L(z)+ (a_{2N-k-1}+1)(v-S_L(z))- S_K(z)\|  - 2\alpha_N\\
    =& \Big((a_{2N-k-1}+1)^2 + 2(a_{2N-k-1}+1)(v-S_L(z),S_L(z) -S_K(z)) \\
    &+ \|S_L(z) - S_K(z)\|^2  \Big)^{1/2} - 2\alpha_N\\
    \ge& a_{2N-k-1}+1 - 2\alpha_N,
\end{align*}
because $(v - S_L(z),S_L(z)-S_K(z))\ge 0$ by definition of $v$. In particular, if $\alpha_N$ is chosen small enough, this means that 
$$
\|S_L(X_{t_N}) - S_K(X_{t_N})\| \ge a_{2N-k-1}.
$$
Moreover, it is clear that for all $K'\in \cK_1$, all $K\in \cK_2\setminus \{L\}$, we have
$$
\inf_{t\in[0,t_N]}\|S_{K'}(X_t) - S_K(X_t)\| \ge \|S_{K'}(z) - S_K(z)\| - 2\alpha_N\ge a_{2N-k} - 2\alpha_N \ge a_{2N-k-1}
$$
if $\alpha_N$ is chosen small enough since $a_{2N-k-1}\le a_{2N-k}/2$. In addition, for all $K\in \cK_1$, we have
\begin{align*}
    \|S_L(X_{t_N}) - S_K(X_{t_N})\|&\ge \|S_L(z) - S_K(z)\| - \alpha_N - a_{2N-k-1} -1 \\
    &\ge a_{2N-k} - a_{2N-k-1} -1 - \alpha_N \\
    &\ge a_{2N-k-1} + 3-\alpha_N \\
    &\ge a_{2N-k-1} ,
\end{align*}
if $\alpha_N$ is chosen small enough, since $a_{2N-k-1} + 2\le a_{2N-k}/2$. Another ingredient is that as in Step 1.1, we also have $ \cK_1\cup\cK_2 \in \cP^{2^{-k-1}\e_N}_{X_{t_N}}$ providing $\alpha_N$ and $\e_N$ are small enough. Finally, according to  Lemma \ref{lem:dispersionformula} we have for all $t\in [0,t_N]$,
\begin{align*}
R_{\ig 1,N\id} (X_t) =& \sum_{K\in \cK_1\cup \cK_2} R_K(X_t) + \frac{1}{2N} \sum_{K\ne K' \in \cK_1\cup\cK_2} |K| |K'| \| S_K(X_t) - S_{K'}(X_t)\|^2\\
\le& |\cK_1\cup\cK_2| ((2^{-k}\e_N)^{1/N }+\alpha_N) + \frac{1}{N} \sum_{K\ne K' \in \cK_1\cup\cK_2} |K||K'| \|S_K(x) - S_{K'}(x)\|^2 \\
&+ 4\alpha_N^2 N + 4N(a_{2N-k-1}+1)^2\\
&\le N((2^{-k}\e_N)^{1/N}+\alpha_N) + 2R_{\ig 1,N\id}(x) + 4\alpha_N^2 N +4N(a_{2N}+1)^2 \\
&\le 2 + 2C_k + 4N(a_{2N}+1)^2 \le C_{k+1}\le C_{2n+2},
\end{align*}
by definition of the sequence $(C_n)_{n\in \ig 0,2N\id}$, if $\alpha_N$ and $\e_N$ are chosen only depending on $N$ and small enough.

\vip

It is thus sufficient to show that there exists $c_N^4>0$ positive and finite such that for all $z\in G^{k}_{\cK_3,\cK_4}$,
$$
\PP^{X,N,\theta,*}_z\Big[ \bigcap_{K\in \cK_1\cup\cK_2} B_K\cap C_K \Big] \ge c^4_N.
$$
By performing the same computation as for \eqref{eq:prepagirsanov}, we get the existence of a positive constant $c^4_N$ such that
$$
\inf_{t\in [0,t_N]} \min_{K\in \cK_1}\min_{i\in K,j\notin K} \|X^i_t - X^j_t\| \ge c^4_N.
$$
Thus we apply Remark \ref{rk:girsanovopen} with $\cK = \cK_1\cup\cK_2$, $T=t_N$ and $\eta = c^4_N/2$, we get,
\begin{align*}
    \PP^{X,N,\theta,*}_z\Big[ \cap_{K\in \cK_1\cup\cK_2} B_K\cap C_K \Big] &\ge \exp (-\tilde C_{\cK,N,\theta,\eta}T)\tilde\QQ^{\cK,T,\eta}_z\Big[\bigcap_{K \in \cK_1\cup\cK_2} B_K\cap C_K \Big]\\
    &= \exp (-\tilde C_{\cK,N,\theta,\eta}T)\prod_{K \in \cK_1\cup\cK_2} \tilde\QQ^{\cK,T,\eta}_z[B_K] \tilde\QQ^{\cK,T,\eta}_z[C_K]\\
    &=\exp (-\tilde C_{\cK,N,\theta,\eta}T)\prod_{K \in \cK_1\cup\cK_2} \tilde \PP[\tilde B_K] \tilde \PP[\tilde C_K],
\end{align*}
where for all $K\in \cK_1\cup\cK_2\setminus \{L\}$, 
\begin{gather*}
     \tilde B_K = \Big\{ \sup_{t\in [0,t_N]} \|S^K_t-S_K(z)\|\le \alpha_N \Big\} \\
     \tilde C_K = \Big\{ \sup_{t\in [0,t_N]}R^K_t\le (2^{-k}\e_N)^{1/|K|} +\alpha_N\quad \mbox{ and } \quad  R^K_{t_N}\le (2^{-k-1}\e_N)^{1/|K|} \Big\},
\end{gather*}
and
\begin{gather*}
     \tilde B_L = \Big\{ \sup_{t\in [0,t_N]} \|S^L_t
     -S_L(z) - t_N^{-1}t(a_{2N-k-1}+1)(v-S_L(z))\|\le \alpha_N \Big\} \\
     \tilde C_L = \Big\{ \sup_{t\in [0,t_N]} R^L_t\le (2^{-k}\e_N)^{1/|L|}+\alpha_N \quad \mbox{ and } \quad R^L_{t_N}\le (2^{-k-1}\e_N)^{1/|L|} \Big\},
\end{gather*}
where for all $K\in \cK_1\cup\cK_2$, $(S^K_t)_{t\ge 0}$ is a $2$-dimensional Brownian motion with diffusivity coefficient $|K|^{-1/2}$ starting from $S_K(z)$, $(R^K_t)_{t\ge 0}$ is a squared Bessel process with dimension $d_{\theta,N}(|K|)$ starting from $R_K(y)\le (2^{-k}\e_N)^{1/|K|}$. The family of processes $\{ (S^K_t)_{t\ge 0}, (R^K_t)_{t\ge 0}, K\in \cK_1\cup\cK_2\}$ is built on some probability space $(\tilde \Omega,\tilde \cF,\tilde \PP)$ and is mutually independent. The result thus holds thanks to classical results on Brownian motion and Bessel processes, see \cite{ry}.  

\vip

{\it Step 2.} Recalling the definition of $\e_N$ and $M_N$ in Step 1, we show that there exists $\e_N'\in (0,\e_N)$ such that
    \begin{align}\label{eq:knb}
   c_2:= \inf_{u\notin E_2^{\e_N}} \PP^{U,N,\theta}_u\Big[\mbox{ it exists } t\in [M_N^{-1},M_N ] \quad \mbox{ such that } \quad U_t \in E^{\e_N'}_2\Big] >0.
    \end{align}   
Indeed, if $(\Omega,\cF,\PP,(\cF_t)_{t\ge 0}, (B_t)_{t\ge 0}, (X_t)_{t\ge 0})$ is a $KS(\theta,N)$-process, we set $(U_t)_{t\ge 0} = (\Phi (X_{\rho_t}))_{t\ge 0}$ with usual notations, which has law $\PP^{U,N,\theta}$ under $\PP$. According to $\mathrm{Decomposition}(\theta,N)$-(i), $(U_t)_{t\ge 0}$ has the law of a $SKS(\theta,N)$-process. On the event where $X_1\in E^1_2$ and $R_{\ig 1,N\id}(X_t) \in [M_N^{-1},M_N]$ for all $t\in [0,1]$, we have $U_{A_1}\in E^{M_N^{-1}}_2$ and $M^{-1}_N\le A_1 \le M_N$, which gives the result with $\e_N' = M_N^{-1}\wedge (\e_N/2)$. 

\vip

{\it Step 3.} We show that there exists $c_3>0$ and $\e_N'' >0$ such that for all $A\in \cB(\sS_N)$,
$$
\inf_{u\in \sS_N} \PP_u^{U,N,\theta} [U_{2M_N} \in A\cap E^{\e_N''}_2] \ge c_3\sigma (A \cap E^{\e_N''}_2), 
$$
where $\sigma$ is the surface measure on $\sS_N$. 

\vip

Let $(U^{\e_N'}_t)_{t\ge 0}$ be the (classical) unique solution of the smoothed SDE
\begin{align*}
    \dd U^{\e_N'}_t = \pi_H\circ \pi_{(U^{\e_N'}_t)^\perp}(\dd B_t + b^{\e_N'}(U^{\e_N'}_t)\dd t) - \frac{2N-3}{2}U^{\e_N'}_t\dd t, 
\end{align*}
where $b^{\e_N'}$ is smooth such that for all $x\in E^{\e_N'/2}_2$, $b^{\e_N'}(x) = b(x)$, recall \eqref{eq:drift}, and there exists $V^{\e_N'} : (\rr^2)^N\to \rr$ smooth such that $b^{\e_N'} = \nabla V^{\e_N'}$. 

\vip

Since $\e_N'/2<\e_N'$, by standard heat kernel theory on weighted manifolds (see, e.g., Grigor'yan \cite[XX, Section 6]{Grigor'yan}) there exists $\eta_N >0$ such that for all $y\in E^{\e_N'}_2$,
\begin{align}\label{eq:minordensity}
\inf_{t\in [(2M_N)^{-1},2M_N]}\inf_{u\in E^{\e_N'}_2} p^{k,\e_N'/2}_t(u,y) \ge \eta_N,
\end{align}
where $p^{k,\e_N'/2}_{M_N}$ is the semigroup at time $M_N$ of the process $(U^{\e_N'}_t)_{t\ge 0}$ killed upon exiting $E^{\e_N'/2}_2$. 
This immediately gives that for $u\in E^{\e_N}_2$, all $A\in \cB( \sS_N)$, 
\begin{align*}
    \PP_u^{U,N,\theta} [U_{2M_N} \in A\cap E^{\e_N'/2}_2] \ge \eta_N \sigma(A\cap E^{\e_N'/2}_2).
\end{align*}
Then, for all $u\notin E^{\e_N}_2$, all $A\in \cB(\sS_N)$, using the strong Markov property,
\begin{align*}
    \PP^{U,N,\theta}_u[ U_{2M_N} \in A\cap E^{\e_N'/2}_2] &\ge \E^{U,N,\theta}_u\Big[ \indiq_{M_N^{-1}\le \tau_2^{\e_N'}\le M_N} \PP_{U_{\tau^{\e'_N}_2}} [U^{\e_N'}_{2M_N - \tau^{\e_N'}_2} \in A\cap E^{\e_N'/2}_2] \Big],\\
    &\ge c_2 \eta_N \sigma(A\cap E^{\e_N'/2}_2 ),
\end{align*}
where we also used Step 2. We conclude by setting $\e_N'' = \e_N'/2$ and $c_3 = \eta_N \wedge (c_2\eta_N)$. 

\vip

{\it Step 4.} We conclude the proof using a classical Harris-Meyn-Tweedie argument. Let $u\in \sS_N$. We build a coupling $(U^1_n, U^2_n)_{n\in \nn}$ as follows : $U^1_0 = u$ and $U^2_0$ follows the law $Z^{-1}\beta^N_\theta$ with $Z$ a normalizing constant. then for every $n\ge 0$, 
\begin{itemize}
    \item if $U^1_n = U^2_n$ then we set $U^1_{n+1} =  U^2_{n+1} = Z_1$ with $Z_1$ following the law of $U_{M_N}$ where $(U_t)_{t\ge 0}$ being a $SKS(\theta,N)$-process starting at $U^1_n$,
    \item if $U^1_n \ne U^2_n$, then with probability $c_3\sigma (E^{\e_N''})$, we sample a random variable $Z_2$ with uniform law on $E^{\e_N''}_2$ and set $U_{n+1} = \tilde U_{n+1} = Z_2$. Otherwise (thus with probability $1-c_3\sigma (E^{\e_N''}_2)$), we draw independently $U^1_{n+1}$ and $U^2_{n+1}$ respectively according to the probabilities $P_{U^1_n}$ and $P_{U^2_n}$ where for al $u'\in \sS_N$, all $A\in \cB(\sS_N)$,
    $$
    P_{u'}(A) = \frac{\PP^{U,N,\theta}_{u'}(U_{M_N} \in A) - c_3\sigma(A\cap E^{\e'}_2)}{1-c_3\sigma (E^{\e_N''}_2)}.
    $$
\end{itemize}
First, we get that for all $u'\in E^{\e_N''}_2$, $P_{u'}$ is a probability so that the preceding coupling is well defined. Indeed, let $u'\in \sS_N$. We first have $P_{u'}(\sS_N) = 1$. Moreover, for all $A \in \cB(\sS_N)$, $P_{u'}(A)\ge 0$ according to Step 3. 

\vip

Now by construction and thanks to the Markov property (see Theorem \ref{mainresult}-(ii) proved in Subsection \ref{subsec:uniquenessSKS}), we get that for all $n\ge 0$ the laws of $U^1_n$ and $U^2_n$ are respectively the laws of a $SKS(\theta,N)$ at time $nM_N$ starting respectively at $u$ and with law $Z^{-1}\beta^N_\theta$. Thus, since $Z^{-1}\beta^N_\theta$ is an invariant probability measure for a $SKS(\theta,N)$ process, and by definition of the coupling, we get
\begin{align*}
    \|\cL_u (U_{2nM_N})
 - Z^{-1}\beta^N_\theta \|_{TV} \le \PP( U^1_n \ne U^2_n)\le \PP(\tau> n),
 \end{align*}
 with 
 \begin{align*}
     \tau = \inf\{ n\ge 0 : U^1_n = U^2_n\}.
 \end{align*}
It is now classical by conditioning on the event $\{\tau>n\} $ to get for all $n\ge 0$, $\PP[\tau> n+1] \le (1-c_3\sigma ( E^{\e_N''}_2))\PP[\tau >n]$, so that $\PP[\tau > n] \le (1-c_3\sigma ( E^{\e_N''}_2))^{n} \le e^{-c'n}$ with $c'>0$.
 
\vip
 
We finally conclude by realizing that for all $t \ge 0$,
\begin{align*}
   \sup_{u\in \sS_N} \|\cL_u (U_t) - Z^{-1}\beta^N_\theta\|_{TV} \le  \sup_{u\in \sS_N} \|\cL_{\mu_u} (U_{2M_N\lfloor t/(2M_N) \rfloor}) - Z^{-1}\beta^N_\theta\|_{TV}, 
\end{align*}
where $\mu_u$ is the law of a $SKS(\theta,N)$-process at time $t-2M_N\lfloor t/(2M_N)\rfloor$ starting from $u$. Of course this directly implies
\begin{align*}
\sup_{u\in \sS_N} \|\cL_u (U_t) - Z^{-1}\beta^N_\theta\|_{TV} &\le \sup_{u\in \sS_N} \|\cL_u (U_{2M_N\lfloor t/(2M_N)\rfloor}) - Z^{-1}\beta^N_\theta\|_{TV} \\
&\le \exp (-c' \lfloor t/(2M_N) \rfloor) \\
&\le e^{c'} \exp (-c'(2M_N)^{-1} t),
\end{align*}
which concludes by setting $C=e^{c'}$ and $c = c' (2M_N)^{-1}$.

\end{proof}

\section{Uniqueness}\label{sec:uniqueness}
In this section we prove the uniqueness part of Theorem \ref{mainresult}. Since the result is already proven in Proposition \ref{prop:girsanoveta} in the cases where $k_0\le N$ (since we don't expect to extend the process after $\tau_{k_0}$) or $k_2>N$ and $X_0\in E_{k_0}\cap E_N$ (since in this case $\tau_N = \infty$ a.s., see Lemma \ref{lem:dimension}), we can suppose that $d_{\theta,N}(N)\in (0,2)$ or $d_{\theta,N}(N)\ge 2$ and $x \in E_{k_0}$ with $R_{\ig 1,N\id}(x) = 0$. However, since these two cases are the same in spirit and that the case where $d_{\theta,N}(N)\ge 2$ is easier since there is only one excursion to handme, we focus on the case where $d_{\theta,N}(N)\in (0,2)$.

\vip

We begin with the following Lemma.
\begin{lemma}\label{lem:besselbarycentretotal}
    Let $(\Omega,\cF,(\cF_t)_{t\ge 0},\PP, (B_t)_{t\ge 0}, (X_t)_{t\ge 0})$ be a $KS(\theta,N)$-process starting at $x\in (\rr^2)^N$. Then, $(S_{\ig 1,N\id}(X_t))_{t\ge 0}$ is a $2$-dimensional Brownian motion with diffusion coefficient $N^{-1/2}$ starting at $S_{\ig 1,N\id}(x)$ and $(R_{\ig 1,N\id}(X_t))_{t\ge 0}$ is a squared Bessel process with dimension $d_{\theta,N}(N)$ starting at $R_{\ig 1,N\id}(x)$. Moreover, $(R_{\ig 1,N\id}(X_t))_{t\ge 0}$ and $(S_{\ig 1,N\id}(X_t))_{t\ge 0}$ are independent.
\end{lemma}
\begin{proof}
The fact that $(S_{\ig 1,N\id}(X_t))_{t\ge 0}$ is as described in the statement is clear. We now treat $(R_{\ig 1,N\id}(X_t))_{t\ge 0}$. We set for $\e >0$, the semi martingale
    \begin{align*}
     \mbox{ for all } t\ge 0, \quad   X^\e_t = x + B_t + \int_0^t b_{\e f(R_{\ig 1,N\id}(X_s))}(X_s)\dd s,
    \end{align*}
    with $f$ given, as in Definition \ref{def:KS/KSK}. Let $\delta >0$ and $\chi :\rr_+\to [0,1]$ smooth such that $\chi(r) = 0$ for all $r\in [0,1]$ and $\chi (r) = 1$ for all $r\ge 2$. We recall that for all $x\in (\rr^2)^N$,
    $$
    \nabla R_{\ig 1,N\id}(x) = 2(x-\gamma (x)) \quad \mbox{ and } \quad \Delta R_{\ig 1,N\id}(x) = 4(N-1).
    $$
    We set for all $x\in (\rr^2)^N$, $F_\delta (x) = R_{\ig 1,N\id}(x) \chi (R_{\ig 1,N\id}(x)/\delta)$ and we immediately have
    \begin{gather*}
\nabla F_\delta (x) = 2(x-\gamma (x)) (\chi(R_{\ig 1,N\id}(x)/\delta) +\delta^{-1} \chi'(R_{\ig 1,N\id}(x)/\delta)R_{\ig 1,N\id}(x)),\\
\Delta F_\delta (x) = 4(N-1)\Big( \chi(R_{\ig 1,N\id}(x)/\delta)+\delta^{-1} \chi'(R_{\ig 1,N\id}(x)/\delta)R_{\ig 1,N\id}(x) \Big)\\
+4R_{\ig 1,N\id} (x) \Big(2\delta^{-1}  \chi'(R_{\ig 1,N\id}(x)/\delta) +  \delta^{-2}\chi''(R_{\ig 1,N\id}(x)/\delta)R_{\ig 1,N\id}(x)) \Big).
    \end{gather*}
    Using the It\^o's formula, and the fact that $b(y)\cdot \gamma (z) = 0$ for all $y\in (\rr^2)^N$, $z\in \rr^2$, we get
    \begin{align*}
        F_\delta (X^\e_t) = F_\delta (x) + M^{\e,\delta}_t + I^{1,\e,\delta}_t + I^{2,\e,\delta}_t +R^{1,\e,\delta}_t + R^{2,\e,\delta}_t + R^{3,\e,\delta}_t +R^{4,\e,\delta}_t,
    \end{align*}
    where
    \begin{gather*}
         M^{\e,\delta}_t = 2\int_0^t \chi(R_{\ig 1,N\id}(X^\e_s)/\delta)(X^\e_s - \gamma (S_{\ig 1,N\id}(X^\e_s)))\cdot \dd B_s \\
         I^{1,\e,\delta}_t = 2\int_0^t \chi(R_{\ig 1,N\id}(X_s)/\delta)\indiq_{R_{\ig 1,N\id}(X_s) \ge \e f(R_{\ig 1,N\id}(X_s))} b(X_s)\cdot X^\e_s  \dd s,\\
         I^{2,\e,\delta}_t = 2(N-1)\int_0^t\chi(R_{\ig 1,N\id}(X^\e_s)/\delta) \dd s,\\
         R^{1,\e,\delta}_t = 2\delta^{-1}\int_0^t \chi'(R_{\ig 1,N\id}(X^\e_s)/\delta)R_{\ig 1,N\id}(X^\e_s))(X^\e_s - \gamma (S_{\ig 1,N\id}(X^\e_s)))\cdot \dd B_s,\\
         R^{2,\e,\delta}_t = 2\delta^{-1}\int_0^t R_{\ig 1,N\id}(X^\e_s) \chi'(R_{\ig 1,N\id}(X^\e_s)/\delta)\indiq_{R_{\ig 1,N\id}(X_s) \ge \e f(R_{\ig 1,N\id}(X_s))} b(X_s)\cdot X^\e_s \dd s,\\
         R^{3,\e,\delta}_t = 2(N+1)\delta^{-1}\int_0^t   R_{\ig 1,N}(X^\e_s) \chi'(R_{\ig 1,N\id}(X^\e_s)/\delta)\dd s,\\
         R^{4,\e,\delta}_t = 2\delta^{-2}\int_0^t  \chi''(R_{\ig 1,N\id}(X^\e_s)/\delta) (R_{\ig 1,N\id}(X^\e_s))^2\dd s. 
    \end{gather*}
    We first aim at letting $\e \to 0.$ Since by definition of a $KS(\theta,N)$-process, $\sup_{t\in [0,T]}\|X^\e_t - X_t\|\to 0$ in probability as $\e \to 0$ for every $T>0$, we get that for all $t\ge 0$,
    \begin{align}\label{eq:besseldelta}
        F_\delta (X_t) = F_\delta (x) + M^{\delta}_t + I^{1,\delta}_t + I^{2,\delta}_t +R^{1,\delta}_t + R^{2,\delta}_t + R^{3,\delta}_t +R^{4,\delta}_t,
    \end{align}
    where
    \begin{gather*}
         M^{\delta}_t = 2\int_0^t \chi(R_{\ig 1,N\id}(X_s)/\delta)(X_s - \gamma (S_{\ig 1,N\id}(X_s)))\cdot \dd B_s \\
         I^{1,\delta}_t = -\theta (N-1)\int_0^t \chi(R_{\ig 1,N\id}(X_s)/\delta)  \dd s,\\
         I^{2,\delta}_t = 2(N-1)\int_0^t\chi(R_{\ig 1,N\id}(X_s)/\delta) \dd s,\\
         R^{1,\delta}_t = 2\delta^{-1}\int_0^t \chi'(R_{\ig 1,N\id}(X_s)/\delta)R_{\ig 1,N\id}(X_s))(X_s - \gamma (S_{\ig 1,N\id}(X_s)))\cdot \dd B_s,\\
         R^{2,\delta}_t = -\theta (N-1)\delta^{-1}\int_0^t R_{\ig 1,N\id}(X_s) \chi'(R_{\ig 1,N\id}(X_s)/\delta) \dd s,\\
         R^{3,\delta}_t = 2(N+1)\delta^{-1}\int_0^t  R_{\ig 1,N}(X_s) \chi'(R_{\ig 1,N\id}(X_s)/\delta)\dd s,\\
         R^{4,\delta}_t = 2\delta^{-2}\int_0^t \chi''(R_{\ig 1,N\id}(X_s)/\delta) (R_{\ig 1,N\id}(X_s))^2\dd s. 
    \end{gather*}
    Indeed, the only non obvious convergences are $I^{1,\e,\delta}\to I^{1,\delta}$ and $R^{2,\e,\delta} \to R^{2,\delta}$ which we now justify. Since we treat them in a similar way, we only justify $I^{1,\e,\delta}\to I^{1,\delta}$. Recalling Step 1 in the proof of Proposition \ref{prop:barycentreetbessel}, we have for all $x\in E_2$, $b(x)\cdot x = -\theta (N-1) /2$. We thus have
    \begin{align*}
        I^{1,\delta}_t = 2\int_0^t \chi( R_{\ig 1,N\id}(X_s) /\delta) b(X_s)\cdot X_s \dd s,
    \end{align*}
since $\int_0^\infty \indiq_{X^N_s\notin E_2}\dd s = 0$ a.s. by definition. From the latest, we deduce,
    \begin{align*}
        \sup_{t\in [0,T]} \|I^{1,\e,\delta}_t - I^{1,\delta}_t\| \le 2 \int_0^T \chi (R_{\ig 1,N\id}(X_s)/\delta) \|b(X_s)\| \dd s \sup_{t\in [0,T]} \|X^\e_t - X_t\|.
    \end{align*}
    By hypothesis, we have $\int_0^T \chi (R_{\ig 1,N\id}(X_s)/\delta) \|b(X_s)\| \dd s<\infty$ a.s. Thus we get the result since $\sup_{t\in [0,T]} \|X^\e_t - X_t\|$ converges to $0$ in probability. 

    \vip

    We now aim to let $\delta \to 0$. It is straightforward that for all $s\ge 0$, a.s. $F_\delta (X_s) \to R_{\ig 1,N\id} (X_s)$ as  $\delta \to 0$. Moreover, we have
    \begin{align*}
         |R^{2,\delta}_t| + |R^{3,\delta}_t| + |R^{4,\delta}_t| \le C\Big(\delta +  \int_0^t \indiq_{ \delta \le R_{\ig 1,N\id} (X_s) \le 2\delta}\dd s \Big),
    \end{align*}
  and this converges a.s. to $0$ as $\delta \to 0$ by the dominated convergence theorem. Besides, by the Ito's isometry, we have 
  \begin{align*}
      \E [(R^{1,\delta}_t)^2] = 4\delta^{-2} \int_0^t (\chi'(R_{\ig 1,N\id}(X_s)/\delta))^2(R_{\ig 1,N\id}(X_s))^3\dd s\le \delta t \|\chi'\|_\infty \to 0 \quad \mbox{ as } \delta \to 0.
  \end{align*}
Another ingredient is that
  \begin{align*}
      I^{1,\delta}_t + I^{2,\delta}_t = d_{\theta,N}(N) \int_0^t \chi ( R_{\ig 1,N\id}(X_s)/\delta)\dd s \to d_{\theta,N}(N) t
  \end{align*}
    a.s. as $\delta \to 0$ by the dominated convergence theorem since $\int_0^\infty\indiq_{R_{\ig 1,N\id}(X_s) = 0}\dd s = 0$ a.s. by assumption. Finally, defining for all $t\ge 0$, 
    $$
    M_t =  2\int_0^t (X_s - \gamma (S_{\ig 1,N\id}(X_s)))\cdot \dd B_s, 
    $$
    we get by the Ito's isometry, 
    \begin{align*}
        \E[(M^\delta_t - M_t)^2] = \E \Big[\int_0^t (1-\chi ( R_{\ig 1,N\id}(X_s)/\delta))^2R_{\ig 1,N\id} (X_s) \dd s \Big] \le \delta t \to 0
    \end{align*}
    as $\delta \to 0$. Combining the previous convergences with \eqref{eq:besseldelta} gives that $(R_{\ig 1,N\id}(X_t))_{t\ge 0}$ is as stated. It only remains to prove independence between $(S_{\ig 1,N\id}(X_t))_{t\ge 0}$ and $(R_{\ig 1,N\id}(X_t))_{t\ge 0}$ but the reasoning is exactly as in Step 3 of the proof of Proposition \ref{prop:barycentreetbessel}. 
 \end{proof}
We are now in position to prove uniqueness in Theorem \ref{mainresult}-(i). According to the strong Markov property and Proposition \ref{prop:girsanoveta}-(i), it only remains to prove that for every $x\in (\rr^2)^N$ such that $R_{\ig 1,N\id}(x)= 0$, there exists a unique in law $KS_x(\theta,N)$-process. To this end it suffices to prove $\mathrm{Decomposition}(\theta,N)$ thanks to Proposition \ref{prop:existuniqSKSeternel}. Let $x\in (\rr^2)^N$ such that $R_{\ig 1,N\id}(x) = 0$ and  $(\Omega,\cF,(\cF_t)_{t\ge 0},\PP, (B_t)_{t\ge 0}, (X_t)_{t\ge 0})$ be a $KS_x(\theta,N)$-process. According to Lemma \ref{lem:besselbarycentretotal}, the process $(R_{\ig 1,N\id}(X_t))_{t\ge 0}$ is a squared Bessel process with dimension $d_{\theta,N}(N)$ starting from $0$. Let $\cZ = \{t\ge 0 : R_{\ig 1,N\id}(X_t) = 0\}$. Since $\cZ^c$ is an open set of $\rr$, we can write $\cZ^c= \cup_{i\ge 0}(a_i,b_i)$ with the intervals $(a_i,b_i)$ disjoints. Moreover, we can always choose an indexation such that $i\in \nn \mapsto (a_i,b_i)$ is measurable (using for example the inverse of the local time of $(R_{\ig 1,N\id}(X_t))_{t\ge 0}$ to describe the sequence $(a_i,b_i)_{i\ge 0}$). We set
$$
\cG =\sigma \Big( (S_{\ig 1,N\id}(X_t))_{t\ge 0}, (R_{\ig 1,N\id}(X_t))_{t\ge 0}\Big).
$$
We set for all $i\ge 0$, $(U^{i}_t)_{t\in \rr} := (\Phi(X_{\rho^{i}_t}))_{t\in \rr}$ where $\rho^{i}_t =(A^i)^{-1}_t$ with 
$$
A^i_t = \int_{(a_i + b_i)/2}^t(R_{\ig 1,N\id}(X_s))^{-1}\dd s.
$$
It is sufficient to show that conditionally on $\cG$, for all $A\subset \nn$ finite, the family $\{ (U^i_t)_{t\in \rr}, i\in A\}$ is i.i.d with for all $i\in A$, $(U^i_t)_{t\in \rr}$ is distributed as an eternal $SKS(\theta,N)$-process. 

\vip

Let $A\subset \nn$. Since the family $\{ (a_i,b_i), i\in A\}$ is finite and composed of disjoint intervals, up to composition with bijections we can always suppose that $A=\ig 0,m\id$ for some $m\in \nn$, that $i\in \ig 0,m\id\mapsto (a_i,b_i)$ is measurable and the sequences $(a_i)_{i\in \ig 0,M\id}$ and $(b_i)_{i\in \ig 0,M\id}$ are increasing. 

\vip

We claim that it is sufficient to show that for all $i\in \ig 0,m\id$, conditionally on $\cG \vee \cF_{b_{i-1}}$, $(U^i_t)_{t\in \rr}$ is distributed as an eternal $SKS(\theta,N)$-process, where we set $b_{-1} = 0$ by convention. Indeed this implies that conditionally on $\cG$, $(U^i_t)_{t\in \rr}$ is distributed as an eternal $SKS(\theta,N)$-process and that for all $i\in \ig 0,m\id$, $(U^i_t)_{t\ge 0}$ is independent from $\{ (U^j_t)_{t\in \rr}, j\le i-1\}$, which implies the result. Let us now fix $i\in \ig 0,m\id$ and conclude. 

\vip

Fix $t_0\in \QQ_+$ and $\sigma^{t_0}_0 = \inf\{t\ge t_0 : R_{\ig 1,N\id} (X_t) = 0\}$. According to the strong Markov property the process $(\indiq_{t<\sigma^{t_0}_0}X_t + \triangle \indiq_{t\ge \sigma^{t_0}_0})_{t\ge 0 }$ is a $KS(\theta,N)$-process killed when exiting $E_N$ starting at $X_{t_0}$. According to $\mathrm{Decomposition}(\theta,N)$-(i), conditionally on $\cF_{t_0}$, there exists a $SKS(\theta,N)$-process $(U^{t_0}_t)_{t\ge 0}$ such that for all $t\in [t_0, \sigma^{t_0}_0)$, 
$$
\Phi (X_{\rho^{t_0}_t}) = U^{t_0}_t, \quad \mbox{ with } \rho^{t_0}_t = (A^{t_0})^{-1}_t \quad \mbox{ and } \quad A^{t_0}_t = \int_{t_0}^t (R_{\ig 1,N\id}(X_s))^{-1}\dd s.
$$
If $t_0\in (a_i,b_i)$, we see that by definition, for all $t\in (a_i,b_i)$, $U^i_{A^i_t} = U^{t_0}_{A^{t_0}_t}$. Since $A^i_t = A^{t_0}_t + T^i_{t_0}$ with $T^i_{t_0} = \int_{(a_i+b_i)/2}^{t_0}(R_{\ig 1,N\id}(X_s))^{-1} \dd s$ and $A^i : (a_i,b_i)\to \rr$ is bijective (see Lemma \ref{lem:horlogeexplosion}), we get that
 \begin{align*}
  \mbox{ for all } s\in \rr, \quad  U^i_{s+T^i_{t_0}} = U^{t_0}_s.
\end{align*}
Since $(U^{t_0}_t)_{t\ge 0}$ solves the SDE \eqref{eq:EDSprobaU} we deduce that there exists a $2N$-dimensional Brownian motion $(B^{t_0}_t)_{t\ge 0}$ such that for all $t\ge 0$,
\begin{align}\label{eq:protoEDSU}
    U^i_{t+T^i_{t_0}} = U^i_{T^i_{t_0}} + \int_{T^i_{t_0}}^{t+T^i_{t_0}}\pi_{H_N}\pi_{(U^{i}_s)^\perp}\dd B^{t_0}_s +\int_{T^i_{t_0}}^{t+T^i_{t_0}}\pi_{(U^{i}_s)^\perp}\pi_{H_N}b_{\e}(U^{i}_s)\dd s +\frac{2N-3}{2} \int_{T^i_{t_0}}^{t+T^i_{t_0}}U^{i}_s\dd s. 
\end{align}
Fix $s_0 \in \rr$. Conditionally on $\cG\vee \cF_{b_{i-1}}$, one can find a sequence $(t_n)_{n\ge 1}$ of elements of $(a_i,b_i)$ decreasing to $a_i$ such that for all $n\ge 1$, $A^i_{t_n}<s_0$. This together with \eqref{eq:protoEDSU} gives the existence of some $2N$-dimensional Brownian motion $(B^{s_0}_u)_{u\ge 0}$ such that for all $t\ge s_0$, 
\begin{align}\label{eq:protoEDSU2}
    U^i_t = U^i_{s_0} + \int_{s_0}^t \pi_{H_N}\pi_{(U^{i}_s)^\perp}\dd B^{s_0}_s+\int_{s_0}^t\pi_{(U^{i}_s)^\perp}\pi_{H_N}b_{\e}(U^{i}_s)\dd s +\frac{2N-3}{2} \int_{s_0}^tU^{i}_s\dd s   
\end{align}
Moreover, by the preceding, we that for all $n\ge 1$, the law of $U^i_{s_0}$ conditionally on $\cG\vee \cF_{b_{i-1}}$ is distributed as a $SKS(\theta,N)$-process at time $s_0-A^i_{t_n}$ starting from some law $\mu$. Thus, according to Proposition \ref{prop:ergodicity},
$$
\| \cL (U^i_{s_0} | \cG\vee \cF_{b_{i-1}} ) - Z^{-1}\beta^N_\theta\|_{TV} \le C'\exp(-c' (s_0 - \rho^i_{t_n})).
$$
However, $A^i_{t_n} \to -\infty$ as $n\to \infty$ (see Lemma \ref{lem:horlogeexplosion}), implying that conditionally on $\cG\vee \cF_{b_{i-1}}$, $U^i_{s_0}$ is $Z^{-1}\beta^N_\theta$-distributed. This together with \eqref{eq:protoEDSU2} gives the result.

\section{Strong Markov property and Feller property}\label{sec:strongMarkov}
 We denote for all $\mu\in (\rr^2)^N$ by $P^N_\mu$ the law of $(X^N_t)_{t\ge0}$ built in Section \ref{sec:constructionKS} when starting with initial law $\mu$.
\begin{proof}[Proof of the strong Markov property in Theorem \ref{mainresult}-(i)]

We show that 
$$
(\Omega,\cF,(\cF_t)_{t\ge 0}, (X^N_t)_{t\ge 0},(P^N_x)_{x\in (\rr^2)^N})
$$
is a strong Markov process such that for all $A\in \cB(C(\rr_+,(E_{k_0})_\triangle))$, $x\in (E_{k_0})_\triangle\mapsto P^N_x(A)$ is measurable. Let $\tau$ a $(\cF_t)_{t\ge 0}$-stopping time and $\mu \in \cP((\rr^2)^N)$. According to the uniqueness result obtained in Section \ref{sec:uniqueness} and a slight adaptation of \cite[Corollary 2]{yw}(since a $KS$-process do not solve rigorously an SDE in the classical sense in general), it is sufficient to show that conditionally on $\cF_\tau$, $(X_{\tau+t})_{t\ge 0}$ is a $KS(\theta,N)$-process starting from $X_\tau$.

\vip

According to Section \ref{sec:constructionKS} we have for all $t\ge 0$,
\begin{align}\label{eq:candiatavectau}
X^N_{\tau+t} =& \gamma (M_{\tau +t}) + (D_{\tau+t})^{1/2}\Big( U^{N}_{A_{\tau+t}} \indiq_{\tau+t\in [0,\tau^0_0)} + \sum_{i\ge 0} U^{N,i}_{A^i_{\tau+t}}\indiq_{t\in (a_i-\tau,b_i-\tau)}\Big).
\end{align}
Conditionally on $\cF_\tau$, $(M_{\tau+t})_{t\ge 0}$ is a $2$-dimensional Brownian motion with diffusivity coefficient $N^{-1/2}$ starting at $M_\tau$ and $(D_{\tau+t})_{t\ge 0}$ is a squared Bessel process with dimension $d_{\theta,N}(N)$ starting at $D_\tau$, these two processes being independent. Moreover, one can write
$$
\{ t\ge 0 : D_{\tau+t} \ne 0\} = [0,\tilde\tau^0_0) \cup \Big( \cup_{i\in J} (\tilde a_i, \tilde b_i )\Big),
$$
with $J = \{i\ge 0: a_i >\tau\}$, $\tilde a_i = (a_i- \tau)\vee 0$, $\tilde b_i = (b_i-\tau)\vee 0$ and $\tilde \tau^0_0 = \inf \{ t\ge 0 : D_{\tau+t} = 0\}$. One can thus rewrite \eqref{eq:candiatavectau} as
\begin{align*}
X^N_{\tau+t} =& \gamma (M_{\tau +t}) + (D_{\tau+t})^{1/2}\Big( \tilde U^{N}_{A^\tau_{t}} \indiq_{t<\tilde\tau^0_0} + \sum_{i\in J}  U^{N,i}_{A^{\tau,i}_{t}}\indiq_{t\in (\tilde a_i,\tilde b_i)}\Big),
\end{align*}
where for all $s\in \rr$, 
$$
\tilde U^{N}_s = \indiq_{\tau < \tau^0_0} U^N_{s+A_\tau} + \sum_{i\ge 0} \indiq_{\tau\in (a_i,b_i)} U^{N,i}_{s + A^i_\tau }
$$
and for all $i\in J$, 
\begin{gather*}
    \mbox{ for all } t\in [0,\tilde \tau^0_0), \quad A^\tau_t = \int_{0}^t \frac{\dd s}{D_{\tau+s}},\\
    \mbox{ for all } i\ge 0, \mbox{ all }t\in (\tilde a_i, \tilde b_i), \quad A^{\tau,i}_t = \int_{(\tilde b_i+\tilde a_i)/2}^t \frac{\dd s}{D_{\tau+s}}. 
\end{gather*}
Conditionally on $\cF_\tau \vee \sigma ( (D_{\tau+t})_{t\ge 0}, (M_{\tau+t})_{t\ge 0})$, the family $\{(\tilde U^N_t)_{t\ge 0}\}\cup\{ (U^{N,i}_t)_{t\in \rr}, i\in J\}$ is mutually independent. For all $i\in J$, $( U^{N,i}_t)_{t\in \rr}$ is distributed as an eternal $SKS(\theta,N)$-process by definition and independence properties.
Finally, if $\tau\in(a_i,b_i)$, then $A^i_\tau$ is $\sigma((D_t)_{t\ge0})$-measurable. Since, conditionally on $(D_t)_{t\ge0}$, $(U^{N,i}_t)_{t\in\rr}$ is a stationary eternal SKS$(\theta,N)$-process, so is $(U^{N,i}_{s+A^i_\tau})_{s\in\rr}$, and the Markov property at time $0$ shows that $(\tilde U^N_s)_{s\ge0}$ is an SKS$(\theta,N)$-process starting from $\tilde U^N_0$, conditionally on $\cF_\tau\vee\sigma((D_{\tau+t},M_{\tau+t})_{t\ge0})$. Putting all together, we conclude that conditionally on $\cF_\tau$, $(X^N_{\tau+t})_{t\ge 0}$ can be decomposed as in Section \ref{sec:constructionKS}, so that it is a $KS(\theta,N)$-process starting from $X_\tau$ according to Lemma \ref{lem:drifttendzero} and Lemma \ref{lem:EDSpv}, which concludes.
\end{proof}
The following Lemma is interesting in itself because it explains that thanks to the symmetry after a global collision between particles, the law of the process is explicit. Indeed, after a global collision, the angle restart with its invariant measure, so that it stays stationary for ever. Moreover this will be crucial to prove the Feller property or the association with the right Dirichlet form. 
\begin{lemma}\label{lem:symdesangles}
Let $x\in (\rr^2)^N$ which we write $x= \gamma (m) +\sqrt{r}u$ with $m\in \rr^2$, $r\ge 0$ and $u\in \sS_N$. For all $t\ge 0$, all $\varphi \in C((\rr^2)^N)$, we have
    \begin{align*}
        E^N_x[\varphi(X_t)\indiq_{t\ge \tau_N}] =\int \E[\varphi (\gamma (M_t) + D_t^{1/2}u') \indiq_{t\ge\tau_0} ]Z^{-1}\beta^N_\theta(\dd u'),
    \end{align*}
    where $(M_t)_{t\ge 0}$ is a $2$-dimensional Brownian motion with diffusivity coefficient $N^{-1/2}$ starting from $m$, $(D_t)_{t\ge 0}$ is a squared Bessel process with dimension $d_{\theta,N}(N)$ starting from $r^2$ and $\tau_0 = \inf \{ t\ge 0 :  D_t = 0\}$.
\end{lemma}
\begin{proof}
    Let $x = m+ru\in (\rr^2)^N$ and $\varphi\in C((\rr^2)^N)$ as in the statement. Let us condition on $\cG = \sigma ( (M_t)_{t\ge 0}, (D_t))_{t\ge 0})$, we get
    \begin{align*}
        E^N_x[\varphi (X^N_t)\indiq_{t\ge \tau_N}|\cG] &= E^N_x\Big[\varphi \Big(\gamma (M_t) + D_t^{1/2}\sum_{i\ge 0} U^{N,i}_{A^i_t}\Big)\indiq_{t\in (a_i,b_i)}|\cG\Big] \\
        &=\sum_{i\ge 0} \indiq_{t\in (a_i,b_i)}E^N_x\Big[\varphi (\gamma (M_t) + D_t^{1/2} U^{N,i}_{A^i_t})|\cG\Big]\\
        &=\sum_{i\ge 0} \indiq_{t\in (a_i,b_i)}\int_{\sS_N}\varphi (\gamma (M_t) + D_t^{1/2} u') Z^{-1}\beta^N_\theta(\dd u'),
    \end{align*}
    where we used the transfer theorem and that conditionally on $\cG$, for all $s\in \rr$, all $i\ge 0$ $U^{N,i}_s$ is distributed as $Z^{-1} \beta^N_\theta$. The last expression factorizes into
    \begin{align*}
        E^N_x[\varphi (X^N_t)\indiq_{t\ge \tau_N}|\cG] =\indiq_{t\ge \tau_0} \int_{\sS_N}\varphi (\gamma (M_t) + D_t^{1/2} u') Z^{-1}\beta^N_\theta(\dd u').
    \end{align*}
   Taking the expectation concludes the proof. 
\end{proof}
\begin{rk}
    If one take $r=0$ in the previous Lemma (meaning all the particles start at the same place), one finds explicitly the law of $X_t$ for all $t\ge 0$, which will lead to Proposition \ref{prop:densité}, proven in Appendix \ref{app:densité}. Moreover, this gives hope to find estimates providing uniqueness of the limit in this precise case when studying Propagation of chaos of the Keller-Segel particle system, even though this is only speculative at this stage. With a similar idea, one can compute explicitly the law of the same quantity in some other models as the Dyson model (the repulsive counterpart of the Keller-Segel model in dimension $1$, see \cite{cepalepingle} or \cite{dyson} for its definition).
\end{rk}
We now conclude this Section with the following:
\begin{proof}[Proof of the Feller Property]
We now show that our process is Feller in the sense that for all $\varphi \in C_0((\rr^2)^N)$, 
\begin{gather*}
\mbox{ for all } t>0, \quad x \mapsto E^N_x[\varphi (X^N_t)] \in C_0((\rr^2)^N),\\
\sup_{x\in (\rr^2)^N} | E^N_x(\varphi (X_t)) - \varphi (x) |\to 0\quad \mbox{ as } t\to 0.
\end{gather*}
Since $C_c^\infty (E_N) + \cC_{\ig 1,N\id}$ is dense in $C_0((\rr^2)^N)$ for the uniform convergence topology (see the paragraph above Theorem \ref{mainresult} for the definition of $\cC_{\ig 1,N\id}$), it is sufficient to prove the latter for all $\psi \in C_c^\infty (E_N) \cup \cC_{\ig 1,N\id}$. Of course it is immediate if $\psi \in \cC_{\ig 1,N\id}$ thanks to Proposition \ref{prop:barycentreetbessel} and properties of Brownian motion and squared Bessel processes, see \cite{ry}. We thus fix $\psi \in C^\infty_c(E_N)$. 

\vip

According to Lemma \ref{lem:symdesangles}, we have for all $x\in (\rr^2)^N$, $E^N_x[\psi (X^N_t)] = P^1_t\psi(x) + P^2_t\psi(x)$ with
\begin{gather*}
     P^1_t\psi(x) = E^N_x[\psi (X^N_t)\indiq_{t<\tau_N}]\\
     P^2_t\psi(x) = \int \E[\psi (\gamma (M_t) + D_t^{1/2}u') \indiq_{t\ge\tau_0} ]Z^{-1}\beta^N_\theta(\dd u'),
\end{gather*}
where we used the notations of Lemma \ref{lem:symdesangles}. According to Proposition \ref{prop:feller}, $P^1_t\psi \in C_0^\infty(E_N) \subset C_0((\rr^2)^N)$ and 
$$
\sup_{x\in E_N} |P^1_t\psi (x) - \psi(x)| = \sup_{x\in (\rr^2)^N} |P^1_t\psi (x) - \psi(x)| \to 0 \quad \mbox{ as } n\to \infty.
$$
It thus only remains to treat $P^2_t\psi$.

\vip

{\it Step 1.} We prove that $P^2_t\psi \in C_0((\rr^2)^N)$. Take a sequence $(x_n)_{n\ge 0}$ of elements of $(\rr^2)^N$ converging to some $x\in (\rr^2)^N$. According to the Feller property of the Brownian motion of the squared Bessel process (see \cite{ry}), the fact that $x\in (\rr^2)^N\mapsto (S_{\ig 1,N\id}(x), R_{\ig 1,N\id}(x))$ is continuous and \cite[Theorem 4.2.5 page 167]{ek}, we get that 
$$
\cL_{x_n}((M_t)_{t\ge 0}, (D_t)_{t\ge 0}) \to \cL_{x}((M_t)_{t\ge 0}, (D_t)_{t\ge 0})\quad \mbox{ weakly as } n\to \infty. 
$$
Moreover, for all $u' \in \sS_N$, the application
$$
(m_t,d_t)_{t\ge 0} \in C(\rr_+,\rr^2\times\rr_+) \to \psi ( \gamma (m_t) + d_t^{1/2} u') \indiq_{t\ge \tau_0(d)}, 
$$
with $\tau_0(d) = \inf\{ t\ge 0 : d_t = 0\}$, is $\cL_x((M_t)_{t\ge 0}, (D_t)_{t\ge 0})$-a.e. continuous by a similar argument as above. Thus for all $u'\in \sS_N$, 
$$
E^N_{x_n}[\psi (\gamma (M_t) + D_t^{1/2}u')\indiq_{t\ge \tau_0}] \to E^N_{x}[\psi (\gamma (M_t) + D_t^{1/2}u')\indiq_{t\ge \tau_0}] \quad \mbox{ as } n\to \infty.
$$
The dominated convergence theorem and the fact that $\beta^N_\theta (\sS_N)<\infty$ gives that $P^2_t\psi (x_n) \to P^2_t\psi (x)$ as $n\to \infty$. Moreover, the same reasoning gives that $P^2_t(x) \to 0$ as $x\to \infty$ which concludes this step. 

\vip

{\it Step 2.} We show that 
$$
\sup_{x\in (\rr^2)^N} | P^2_t\psi (x)| \to 0\quad \mbox{ as } t\to 0.
$$
Fix $\e >0$. Since $\psi \in C^\infty_c(E_N)$, there exists $\eta >0$ such that $\mathrm{Supp}\psi \subset E^\eta_N$. Using the strong Markov property, we have,
\begin{align*}
    \sup_{x\in (\rr^2)^N}|P^2_t\psi(x)| \le \|\psi\|_\infty \PP_0[t\ge \tau_\eta],
\end{align*}
with $\PP_0$ the law of a squared Bessel process with dimension $d_{\theta,N}(N)$ starting at $0$ and $\tau_\eta$ its first hitting time of $\eta$. We get the result since $\tau_\eta >0$ a.s. 
\end{proof}

\section{Ito's formula}\label{sec:ito}
We now prove $\mathrm{Ito}(\theta,N)$. The difficulty here is that potentially $\cB^N_\theta = \{\ig 1,N\id\}$, so that the classical Ito's formula does not hold. We will rely on the compensations of the angles somehow in the same way the Ito's integral relies on the compensations of the increments of the Brownian motion. 
\begin{proof}[Proof of $\mathrm{Ito}(\theta,N)$]
Let $(\Omega,\cF,(\cF_t)_{t\ge 0},\PP, (B_t)_{t\ge 0}, (X_t)_{t\ge 0})$ a $KS(\theta,N)$-process. 

\vip

Set $(D_t)_{t\ge 0} := (R_{\ig 1,N\id}(X_t))_{t\ge 0}$ which is a squared Bessel process with dimension $d_{\theta,N}(N)$. We define the following sequence of stopping time $(\tilde T^\eta_k)_{k\ge 0}$ and $( T^\eta_k)_{k\ge 1}$ for all $\eta>0$ by induction by setting $\tilde T^\eta_0 = 0$ and for all $k\ge 0$,
\begin{gather*}
    T^\eta_{k+1} = \inf \{ t\ge \tilde T^\eta_{k} :  D_t \ge \eta \},\\
     \tilde T^\eta_{k+1} = \inf \{ t\ge  T^\eta_{k+1} :  D_t =0 \}. 
\end{gather*}
According to the strong Markov property proven in Section \ref{sec:strongMarkov} and Proposition \ref{prop:girsanoveta}-(ii)-(b), for all $\varphi \in C_c^\infty ((\rr^2)^N)$ such that $\sup_{\alpha \in (0,1]}\sup_{x'\in (\rr^2)^N}|\cL_\alpha \varphi (x')| <\infty$, all $x\in (\rr^2)^N$, all $\eta >0$, all $k\ge 0$, 
\begin{align*}
 \mbox{for all } t\in [ T^\eta_k, \tilde T^\eta_k], \quad \varphi (X_t) = \varphi (X_{T^\eta_k}) + \int_{T^\eta_k}^t \nabla \varphi (X_s)\cdot \dd B_s + \int_{T^\eta_k}^t \cL^X\varphi (X_s)\dd s.
\end{align*}
Fix $t\ge 0$. By summing we find
\begin{align*}
    \varphi (X_t) = \varphi (X_0) + \int_0^t \nabla \varphi (X_s)\cdot \dd B_s + \int_0^t\cL^X\varphi (X_s)\dd s + R^{1,\eta}_t + R^{2,\eta}_t + R^{3,\eta}_t,
\end{align*}
where 
\begin{gather*}
    R^{1,\eta}_t = \sum_{k :\tilde T^\eta_k \le t} \varphi (X_{t\wedge T^\eta_{k+1}}) - \varphi ( X_{t\wedge \tilde T^\eta_{k}}) ,\\
    R^{2,\eta}_t = - \sum_{k : \tilde T^\eta_k \le t} \int_{t\wedge \tilde T^\eta_k}^{t\wedge T^\eta_{k+1}} \nabla \varphi (X_s) \cdot \dd B_s,\\
    R^{3,\eta}_t = - \sum_{k : \tilde T^\eta_k \le t} \int_{t\wedge \tilde T^\eta_k}^{t\wedge T^\eta_{k+1}} \cL^X\varphi (X_s)\dd s. 
\end{gather*}
Since $\cL^X\varphi$ is bounded by hypothesis, it is clear that $R^{3,\eta}_t \to 0$ a.s. as $\eta\to 0$ since 
\begin{align}\label{eq:intindicpetit}
\sum_{k : \tilde T^\eta_k \le t} \int_{t\wedge \tilde T^\eta_k}^{t\wedge T^\eta_{k+1}}\dd s \le \int_0^t \indiq_{D_s\le \eta}\dd s \to 0 \quad \mbox{ a.s. as }\eta \to 0,
\end{align}
according to the dominated convergence theorem, since $\int_0^t \indiq_{D_s = 0}\dd s = 0$ a.s. by hypothesis. Moreover, since $\varphi$ is smooth and using the Ito's isometry, we get that $R^{2,\eta}_t \to 0$ in $L^2$ as $\eta \to 0$. It only remains to show that $R^{1,\eta}_t \to 0$ in $L^2$ as $\eta \to 0$ in order to conclude. However, setting $(M_t)_{t\ge 0} := (S_{\ig 1,N\id}(X_t))_{t\ge 0}$, for all $k\ge 1$, we can write $\indiq_{\tilde T^\eta_{k}\le t}  (\varphi (X_{\tilde T^\eta_k}) - \varphi ( X_{ T^\eta_{k+1}\wedge t}))$ as
\begin{align*}
  \indiq_{\tilde T^\eta_{k}\le t}( \varphi (\gamma (M_{\tilde T^\eta_k})) -  \varphi (\gamma (M_{ T^\eta_{k+1}\wedge t}) + \eta^{1/2}\Theta_k)),
\end{align*}
where $\Theta_k = \Phi(X_{T^\eta_{k+1}\wedge t})$.

Summing in $k$ and using Taylor expansion with integral remainder and the fact that $\varphi \in C^\infty_c((\rr^2)^N)$, we get 
\begin{align}\label{eq:RegalS}
    (R^{1,\eta}_t)^2\le C(S^{1,\eta} + S^{2,\eta} + S^{3,\eta})
\end{align}
where
\begin{gather*}
    S^{1,\eta} = \Big(\sum_{k \ge 0} \indiq_{\tilde T_k^\eta\le t}\eta^{1/2} \Theta_k  \cdot \nabla \varphi (\gamma (M_{T^\eta_{k+1}\wedge t})) \Big)^2,\\
    S^{2,\eta}=\Big(\sum_{k \ge 0} \indiq_{\tilde T_k^\eta\le t} \varphi (\gamma ( M_{\tilde T^\eta_k})) - \varphi ( \gamma (M_{T^\eta_{k+1}\wedge t}))\Big)^2 ,\\
    S^{3,\eta} = \eta^2 |\cN^\eta_t|^2 ,
\end{gather*}
where $\cN^\eta_t = \{ k\ge 0 : \tilde T^\eta_k \le t\}$. The sequence $(\Theta_k)_{k\ge 0}$ is i.i.d., independent from $(M_t,D_t)_{t\ge 0}$ and $Z^{-1}\beta^N_\theta$-distributed thanks to Lemma \ref{lem:symdesangles} and $\mathrm{Decomposition}(\theta,N)$. In particular, for all $k\ge 0$, the expectation of $\Theta_k$ is $0$. We thus have the existence of a positive constant $C$ such that
\begin{gather}\label{eq:boundS1}
    \E[S^{1,\eta}] = \eta\sum_{k\ge 0} \E [\indiq_{\tilde T^\eta_k\le t}(\Theta_k  \cdot \nabla \varphi (\gamma (M_{T^\eta_{k+1}\wedge t})))^2]\le C\|\nabla \varphi\|_\infty \eta \E[|\cN^\eta_t|]. 
\end{gather}
Using the Ito's formula, one gets 
\begin{align*}
    \E[S^{2,\eta}] =  \E \Big[\Big(\sum_{k : \tilde T^\eta_k \le t} \int_{\tilde T^\eta_k}^{T^\eta_{k+1}\wedge t} \frac1N \sum_{i=1}^N\nabla_{x^i} \varphi (\gamma (M_s)) \cdot \dd M_s + \frac1N\int _{\tilde T^\eta_k}^{T^\eta_{k+1}\wedge t} \Delta \varphi (\gamma (M_s)) \dd s\Big)^2 \Big].
\end{align*}
Thanks to the Ito's isometry and the fact that $(x+y)^2\le 2x^2+2y^2$ for all $x,y\in \rr$, one obtains the existence of $C\in (0,\infty)$ such that
\begin{align*}\label{eq:boundS2}
    \E[S^{2,\eta}] \le C\Big( \|\nabla \varphi\|_\infty\E \Big[\sum_{k : \tilde T^\eta_k \le t}\int_{\tilde T^\eta_k}^{T^\eta_{k+1}\wedge t} \dd s\Big] +  \|\Delta \varphi\|_\infty \E \Big[\Big(\sum_{k : \tilde T^\eta_k \le t}\int_{\tilde T^\eta_k}^{T^\eta_{k+1}\wedge t} \dd s\Big)^2\Big] \Big].
\end{align*}
According to \eqref{eq:intindicpetit} we deduce that $\E[S^{2,\eta}]\to 0$ as $\eta \to 0$. According to \eqref{eq:RegalS}, is sufficient to conclude to show that $\eta^2\E[|\cN^\eta_t|^2] \to 0$ as $\eta \to 0$ since by the Cauchy-Schwarz inequality and \eqref{eq:boundS1}, 
$$
\E[ S^{1,\eta} + S^{3,\eta}] \le C\Big( \eta^2\E[|\cN^\eta_t|^2] + \Big( \eta^2\E[|\cN^\eta_t|^2]\Big)^{1/2}\Big).
$$

\vip

Using the Ito's formula, we get
\begin{align*}
  \eta |\cN^\eta_t|\le  \sum_{k : \tilde T^\eta_k \le t} D_{T^\eta_{k+1}} - D_{\tilde T^\eta_k} = 2 \sum_{k : \tilde T^\eta_k \le t} \int_{\tilde T^\eta_k}^{T^\eta_{k+1}} D_s^{1/2} \dd W_s + d_{\theta,N}(N) \int_{\tilde T^\eta_k}^{T^\eta_{k+1}}\dd s,
\end{align*}
where $(W_t)_{t\ge 0}$ is a standard Brownian motion. Thus we get, using that $(x+y)\le 2x^2 + 2y^2$ for all $x,y\in \rr$ and the Ito's isometry,
\begin{align*}
    \E[\eta^2 |\cN^\eta_t|^2]\le 8 \E \Big[ \sum_{k : \tilde T^\eta_k \le t} \int_{\tilde T^\eta_k}^{T^\eta_{k+1}} D_s \dd s \Big] + 2(d_{\theta,N}(N))^2 \E \Big[ \Big(\sum_{k : \tilde T^\eta_k \le t} \int_{\tilde T^\eta_k}^{T^\eta_{k+1}} \dd s \Big)^2\Big],
\end{align*}
which tends to $0$ as $\eta \to 0$ according to \eqref{eq:intindicpetit} and the dominated convergence theorem.
\end{proof}

\section{Dirichlet form}\label{sec:dirichlet}
In this section, we prove the association between the process we built in Section \ref{sec:constructionKS} and the Dirichlet form described in Theorem \ref{mainresult}-(i). 
\begin{proof}[Proof of the association with the Dirichlet form.]
We define for all $\varphi : (\rr^2)^N\to \rr_+$ measurable, all $t\ge 0$ and all $x\in (\rr^2)^N$,
$$
P^N_t\varphi (x) =E^N_x[\varphi (X_t)]. 
$$

\vip

{\it Step 1.} We show that for all $\varphi,\psi : (\rr^2)^N \to \rr_+$ measurable, 
$$
\int_{(\rr^2)^N} P^N_t\varphi (x) \psi (x) \mu^N_\theta (\dd x) = \int_{(\rr^2)^N} \varphi (x) P^N_t\psi (x)\mu^N_\theta (\dd x).
$$
Let $\varphi,\psi : (\rr^2)^N \to \rr_+$ measurable and $t>0$. We first have
$$
\int_{(\rr^2)^N} \tilde P^N_t\varphi (x) \psi (x) \mu^N_\theta (\dd x) = \int_{(\rr^2)^N} \varphi (x) \tilde P^N_t\psi (x)\mu^N_\theta (\dd x),
$$
where for all $f :(\rr^2)^N\to \rr_+$ measurable, $\tilde P^N_tf (x) = E^N_x[f (X^N_t)\indiq_{t<\tau_N}]$, according to Proposition \ref{prop:dirichlettué}. It only remains to show that
$$
\int_{(\rr^2)^N}  E^N_x[\varphi (X_t)\indiq_{t\ge \tau_N}] \psi (x) \mu^N_\theta (\dd x) = \int_{(\rr^2)^N} \varphi (x) E^N_x[\psi (X_t)\indiq_{t\ge \tau_N}]\mu^N_\theta (\dd x).
$$
According to Lemma \ref{lem:symdesangles}, we have
\begin{align*}
    \int_{(\rr^2)^N}  E^N_x[\varphi (X_t)\indiq_{t\ge \tau_N}]& \psi (x) \mu^N_\theta (\dd x) \\
    &= \int_{(\rr^2)^N} \int_{\sS_N} E^N_x[\varphi (\gamma (M_t) + D_t^{1/2}u')\indiq_{t\ge \tau_N}] \psi (x) Z^{-1}\beta^N_\theta(\dd u')\mu^N_\theta (\dd x).
\end{align*}
Using \cite[Lemma A.2 page 56]{ft}, we get
\begin{align*}
   &\int_{(\rr^2)^N}  E^N_x[\varphi (X_t)\indiq_{t\ge \tau_N}] \psi (x) \mu^N_\theta (\dd x) \\
   &= \frac12\int_{\rr^2}\int_0^\infty \int_{\sS_N}\int_{\sS_N} \E_{m,r}[\varphi (\gamma (M_t) + D_t^{1/2} u')  \indiq_{t\ge \tau_0}] \psi (\gamma (m) + \sqrt{r}u)r^\nu\beta^N_\theta (\dd u) \beta^N_\theta (\dd u') \dd m r\dd r,
\end{align*}
where $\nu = d_{\theta,N}(N)/2-1$ and $\E_{m,r}$ is the expectation with respect to the process $(M_t,D_t)_{t\ge 0}$ starting from $(m,r)$. We set for all $u\in \sS_N$, 
\begin{gather*}
    f_u : (m,r)\in \rr^2\times\rr_+ \mapsto \varphi (\gamma (m) + \sqrt r u),\\
    g_u : (m,r)\in \rr^2\times\rr_+ \mapsto \psi (\gamma (m) + \sqrt r u).
\end{gather*} 
We see that in order to conclude, it is sufficient to show that for all $u,u'\in \sS_N$, 
\begin{align*}
\int_{\rr^2}\int_0^\infty&  \E_{m,r}[f_u (M_t, D_t) \indiq_{t\ge \tau_0}] g_{u'}(m,r)r^\nu \dd m r\dd r \\
&= \int_{\rr^2}\int_0^\infty   f_u(m,r)\E_{m,r}[g_{u'}(M_t,D_t) \indiq_{t\ge \tau_0}] r^\nu \dd m r\dd r.
\end{align*}
However this is immediate since both semigroup of $(M_t,D_t)_{t\ge 0}$ and $(M_t,D_t)_{t\ge 0}$ killed at $\tau_0$ are symmetric with respect to $r^\nu \dd m \dd r$ (see for example \cite[Step 1 of Proposition 10, Lemma B.5, Lemma B.6]{ft}) and for all $u,u'\in \sS_N$,
\begin{align*}
    \int_{\rr^2}\int_0^\infty  &\E_{m,r}[f_u (M_t, D_t) \indiq_{t\ge \tau_0}]g_{u'}(m,r)r^\nu \dd m r\dd r\\
    =& \int_{\rr^2}\int_0^\infty  \E_{m,r}[f_u (M_t, D_t) ] g_{u'}(m,r)r^\nu \dd m r\dd r\\
    &- \int_{\rr^2}\int_0^\infty  \E_{m,r}[f_u (M_t, D_t) \indiq_{t< \tau_0}] g_{u'}(m,r)r^\nu \dd m r\dd r.
\end{align*}

{\it Step 2.} According to the strong Markov property (see Section \ref{sec:strongMarkov}) and the fact that $(X^N_t)_{t\ge 0}$ is a.s. continuous by definition, Step 1 implies that $\mM^{X,N,\theta}$ is associated with a Dirichlet space $(\cE,\cF)$ on $L^2((\rr^2)^N,\mu^N_\theta)$. We show that $C^\infty_c(E_N)+ \cC_{\ig 1,N\id} \subset \cF$ and for all $\varphi \in C^\infty_c(E_N)\cup \cC_{\ig 1,N\id}$,
$$
\cE (\varphi, \varphi) = \frac12 \int_{(\rr^2)^N} \|\nabla \varphi (x)\|^2 \mu^N_\theta(x)\dd x.
$$
According to \cite[Lemma 1.3.4 page 23]{f}, it is sufficient to show that for all $\varphi \in C^\infty_c(E_N) \cup \cC_{\ig 1,N\id}$,
\begin{align*}
    \int_{(\rr^2)^N} \frac{\varphi(x) - P^N_t\varphi(x)}{t}\varphi(x) \mu^N_\theta(x)\dd x \to \frac12\int_{(\rr^2)^N} \|\nabla\varphi(x)\|^2 \mu^N_\theta(x)\dd x,
\end{align*}
as $t\to 0$.

\vip

{\it Step 2.1} We prove the statement if $\varphi \in C_c^\infty (E_N)$. A straightforward consequence of Proposition \ref{prop:dirichlettué} and \cite[Lemma 1.3.4 page 23]{f} gives that
\begin{align*}
    \int_{(\rr^2)^N} \frac{\varphi(x) - \tilde P^N_t\varphi(x)}{t}\varphi(x) \mu^N_\theta(x)\dd x \to \frac12\int_{E_N} \|\nabla\varphi(x)\|^2 \mu^N_\theta(x)\dd x = \frac12\int_{(\rr^2)^N} \|\nabla\varphi(x)\|^2 \mu^N_\theta(x)\dd x.
\end{align*}
It thus only remains to show that
\begin{align*}
    I_t := \int_{(\rr^2)^N} \frac{E^N_x[\varphi(X_t)\indiq_{t\ge \tau_N}]}{t}\varphi(x) \mu^N_\theta(x)\dd x \to 0 \quad \mbox{ as } t\to 0.
\end{align*}
Since $\varphi \in C_c^\infty (E_N)$, there exists $\eta >0$ such that $\mathrm{Supp} \varphi \subset E^\eta_N$. This implies that
\begin{align*}
     |I_t| \le \|\varphi\|_\infty\int_{(\rr^2)^N} \frac{\E_{R_{\ig 1,N\id}(x)}[ t\ge \tau^2_\eta]}{t} |\varphi (x)|\mu^N_\theta(x)\dd x,
\end{align*}
where $\E_{R_{\ig 1,N\id}(x)}$ is the law of a squared Bessel process with dimension $\delta$, starting at $R_{\ig 1,N\id}(x)$ and $\tau_0$ its first hitting time of $0$ and $\tau^2_\eta$ its first hitting time of $\eta$ after time $\tau_0$. By comparison theorem (see \cite[Theorem 2.7 page 394]{ry}) and the strong Markov property at time $\tau_0$, we get
\begin{align}\label{eq:majoI}
|I_t| \le \|\varphi\|_\infty \frac{\PP_0[t\ge \tau_\eta]}{t} \int_{(\rr^2)^N}  |\varphi|\mu^N_\theta,
\end{align}
where $\tau_\eta$ is the first hitting time of $\eta$. This together with Lemma \ref{lem:estimetuagebessel} implies that $|I_t|\to 0$ as $t\to0$.

\vip

{\it Step 2.2} Let us treat the case where $\varphi \in \cC_{\ig 1,N\id}$. This means that there exists $\psi \in C^\infty_c(\rr^2\times\rr_+)$ such that $\varphi (x) = \psi (S_{\ig 1,N\id}(x),R_{\ig 1,N\id}(x))$ for all $x\in (\rr^2)^N$. Using \cite[Lemma A.2]{ft} we get
\begin{align}\label{eq:convdirichlet1}
    \int_{(\rr^2)^N} \frac{\varphi(x) - P^N_t\varphi(x)}{t} \varphi (x) \mu^N_\theta (x) \dd x &= \frac12 \int_{\rr^2}\int_0^\infty \int_{\sS_N} \frac{\psi (m,r) - \E_{m,r}[\psi (M_t,D_t)]}{t} r^\nu \dd m \dd r\beta^N_\theta (\dd u)\notag\\
    &=\frac{\beta^N_\theta (\sS_N)}{2}\int_{\rr^2}\int_0^\infty  \frac{\psi (m,r) - \E_{m,r}[\psi (M_t,D_t)]}{t} r^\nu \dd m \dd r.
\end{align}
We can now apply \cite[Lemma 1.3.4 page 23]{f} again to get
\begin{align*}
    \int_{\rr^2}\int_0^\infty& \frac{\psi(m,r) - \E_{(m,r)}[\psi (M_t,D_t)]}{t} \psi (m,r) r^\nu \dd m \dd r \\
    &\to \frac12\int_{\rr^2}\int_0^\infty \Big(\frac1N\|\nabla_m \psi (m,r)\|^2 + 4r|\partial_r\psi (m,r)|^2 \Big) r^\nu \dd m\dd r,
\end{align*}
as $t\to 0$. This together with \eqref{eq:convdirichlet1} implies that $\varphi \in \cF$ and 
\begin{align}\label{eq:calculdirichlet}
  \cE (\varphi,\varphi) &= \frac{\beta^N_\theta(\sS_N)}{4}\int_{\rr^2}\int_0^\infty \Big(\frac1N \|\nabla_m \psi (m,r)\|^2 +4 r|\partial_r\psi (m,r)|^2 \Big) r^\nu \dd m\dd r\notag \\
  &= \frac{1}{4}\int_{\rr^2}\int_0^\infty\int_{\sS_N} \Big(\frac1N\|\nabla_m \psi (m,r)\|^2 + 4r|\partial_r\psi (m,r)|^2 \Big) r^\nu \dd m\dd r \beta^N_\theta (\dd u).
\end{align}
However, for all $x\in (\rr^2)^N$, $\varphi (x) = \psi (S_{\ig 1,N\id}(x), R_{\ig 1,N\id}(x))$ so that
\begin{align*}
    \nabla \varphi (x) = \frac1N\gamma(\nabla_m\psi(S_{\ig 1,N\id} (x),R_{\ig 1,N\id}(x))) +\partial_r \psi (S_{\ig 1,N\id}(x),R_{\ig 1,N\id}(x)) \nabla R_{\ig 1,N\id}(x).
\end{align*}
Since for all $y\in \rr^2$, $\gamma (y)\cdot \nabla R_{\ig 1,N\id}(x) = 0$ by symmetry, we get
\begin{align*}
    \|\nabla \varphi (x)\|^2 = \frac1N \| \nabla_m\psi ((S_{\ig 1,N\id} (x),R_{\ig 1,N\id}(x))\|^2 + 4R_{\ig 1,N\id}(x)|\partial_r \psi (S_{\ig 1,N\id}(x),R_{\ig 1,N\id}(x))|^2,
\end{align*}
since $\|\nabla R_{\ig 1,N\id}(x)\|^2 = 4 R_{\ig 1,N\id}(x)$. By \eqref{eq:calculdirichlet} and a reverse application of \cite[Lemma A.2]{ft} we conclude this step.

\vip

{\it Step 3.} We conclude. According to Step 2, we have $C_c^\infty (E_N)+\cC_{\ig 1,N\id} \subset \cF$. According to Lemma \ref{lem:core}, we get that $C_c^\infty (E_N)+\cC_{\ig 1,N\id}$ is dense in $C_c^\infty ((\rr^2)^N)$ for the $\cE^N_{X,1}$ norm, see Theorem \ref{mainresult}. Since $\cE^N_{X,1}$ and $\cE_1 := \cE + (\cdot ,\cdot )_{L^2((\rr^2)^N,\mu^N_\theta)}$ coincide on $C_c^\infty (E_N)+\cC_{\ig 1,N\id}$, and $\cF$ is complete for the $\cE_1$ norm, we get that 
\begin{gather*}
\overline{C^\infty_c((\rr^2)^N)}^{\cE^N_{X,1}} \subset \cF,\\
\mbox{ for all } \varphi \in C_c^\infty ((\rr^2)^N), \quad \cE(\varphi, \varphi) = \frac12 \int_{(\rr^2)^N} \|\nabla \varphi (x)\|^2 \mu^N_\theta (x) \dd x. 
\end{gather*}
It is thus sufficient to show that $\cF\subset \overline{C^\infty_c((\rr^2)^N)}^{\cE^N_{X,1}}$ which we now do. We set for all $\cK \in \gK_N$, all $\varphi = \psi ( (S_K(x),R_K(x))_{K\in \cK}) \in \cC_\cK$, all $x\in (\rr^2)^N$,
\begin{align}\label{eq:operateurcluster}
\cL\varphi (x) = \sum_{K\in \cK} \cL^{S_K}\varphi (x) + \cL^{R_K}\varphi (x),
\end{align}
where for all $K\in \cK$,
\begin{align*}
\cL^{S_K}\varphi (x)=& \frac{1}{2|K|}\Delta_{K,1} \psi ((S_K(x),R_K(x))_{K\in \cK}) \\
&-\frac{\theta}{N|K|} \sum_{i\in K} \sum_{j\notin K} \frac{x^i-x^j}{\|x^i-x^j\|^2} \cdot \nabla_{K,1}\psi ((S_K(x),R_K(x))_{K\in \cK}),\\
\cL^{R_K}\varphi (x) =& 2R_K(x)\partial_{K,2}^2\psi ((S_K(x),R_K(x))_{K\in \cK}) \\
&+ \Big(d_{\theta,N}(|K|)  - \frac{2\theta}{N}\sum_{i\in K}\sum_{j\notin K}(x^i-S_K(x))\cdot \frac{x^i-x^j}{\|x^i-x^j\|^2}\Big)\partial_{K,2}\psi ((S_K(x),R_K(x))_{K\in \cK}),
\end{align*}
where $\Delta_{K,1}$, $\nabla_{K,1}$ act on the variable corresponding to $S_K$ and $\partial_{K,1}$ and $\partial^2_{K,2}$ act on the variable corresponding to $R_K$.
The differential operator defined in \eqref{eq:operateurcluster} is in fact the $KS(\theta,N)$-generator applied on particular test function which group particles according to clusters by replacing the variable of each subgroup by the variables $S_K$ and $R_K$. Morally speaking, this allows us to define the operator when collision occurs in a configuration. 

\vip

Observe that $(\cL, \sum_{\cK\in \gK_N}\cC_\cK)$ is crucially a densely defined symmetric diffusion operator on $L^2((\rr^2)^N,\mu^N_\theta)$ in the sense of \cite{Eberle} because $\sum_{\cK\in \gK_N}\cC_\cK$ is dense in $C_c^\infty ((\rr^2)^N)$ for both the uniform convergence topology and $\cE^N_{X,1}$. 

\vip

Let $(\tilde \cL,\tilde \cD)$ the generator associated with $(\cE,\cF)$ in the sense of \cite[Theorem 1.3.1 page 20]{f}. It only remains to prove that

\vip

(a) $(\tilde \cL,\tilde \cD)$ extends $(\cL,\sum_{\cK \in \gK_N} \cC_\cK)$ in the sense that for all $\varphi \in \sum_{\cK\in \gK_N}\cC_\cK$, $\varphi \in \tilde \cD$ and $\tilde \cL \varphi = \cL \varphi$,

\vip

(b) $\overline{W^{1,2}\cap L^\infty ((\rr^2)^N)}^{\cE^N_{X,1}} \subset \overline{C_c^\infty ((\rr^2)^N)}^{\cE^N_{X,1}}$ where $L^\infty ((\rr^2)^N)$ is the set of bounded function defined on $(\rr^2)^N$ and
\begin{align*}
W^{1,2} =& \Big\{ u \in L^2((\rr^2)^N, \mu^N_\theta)) : \mbox{for all } i\in \ig 1,N\id, \mbox{ all } \e\in \{1,2\}, \exists v^{i,\e}\in L^2( (\rr^2)^N, \mu^N_\theta)) \mbox{ such that} \\&
\mbox{ for all }\varphi \in C_c^\infty ((\rr^2)^N),
\int_{(\rr^2)^N} v^{i,\e}(x) \varphi(x)  \mu^N_\theta(x) \dd x = - \int_{(\rr^2)^N} u(x)\partial_{x^i,\e}(\varphi \mu^N_\theta)(x) \dd x \Big\}.
\end{align*}

This will be sufficient to conclude since according to (a) and \cite[Theorem 3.1 page 114]{Eberle} we get $\cF \subset \overline{W^{1,2}\cap L^\infty ((\rr^2)^N)}^{\cE^N_{X,1}} $. Combined with (b) gives the result. We proceed to the proofs of (a) and (b) split in several steps.

\vip

{\it Step 3.1} We prove (a). Fix $\cK\in \gK_N$ and $\varphi = \psi ((S_K,R_K)_{K\in \cK}) \in \cC_\cK$. We need to show that $\varphi \in \tilde \cD$ and $\tilde \cL\varphi = \cL\varphi$. According to \cite[Appendix B.]{ft}, it is sufficient to prove that 
\begin{align*}
    \frac{P^N_t\varphi - \varphi }{t} \to \cL\varphi \quad \mbox{ in } L^2((\rr^2)^N,\mu^N_\theta).
\end{align*}
Observe that this is direct in the case where $\cK = \{\ig 1,N\id\}$ since $(S_{\ig 1,N\id}(X_t))_{t\ge 0}$ and $(R_{\ig 1,N\id}(X_t))_{t\ge 0}$ are respectively Brownian motion and squared Bessel process. Thus we can suppose that $\cK\ne \{ \ig 1,N\id\}$. 

\vip

We have for some constant $C>0$,
\begin{align}
    \Big\| \frac{P^N_t\varphi -\varphi}{t} - \cL\varphi \Big\|^2_{L^2(\mu^N_\theta)} \le C(I^1_t + I^2_t + I^3_t),
\end{align}
where
\begin{gather*}
    I^1_t = \int_{(\rr^2)^N} \Big| \frac{\E^{X,N,\theta}_x[\varphi(X_t)\indiq_{t<\tau_N}] - \varphi (x)}{t} - \cL\varphi(x)\E^{X,N,\theta}_x[(t\wedge\tau_N)/t]\Big|^2 \mu^N_\theta(\dd x),\\
    I^2_t = \int_{(\rr^2)^N}  \frac{\E^{X,N,\theta}_x[\varphi (X_t) \indiq_{t\ge \tau_N}]^2 }{t^2}\mu^N_\theta(\dd x),\\
    I^3_t = \int_{(\rr^2)^N} ( \cL\varphi (x) )^2\E^{X,N,\theta}_x[(t-t\wedge\tau_N)/t]^2 \mu^N_\theta(\dd x).
\end{gather*}
We prove that $I^1_t \to 0$ as $t\to 0$. According to Remark \ref{rk:dimensions}-(iii), $(X_t)_{t\in [0,\tau_N)}$ solves the SDE \eqref{eq:EDS} in the classical sense. Using the Ito's formula, we get for all $t\in [0,\tau_N]$,
\begin{align*}
    \varphi (X_t) =\varphi (x) +  N_t + \int_0^t \cL\varphi (X_s)\dd s,
\end{align*}
where $(N_t)_{t\ge 0}$ is a $L^2$-martingale. Thus we have
\begin{align*}
\E^{X,N,\theta}_x[\varphi(X_{t\wedge \tau_N})] = \varphi (x)  + \E^{X,N,\theta}\Big[\int_0^{t\wedge \tau_N} \cL\varphi (X_s)\dd s \Big].
\end{align*}
Since $\varphi \in \cC_\cK$, $\cL\varphi$ is bounded and continuous. The fact that $\varphi (X_{\tau_N}) = 0$ a.s. and the Jensen inequality gives that
\begin{align*}
    \Big| \frac{\E^{X,N,\theta}_x[(\varphi(X_t))\indiq_{t<\tau_N}] - \varphi (x)}{t} - \cL\varphi(x)\E^{X,N,\theta}_x[(t\wedge\tau_N)/t]\Big|^2 \le \frac1t \E^{X,N,\theta}_x\Big[\int_0^{t\wedge \tau_N} |\cL\varphi (X_s) - \cL\varphi(x)|^2 \dd s\Big].
\end{align*}
Thus, we get
$$
I^1_t \le \frac{1}{t} \int_0^t \int_{(\rr^2)^N} \E^{X,N,\theta}_x[|\cL\varphi (X_s) - \cL\varphi(x)|^2 ]\mu^N_\theta (\dd x)\dd s.
$$
Expanding the square, using the fact that $P^N_t$ is $\mu^N_\theta$-reversible and for all $s\ge 0$, $P^N_s1 = 1$, we obtain
\begin{align*}
I^1_t \le& \frac{1}{t} \int_0^t \int_{(\rr^2)^N} \Big(\E^{X,N,\theta}_x[|\cL\varphi (X_s)|^2] - 2P^N_s\cL\varphi (x)\cL\varphi(x) + |\cL\varphi (x)|^2\Big)\mu^N_\theta (\dd x)\dd s\\
=& \frac{2}{t} \int_0^t \int_{(\rr^2)^N} \Big(|\cL\varphi (x)|^2 - P^N_s\cL\varphi (x)\cL\varphi(x)\Big)  \mu^N_\theta (\dd x)\dd s.
\end{align*}
According to \cite[Lemma 1.4.3 page 30]{f}, since $\cL\varphi$ is smooth and compactly supported, $P^N_s\cL\varphi \to \cL\varphi$ in $L^2((\rr^2)^N, \mu^N_\theta)$ as $s\to 0$, implying that $I^1_t \to 0$ as $t\to 0$. 

\vip

We now prove that $I^2_t \to 0$ as $t\to 0$. By definition of $\varphi$, there exists $\eta>0$ such that $\mathrm{Supp}\varphi \subset E^\eta_N$. Thus, by the strong Markov property,
\begin{align*}
    I^2_t \le  \frac{\|\varphi\|_\infty}{t^2} \PP_0[t\ge \tau_{\eta}] \mu^N_\theta(\mathrm{Supp}\varphi),
\end{align*}
where $\PP_0$ is the law of a squared Bessel process with dimension $d_{\theta,N}(N)$ starting at $0$ and $\tau_\eta$ is the first hitting time of $\eta$ of $(R_{\ig 1,N\id}(X_t))_{t\ge 0}$. We deduce that $I^2_t \to 0$ as $t\to 0$ according to Lemma \ref{lem:estimetuagebessel}. 

\vip

Finally, we prove that $I^3_t\to 0$ as $t\to 0$. We Immediately have for $t\in (0,1]$, thanks to \cite[Lemma A.2]{ft},
\begin{align*}
    I^3_t \le C_\varphi \int_{a}^b\PP_r[t\ge \tau_0]^2r^\nu \dd r,
\end{align*}
where $\PP_r$ is the law of a squared Bessel process with dimension $d_{\theta,N}(N)$ starting at $r$, $\tau_0$ its first hitting time of $0$ and $C_\varphi$ is a positive constant depending on $\mathrm{Supp}\varphi$, $\|\varphi\|_\infty$ and $\|\cL\varphi\|_\infty$, and $a,b$ are also positive constant depending on $\mathrm{Supp}\varphi$. We conclude according to Lemma \ref{lem:estimetuagebessel}. 

\vip

{\it Step 3.2} We end the proof by showing (b). According to \cite[Theorem 3 page 31]{Kilpelainen1997}, it is sufficient to show that the weight $\mu^N_\theta$ is $2$-admissible, see \cite{Kilpelainen1997} for more details. Moreover, according to \cite[page 10]{HeinonenKilpelainenMartio1993}, this is the case if $\mu^N_\theta$ belongs to the Muckenhoupt class $A_2$, i.e.
$$
\sup_{x_0\in (\rr^2)^N, r>0} \Big( \frac{1}{|B(x_0,r)|} \int_{B(x_0,r)} \mu^N_\theta (x) \dd x \Big) \Big( \frac{1}{|B(x_0,r)|} \int_{B(x_0,r)} (\mu^N_\theta (x) )^{-1}\dd x\Big) <\infty, 
$$
which we prove now. First of all, by homogeneity of $\mu^N_\theta$, we get thanks to a change of variable that it is equivalent to show that
$$
\sup_{y_0\in (\rr^2)^N} \Big(  \int_{B(y_0,1)} \mu^N_\theta (x) \dd x \Big) \Big( \int_{B(y_0,1)} (\mu^N_\theta (x) )^{-1}\dd x\Big) <\infty.
$$
For any $y_0\in (\rr^2)^N$, we have
$$
  \int_{B(y_0,1)} (\mu^N_\theta (x) )^{-1}\dd x \le \tilde C \prod_{i\ne j} (\|y^i_0 - y^j_0\| + 1)^{\theta/N}.
$$
Thus, a change of variable gives for all $y_0 \in (\rr^2)^N$,
\begin{align*}
 \Big(  \int_{B(y_0,1)} \mu^N_\theta (x) \dd x \Big) \Big( \int_{B(y_0,1)} (\mu^N_\theta (x) )^{-1}\dd x\Big)&\le \int_{B(0,1)} \prod_{i\ne j} \Big( \frac{\|y^i_0-y^j_0\|+1}{\|y^i_0-y^j_0 + y^i-y^j\|}\Big)^{\theta/N} \dd y\\
 &\le \int_{B(0,1)} \prod_{i\ne j} \Big( 4 + \frac{4\indiq_{\|y^i_0-y^j_0\|\le 3}}{\|y^i_0-y^j_0 + y^i-y^j\|}\Big)^{\theta/N} \dd y,
\end{align*}
where in the last line we use that for every $i\ne j$,
$$
\frac{\|y^i_0-y^j_0\|+1}{\|y^i_0-y^j_0 + y^i-y^j\|} \le 4 + \frac{4\indiq_{\|y^i_0-y^j_0\|\le 3}}{\|y^i_0-y^j_0 + y^i-y^j\|},
$$
by distinguishing whether $\|y^i_0-y^j_0\|$ is smaller than $3$ or not. Lemma \ref{lem:radon} concludes the proof by distinguishing on if $\|y^i_0-y^j_0\|$ is smaller than $3$ or not for every $i\ne j \in \ig 1,N\id$. 
\end{proof}

\appendix

\section{Proofs of technical results}\label{app:estimations}

\subsection{Proof of existence and properties of $\eta$-pre-$KS$-process}\label{app:preuvepropetapreKS}

\begin{proof}[Proof of Lemma \ref{lem:etapreKS}] Let us fix $\theta >0$ and $\eta >0$.

\vip

{\it Step 1.} We first show (i). Fix $\mu \in \cP( (E_{k_0}\cap E^\eta_N)_\triangle)$.

\vip

{\it Step 1.1} We build by induction a $\eta$-$pre$-$KS_\mu(\theta,N)$ process. We set $\tilde X^N_0$ with law $\mu$, $\sigma^\eta_0=0$. Now take $k\ge 0$ and suppose that there exists a filtered probability space $(\Omega^k,\cF^k,(\cF^k_t)_{t\ge 0},\PP^k)$ on which are defined a $(\rr^2)^N$-valued Brownian motion $(B^k_t)_{t\ge 0}$ and $(\tilde X^N_t)_{t\in [0, \sigma^\eta_k]}$ a $\eta$-$pre$-$KS_x(\theta,N)$-process driven by $(B^k_t)_{t\ge 0}$ and stopped at $\sigma^\eta_k$ such that $(B^k_t)_{t\ge 0}$ is independent from $\tilde X^N_0$. We now build a $\eta$-$pre$-$KS_\mu(\theta,N)$-process stopped at $\sigma^\eta_{k+1}$.

\vip

We can suppose that $\cK^\eta_k\notin \{\triangle, \ig 1,N\id\}$ since the result is clear otherwise. According to $\cP(m)$ for all $m\in \ig 1,N-1\id$, enlarging the probability space to some $(\Omega^{k+1},\cF^{k+1}, (\cF^{k+1}_t)_{t\ge 0}, \PP^{k+1})$ if necessary, there exists a family of mutually independent processes knowing $\cF^{k+1}_{\sigma^\eta_k}$,
$$
\Big\{((\hat X^{N,i}_t)_{t\ge 0})_{i\in K}), K\in \cK^\eta_k(( \tilde X^{N}_{t\wedge \sigma^\eta_k})_{t\ge 0}) \Big\},
$$
such that for all $K\in \cK^\eta_k(\tilde X^N)$, $((\hat X^{N,i}_t)_{t\ge 0})_{i\in K})$ is a $KS(\theta |K|/N, |K|)$-process with initial configuration $(\tilde X^{N,i}_{\sigma^\eta_k})_{i\in K}$.

\vip

We set for all $t\ge 0$, $\hat X^N_t = (\hat X^{N,i}_t)_{i\in \ig 1,N\id}$ and write $\hat \sigma^\eta_1$ the time of the first change of $\eta$-configuration of $(\hat X^N_t)_{t\ge 0}$. We set
$$
\mbox{ for all } t \in [\sigma^\eta_k,\sigma^\eta_{k} + \hat \sigma^\eta_1], \quad \tilde X^N_t = \hat X^N_{t-\sigma^\eta_k}.
$$
It is then clear that $\sigma^\eta_{k+1}((\tilde X^N_t)_{t\ge 0}) = \sigma^\eta_k + \hat \sigma^\eta_1$. and that $(\tilde X^N_t)_{t\in [0, \sigma^\eta_{k+1}]}$ is a $\eta$-$pre$-$KS_\mu(\theta,N)$-process on $[0,\sigma^\eta_{k+1}]$, proving the induction.

\vip

We conclude using the fact that $\lim_{k\to \infty} \sigma^\eta_k\wedge \sigma^\eta_* = \sigma^\eta_*$, see Definition \ref{def:sigmastar} and the Kolmogorov's extension theorem.
 
\vip
 
{\it Step 1.2} We conclude this step by showing uniqueness in law. Since $\lim_{k\to \infty} \sigma^\eta_k\wedge \sigma^\eta_* = \sigma^\eta_*$ a.s., it is sufficient to show that for all $k\ge 0$, the law of a $\eta$-$pre$-$KS_\mu(\theta,N)$-process stopped at $\sigma^\eta_k\wedge \sigma^\eta_*$ is unique. We prove this by induction on $k\ge 0$. The result being clear for $k=0$, we suppose this is true for some $k\ge 0$ and show this still holds for $k+1$. According to $\cP(N-1)$-(a) concerning the uniqueness in law of $KS(\theta,k)$-processes with $k\le N-1$ and Definition \ref{def:etapreKS}, it is sufficient to prove that setting $m\in \nn^*$ and for all $i\in \ig 1,m\id$, $(X^i_t)_{t\ge 0}$ a $KS_{x_i}(\theta_i,k_i)$-process driven by a standard Brownian motion $B^i$ with $x_i \in E_{k_0(\theta_i,k_i)}$, $\theta_i >0$, $k_i\le N-1$ and $(B^1,\dots,B^m)$ being mutually independent, we have that the family $((X^1_t)_{t\ge 0}, \dots, (X^m_t)_{t\ge 0})$ is mutually independent. This is clear thanks to the uniqueness part in $\cP(N-1)$-(a) and Proposition \ref{prop:indepuniciteloi} (property $\mathrm{(U)}$ clearly holds by construction).

\vip

{\it Step 2.} We show (ii).

\vip

{\it Step 2.1.} We first show that for all $\mu\in \cP((E_{k_0}\cap E^\eta_N)_\triangle)$, $\sigma^\eta_1$ and  $\indiq_{\sigma^\eta_1<\sigma^\eta_*}$ are continuous applications for the one-point compactification topology $\tilde \PP^{\eta,N,\theta}_\mu$-a.e. and for all $t>0$, $\tilde \PP^{\eta,N,\theta}_\mu(\sigma^\eta_1 = t) =0$.

\vip

Fix $\mu\in \cP(E_{k_0})$. According to Step 1, the definition of a $\eta$-$pre$-$KS$-process, $\mathrm{Decomposition}(\theta,N)$ in $\cP(N-1)$-(a) and Lemma \ref{lem:dispersionformula}, there exists a filtered probability space $(\Omega,\cF,(\cF_t)_{t\ge 0},\PP)$ on which are defined a $\eta$-$pre$-$KS_\mu(\theta,N)$-process $(\tilde X^N_t)_{t\ge 0}$ and a $(E_{k_0})_\triangle$-valued process $(\hat X^N_t)_{t\ge 0}$ and a family of real valued process $\{ (S^K_t)_{t\ge 0}, K\in \cK^\eta_0\}\cup\{ (R^K_t)_{t\ge 0}, K\in \cK^\eta_0\}$ such that 
\begin{itemize}
\item for all $K\in \cK^\eta_0$, $((\hat X^{N,i}_t)_{t\ge 0})_{i\in K}$ is a $KS(\theta |K| /N, |K|)$-process starting at $(\tilde X^{N,i}_0)_{i\in K}$,
\item the family $\{ (S^K_t)_{t\ge 0}, K\in \cK^\eta_0\}\cup\{ (R^K_t)_{t\ge 0}, K\in \cK^\eta_0\}$ is mutually independent and constituted of Brownian motions and squared Bessel processes,
\item for all $t\in [0,\sigma^\eta_1]$, $\tilde X^N_t = \hat X^N_t$, 
\item for all $K\in \cK^\eta_0$, all $t\in [0,\sigma^\eta_1]$, $S_K(\tilde X^N_t) = S^K_t$ and $R_K(\tilde X^N_t) = R^K_t$,
\item we can write $\sigma^\eta_1 := \sigma^\eta_1((\tilde X^N_t)_{t\ge 0}) = \tau_1\wedge \tau_2\wedge \tau_3$ where
\begin{gather*}
\tau_1 = \inf \Big\{ t\ge 0 : R^K_t\ge 2\eta^{1/|K|}\mbox{ for some } K\in \cK^\eta(x)\Big\},\\
\tau_2 = \inf \Big\{ t\ge 0 : \sum_{i=1}^pR^{K_i}_t + \frac{1}{2\tilde N} \sum_{1\le i\ne j\le p} |K_i||K_j| \| S^{|K_i|}_t - S^{|K_j|}_t\|^2  \le \frac12 \eta^{1/\tilde N}\\
 \mbox{ for some subfamily } (K_i)_{i\in \ig 1,p\id} \mbox{ of } \cK^\eta_0(x) \mbox{ with } p\ge 2\Big\},\\
\tau_3 = \inf \{ t\ge 0 : \hat X^N_t = \triangle\},
\end{gather*}
with $\tilde N := \sum_{i=1}^p |K_i|$ and we adopt the convention $\inf \emptyset = \infty$,
\item we can write $\indiq_{\{\sigma^\eta_1<\sigma^\eta_*\}} = \indiq_{A_1} \indiq_{A_2}\indiq_{\tau_1\wedge \tau_2<\tau_3}$, where
\begin{gather*}
A_1 =\Big\{ \sum_{K\in \cK^\eta(x)} R^K_{\tau_1\wedge \tau_2} + \frac{1}{2N}\sum_{K\ne L\in \cK^\eta(x)} |K||L| \|S^K_{\tau_1\wedge \tau_2} - S^L_{\tau_1\wedge \tau_2}\|^2>\eta^{1/N}\Big\},\\
 A_2= \Big\{ \sum_{K\in \cK^\eta(x)} R^K_{\tau_1\wedge \tau_2} + \frac{1}{2N}\sum_{K\ne L\in \cK^\eta(x)} |K||L| \|S^K_{\tau_1\wedge \tau_2} - S^L_{\tau_1\wedge \tau_2}\|^2 \\
 + N\Big(\sum_{K\in \cK^\eta(x)} \frac{|K|}{N}S^K_{\tau_1\wedge \tau_2}\Big)^2 < \eta^{-2/N}\Big\}.
\end{gather*}
\end{itemize}
Observe that according to Lemma \ref{lem:dispersionformula}, 
\begin{gather*}
A_1 \cap \{\tau_1\wedge \tau_2<\tau_3\} = \{ R_{\ig 1,N\id}( \hat X^N_{\tau_1\wedge \tau_2}) > \eta^{1/N}\}\cap \{\tau_1\wedge \tau_2<\tau_3\},\\
    A_2 \cap \{\tau_1\wedge \tau_2<\tau_3\} = \{ \| \hat X^N_{\tau_1\wedge \tau_2}\|^2 < \eta^{-2/N}\}\cap \{\tau_1\wedge \tau_2<\tau_3\}.
\end{gather*}
It is worth emphasizing that $(\hat X^N_t)_{t\ge 0}$ is certainly not independent from the family 
$$
\{ (S^K_t)_{t\ge 0}, K\in \cK^\eta_0\}\cup\{ (R^K_t)_{t\ge 0}, K\in \cK^\eta_0\}.
$$

\vip

According to $\cP(N-1)$-(d), we know that $\tau_3$ is continuous for the one-point compactification topology $\tilde \PP^{\eta,N,\theta}_\mu$-a.e. and $\tilde \PP^{\eta,N,\theta}_\mu(\tau_3 = t) = 0$ for all $t>0$. Moreover it is classical that $\tau_1$ is continuous $\tilde \PP^{\eta,N,\theta}_\mu$-a.e since a Bessel process a.s. does not have any local extremum and $\tilde \PP^{\eta,N,\theta}_\mu(\tau_1 = t) = 0$ for all $t>0$. It is also clear that $\tau_2$ is continuous $\tilde \PP^{\eta,N,\theta}_\mu$-a.e. and that $\tilde \PP^{\eta,N,\theta}_\mu(\tau_2 = t) = 0$. Indeed, $(Z_t)_{t\ge 0} := (\sum_{i=1}^pR^{K_i}_t + \frac{1}{2N} \sum_{1\le i\ne j\le p} |K_i||K_j| \| S^{|K_i|}_t - S^{|K_j|}_t\|^2)_{t\ge 0}$ does not have local extremum a.s that equals some positive value. (For example, one can check that thanks to the Ito's formula and the Dubins-Schwarz theorem, we can write a.s., 
$$
\mbox{ for all } t\in [0,\sigma^\eta_1], \quad Z_t = \tilde B_{F(t)} + A_t,
$$
where $(\tilde B_t)_{t\ge 0}$ is a Brownian motion and $t\mapsto A_t$ is $C^1$ and $F$ an increasing function. The result is clear since a Brownian motion does not have local extremum a.s.)

\vip

Since we have already proven that $\tau_1\wedge \tau_2$ and $\tau_3$ were continuous $\tilde \PP^{\eta,N,\theta}_\mu$-a.e., it remains to show that $P^1 = P^2 = 0$ where
\begin{gather*}
P^1 =\PP\Big(\sum_{K\in \cK^\eta_0} R^K_{\tau_1\wedge \tau_2} + \frac{1}{2N}\sum_{K\ne L\in \cK^\eta_0}|K||L| \|S^K_{\tau_1\wedge \tau_2} - S^L_{\tau_1\wedge \tau_2}\|^2 = \eta^{1/N}\Big),\\
P^2=\PP\Big(\sum_{K\in \cK^\eta_0} R^K_{\tau_1\wedge \tau_2} + \frac{1}{2N}\sum_{K\ne L\in \cK^\eta_0}|K||L| \|S^K_{\tau_1\wedge \tau_2} - S^L_{\tau_1\wedge \tau_2}\|^2 + N\Big(\sum_{K\in \cK^\eta_0} \frac{|K|}{N}S^K_{\tau_1\wedge \tau_2}\Big)^2 = \eta^{-2/N}\Big),\\
 P^3:=\PP (\tau_1\wedge \tau_2 = \tau_3) = 0.
\end{gather*}
We have $P^1 \le P^{1,1} + P^{1,2}$ and $P^2 \le P^{2,1} + P^{2,2}$ with
\begin{gather*}
P^{1,1} = \PP\Big(\sum_{K\in \cK^\eta_0} R^K_{\tau_1} + \frac{1}{2N}\sum_{K\ne L\in \cK^\eta_0}|K||L| \|S^K_{\tau_1} - S^L_{\tau_1}\|^2 = \eta^{1/N}\Big)\\
P^{1,2} = \PP\Big(\sum_{K\in \cK^\eta_0} R^K_{\tau_2} + \frac{1}{2N}\sum_{K\ne L\in \cK^\eta_0}|K||L| \|S^K_{ \tau_2} - S^L_{ \tau_2}\|^2= \eta^{1/N}\Big),\\
P^{2,1} = \PP\Big(\sum_{K\in \cK^\eta_0} R^K_{\tau_1} + \frac{1}{2N}\sum_{K\ne L\in \cK^\eta_0}|K||L| \|S^K_{\tau_1} - S^L_{\tau_1}\|^2 + N\Big(\sum_{K\in \cK^\eta_0} \frac{|K|}{N}S^K_{\tau_1}\Big)^2 = \eta^{-2/N}\Big),\\
P^{2,2} = \PP\Big(\sum_{K\in \cK^\eta_0} R^K_{\tau_2} + \frac{1}{2N}\sum_{K\ne L\in \cK^\eta_0}|K||L| \|S^K_{ \tau_2} - S^L_{ \tau_2}\|^2 + N\Big(\sum_{K\in \cK^\eta_0} \frac{|K|}{N}S^K_{\tau_2}\Big)^2 = \eta^{-2/N}\Big).
\end{gather*}
Clearly $P^{1,1} + P^{2,1} = 0$ since $\tau_1>0$ a.s., the family $\{(S^K_t)_{t\ge 0}, K\in \cK^\eta_0\}\cup\{(R^K_t)_{t\ge 0}, K\in \cK^\eta_0\}$ is mutually independent and, for all $t>0$, all $K\in \cK^\eta_0$, the law of $S^K_t$ has a density. Then, $P^{2,2} = 0$ since $\tau_2 >0$ a.s., for all $t>0$ the law of $\sum_{K\in \cK^\eta_0} |K| S^K_t$ has a density and $(\sum_{K\in \cK^\eta_0} |K| S^K_t)_{t\ge 0}$ is independent from
$$
\Big\{\Big( \sum_{K\in \cK} R^K_t + \frac{1}{2\sum_{K\in \cK}|K|} \sum_{K\ne L \in \cK} |K||L| \|S^K_t - S^L_t\|^2\Big)_{t\ge 0}, \cK \subset  \cK^\eta_0(x) \Big\}.
$$
To see the last point, it is sufficient to remember that the family $$
\{ (R^K_t)_{t\ge 0}, K\in \cK^\eta_0\}\cup\{ (S^K_t)_{t\ge 0}, K\in \cK^\eta_0\}
$$
is mutually independent and that the family $\{ (S^K_t)_{t\ge 0}, K\in \cK^\eta_0\}$ has the same law as the family $\{(S_K(B_t))_{t\ge 0}, K\in \cK^\eta_0\}$ where $(B_t)_{t\ge 0}$ is a $2N$-dimensional standard Brownian motion. We can thus write for any $\cK \subset \cK^\eta_0$,
\begin{align*}
\sum_{K\ne L\in \cK} |K||L| \|S_K(B_t)-S_L(B_t)\|^2 = \sum_{K\ne L\in \cK} (|K||L|)^{-1} \Big\|\sum_{k\in K,\ell\in L}(B^k_t - B^\ell_t)\Big\|^2,
\end{align*}
which is of course independent from $(\sum_{i=1}^NB^i_t)_{t\ge 0}$.

\vip

We conclude that $P^1 + P^2 = 0$ by showing that $P^{1,2} = 0$. We have
\begin{align*}
P^{1,2} \le \E \Big[ \sum_{\substack{\cK \subset\cK^\eta_0\\ |\cK|\ge 2}} \PP \Big( \sum_{K\in \cK^\eta_0} R^K_{\tau_\cK} + \frac{1}{2N} \sum_{K\in \cK^\eta_0} |K||L| \|S^K_{\tau_\cK} - S^L_{\tau_\cK}\|^2 = \eta^{1/N}|\cF_0 \Big)\Big],
\end{align*}
where for all $\cK\in \cP(\cK^\eta_0)$ with $|\cK|\ge 2$,
\begin{align*}
\tau_\cK = \inf\Big\{ t\ge 0 : \sum_{K\in \cK} R^{K}_t + \frac{1}{2\sum_{K\in \cK} |K|} \sum_{K\ne L \in \cK} |K||L| \|S^{K}_t - S^{L}_t\|^2 = 2^{-1} \eta^{1/(\sum_{K\in \cK} |K|)}\Big\}.
\end{align*}
If $\cK\ne \cK^\eta_0$, then the same kind of reasoning as above works. If $\cK = \cK^\eta_0$, then the probability is obviously $0$ since $\eta^{1/N} \ne 2^{-1}\eta^{1/N}$. This leads to $P^{1,2}=0$. 

\vip

To end this Step, it remains to show that $P^3 = 0$.
First we have 
$$
\PP(\tau_1 = \tau_3)\le \E\Big[ \sum_{K,L\in \cK^\eta_0} \PP( \tau_K = \tau^L_3 |\cF_0)\Big],
$$
where $\tau_K = \inf \{ t\ge 0 : R^K_t \ge \eta^{1/|K|} \}$ and $\tau^L_3$ is the life time associated to $((\hat X^{N,i}_t)_{t\ge 0})_{i\in L}$. If $K = L$ then a.s. $\PP( \tau_K = \tau^L_3 |\cF_0) = 0$ thanks to $\cP(N-1)$-(d). If $|K|\ne |L|$, then a.s.
$$
\PP( \tau_K = \tau^L_3 |\cF_0) \le \PP(R^K_{\tau^L_3} = \eta^{1/|K|} |\cF_0) = 0,
$$
since $\tau^L_3$ is independent from $(R^K_t)_{t\ge 0}$ and for all $t\ge 0$, all $c>0$, $\PP( R^K_t = c |\cF_0) = 0$, we then deduce that $\PP(\tau_1 = \tau_3) = 0$. Finally, 
$$
\PP(\tau_2 = \tau_3) \le \E\Big[ \sum_{L\in \cK^\eta_0}\sum_{\substack{\cK\subset \cK^\eta_0\\|\cK|\ge 2}} P^3_{\cK,L}\Big].
$$
where
$$
P^3_{\cK,L} = \PP \Big(\sum_{K\in \cK} R^K_{\tau^L_3} + \frac{1}{2\sum_{K\in \cK} |K|} \sum_{K\ne K' \in \cK} |K||K'| \|S^K_{\tau^L_3}-S^{K'}_{\tau^L_3}\|^2 = \frac12 \eta^{1/\sum_{K\in \cK} |K|} |\cF_0\Big).
$$
Since $\cK\ne \{ L\}$ because $|\cK|\ge 2$, then $P^3_{\cK,L} = 0$ a.s. by the same reasoning as before.

\vip

{\it Step 2.2.} We show that the application $\cK^\eta_1$ is a  continuous application for the one-point compactification topology $\tilde \PP^{\eta,N,\theta}_\mu$-a.e.

\vip

We have that $\cK^\eta_1$ is continuous on the event 
$$
\bigcap_{K\subset \ig 1,N\id}\Big\{ R_K(X_{\tau_1\wedge\tau_2} )\ne \eta^{1/|K|} \quad \mbox{ and } \quad \min_{i\notin K} R_{K\cup \{i\}}(X_{\tau_1\wedge \tau_2}) \ne \eta^{1/(|K|+1)} \Big\},
$$
thus it suffices to show that for all $K\subset \ig 1,N\id$, 
$$
\PP \Big(R_K(X_{\tau_1\wedge\tau_2} )= \eta^{1/|K|}\Big)=0.
$$
and this is clear using the same reasoning as in Step 2.1.

\vip

{\it Step 2.3.} We conclude Step 2 here by induction on $k\ge 1$. The result has been proven in Step 2.1 for $k=1$. Take $k\ge 1$ such that the result holds for $k$ and we prove that the result still holds for $k+1$.

\vip

Let $(\Omega,\cF,(\cF_t)_{t\ge 0},\PP)$ a filtered probability space on which is defined a $2N$-dimensional Brownian motion $(B_t)_{t\ge 0}$ driving an $\eta$-$pre$-$KS_\mu$-process $(\tilde X^N_t)_{t\ge 0}$ with $\mu \in \cP(E_{k_0})$. According to the definition of $(\sigma^\eta_k)_{k\ge 0}$ and $(\cK^\eta_k)_{k\ge 0}$, we have
\begin{align}
\sigma^\eta_{k+1}\indiq_{\sigma^\eta_k<\sigma^\eta_*}((\tilde X^N_t)_{t\ge 0}) = \indiq_{\sigma^\eta_k<\sigma^\eta_*} ( (\tilde X^N_t)_{t\ge 0}) \sigma^\eta_1((\tilde X^N_{\sigma^\eta_k + t})_{t\ge 0}),\label{eq:recstoppingtime1}\\
 \indiq_{\sigma^\eta_{k+1}<\sigma^\eta_*}((\tilde X^N_t)_{t\ge 0}) = \indiq_{\sigma^\eta_{k}<\sigma^\eta_*}((\tilde X^N_t)_{t\ge 0})\indiq_{\sigma^\eta_1<\sigma^\eta_*} ((\tilde X^N_{\sigma^\eta_k + t})_{t\ge 0}),\label{eq:recstoppingtime2}\\
 \indiq_{\sigma^\eta_k<\sigma^\eta_*}\cK^\eta_{k+1} ( (\tilde X^N_t)_{t\ge 0}) = \indiq_{\sigma^\eta_{k}<\sigma^\eta_*}((\tilde X^N_t)_{t\ge 0})\cK^\eta_1 ( (\tilde X^N_{\sigma^\eta_k+t})_{t\ge 0}) \label{eq:recstoppingtime3}
\end{align}
According to the induction hypothesis, can write $(\tilde X^N_{\sigma^\eta_k + t})_{t\ge 0} = f((\tilde X^N_t)_{t\ge 0})$ with $f$ continuous at $\omega$ for all $\omega \in A$ with $A\in \cF_{\sigma^\eta_k}$ such that $\PP (A) = 1$. Moreover, thanks to the induction hypothesis, we can also suppose that $\indiq_{\sigma^\eta_k<\sigma^\eta_*}$ is continuous on $A$. By definition of a $\eta$-$pre$-$KS_\mu$-process, $(\tilde X^N_{\sigma^\eta_k + t})_{t\ge 0}$ is also a $\eta$-$pre$-$KS_\nu$-process on $(\Omega,\cF,(\cF_{\sigma^\eta_k + t})_{t\ge 0}, \PP(\cdot | \sigma^\eta_k<\sigma^\eta_*))$ driven by the Brownian motion $(B_{\sigma^\eta_k+t} - B_{\sigma^\eta_k})_{t\ge 0}$ with some initial law $\nu \in \cP(E_{k_0})$ (see Remark \ref{rk:unicitepreKS}-(iii)). Thus, according to Step 2.1, there exists $B\in \cF$ such that $ \PP(B|\sigma^\eta_k<\sigma^\eta_*) = 1$ and for all $\omega' \in B$, $\sigma^\eta_1\indiq_{0<\sigma^\eta_*}$ and $\indiq_{\sigma^\eta_1<\sigma^\eta_*}$ are continuous at $\omega$. 

\vip

Thanks to \eqref{eq:recstoppingtime1}, \eqref{eq:recstoppingtime2} and \eqref{eq:recstoppingtime3}, we have that $\sigma^\eta_{k+1}$ and $\cK^\eta_{k+1}$ are continuous at $\omega$ for all $\omega \in A\cap f^{-1}(B)$. Since $\PP (A\cap f^{-1}(B)) = 1$, we get that $\sigma^\eta_{k+1}$ and $\cK^\eta_{k+1}$ are continuous $\tilde\PP^{\eta,N,\theta}_\mu$-a.e. We conclude thanks to the fact that according to \eqref{eq:recstoppingtime1} and the fact that $(\tilde X^N_{\sigma^\eta_k + t})_{t\ge 0}$ is a $\eta$-$pre$-$KS(\theta,N)$-process, we have that for all $t>0$, $\PP(\sigma^\eta_{k+1} = t) = 0$.

\vip

{\it Step 3.} We show (iii). Let the $\nn\cup\{\infty\}$-valued random variable $p_0 := \inf \{ k\ge 1 : \cK^\eta_k = \{\ig 1,N\id\}\}$, with the convention that $\inf \emptyset = \infty$. According to (ii), we have that $p_0$, viewed as a function of the trajectory, is continuous $\tilde \PP^{\eta,N,\theta}_\mu$-a.e. This together with (ii) and the fact that $\sigma^\eta_* = \indiq_{p_0<\infty}\indiq_{\sigma^\eta_{p_0-1} < \sigma^\eta_*} \sigma^\eta_{p_0} + \infty\indiq_{p_0=\infty}$ gives the result. 

\vip

{\it Step 4.} We show (iv).

\vip

{\it Step 4.1.} Let for any $x\in (E_{k_0}\cap E^\eta_N)_\triangle$, $(X^x_t)_{t\ge 0}$ a $\eta$-$pre$-$KS(\theta,N)$-process driven by $(B^x_t)_{t\ge 0}$ on a filtered probability space $(\Omega^x,\cF^x,(\cF^x_t)_{t\ge 0},\PP^x)$. We show that for all $\varphi : C(\rr_+,(E_{k_0}\cap E^\eta_N \times (\rr^2)^N)_\triangle)\to \rr$ bounded and Lipschitz, all $k\ge 0$, the application
$$
x \mapsto \E^x[\varphi ( (X^x_{t\wedge \sigma^\eta_k})_{t\ge 0})]
$$
is continuous for the one-point compactification topology. The result is clear if $k=0$, we thus take $k\ge 0$  such that the property holds and we show it is still true for $k+1$. 

\vip

Take $\varphi : C(\rr_+,(E_{k_0}\cap E^\eta_N)_\triangle) \to \rr$ Lipschitz for the topology induced by the product topology and some sequence $(x_n)_{n\ge 0}$ of elements of $(E_{k_0}\cap E^\eta_N)_\triangle$ such that $x_n \to x \in (E_{k_0}\cap E^\eta_N)_\triangle$ as $n\to \infty$. Then we have for all $y\in (E_{k_0}\cap E^\eta_N)_\triangle$,
\begin{align*}
\varphi (( X^y_{t\wedge \sigma^\eta_{k+1}})_{t\ge 0}) =  \psi \Big( ( X^y_{t\wedge \sigma^\eta_{k}})_{t\ge 0},( X^y_{(t + \sigma^\eta_k)\wedge \sigma^\eta_{k+1}}
)_{t\ge 0}\Big),
\end{align*}
where $\psi$ is an $\PP^x$-a.e. continuous and bounded application from $(C(\rr_+, (E_{k_0}\cap E^\eta_N\times (\rr^2)^N)_\triangle)^2$ to $\rr$, for all $x\in (E_{k_0}\cap E^\eta_N)_\triangle$ according to Step 2.  According to the induction hypothesis, we get that $\cL( (X^{x_n}_{t\wedge \sigma^\eta_k})_{t\ge 0})$ converges weakly to $\cL((X^x_{t\wedge \sigma^\eta_k})_{t\ge 0})$ as $n\to \infty$, so that, according to the Skorokhod representation theorem, there exists a probability space $(\tilde \Omega,\tilde \cF,\tilde \PP)$ on which are defined a sequence $((X^{k,n}_t)_{t\ge 0})_{n\ge 0}$ of random variable and a random variable $(X^k_t)_{t\ge 0}$ such that for all $n\ge 0$, $(X^{k,n}_t)_{t\ge 0}$ is distributed as $(X^{x_n}_{t\wedge \sigma^\eta_k})_{t\ge 0}$ under $\tilde \PP$,  $(X^k_t)_{t\ge 0}$ is distributed as $(X^x_{t\wedge \sigma^\eta_k})_{t\ge 0}$ under $\tilde\PP$ and there exists $A\in \tilde \cF$ such that $\tilde \PP(A) = 1$ and for all $\tilde \omega \in \tilde \cF$, $(X^{k,n}_t)_{t\ge 0}(\tilde \omega) \to (X^k_t)_{t\ge 0}(\tilde \omega)$ as $n\to \infty$. 

\vip

Then, conditionning on $\tilde\cF_{X^{k,n}_{\sigma^\eta_k}}$, using Remark \ref{rk:unicitepreKS}-(iii) and writing for simplicity for both expectation $(X_t)_{t\ge 0}$ the associated $\eta$-pre-$KS$-process and driving Brownian motion in the second argument, we have
\begin{align*}
   \E^{x_n}\Big[\psi ((X^{k,n}_{t\wedge \sigma^\eta_k})_{t\ge 0},& (X^{k,n}_{t\wedge \sigma^\eta_1})_{t\ge 0})|\tilde\cF_{X^{k,n}_{\sigma^\eta_k}}\Big]=\E^{X^{k,n}_{\sigma^\eta_k}}[\psi ((X^{k,n}_{t\wedge \sigma^\eta_k})_{t\ge 0},(X_{t\wedge \sigma^\eta_1})_{t\ge0}],\\
    \E^{x}\Big[\psi ((X^k_t)_{t\ge 0}, &(X^k_{t\wedge\sigma^\eta_1})_{t\ge0})|\tilde\cF_{X^k_{\sigma^\eta_k}}\Big]=\E^{X^k_{\sigma^\eta_k}}[\psi ((X^k_t)_{t\ge 0}, (X_{t\wedge\sigma^\eta_1})_{t\ge 0})].
\end{align*}
We deduce
\begin{align*}
\Delta_n :=\Big|\E^{x_n}&\Big[\psi ((X^{k,n}_{t\wedge \sigma^\eta_k})_{t\ge 0}, (X^{k,n}_{t\wedge \sigma^\eta_1})_{t\ge 0})|\tilde \cF_{X^{k,n}_{\sigma^\eta_k}}\Big]\\
&-\E^{x}\Big[\psi ((X^k_t)_{t\ge 0}, (X^k_{t\wedge\sigma^\eta_1})_{t\ge 0}) | \cF_{X^k_{\sigma^\eta_k}}\Big]  \Big|\\
=\Big|\E^{X^{k,n}_{\sigma^\eta_k}}&[\psi ((X^{k,n}_{t\wedge \sigma^\eta_k})_{t\ge 0}, (X_{t\wedge \sigma^\eta_1})_{t\ge 0})]-\E^{X^{k}_{\sigma^\eta_k}}[\psi ((X^k_t)_{t\ge 0}, (X_{t\wedge\sigma^\eta_1})_{t\ge 0})]  \Big|.
\end{align*}
Fix $\tilde \omega \in A$ Using the of the fact that any $KS(\theta',k)$-process with $\theta' >0$ and $k\le N$ is Feller thanks to $\cP(N-1)$-(a) and \cite[Theorem 4.2.5 page 167]{ek}, and again the Skorokhod representation theorem, there exists a probability space $(\hat\Omega,\hat\cF,\hat\PP)$ on which are defined a sequence $(\hat X^n_t)_{t\ge 0}$ and a random variable $(\hat X_t)_{t\ge 0}$ such that for all $n\ge 0$, $(\hat X^n_t)_{t\ge 0}$ is $\PP^{X^{k,n}_{\sigma^\eta_k}(\tilde \omega)}$-distributed, $(\hat X_t)_{t\ge 0}$ is $\PP^{X^{k}_{\sigma^\eta_k}(\tilde \omega)}$-distributed and $(\hat X^n_t)_{t\ge 0} \to (\hat X_t)_{t\ge 0}$ as $n\to \infty$ $\hat \PP$-a.s. Thus, we get by the dominated convergence theorem that $\Delta_n(\tilde \omega) \to 0$ as $n\to \infty$. We conclude that $\tilde\E[\Delta_n] \to 0$ as $n\to \infty$ by the dominated convergence theorem, ending the proof by triangular inequality. 

\vip

{\it Step 4.2.} We conclude the Step 4 here. According to the Portemanteau theorem, it is sufficient to show that fixing $\varphi :C(\rr_+,(E_{k_0}\cap E_N^\eta)_\triangle) \to \rr$ bounded and Lipschitz and a sequence $(x_n)_{n\ge 0}$ of elements of $(E_{k_0}\cap E_N^\eta)_\triangle$ converging to some $x\in (E_{k_0}\cap E_N^\eta)_\triangle$, we get
\begin{align*}
 \E^{x_n}[\varphi ( (X^n_t)_{t\ge 0})] \to \E^x[\varphi (( X_t)_{t\ge 0})]\quad \mbox{ as } n \to \infty. 
\end{align*}
We have for all $k\ge 0$, 
\begin{align*}
\Big| \E^{x_n}[\varphi ((X^n_t)_{t\ge 0})] &-\E^x[\varphi (( X_t)_{t\ge 0})] \Big|\le 2\|\varphi\|_\infty \sup_{z\in (E_{k_0}\cap E^\eta_N)_\triangle} \tilde \PP^{\eta,N,\theta}_z ( \sigma^\eta_k<\sigma^\eta_*) \\
&+ \Big| \E^{x_n}[\varphi ((X^n_{t\wedge \sigma^\eta_k})_{t\ge 0})] -\E^x[\varphi (( X_{t\wedge \sigma^\eta_k})_{t\ge 0})] \Big|.
\end{align*}
According to Step 4.1, it is sufficient to show that $ \sup_{z\in (E_{k_0}\cap E^\eta_N)_\triangle} \tilde \PP^{\eta,N,\theta}_z ( \sigma^\eta_k<\sigma^\eta_*)\to 0$ as $k\to \infty$. But this is clear thanks to the first Dini's theorem since for all fixed $z \in (E_{k_0}\cap E^\eta_N)_\triangle$, $\tilde \PP^{\eta,N,\theta}_z(\sigma^\eta_k<\sigma^\eta_*)$ decrease to $0$ as $k\to \infty$, $(E_{k_0}\cap E^\eta_N)_\triangle$ is a compact space for the one-point compactification topology and the application $z \mapsto \tilde \PP^{\eta,N,\theta}_z(\sigma^\eta_k<\sigma^\eta_*)$ is continuous for the one-point compactification topology according to (ii) and Step 4.1.

\vip

{\it Step 5.} We prove (v). Fix $\varphi \in C_0(E_{k_0}\cap E^\eta_N)$. Thanks to Step 4, it only remains to show that
$$
\sup_{x \in E_{k_0}\cap E^\eta_N}\Big|\tilde\E^{\eta,N,\theta}_x[\varphi(X_t)\indiq_{t<\sigma^\eta_*}] - \varphi (x)\Big| \to 0,
$$
as $t \to 0$. 

\vip

We have for $t\ge 0$,
\begin{align}\label{eq:felleretabb}
\sup_{x \in E_{k_0}\cap E^\eta_N}\Big|\tilde\E^{\eta,N,\theta}_x[\varphi(X_t)\indiq_{t<\sigma^\eta_*}] - \varphi (x)\Big| \le E^1_t + E^2_t,
\end{align}
where
\begin{gather*}
E^1_t = \sup_{x \in E_{k_0}\cap E^\eta_N}\Big|\tilde\E^{\eta,N,\theta}_x[\varphi(X_t) \indiq_{t< \sigma^\eta_1}] - \varphi (x)\Big|,\\
E^2_t = \sup_{x \in E_{k_0}\cap E^\eta_N}\tilde\E^{\eta,N,\theta}_x\Big[\sup_{s\in [0,t]}|\varphi(X_s)|\indiq_{\sigma^\eta_1 \le t<\sigma^\eta_*}\Big].
\end{gather*} 
Firstly, one has
$$
 E^1_t \le \sum_{\cK \in \gK_N} E^{1,\cK}_t \quad \mbox{ with }\quad E^{1,\cK}_t = \sup_{x : \cK^\eta (x) = \cK}|\tilde\E^{\eta,N,\theta}_x[\varphi(X_t) \indiq_{t< \sigma^\eta_1}] - \varphi (x)|.
$$
By definition \ref{def:etapreKS}, for each $\cK\in \gK_N$, each $x\in E_{k_0}\cap E^\eta_N$ such that $\cK^\eta(x) = \cK$ there exists a family of mutually independent processes 
$$
\{(X^K_t)_{t\ge 0}, K\in \cK\}
$$ 
such that for all $K\in \cK$, $(X^K_t)_{t\ge 0}$ is a $KS_x(\theta |K|/N,|K|)$-process and for all $t\in [0,\sigma^\eta_1]$, $X^K_t = (X^i_t)_{i\in K}$. Then
\begin{align}\label{eq:decompoenpetitconfig}
E^{1,\cK}_t &\le \sup_{x : \cK^\eta(x) = \cK}|\E_x[\varphi ((X^K_t)_{K\in \cK})] - \varphi (x)| + E^2_t,
\end{align}
where $\PP_x$ is the law of independent $KS$-processes as above and we do the abuse of notation writing again $((X^K_t)_{t\ge 0})_{K\in \cK}$ for the canonical process. According to $\cP(N-1)$-(a), the first term of the right hand side of \eqref{eq:decompoenpetitconfig} tends to $0$ as $t\to 0$, this implies, combined with \eqref{eq:decompoenpetitconfig} that
$$
\limsup_{t\to 0} E^1_t\le \Big(|\gK_N|+1\Big) \limsup_{t\to 0} E^2_t.
$$
All in all it is sufficient to conclude to show that $E^2_t \to 0$ as $t\to 0$, but this is a straight application of the first Dini's theorem since for all $x\in (E_{k_0}\cap E^\eta_N)_\triangle$, $\tilde\E^{\eta,N,\theta}_x[\sup_{s\in [0,t]}|\varphi(X_s)|\indiq_{\sigma^\eta_1 \le t}]$ non increases to $0$ as $t\to 0$ (by the dominated convergence theorem for example) and for all $t>0$, the application
$$
x \in E_{k_0} \cap E^\eta_N \mapsto \tilde\E^{\eta,N,\theta}_x\Big[\sup_{s\in [0,t]}|\varphi(X_s)|\indiq_{\sigma^\eta_1 \le t}\Big]
$$
is continuous for the one-point compactification topology according to Lemma \ref{lem:etapreKS}-(ii).
\end{proof}

\subsection{Estimations for the excursion argument}\label{app:estimation}

The proof of the following is given in \cite[page 37]{G.F.LawlerNotes}
\begin{lemma}\label{lem:EDSexcursion}
Let $(B^r_t)_{t\ge 0}$ be a Bessel process of dimension $\delta \in (0,2)$ starting at $r>0$. Setting $\tau^r_0 := \inf \{ t\ge 0 : B^r_t = 0\}$, and $\gge^r_t = (\tau^r_0)^{-1/2}B^r_{t\tau^r_0}$ for all $t\ge 0$, then $(\gge^r_t)_{t\ge 0}$ solves the following SDE
\begin{align}\label{eq:SDEexcursionbessel} \gge^r_t = \frac{r}{\sqrt{\tau^r_0}} + W_t + \int_0^t\frac{3-\delta}{2\gge^r_s}\dd s - \int_0^t\frac{\gge^r_s}{1-s}\dd s,
\end{align}
where $(W_t)_{t\ge 0}$ is a $1$-dimensional Brownian motion independent of $\tau^r_0$.
\end{lemma}

\begin{lemma}\label{lem:estimexcursionbessel}
(i) For all $r>0$ we define for all $t\ge 0$, $Z^r_t = \gge^r_{\tilde\rho^r_t}$ where $(\gge^r_t)_{t\in [0,1]}$ is as in Lemma \ref{lem:EDSexcursion} and $\tilde\rho^r_t = \inf \{ s\ge 0 :  \tilde A^r_s\ge t\}$ with $\tilde A^r_s = \int_0^s(\gge^r_u)^{-2} \dd u$. Then, there exists $c_\delta,C_\delta \in \rr_+^*$ such that
\begin{align*}
\mbox{ for all } t \ge 0, \quad \sup_{r\in \rr_+} \E [ (Z^r_t)^2] \le C_\delta e^{-c_\delta \sqrt{t}}.
\end{align*}

\vip

(ii) Let $(\gge_t)_{t\in [0,1]}$ be a normalized Bessel excursion and for all $t\ge 0$, $\tilde\rho_t = \inf \{ s\ge 0 : \tilde A_s \ge t\}$ with $\tilde A_s = \int_{1/2}^s \gge_u^{-2}\dd u$. Setting for all $t\ge 0$, $Z_t = \gge_{\tilde\rho_t}$, then there exists $c_\delta,C_\delta\in \rr_+^*$ such that 
\begin{align*}
\mbox{ for all }t \ge 0,\quad \E [ (Z_t)^2 ] \le C_\delta e^{-c_\delta\sqrt{|t|}}.
\end{align*} 
\end{lemma}
\begin{proof}
{\it Step 1.} Fix $r>0$. We show that 
\begin{align*}
\mbox{ for all } t\ge 0, \quad Z^r_t = \frac{\frac{r}{\sqrt{\tau^r_0}}e^{W_t + \nu t}}{1+\frac{r^2}{\tau^r_0} \int_0^t e^{2W_s + 2\nu s}\dd s},
\end{align*}
where $(W_t)_{t\ge 0}$ is a standard Brownian motion and $\nu = 1-\delta/2$. By Lemma \ref{lem:EDSexcursion} and the definition of $(\tilde\rho^r_t)_{t\ge 0}$, we classically have that for all $t\ge 0$,
\begin{align*}
Z^r_t = \frac{r}{\sqrt{\tau^r_0}} + \int_0^t Z^r_s \dd \tilde W_s + \frac{3-\delta}{2}\int_0^t Z^r_s \dd s - \int_0^t \frac{(Z^r_s)^3}{1-\tilde \rho^r_s}\dd s,
\end{align*}
where $(\tilde W_t)_{t\ge 0}$ is a standard Brownian motion. We deduce that
\begin{align}\label{bbb}
Z^r_t = \frac{r}{\sqrt{\tau^r_0}}\exp \Big( \tilde W_t + \nu t - \int_0^t \frac{(Z^r_s)^2}{1-\tilde\rho^r_s}\dd s \Big).
\end{align}
However, thanks to the change of variable $u=s\tau^r_0$ and to Lemma \ref{lem:horlogeexplosion}, we get that a.s., $\tilde A^r_1 = \tau^r_0 \int_0^1 \frac{\dd s}{(\gge^r_{s\tau^r_0})^2}= \int_0^{\tau^r_0} \frac{\dd u}{(\gge^r_u)^2} = \infty$ a.s. and $\tilde A^r_t <\infty$ for all $t<1$, so 
\begin{equation}\label{ddd}
1 = \lim_{t\to \infty} \tilde\rho^r_t = \int_0^\infty (Z^r_s)^2 \dd s,
\end{equation}
since for all $t\ge 0$, $\tilde\rho^r_t = \int_0^t (Z^r_s)^2\dd s$. Thus,
\begin{align}\label{aaa}
- \int_0^t \frac{(Z^r_s)^2}{1-\tilde\rho^r_s}\dd s = - \int_0^t \frac{(Z^r_s)^2}{\int_s^\infty (Z^r_u)^2\dd u}\dd s = \log \Big( \int_t^\infty (Z^r_s)^2 \dd s \Big).
\end{align}
Combining  \eqref{bbb} and \eqref{aaa}, we deduce, 
\begin{align}\label{ccc}
Z^r_t = \frac{r}{\sqrt{\tau^r_0}} \int_t^\infty (Z^r_s)^2\dd s \exp \Big( \tilde W_t + \nu t \Big).
\end{align}
We set $f(t) = \int_t^\infty (Z^x_s)^2\dd s$, so that \eqref{ccc} gives
\begin{align*}
-\frac{f'(t)}{f(t)^2} = \frac{r^2}{\tau^r_0} \exp \Big( 2\tilde W_t + 2\nu t \Big).
\end{align*}
By integrating the last equality and using that $f(0) = \int_0^\infty (Z^r_s)^2\dd s = 1$ thanks to \eqref{ddd}, we get that
\begin{align*}
\frac{1}{f(t)} = 1 + \frac{r^2}{\tau^r_0} \int_0^t \exp \Big( 2 \tilde W_s + 2\nu s\Big) \dd s,
\end{align*}
which implies
\begin{align*}
f(t) = \frac{1}{1 + \frac{r^2}{\tau^r_0} \int_0^t \exp \Big( 2 \tilde W_s + 2\nu s\Big) \dd s}.
\end{align*}
Since for all $t\ge 0$, $Z^r_t = \sqrt{-f'(t)}$ by definition of $f$, this concludes the step.

\vip

{\it Step 2.} We conclude point (i) here. According to Step 1, we have for all $t\ge 0$, all $r>0$, 
\begin{align*}
\E[ (Z^r_t)^2] = \E \Big[ \Big(\frac{\frac{r}{\sqrt{\tau^r_0}}e^{\tilde W_t + \nu t}}{1+\frac{r^2}{\tau^r_0}e^{2\tilde W_t + 2\nu t} Y_t}\Big)^2\Big],
\end{align*}
where $Y_t = \int_0^t \exp ( 2 B_s - 2\nu s )\dd s$ with $(B_s)_{s\in [0,t]} = (\tilde W_{t-s} - \tilde W_t)_{s\in [0,t]}$ a $1$-dimensional Brownian motion. In particular, $Y_t$ is independent of $\tau^r_0$.

\vip

Since $r/\sqrt{\tau^r_0}$ has the same law as $1/ \sqrt{\tau^1_0}$ by scaling, and since $\tau^r_0$ is independent from $(\tilde W_t, Y_t)$, it suffices to estimate $\E [(Z^1_t)^2]$.

\vip

According to \cite[page 36]{G.F.LawlerNotes}, the density on $\rr$ of the law of $\tau^1_0$ is given by
\begin{align}\label{eq:densitetau}
    g^1_0(a) := Ca^{\delta/2-2}e^{-1/(2a)}\indiq_{a>0},
\end{align}
with $C \in (0,\infty)$ a constant, implying that $\E[(\tau^1_0)^{-1}]<\infty$.
This together with independence and the fact that the law of $\sup_{t\in [0,1]} \tilde W_t$ equals the law of $|\tilde W_t|$ which has a gaussian tail we get,
\begin{align*}
\sup_{t\in [0,1]}\E [(Z^1_t)^2] \le \E [ (\tau^1_0)^{-1} ] \E[ e^{2\sup_{t\in [0,1]}\tilde W_t + 2\nu}] <\infty.
\end{align*}

\vip

Let us fix $t\ge 1$. Writing $Z^1_t = \frac{\exp \Big(\tilde W_t + \nu t - \frac{1}{2}\log(\tau^1_0)\Big)}{1+\exp \Big(2\tilde W_t + 2\nu t -\log ( \tau^1_0 ) \Big) Y_t}$, we have
\begin{align}\label{eee}
(Z^1_t)^2 \le \indiq_{\tilde W_t \ge -\nu t/2 + \log (\tau^1_0)/2} \frac{e^{-\nu t }}{(Y_t)^2} + \indiq_{\tilde W_t \le -\nu t /2 + \log (\tau^1_0)/2} \frac{1}{4Y_t}.
\end{align}
Indeed, the first term of the left hand side comes from the fact that $1 + (\tau^1_0)^{-1} e^{2\tilde W_t + 2\nu t}Y_t \ge (\tau^1_0)^{-1} e^{2\tilde W_t + 2\nu t}Y_t$ 
and the second term of the left hand side comes from the fact that for all $b>0$, $\sup_{x\ge 0} x(1+x^2 b)^{-1} = 2^{-1}b^{-1/2}$. Taking the expectation in \eqref{eee}, the Cauchy-Schwarz inequality gives
\begin{align}\label{fff}
\E [(Z^1_t)^2] \le \E[ Y_t^{-2}] e^{-\nu t } + \frac{1}{4} \Big(\E [ Y_t^{-2}] \PP [ \tilde W_t \le -\nu t/2  + \log (\tau^1_0)/2]\Big)^{1/2},
\end{align}
so that, in order to conclude, it suffices to show that there exists a finite constants $C_\delta >0$ such that
\begin{itemize}
\item (a) $\sup_{t\ge 1} \E [ Y_t^{-2} ] \le C_\delta$, 
\item (b) for all $t\ge 1, \quad \PP [ \tilde W_t \le -\nu t /2 + \log (\tau^1_0)/2] \le C_\delta e^{-C_\delta^{-1}\sqrt{t}}.$
\end{itemize}

\vip

We first prove (a). For all $t\ge 1$, 
\begin{align*}
\E[Y_t^{-2}] \le \E \Big[  \Big(\int_0^{1} e^{2B_s - 2\nu s } \dd s \Big)^{-2} \Big] \le \E \Big[ e^{4\sup_{s\in [0,1]}( -B_s) -4\nu} \Big] <\infty,
\end{align*}
where we used again that the law of $\sup_{s\in [0,1]} -B_s$ equals the law of $|B_1|$. 

\vip

We now prove (b). We have $\PP [ \tilde W_t \le - \nu t/2 + \log (\tau^1_0)/2 ] \le P^1_t + P^2_t$ where
\begin{gather*}
P^1_t = \PP [\tilde W_t \le - t^{3/4}], \\
P^2_t = \PP [ \tau^1_0 \ge e^{\nu t-2t^{3/4}} ],
\end{gather*}
by discussing if whether $-\nu t / 2 + \log (\tau^1_0)/2 \le -t^{3/4}$ or not. Since $-t^{-1/2}\tilde W_t$ is a standard Gaussian variable, there exists a universal constants $c',C'>0$ such that $P^1_t \le C' \exp (-c'\sqrt{t})$.

\vip

According to \eqref{eq:densitetau} and the fact that $\delta <2$, we conclude that 
$$
P^2_t  = \tilde C\int_{e^{\nu t -2t^{3/4}}} e^{-1/a} a^{\delta/2-2}\dd a \lesssim   e^{-(\delta/2-1)\nu t}.
$$
which ends the proof of (i).

\vip

{\it Step 3.} We now prove (ii). According the explicit construction of a Bessel excursion given in \cite[Theorem 3.1]{Pitman-Yor1981} we get that $(\gge_t)_{t\in [0,1]}$ equals in law $(\gge_{1-t})_{t\in [0,1]}$, implying that $(Z_t)_{t\ge 0}$ and $(Z_{-t})_{t\ge 0}$ has the same law. So we only need to show that there exists $C>0$ such that for all $t\ge 0$, $\E [ Z_t ] \le Ce^{-C^{-1}t}$. However since the process $(\gge_t)_{t\in [1/2,1]}$ solves the SDE \eqref{eq:SDEexcursionbessel} (see \cite{G.F.LawlerNotes}) with initial law $\gge_{1/2}$ almost surely positive, Step 2. concludes. 
\end{proof}

\begin{lemma}\label{lem:estimeduU}
(i) If $d_{\theta,N}(N-1) >1$, then 
$$
\sup_{u\in \sS_N} \E^{U,N,\theta}_u \Big[ \int_0^1 \| b(U_s)\| \dd s \Big] < \infty\quad \mbox{and for all } K>0 , \quad \sup_{x\in B(0,K)} \E^{X,N,\theta}_x \Big[ \int_0^{\tau_N \wedge 1} \| b(X_s)\|\dd s \Big] <\infty,
$$
where $(X_t)_{t\ge 0}$ and $(U_t)_{t\ge 0}$ are the canonical processes.

\vip

(ii) If $d_{\theta,N}(N)>1$, for all $K>0$, there exists a finite $C_K>0$ such that
$$
\sup_{x\in B(0,K)} \E^{X,N,\theta}_x \Big[ \int_0^{1} \| b(X_s)\|\dd s \Big] \le C_K,
$$
where $(X_t)_{t\ge 0}$ is the canonical process.
\end{lemma}
\begin{proof}
{\it Step 1}. Let $x\in E_N$ and $(\Omega, \cF,(\cF_t)_{t\ge 0}, \PP, (B_t)_{t\ge 0},(X_t)_{t\ge 0})$ be a $KS_x(\theta,N)$-process. We first show that for all stopping time $\tau$ with respect to the filtration $(\cF_t)_{t\ge 0}$, 
\begin{align*}
\int_0^{\tau\wedge \tau_N} \| b(X_s)\| \dd s = \sqrt{\tau_N} \int_0 ^{\tilde A^r_{(\tau/\tau_N)\wedge 1}} Z^r_s \|b(U_s)\|\dd s. \quad \mbox{ a.s.},
\end{align*}
where $r = \sqrt{R_{\ig 1,N\id}(x)}$ and $(U_t)_{t\ge 0} = (\Phi (X_{\rho_t}))_{t\ge 0}$ with the same notation as in  Lemma \ref{lem:estimexcursionbessel} and $\mathrm{Decomposition}(\theta,N)$.

\vip

First of all, performing the change of variable $v = s/ \tau_N$ and recalling Lemma \ref{lem:estimexcursionbessel} for the definition of $(\gge^r_t)_{t\in [0,1]}$, we have
\begin{align*}
\int_0^{\tau\wedge \tau_N} \| b(X_s)\| \dd s = \sqrt{\tau_N} \int_0^{(\tau/\tau_N)\wedge 1} (\gge^r_v)^{-1} \| b(U_{A_{v\tau_N}}) \|\dd v.
\end{align*} 
Moreover, we have for all $v\in [0,1]$, 
$$
A_{v\tau_N} = \int_0^{v\tau_N} \frac{\dd s}{R_{\ig 1,N\id}(X_s)} = (\tau_N)^{-1} \int_0^v \frac{\dd s}{(\gge^r_{v/\tau_N})^2} = \int_0^v \frac{\dd s}{(\gge^r_s)^2} =\tilde A^r_v,
$$ 
which concludes together with the change of variable $v = \tilde \rho^r_s$.

\vip

{\it Step 2} Fix $\eta >0$ and $K>0$. We prove that 
$$
\sup_{x\in B(0,K)} \E_x^{X,N,\theta} \Big[ \int_0^{\tau^\eta_N\wedge 1} \| b(X_s)\|\dd s\Big] <\infty.
$$
Recalling \eqref{config} for the definition of the sequence $(\sigma^\eta_k)_{k\ge 0}$, we have for $x\in B(0,K)$,
\begin{align*}
\E_x^{X,N,\theta}\Big[ \int_0^{\tau^\eta_N \wedge 1} \| b(X_s)\| \dd s \Big] = \sum_{k\ge 0} \E_x^{X,N,\theta}\Big[\int_{\tau^\eta_N \wedge\sigma_k^\eta \wedge 1}^{\tau^\eta_N \wedge\sigma_{k+1}^\eta\wedge 1} \| b(X_s)\| \dd s \Big].
\end{align*}
By the strong Markov property, Proposition \ref{prop:girsanoveta}-(ii)-(c) and $\cP(N-1)$-(b), we get for all $k\ge 0$,
\begin{align*}
\E_x^{X,N,\theta}\Big[\int_{\tau^\eta_N \wedge\sigma_k^\eta \wedge 1}^{\tau^\eta_N \wedge\sigma_{k+1}^\eta\wedge 1} \| b(X_s)\| \dd s \Big] &\le \E_x^{X,N,\theta}\Big[\E_{X_{\tau^\eta_N \wedge\sigma_k^\eta\wedge 1}}^{X,N,\theta} \Big[\int_0^{\tau^\eta_N \wedge\sigma_{1}^\eta\wedge 1} \| b(X_s)\| \dd s \Big]\Big] \le C_{\eta}.
\end{align*}
with $C_\eta$ a positive finite constant only depending on $\eta$. Thus we get for all $x\in B(0,K)$,
\begin{align*}
\E_x^{X,N,\theta}\Big[\int_0^{\tau^\eta_N\wedge 1} \| b(X^N_s)\| \dd s \Big] \le C_\eta \E_x^{X,N,\theta}[|\cN^\eta|].
\end{align*}
where $\cN^\eta := \{ k\ge 0 : \sigma^\eta_k \le 1\}$. It suffices to show that $\sup_{x\in B(0,K)} \E_x^{X,N,\theta}[|\cN^\eta|] < \infty$, which we now do.

\vip

By definition of the sequence $(\sigma^\eta_k)_{k\ge 0}$, and since there is only a finite number of partitions of $\ig 1,N\id$, there exists a constant $c_{N,\eta} >0$ such that for all $k\ge 0$, 
\begin{align*}
\sum_{K\in \cK^\eta_k} \Big[(R_K(X_{\sigma^\eta_{k+1}}) - R_K(X_{\sigma^\eta_k}))^2 + \|S_K(X_{\sigma^\eta_{k+1}}) - S_K(X_{\sigma^\eta_k})\|^2\Big] \ge c_{N,\eta}. 
\end{align*}
Indeed, this comes from the fact that either there exists $K\in \cK^\eta_k$ such that $R_K(X_{\sigma^\eta_{k+1}}) - R_K(X_{\sigma^\eta_k}) \ge 2\eta^{1/|K|} - \eta^{1/|K|} =\eta^{1/|K|} $, either there exists a subfamily $(K_i)_{i\in \ig 1,p\id}$ of $\cK^\eta_k$ with $p\ge 2$ such that 
\begin{align*}
R_{K_1\cup\dots \cup K_p}(X_{\sigma^\eta_{k}}) =& \sum_{i=1}^p R_{K_i}(X_{\sigma^\eta_{k}}) + \frac{1}{2\tilde N} \sum_{1\le i\ne j\le p}|K_i||K_j|\|S_{K_i}(X_{\sigma^\eta_{k}}) - S_{K_i}(X_{\sigma^\eta_{k}})\|^2 \\
&\ge  2 \eta^{1/(|K_1|+\dots |K_p|)},\\
R_{K_1\cup\dots \cup K_p}(X_{\sigma^\eta_{k+1}}) &= \sum_{i=1}^p R_{K_i}(X_{\sigma^\eta_{k+1}}) + \frac{1}{2\tilde N} \sum_{1\le i\ne j\le p}|K_i||K_j|\|S_{K_i}(X_{\sigma^\eta_{k+1}}) - S_{K_i}(X_{\sigma^\eta_{k+1}})\|^2 \\
&\le \frac12 \eta^{1/(|K_1|+\dots |K_p|)},
\end{align*}
where $\tilde N = \sum_{i=1}^p |K_i|$. Thus, we get 
\begin{align*}
    \sum_{k\ge 0} \sum_{K\in \cK^\eta_k} \Big[(R_K(X_{\sigma^\eta_{k+1}\wedge 1}) - R_K(X_{\sigma^\eta_k\wedge 1}))^2 + \|S_K(X_{\sigma^\eta_{k+1}\wedge 1}) - S_K(X_{\sigma^\eta_k\wedge 1})\|^2\Big] \ge c_{N,\eta} (|\cN^\eta|-1). 
\end{align*}
Thus, we get for all $x\in B(0,K)$, by Proposition \ref{prop:girsanoveta}-(ii)-(c) that 
\begin{align*}
  \E_x^{X,N,\theta}[|\cN^\eta|]\le&  1+c_{N,\eta}^{-1} \exp (C_{N,\theta,\eta}T) \E_{\QQ^{T,\eta}_x}\Big[ \sum_{k\ge 0} \sum_{K\in \cK^\eta_k} \Big[(R_K(X_{\sigma^\eta_{k+1}\wedge 1}) - R_K(X_{\sigma^\eta_k\wedge 1}))^2 \\
  &+ (S_K(X_{\sigma^\eta_{k+1}\wedge 1}) - S_K(X_{\sigma^\eta_k\wedge 1}))^2\Big]\Big].
\end{align*}
By definition of $\QQ^{T,\eta}_x$, the SDE of squared Bessel processes and the Itô's isometry, we get
\begin{align*}
    \E_x^{X,N,\theta}[|\cN^\eta|]\le &1+c_{N,\eta}^{-1}\exp (C_{N,\theta,\eta}T)\E_{\QQ^{T,\eta}_x} \Big[ \sum_{k\ge 0}\sum_{K\in \cK^\eta_k} 8\int_{\sigma^\eta_k\wedge 1}^{\sigma^\eta_{k+1}\wedge 1} R_K(X_s) \dd s\\
    &+ 2d_{\theta,N}(|K|)^2 (\sigma^\eta_{k+1}\wedge 1 - \sigma^\eta_k\wedge 1)^2 + |K|^{-1}(\sigma^\eta_{k+1}\wedge 1 - \sigma^\eta_k\wedge 1)  \Big].
\end{align*}
Since the family of non empty subsets of $\ig 1,N\id$ is finite and for all $K\subset \ig 1,N\id$, there exists $C_K>0$ such that for all $x\in (\rr^2)^N$, $R_K(x) \le C_K\|x\|^2$, we deduce the existence of a constant $C\in (0,\infty)$ such that
\begin{align*}
    E_x^{X,N,\theta}[|\cN^\eta|] \le C\Big(1+\E_{\QQ^{T,\eta}_x} \Big[ \int_0^1 \|X_s\|^2 \dd s \Big] \Big).
\end{align*}
Using Proposition \ref{prop:girsanoveta}(ii)-(c) again and Lemma \ref{lem:dispersionformula}, we obtain up to changing the constant $C$,
\begin{align*}
    E_x^{X,N,\theta}[|\cN^\eta|] \le C\Big(1+\E_x^{X,N,\theta} \Big[ \int_0^1 (R_{\ig 1,N\id}(X_s) + \|S_{\ig 1,N\id}(X_s)\|^2 ) \dd s \Big] \Big),
\end{align*}
which is classically finite since $(R_{\ig 1,N\id}(X_t))_{t\ge 0}$ is a squared Bessel process and $(S_{\ig 1,N\id}(X_t))_{t\ge 0}$ is a Brownian motion, see for example \cite{ry}. The uniformity with respect to $x\in B(0,K)$ comes from the fact that $B(0,K)$ is bounded.  

\vip

{\it Step 3.} We prove that
$$
\sup_{u\in \sS_N} \E_u^{U,N,\theta} \Big[ \int_0^1 \|b(U_s)\|\dd s \Big] <\infty.
$$
Let $u\in \sS_N $. Fix $\eta >0$. According to Step 1, Step 2, and since $u \in B(0,1)$, there exists $C_\eta >0$ such that 
\begin{align*}
\E_u^{X,N,\theta}\Big[ \sqrt{\tau_N}\int_0^{\tilde A^1_{(\tau^\eta_N/\tau_N) \wedge 1}} Z^1_s\| b(U_s) \|\dd s \Big] \le C_\eta.
\end{align*}
According to the definition of $\tau^\eta_N$, we have
\begin{align*}
\eta \E_u^{X,N,\theta} \Big[ \indiq_{\tilde A^1_{(\tau^\eta_N/\tau_N) \wedge 1} \ge 1} \int_0^1 \| b(U_s) \| \dd s\Big] \le C_\eta.
\end{align*}
According to the independence from $\mathrm{Decomposition}(\theta,N)$, we get that 
\begin{align*}
\E_u^{U,N,\theta} \Big[\int_0^1 \| b(U_s) \| \dd s\Big] \le \frac{C_\eta}{\eta \PP^{X,N,\theta}_u[\tilde A^1_{(\tau^\eta_N/\tau_N\wedge 1)} \ge 1]}, 
\end{align*}
which is a uniform bound in $u\in \sS$ since the event $ \{\tilde A^1_{(\tau^\eta/\tau_N) \wedge 1}\ge 1\}$ does not depend on $u$ and has a positive probability.

\vip

{\it Step 4.} We conclude (i) by proving here that if $d_{\theta,N}(N-1)>1$, then for all $K>0$,
$$
\sup_{x\in B(0,K)}\E_x^{X,N,\theta}\Big[\int_0^{\tau_N\wedge 1} \|b(X_s)\|\dd s \Big] <\infty.
$$
Since $\tau_N = 0$ a.s. if $R_{\ig 1,N\id}(x) = 0$, thus we can focus on initial configurations in $E_N$. Fix $M>1$ and $x\in E_N$. We have 
$$
\E_x^{X,N,\theta} \Big[ \int_0^{\tau_N\wedge 1} \| b(X_s)\| \dd s\Big] = E_1+E_2,
$$
where
\begin{gather*}
E_1 = \E_x^{X,N,\theta} \Big[ \indiq_{\tau_N\le M} \int_0^{\tau_N\wedge 1} \| b(X_s)\| \dd s \Big],\\
E_2 = \E_x^{X,N,\theta} \Big[ \indiq_{\tau_N> M} \int_0^{\tau_N\wedge 1} \| b(X_s)\| \dd s \Big].
\end{gather*}
We set $r = (R_{\ig 1,N\id}(x))^{1/2}>0$. According to Step 1, the independence between $(Z^r_s)_{s\ge 0}$ and $(U_s)_{s\ge 0}$ (see $\mathrm{Decomposition}(\theta,N)$) and Lemma \ref{lem:estimexcursionbessel}, we get 
\begin{align*}
E_1 &\le \sqrt M\E_x^{X,N,\theta}\Big[ \indiq_{\tau_N \le M} \int_0 ^\infty Z^r_s \|b(U_s)\|\dd s \Big]\\
&\le C\sqrt M\E_u^{U,N,\theta}\Big[ \int_0 ^\infty e^{-c\sqrt{s}} \|b(U_s)\|\dd s \Big],
\end{align*}
where $u = \Phi (x)$. This implies by Fubini that 
\begin{align*}
E_1 &\le C \sqrt{M} \sum_{n=0}^\infty e^{-c\sqrt{n}}\E_u^{U,N,\theta} \Big[ \int_n^{n+1} \| b(U_s) \|\dd s \Big]\\
&\le  C \sqrt{M} \sup_{u\in \sS_N} \E_u^{U,N,\theta}\Big[ \int_0^{1} \| b(U_s) \|\dd s \Big] \sum_{n=0}^\infty e^{-c\sqrt{n}} < \infty,
\end{align*}
by Step 3 and the Markov property of $SKS(\theta,N)$-process.

\vip

Then, using Step 1 again, we get
\begin{align*}
E_2 &=  \E_x^{X,N,\theta}\Big[ \indiq_{\tau_N > M} \sqrt{\tau_N} \int_0 ^{\tilde A^r_{(\tau_N)^{-1} \wedge 1}} Z^r_s \|b(U_s)\|\dd s \Big]\\
&\le \E_x^{X,N,\theta}\Big[\indiq_{\tau_N > M}\sqrt{\tau_N} \sup_{s\in [0, \tilde A^r_{(\tau_N)^{-1}\wedge 1}]}Z^r_s \int_0 ^{\tilde A^r_{(\tau_N)^{-1} \wedge 1}} \|b(U_s)\|\dd s \Big]\\
&\le  \E_x^{X,N,\theta}\Big[ \indiq_{\tau_N > M}\sup_{s\in [0,1]}\sqrt{R_{\ig 1,N\id}(X_s)} \int_0 ^{\tilde A^r_{(\tau_N)^{-1} \wedge 1}} \|b(U_s)\|\dd s \Big].
\end{align*}
Then, using the fact that on $\{ \tau_N>M\}$, we have $(\tau_N)^{-1} \wedge 1 \le M^{-1}$, we get that
\begin{align*}
E_2&\le  \E_x^{X,N,\theta}\Big[\sup_{s\in [0,1]}\sqrt{R_{\ig 1,N\id}(X_s)}\sum_{k=0}^{\lfloor\tilde A^r_{M^{-1}}\rfloor}\int_k^{k+1} \|b(U_s)\|\dd s\Big]. 
\end{align*}
Then, using Step 3 and the strong Markov property for $SKS(\theta,N)$-process, there exists a constant $C>0$ such that
\begin{align*}
E_2&\le  C\E_x^{X,N,\theta}\Big[\sup_{s\in [0,1]}\sqrt{R_{\ig 1,N\id}(X_s)}(\lfloor\tilde A^r_{M^{-1}}\rfloor+1)\Big]\\
&\le \Big( \E_x^{X,N,\theta} \Big[ \sup_{s\in [0,1]} R_{\ig 1,N\id}(X_s) \Big]\Big)^{1/2} \Big( \E_x^{X,N,\theta} \Big[ \big(\tilde A^r_{M^{-1}} +1\big)^2 \Big]\Big)^{1/2}, 
\end{align*}
where we used the Cauchy-Schwarz inequality in the last line. The first term of the product of the right-hand-side is clearly finite according to $\mathrm{Decomposition}(\theta,N)$ by classical properties of Bessel processes (for example, one can use the comparison theorem \cite[Theorem 2.7 page 394]{ry} to show that there exists a $2$-dimensional Brownian motion $(\tilde B_t)_{t\ge 0}$ such that a.s. $R_{\ig 1,N\id}(X_t) \le \|\tilde B_t\|^2$ for all $t\ge 0$). The second term is also finite since for all $\alpha >0$,
\begin{align*}
\PP_x^{X,N,\theta}[ \tilde A^r_{M^{-1}} \ge \alpha ] = \PP^{X,N,\theta}_x [\tilde\rho^r_\alpha \le M^{-1}] = \PP_x^{X,N,\theta} [ 1-\tilde\rho^r_\alpha > 1-M^{-1}]
\end{align*}
According to the proof of \ref{lem:estimexcursionbessel} (keeping the same notations),
$$
1-\tilde\rho^r_\alpha =\Big(1+(r^2/\tau_N)\int_0^\alpha \exp (2\tilde W_s + 2\nu s) \dd s \Big)^{-1},
$$
we get that
\begin{align*}
\PP_x^{X,N,\theta} [\tilde A^r_{M^{-1}} \ge \alpha ] &= \PP_x^{X,N,\theta} \Big[ \int_0^\alpha \exp (2\tilde W_s + 2\nu s) \dd s < \frac{\tau_N}{r^2(M-1)}\Big]\\
&\le \PP_x^{X,N,\theta} \Big[ \exp (2\inf_{s\in [0, \alpha]}\tilde W_s)(e^{2\nu \alpha}-1) < 2\nu \frac{\tau_N}{r^2(M-1)}\Big]\\
&= \PP_x^{X,N,\theta} \Big[ \exp (-2\sqrt{\alpha}|G|)(e^{2\nu \alpha}-1) <2\nu \frac{\tau_N}{r^2(M-1)}\Big],
\end{align*}
where $G$ is a standard gaussian random variable. We used in the last inequality that $(\tilde W_t)_{t\ge 0}$ is independent of $\tau_N$, that $\sup_{s\in [0,\alpha]} \tilde W_s$ has the law of $|\tilde W_\alpha|$ and finally a scaling argument. The random variable $\tau_N/r^2,G$ are independent and by scaling again, $(\tau_N/r^2,G)$ is distributed as $(\tau^1_0,G)$ defined on some probability space $(\tilde \Omega, \tilde \cF, \tilde \PP)$ where $G$ is a standard gaussian law and $\tau^1_0$ is the first hitting time of $0$ of a Bessel process with dimension $\delta$ starting from $1$ independent from $G$. We deduce that 
\begin{align*}
\PP_x^{X,N,\theta}[ \tilde A^r_{M^{-1}} &\ge \alpha ] \le \tilde \PP \Big[ \tau^1_0e^{2\sqrt{\alpha}|G|} > \frac{(M-1)}{2\nu}(e^{2\nu\alpha}-1)\Big]\le P_1 + P_2,
\end{align*}
where
\begin{align*}
    P_1 = \tilde \PP [ |G| > \alpha^{1/4}]\quad \mbox{ and } \quad P_2 = \tilde \PP \Big[\tau^1_0 >\frac{(M-1)}{2\nu} (e^{2\nu \alpha - 2\alpha^{3/4}} -e^{-2\alpha^{3/4}})\Big].
\end{align*}
Since $\tau^1_0$ has an explicit density given in \eqref{eq:densitetau}, we easily get that $\PP_x^{X,N,\theta}[\tilde A^r_{M^{-1}} \ge \alpha] = O(\alpha^{-3})$ as $\alpha$ goes to infinity, which gives the result. 

\vip

{\it Step 5.} We conclude the proof of (ii) by showing that for all $K>0$,
$$
\sup_{x\in B(0,K)} \E_x^{X,N,\theta} \Big[ \int_{\tau_N\wedge 1}^{1} \|b(X_s)\|\dd s\Big] <\infty.
$$
Indeed, this together with Step 4 will conclude according to Lemma \ref{lem:dimension}, giving that for all $N\ge 3$, $d_{\theta,N}(N)>1$ implies that $d_{\theta,N}(N-1)>1$ (the case $N=2$ being trivial).

\vip

Fix $K>0$. According to the Strong Markov property, Lemma \ref{lem:symdesangles} and the Fubini's theorem, we have for all $x\in B(0,K)$,
\begin{align*}
     \E_x^{X,N,\theta} \Big[ \int_{\tau_N\wedge 1}^{1} \|b(X_s)\|\dd s\Big] &=  \E_x^{X,N,\theta} \Big[\indiq_{\tau_N<1} \E_{X_{\tau_N}}^{X,N,\theta} \Big[ \int_0^{1-\tau_N} \|b(X_s)\|\dd s\Big]\Big]\\
     &= \E_x^{X,N,\theta}\Big[ \indiq_{\tau_N<1} \E_{X_{\tau_N}}^{X,N,\theta} \Big[ \int_0^{1-\tau_N} \frac{\dd s}{\sqrt{R_{\ig 1,N\id}(X_s)}}\Big] \int_{\sS_N} \|b(u)\| \beta^N_\theta(\dd u)\Big]\\
     &\le \tilde \E \Big[ \int_0^1 \frac{\dd s}{\sqrt{Z_s}} \Big]\int_{\sS_N} \|b(u)\| \beta^N_\theta(\dd u),
\end{align*}
where $(Z_t)_{t\ge 0}$ is a squared Bessel process with dimension $d_{\theta,N}(N)>1$ starting from $0$ built on some probability space $(\tilde \Omega, \tilde \cF, \tilde \PP)$. According to the triangular inequality, Lemma \ref{lem:radon}-(iii) and \cite[Corollary 1.4 page 441]{ry}, we conclude. 
\end{proof}

\begin{lemma}\label{lem:estimedun}
(i) Setting for all $\gge \in \cE = \{ f \in C(\rr_+,\rr_+) : \exists t_0\ge 0 : \forall t\ge t_0, f(t)=0\} $, all $\gu \in \cU = C(\rr_+,\sS_N)$, $\tilde F(\gge,\gu) = \int_0^{\ell (\gge)} \gge ^{-1}_s\| b(\gu_{A_s})\|\dd s$, we have
$$
\int_\cE \int_\cU 1 \wedge \tilde F(\gge,\gu)^2 \gn(\dd \gge )\Lambda(\dd \gu) < \infty.
$$
(ii) In the case where $k_1\ge N+1$ (and thus $d_{\theta,N}(N)>1$ according to Lemma \ref{lem:dimension}), we have
$$
\int_\cE \int_\cU 1 \wedge \tilde F(\gge,\gu)\gn(\dd \gge )\Lambda(\dd \gu) < \infty.
$$
\end{lemma}
\begin{proof}
{\it Step 1.} Let $\alpha \in \{1,2\}$. According to \eqref{eq:mesureexcursionformule}, we have
\begin{align*}
\int_\cE\int_\cU 1 \wedge \tilde F(\gge,\gu)^\alpha \gn(\dd \gge)\Lambda (\gu)  &= c_\delta\int_0^\infty \frac{\dd \ell}{\ell ^{2-\delta/2}}1\wedge \E\Big[\Big( \int_0^\ell \| b(\ell^{1/2} \gge_{\ell^{-1} s} U_{A_s})\|\dd s\Big)^\alpha \Big]\\
&= c_\delta\int_0^\infty \frac{\dd \ell}{\ell ^{2-\delta/2}} 1\wedge\E\Big[ \Big(\ell^{1/2} \int_0^1 \| b(\gge_{u}U_{A_u})\|\dd u\Big)^\alpha\Big] ,
\end{align*}
where $(\gge_t)_{t\in [0,1]}$ is an excursion of a Bessel process with dimension $d_{\theta,N}(N)$, $(U_t)_{t\in \rr}$ an eternal $SKS(\theta,N)$-process independent from $(\gge_t)_{t\in [0,1]}$, $\E$ the expectation relative to $(\gge_t)_{t\in [0,1]}$ and $(U_t)_{t\in \rr}$ and where we did the change of variable $u = \ell^{-1}s$.
 Then we perform the change of variable $v = A_u$ to get
\begin{align*}
\int_\cE\int_\cU 1 \wedge \tilde F(\gge,\gu)^2 \gn(\dd\gge)\Lambda (\gu) = c_\delta\int_0^\infty \frac{\dd \ell}{\ell ^{2-\delta/2}} 1\wedge \Big(\ell^{1/2} \int_{-\infty}^\infty Z_v\| b(U_{v})\|\dd v\Big)^\alpha\Big],
\end{align*}
where $Z_v = \gge_{\rho_v}$, with $\rho_v = \inf\{t\ge 0 : A_t\ge v\}$, the generalized inverse of $A$. Using Fubini and the change of variable 
$$
u = \ell \E\Big[\Big(\int_{-\infty}^\infty Z_v \|b(U_v)\|\dd v\Big)^2\Big],
$$we get that 
\begin{align*}
\int_\cE\int_\cU1 \wedge \tilde F(\gge,\gu)^2 \gn_1(\dd \gge)\Lambda (\gu) &= c_\delta\E \Big[ \Big(\int_{-\infty}^\infty Z_v \|b(U_v)\|\dd v\Big)^{2}\Big]^{1-\delta/2}\int_0^\infty \frac{1\wedge u^{\alpha/2}}{u^{2-\delta/2}}\dd u.
\end{align*}
The integral $\int_0^\infty(1\wedge u^{\alpha/2 )}) u^{\delta/2-2}\dd u$ is finite in any case because either in case (i) $\alpha =2$ and then $2-\delta/2 \in (1,2)$, or in case (ii) $\alpha =1$ and we have $\delta >1$. It suffices to show that
$$
\E \Big[ \Big(\int_{-\infty}^\infty Z_v \|b(U_v)\|\dd v\Big)^{2}\Big]<\infty,
$$
to get the result. According to Fubini and Lemma \ref{lem:estimexcursionbessel}, we have 
\begin{align*}
\E \Big[ \Big(\int_{-\infty}^\infty Z_v \|b(U_v)\|\dd v\Big)^{2}\Big] &= \int_{-\infty}^\infty\int_{-\infty}^\infty \E [ Z_s Z_t ] \E[\|b(U_s)\|\|b(U_t)\|] \dd s \dd t\\
&\le C^2 \int_{-\infty}^\infty\int_{-\infty}^\infty e^{-c(\sqrt{|s|}+\sqrt{|t|})} \E[\|b(U_s)\|\|b(U_t)\|] \dd s \dd t.
\end{align*}
Moreover, using the Markov property and the stationarity of $(U_t)_{t\in \rr}$, we have for all $s,t\in \rr$,
$$
\E[\|b(U_s)\|\|b(U_t)\|] = \int_{\sS_N} \varphi_u ( |s-t|) \| b(u) \| \beta^N_\theta (\dd u),
$$
where $\varphi_u(\cdot) = \E_u[\|b(U_\cdot)\|]$. We deduce that
\begin{align*}
\E \Big[ \Big(\int_{-\infty}^\infty Z_v \|b(U_v)\|\dd v\Big)^{2}\Big] &\le C^2 \int_{\sS_N}\beta^N_\theta(\dd u)\|b(u)\| \int_{-\infty}^\infty\dd s\int_{-\infty}^\infty\dd t e^{-c(\sqrt{|s|}+\sqrt{|t|})} \varphi_u(|s-t|)\\
&=  C^2 \int_{\sS_N} \beta^N_\theta(\dd u)\|b(u)\| \int_{-\infty}^\infty\dd w\int_{-\infty}^\infty\dd t e^{-c(\sqrt{|w+t|}+\sqrt{|t|})} \varphi_u(|w|) ,
\end{align*}
where we made the change of variable $w = s-t$ in the last integral. We have for all $t\in \rr$, all $u\in \sS_N$,
\begin{align*}
\int_{-\infty}^\infty e^{-c\sqrt{|w+t|}} \varphi_u(|w|)\dd w &= \sum_{k\in \zz} \int_k^{k+1} e^{-c\sqrt{|w+t|}} \varphi_u(|w|)\dd w\\
&\le \sum_{k\in \zz} e^{-c\sqrt{|k+t|}+c} \int_k^{k+1}\varphi_u(|w|)\dd w\\
&\le C'',
\end{align*}
where $C''$ is a finite constant independent of $t$ and $u$, according to Lemma \ref{lem:estimeduU}-(i). We deduce that 
\begin{align*}
\E \Big[ \Big(\int_{-\infty}^\infty Z_v \|b(U_v)\|\dd v\Big)^{2}\Big] \le C^2C''\int_{-\infty}^\infty e^{-c\sqrt{|t|}}\dd t < \infty,
\end{align*}
which ends the step.
\end{proof}

\section{Independence and weak uniqueness of S.D.E}\label{sec:unicitefaible}
This section is devoted to the proof of Propositions \ref{prop:unicitechelou} and \ref{prop:indepuniciteloi}. Although these results are crucial for establishing the independence of Keller-Segel processes driven by independent noises, we believe that they are of independent interest. Moreover, their proofs are rather simple, which motivates presenting them separately in a dedicated section. For any $d_1,d_2\in \nn^*$, we say the families of measurable functions from $\rr^{d_1}$ to respectively the space of $d_1\times d_2$ matrices and $\rr^{d_1}$ $(\sigma_\e)_{\e>0}$, $(b_\e)_{\e>0}$ satisfy the property $\mathrm{(U)}$ if 
\begin{itemize}
    \item for every probability measure $\mu$ on $\rr^{d_1}\cup \{\triangle\}$ (where $\triangle$ is a cemetery point), there exists a tuple $(\Omega,\cF,(\cF_t)_{t\ge 0},(B_t)_{t\ge 0}, (X_t)_{t\ge 0},\PP)$ such that $(\Omega,\cF,(\cF_t)_{t\ge 0},\PP)$ is a filtered probability space, $X_0$ is $\mu$-distributed,  $(B_t)_{t\ge 0}$ a $d_2$-dimensional Brownian motion $(\cF_t)_{t\ge 0}$-adapte and that for all $T\ge0$, 
\begin{align}\label{edsclassique}
\sup_{t\in [0,T]}\Big\|X_t - X_0 - \int_0^t \sigma_{\e} (X_s)\dd B_s - \int_0^t b_{\e}(X_s)\dd s\Big\| \to 0 \mbox{ in probability as } \e \to 0.
\end{align}
We call the tuple a weak solution of \eqref{edsclassique}.
\item The law of $(X_t)_{t\ge0}$ is the same for all weak solutions of \eqref{edsclassique} such that $X_0$ is $\mu$-distributed; we denote it by $\PP_\mu$ and write $\PP_x=\PP_{\delta_x}$. Moreover, for every measurable set $\Gamma$ of $C(\rr_+,\rr^{d_1})$, the map $x\mapsto\PP_x(\Gamma)$ is measurable and $\PP_\mu(\Gamma)=\int_{\rr^{d_1}}\PP_x(\Gamma)\,\mu(\dd x)$.
\end{itemize}

\begin{prop}\label{prop:unicitechelou}
Let $d_1$, $d_2\in \nn^*$ and $(\sigma_\e)_{\e>0}$, $(b_\e)_{\e>0}$ two families of measurable applications from $\rr^{d_1}$ to respectively the space of $d_1\times d_2$ matrices and $\rr^{d_1}$ such that $\mathrm{(U)}$ is satisfied. Then, if 
$$
(\Omega,\cF,(\cF_t)_{t\ge 0},(B_t)_{t\ge 0}, (X_t)_{t\ge 0},\PP)
$$
is such a solution and $\cG$ is a sub $\sigma$-algebra of $\cF$ such that for all $t\ge 0$, $\cG\vee \cF_t$ is independent from $(B_{s+t} - B_t)_{s\ge 0}$, then for every $A\in \cB(C(\rr_+,\rr^{d_1}\cup\{\triangle\})$, $\PP((X_t)_{t\ge 0}\in A|\cG\vee \sigma (X_0)) = \PP_{X_0}(A)$ a.s. In particular, if $X_0$ is independent from $\cG$, then $(X_t)_{t\ge 0}$ is independent from $\cG$. 
\end{prop}
\begin{proof}
Let $d_1$, $d_2$, $(\sigma_\e)_{\e >0}$, $(b_\e)_{\e >0}$ as above and fix $\mu \in \cP(\rr^{d_1}\cup \{\triangle\})$ and $(\Omega,\cF,(\cF_t)_{t\ge 0},(B_t)_{t\ge 0}, (X_t)_{t\ge 0},\PP)$ a weak solution to \eqref{edsclassique} with $X_0$ having law $\mu$ and $\cG$ a sub $\sigma$-algebra of $\cF$ as in the statement.

\vip

Take $G\in \cG$ and $C$ a borelian of $\rr^{d_1}$ with $\PP(\tilde G)>0$ where $\tilde G = G \cap \{X_0\in C\}$. We define the filtered probability space $(\Omega, \cF, (\cF_t)_{t\ge 0}, \PP(\cdot|\tilde G))$. According to the hypothesis on $\cG$, $(B_t)_{t\ge 0}$ is still a $d_2$-dimensional Brownian motion in the probability space $(\Omega, \cF, (\cF_t)_{t\ge 0},\PP(\cdot|\tilde G))$. Moreover we have for all $T\ge 0$, all $\eta >0$,
\begin{align*}
\PP\Big[ \sup_{t\in [0,T]}\Big\|&X_t - X_0 - \int_0^t \sigma_\e (X_s)\dd B_s - \int_0^t b_\e(X_s)\dd s\Big\| \ge \eta\Big|\tilde G\Big] \\
&\le (\PP(\tilde G))^{-1}\PP\Big[ \sup_{t\in [0,T]}\Big\|X_t - X_0 - \int_0^t \sigma_\e (X_s)\dd B_s - \int_0^t b_\e(X_s)\dd s\Big\| \ge \eta\Big] \to 0 \mbox{ as } \e \to 0.
\end{align*}
Denote by $\mu_{\tilde G}$ the law of $X_0$ under $\PP(\cdot |\tilde G)$. We deduce that for all $A\in \cB(C([0,T],\rr^{d_1}\cup \{\triangle\})$,
\begin{align*}
\PP ( \{ (X_t)_{t\in [0,T]} \in A \} \cap \tilde G) &= \E [\PP_{X_0} (\{ (X_t)_{t\in [0,T]} \in A \}|\tilde G)\indiq_{\tilde G} ] \\
&= \E [\int_{\rr^{d_1}} \PP_x ( \{ (X_t)_{t\in [0,T]} \in A \})\mu_{\tilde G}(\dd x) \indiq_{\tilde G}] \\
&= \int_{\rr^{d_1}} \PP_x ( \{ (X_t)_{t\in [0,T]} \in A \})\mu_{\tilde G}(\dd x) \PP(\tilde G).
\end{align*}
Since this equality trivially holds if $\PP(\tilde G) = 0$, we get that $(X_t)_{t\in [0,T]}$ is independent of $\cG$ for all $T>0$, and is distributed as stated, which concludes the proof.  
\end{proof}
\begin{prop}\label{prop:indepuniciteloi}
Let $d_1$, $d_2\in \nn^*$, $n\in \nn^*$, $(x^p)_{p\in \ig 1,n\id}\in (\rr^{d_1})^n$ and $(\sigma^p_\e)_{\e>0,p\in \ig 1,n\id}$, $(b^p_\e)_{\e>0,p\in \ig 1,n\id}$ two families of measurable applications from $\rr^{d_1}$ to respectively the space of $d_1\times d_2$ matrices and $\rr^{d_1}$ satisfying $\mathrm{(U)}$. Suppose that there exists a tuple $(\Omega, \cF,(\cF_t)_{t\ge 0}, ((B^p_t)_{t\ge 0})_{p\in \ig 1,n\id}, ((X^p_t)_{t\ge 0})_{p\in \ig 1,n\id})$ where $(\Omega,\cF,(\cF_t)_{t\ge 0},\PP)$ is a filtered probability space on which the family $((B^p_t)_{t\ge 0})_{p\in \ig 1,n\id}$ is a family of i.i.d $d_2$-dimensional Brownian motion and for all $p \in \ig 1,n\id$, $(X^p_t)_{t\ge 0}$ solves \eqref{edsclassique} with coefficients $(\sigma^p_\e)_{\e>0}, (b^p_\e)_{\e>0}$, $(B^p_t)_{t\ge 0}$ as driving Brownian motion and starting at $x^p$. Then the family $((X^p_t)_{t\ge 0})_{p\in \ig 1,n\id}$ is a family of mutually independent random variable. 
\end{prop}
\begin{proof}
Let $(\Omega, \cF,(\cF_t)_{t\ge 0}, ((B^p_t)_{t\ge 0})_{p\in \ig 1,n\id}, ((X^p_t)_{t\ge 0})_{p\in \ig 1,n\id})$ as in the statement. Fix $p\in \ig 1,n\id$. We will prove by induction that for all $p\in \ig 0,n\id$, 
\begin{align*}
\cP(p) :''\mbox{ the family } \{ (X^i_t)_{t\ge 0}, i\in \ig 1,p\id\} \cup \{ (B^i_{t})_{t\ge 0}, i\in \ig p+1, n\id \} \mbox{ is mutually independent }\\
\mbox{and for all } t\ge 0, \cH^p_t := \cF_t \vee \sigma ((X^i_t)_{t\ge 0}, i\le p) \mbox{ is independent from } \sigma ( (B^i_{t+s}-B^i_t)_{s\ge 0}, i\ge p+1)  .''
\end{align*}
$\cP(0)$ is clear since $((B^i_t)_{t\ge 0})_{i\in \ig 1,n\id}$ is mutually independent by hypothesis. 

\vip

Suppose that $\cP(p)$ is true for some $p\in \ig 0,n-1\id$. We show that $\cP(p+1)$ is true. We set
$$
\cG_{p+1} = \sigma \Big( \{ (X^i_t)_{t\ge 0}, i\in \ig 1,p\id\}\cup \{ (B^i_t)_{t\ge 0}, i\in \ig p+2,n\id\}\Big).
$$
We have that for all $t\ge 0$, $(B^{p+1}_{s+t} - B^{p+1}_t)_{s\ge 0}$ is independent from $\cG_{p+1}\vee \cF_t$ thanks to $\cP(p)$, so that, using Proposition \ref{prop:unicitechelou}, we get that $(X^{p+1}_t)_{t\ge 0}$ is independent from $\cG_{p+1}$. This implies that $(X^{p+1}_t)_{t\ge 0}$ is independent from the family of random variables
$$
\cI =\{ (X^i_t)_{t\ge 0}, i\in \ig 1,p\id\}\cup \{ (B^i_t)_{t\ge 0}, i\in \ig p+2,n\id\}.
$$
Since $\cI$ is a family of mutually independent random variable according to $\cP(p)$, we conclude the first part of $\cP(p+1)$. Fix $t\ge 0$. It only remains to prove that $\cH^{p+1}_t$ is independent from $\sigma ((B^i_{s+t}-B^i_s)_{s\ge 0}, i\ge p+2)$. However this is direct using the first part of $\cP(p)$ and the property (U). 
\end{proof}
This result implies directly the following, which will not be useful in this paper but seems to be interesting on its own since it gives an additional flavor about the meaning of uniqueness in law.
\begin{corollary}\label{cor:uniciteindep}
Let $d_1$, $d_2\in \nn$, $n\in \nn^*$ and $(\sigma^p)_{p\in \ig 1,n\id}$, $(b^p)_{p\in \ig 1,n\id}$ two families of measurable applications from $\rr^{d_1}$ to respectively the space of $d_1\times d_2$ matrices  and $\rr^{d_1}$ such that for all $p\in \ig 1,n\id$, all $\mu \in \cP(\rr^{d_1}\cup\{\triangle\})$, all tuple $(\Omega,\cF,(\cF_t)_{t\ge 0},\PP,(B_t)_{t\ge 0},(X_t)_{t\ge 0})$ such that $(\Omega,\cF,(\cF_t)_{t\ge 0},\PP)$ is a filtered probability space and 
$$
\mbox{ for all } t\ge 0, \quad X_t = X_0 + \int_0^t \sigma^p(X_s)\dd B_s + \int_0^t b^p(X_s)\dd s,
$$
with $X_0$ having law $\mu$,
the law of $(X_t)_{t\ge0}$ is characterized and denoted by $\PP_\mu$ (we write $\PP_x=\PP_{\delta_x}$). Suppose moreover that for every measurable set $\Gamma$ of $C(\rr_+,\rr^{d_1}\cup\{\triangle\})$, the map $x\mapsto\PP_x(\Gamma)$ is measurable and $\PP_\mu(\Gamma)=\int_{\rr^{d_1}\cup\{\triangle\}}\PP_x(\Gamma)\,\mu(\dd x)$. Then, if there exists $(x^p)_{p\in \ig 1,n\id}$ a tuple 
$$
(\Omega, \cF,(\cF_t)_{t\ge 0},\PP, ((B^p_t)_{t\ge 0})_{p\in \ig 1,n\id}, ((X^p_t)_{t\ge 0})_{k\in \ig 1,n\id})
$$
where $(\Omega,\cF,(\cF_t)_{t\ge 0},\PP)$ is a filtered probability space on which the family $((B^p_t)_{t\ge 0})_{p\in \ig 1,n\id}$ is a family of i.i.d $d_2$-dimensional Brownian motion and for all $p \in \ig 1,n\id$, $(X^p_t)_{t\ge 0}$ solves the SDE with coefficients $\sigma^p,b^p$, initial condition $x^p$ and $(B^p_t)_{t\ge 0}$ as driving Brownian motion, therefore the family $((X^p_t)_{t\ge 0})_{p\in \ig 1,n\id}$ is a family of mutually independent random variables. 
\end{corollary}

\section{Other results}\label{app:otherresult}
\subsection{Some results about Bessel processes}
We show the largely suboptimal following result.:
\begin{lemma}\label{lem:estimetuagebessel}
    Let $\delta \in (0,2)$ and $(Z_t)_{t\ge 0}$ a squared Bessel process with dimension $\delta$ defined on some probability space $(\Omega,\cF,(\cF_t)_{t\ge 0},\PP)$ driven by a standard Brownian motion $(B_t)_{t\ge 0}$. We note $\PP_r$ the law of $(Z_t)_{t\ge 0}$ starting from $r\ge 0$. We have for all $r> 0$,
    $$
    \sup_{a\ge r}\PP_a[t\ge \tau_0] \le C_rt^2 \quad \mbox{ and }\quad \PP_0[t\ge \tau_r] \le C_rt^2,
    $$
    where for all $\eta\ge 0$, $\tau_\eta = \inf \{ t\ge 0 : Z_t =\eta\}$ and $C_r$ is a positive constant only depending on $r$. 
\end{lemma}
\begin{proof}
   We begin with the first estimate. Using \eqref{eq:densitetau} and a scaling property, we get that
   \begin{align*}
   \PP_r[t\ge \tau_0] = \PP_1[t/r\ge \tau_0] &= C\int_0^{t/r} a^{\delta/2-2} e^{-1/(2a)} \dd a\\
   &\le \tilde Ce^{-r/(4t)}.
   \end{align*}
This clearly gives the result if $r>0$ is away from $0$. 

\vip

Let us turn to the second estimate. Fix $r\ge 0$. Applying the Ito's formula, we get
\begin{align*}
   r =  D_{\tau_r} = 2\int_0^{\tau_r}\sqrt{D_s}\dd B_s + \delta \tau_r.
\end{align*}
  Thus, on the event $\{t\ge \tau_r\}$, we have $\int_0^{\tau_r}\sqrt{D_s}\dd B_s\ge (r - \delta t)/2$, implying for $t$ small enough,
  $$
  \{ t\ge \tau_r\}\subset \Big\{ \sup_{s\in [0,t]} \int_0^s\sqrt{D_s}\dd B_s \ge \frac{r}{4}\Big\}.
  $$
Using the Markov inequality and the Burkholder-Davis Gundy inequality with power $4$, we get the existence a constant $\tilde C_r>0$ such that
\begin{align*}
    \PP_0[t\ge \tau_r] \le \tilde C_r\E_0 \Big[ \Big( \int_0^t D_s\dd s\Big)^2\Big].
\end{align*}
For all $t\ge 0$ smaller than $1$ we get
$$
\PP_0[t\ge \tau_r] \le \tilde C_r \E_0\Big[\sup_{s\in [0,1]}D_s^2\Big] t^2.
$$
Here again, it is classical (using a Gronwall argument) that $\E\Big[\sup_{s\in [0,1]}D_s^2\Big] <\infty$, which concludes. 
\end{proof}

\begin{lemma}\label{lem:horlogeexplosion}
Let $(\Omega,\cF,(\cF_t)_{t\ge 0},\PP)$ be a probability space on which is defined $(D_t)_{t\ge 0}$ a squared Bessel process with dimension $\delta \in \rr$ driven by a standard Brownian motion $(B_t)_{t\ge 0}$.

\vip

(i)-(a) If $\delta\ge 2$ and $D_0 >0$, then setting for all $t\ge 0$, $A_t =\int_0^t D_s^{-1}\dd s$ we have $A_\infty := \lim_{t\to\infty} A_t = \infty$. Setting for all $t\in \rr_+$, $\rho_t = \inf \{ s\in \rr : A_s \ge t\}$, we have that $\rho : \rr_+ \to \rr_+$ is a continuous bijection.

\vip

(i)-(b) If $\delta\ge 2$ and $D_0 =0$, then setting for all $t_0>0$, all $t>0$, $A^{t_0}_t = \int_{t_0}^t D_s^{-1}\dd s$ we have $A^{t_0}_\infty := \lim_{t\to\infty} A^{t_0}_t = \infty$ and $A^{t_0}_0 := \lim_{t\to 0} A^{t_0}_t = -\infty$ a.s. Setting for all $t\in \rr$, $\rho^{t_0}_t = \inf \{ s\in \rr : A^{t_0}_s \ge t\}$, we have that $\rho^{t_0} : \rr \to \rr_+$ is a continuous bijection.  

\vip

(ii) If $\delta \in (0,2)$, let $\cZ = \{ t\ge 0 : D_t = 0\}$. We have $\rr_+\setminus \cZ=[0,\tau^0_0)\cup_{i\ge 0}(a_i,b_i)$ with $\tau_0$ the first hitting time of $0$ and the family $\{(a_i,b_i),i\ge 0\}$ composed of disjoint intervals (since $\rr_+\setminus \cZ$ is open in $\rr_+$) with $i\mapsto (a_i,b_i)$ measurable. Define for every $i\ge 0$,
\begin{gather*}
\mbox{ for all } t\in [0,\tau^0_0), \quad A_{t} = \int_{0}^t \frac{\dd s}{D_s},\\
    \mbox{ for all } t\in (a_i,b_i), \quad A^i_{t} = \int_{(a_i+b_i)/2}^t \frac{\dd s}{D_s},
\end{gather*}
and 
$$
A_{\tau^0_0} = \lim_{t\to \tau^0_0} A_t,\quad A^i_{a_i} = \lim_{t\to a_i} A^i_t \quad \mbox{ and } A^i_{b_i} = \lim_{t\to b_i} A^i_t
$$
which all exist by monotonicity. Then, $A_{\tau^0_0}=A^i_{b_i}=- A^i_{a_i} = \infty$ a.s. Thus, defining for all $t\ge 0$,
\begin{gather*}
\rho_t = \inf\{ s\ge 0 : A_s\ge t\},\\
    \rho^i_t = \inf\{ s\ge 0 : A^i_s \ge t\},
\end{gather*}
we get that $\rho : \rr_+\to [0,\tau^0_0)$ and $\rho^i : \rr \to (a_i,b_i)$ are continuous bijection.

\vip

(iii) If $\delta \le 0$, we get analogously $\tau^0_0$, $(A_t)_{t\in [0,\tau^0_0)}$, $(\rho_t)_{t\ge 0}$ as above. Then, $A_{\tau^0_0} = \infty$ a.s. and $\rho : \rr_+ \to [0,\tau^0_0)$ is a continuous bijection. 
\end{lemma}
\begin{proof}
{\it Step 1.} We show (i)-(a). Since $D_0>0$ and $\delta \ge 2$, we have $\tau^0_0 = \infty$ a.s., and one can apply the Ito's formula to get for all $t\ge 0$,
\begin{align}\label{eq:itobessel}
    \log (D_t) = \log(D_0) + 2\int_0^t \frac{\dd W_s}{\sqrt{D_s}} + (\delta - 2) \int_0^t \frac{\dd s}{D_s}. 
\end{align}
Using the Dubins-Schwarz theorem, we have the existence of a standard brownian motion $(B_t)_{t\ge 0}$ such that
\begin{align}\label{eq:logbessel}
    \log (D_t) = \log(D_0) + 2B_{A_t} + (\delta - 2) A_t. 
\end{align}
Since $(D_t)_{t\ge 0}$ is transient if $\delta >2$ and recurrent if $\delta =2$, in any case there exists a.s. a sequence $(t_n)_{n\ge 0}$ such that  $D_{t_n} \to \infty$ as $n\to \infty$. This together with \eqref{eq:logbessel} gives the result.

\vip

{\it Step 2.} Let $\delta \ge 2$ and suppose that $D_0 = 0$. Fix $t_0>0$. Let us show that $A^{t_0}_\infty = \infty$ a.s., but this is the same proof as in Step 1 providing that if $t_0>0$, one has for all $t\ge t_0$,
$$
\log (D_t) = \log(D_{t_0}) + 2\int_{t_0}^t \frac{\dd W_s}{\sqrt{D_s}} + (\delta - 2) \int_{t_0}^t \frac{\dd s}{D_s},
$$
which is true since for all $\e >0$, one can apply Ito's formula on $[\tau_\e,\infty)$ with $\tau_\e = \inf \{ t\ge 0 : D_t \ge \e\}$ and a.s. there exists $\e \in \QQ_+^*$ such that $t_0\in [\tau_\e,\infty)$. 

\vip

{\it Step 3.} Suppose that $\delta\in (0,2)$. We show that $A^i_{b_i} = \infty$ a.s. for all $i\ge 0$. Using the same trick as in Step 2, one can write for all $t\in (a_i,b_i)$,
\begin{align*}
    \log (D_t) = \log (D_{(a_i+b_i)/2}) + 2\int_{(a_i+b_i)/2}^t\frac{\dd W_s}{\sqrt{D_s}} + (\delta-2)\int_{(a_i+b_i)/2}^t\frac{\dd s}{D_s}. 
\end{align*}
According to the Dubins-Schwarz theorem, one gets the existence of a standard Brownian motion $(B_t)_{t\ge 0}$ such that
\begin{align}\label{eq:itoexcursion}
    \log (D_t) = \log (D_{(a_i+b_i)/2}) + 2B_{A^i_t} + (\delta-2)A^i_t. 
\end{align}
We have $\log(D_t)\to -\infty$ as $t\to b_i$. This implies that $A^i_{b_i} = \infty$ a.s. because otherwise, the RHS of \eqref{eq:itoexcursion} would remain a.s. bounded. 

\vip

{\it Step 4.} Point (iii) can be proved the same way as Step 3. We omit the details. 

\vip

{\it Step 5.} We prove (ii). It only remains to prove that for all $i\ge 0$, $A^i_{a_i} = -\infty$ a.s. According to the strong Markov property, it is sufficient to show the result for $(D_t)_{t\ge 0}$ a squared Bessel process with dimension $\delta$ starting from $0$. However it is clear since for every $\delta$-dimensional Bessel excursion $(\gge_t)_{t\in [0,1]}$, its law equals $(\gge_{1-t})_{t\in [0,1]}$ according to the explicit construction in \cite[Theorem 3.1 page 302]{Pitman-Yor1981}.

\vip

{\it Step 6.} Point (i)-(b) can be proved using a similar, albeit simpler, argument as in Step 5. We omit the details.
\end{proof}

\subsection{Proof of Lemma \ref{lem:dimension}}\label{app:preuveLemma1}
\begin{proof}
    We follow the proof of \cite[Lemma 1]{ft}. For \eqref{eq:dimpositive}, the same proof as in \cite{ft} concludes.

    \vip

    We now fix $N\ge 2$ and $\theta >0$ such that $N>3\theta$. If $k_2\ge 4$, then $2N/\theta >4$, which imposes $d_{\theta,N}(2) = 2-2\theta/N\in (1,2)$. If $k_2= 3$ and $k_2=k_1-1$, then we conclude as before, and it thus remains the case where $k_2 = 3$ and $d_{\theta,N}(k_2)\ge 2$. However, in this case, $2(2-3\theta/N)\ge 2$ i.e. $\theta/N\le 1/3$, implying that $2-2\theta/N\ge 4/3$. 

    \vip 

    We put ourselves in the case where $k_2\ge 4$ and prove that for all $k\in \ig 3,k_2-1\id$, $d_{\theta,N}(k) \ge 2$. By concavity of $x\mapsto (x-1)(2-\theta x/N)$, it is sufficient to show that (i) $d_{\theta,N}(3)\ge 2$, (ii) $d_{\theta,N}(k_2-1)\ge 2$. We introduce $a = 2N/\theta$ and observe that $d_{\theta,N}(k) = 2a^{-1}(k-1)(a-k)$ and $k_0 = \lceil a\rceil$.

    \vip

    For (i), in the case where $k_2\ge 6$, we have $d_{\theta,N}(3) = 4(1-3a^{-1}) \ge 2$ since $a\ge 6$. 
    
    In the case where $k_2 = 5$ and $d_{\theta,N}(k_1-1)<2$, the reasoning is the same. 
    
    If now $k_2 = 5$ and $d_{\theta,N}(k_1-1)\ge 2$, we have $\theta/N \le 3/10$, which implies that $d_{\theta,N}(3) = 2(2-3\theta/N)\ge2(2-9/10)\ge 2$. 
    
    In the case where $k_2 = 4$ and $d_{\theta,N}(k_1-1)<2$, we have $3(2-4\theta/N)<2$ so that $N<3\theta$, which contradicts the hypothesis $N>3\theta$. 

    In the case where $k_2 = 4$ and $d_{\theta,N}(k_1-1)\ge 2$, we have $\theta/N \le 1/3$, implying that $d_{\theta,N}(3)=2(2-3\theta/N) \ge 2(2-1)=2$.

    \vip

    For (ii), we can suppose $d_{\theta,N}(k_1-1)<2$ because it is immediate in the other scenario. In the case where $k_2\ge 5$, we have $d_{\theta,N}(k_0-3) = 2a^{-1}(\lceil a\rceil - 4)(a- \lceil a\rceil +3)$. It amounts to prove that $(\lceil a\rceil - 4)(a- \lceil a\rceil +3) \ge a$. Writing $a = n-r$ with $n\in \nn^*$ and $r\in [0,1)$, we need that $(n-4)(3-r) \ge n-r$ i.e. $(2-r)n + 5r \ge 12$. However, $n\ge 7$ so that $(2-r)n + 5r \ge 14 - 2r \ge 12$. In the case where $k_2 = 4$, we have $k_2-1 = 3$ which amounts to point (i). 

    \vip

    To conclude, it is sufficient to prove that (a) $d_{\theta,N}(k_1-1)>1$ if $N/\theta\ne 2$, (b) $d_{\theta,N}(k_1)<2$ and (c) $d_{\theta,N}(k_0-1) = 1$ iff there exists $n\in \nn \setminus\{ 0,1,2,3\}$ such that $2N/\theta = n - (n-4)/(2n-5)$.

    \vip

    For (a) we can of course focus on the case where $d_{\theta,N}(k_1-1)<2$. We write $a = n-r$ with $n\in \nn^*$ and $r\in [0,1)$. In the case where $n\ge 5$ we have $d_{\theta,N}(k_0-2) = 2(n-3)(2-r)/(n-r) > 2(n-3)/(n-1) \ge 1$. In the case where $n=4$, we get $d_{\theta,N}(k_0-2) = d_{\theta,N}(2) = 2-4a^{-1} = 2- 4/(4-r) >1$ since $r>0$ because $N/\theta \ne 2$ by assumption. 

    \vip

    For (b), one has $d_{\theta,N}(k_1) = 2(n-2)(1-r)/(n-r) \le 2(n-2)/n < 2$.

    \vip

    Finally for (c), if $d_{\theta,N}(k_0-1) = 1$ then $2(n-2)(1-r)/(n-r) =1$ iff $(2n-5)r = n-4$ iff $r = (n-4)/(2n-5)$, which is positive only if $n\ge 4$.    
\end{proof}

\subsection{Proof of Lemma \ref{lem:radon}-(iii)}\label{app:preuveradon}
\begin{proof}[Proof of Lemma \ref{lem:radon}-(iii)]
We prove by induction that $\cP(N)$:'' for all $\theta\in (0,N)$, all $k\ge 2$ such that $d_{\theta,N}(k)>1$, for all $\cC$ compact of $E_{k+1}$,
$$
\int_{\cC}\frac{1}{\|x^1-x^2\|}\mu^N_\theta (\dd x) < \infty \quad \mbox{ and } \quad \int_{\sS_N\cap \cC}\frac{1}{\|u^1-u^2\|}\beta^N_\theta (\dd u) < \infty''
$$
holds for all $N\ge 2$. First $\cP(2)$ clearly holds. Indeed, $d_{\theta,2}(2) = 2-\theta > 1$ iff $\theta <1$, and for every compact set $\cC$ of $E_3 = (\rr^2)^N$, if $\theta<1$,
$$
\int_{\cC} \|x^1-x^2\|^{-1-2\theta/2} \dd x <\infty. 
$$
Moreover the integral over $\sS_2$ is trivially finite since for all $u\in \sS_2$, $\|u^1-u^2\|^2 = 4\|u^1\|^2 = 4\|u^2\|^2 = 2$. 

\vip 

Let $N\ge 3$ such that $\cP(k)$ holds true for every $k\in \ig 2,N-1\id$. We prove $\cP(N)$. Let $\theta \in (0,N)$ and $k\ge 2$ such that $d_{\theta,N}(k) >1$.

\vip

{\it Step 1.} Let $p\in \ig 2,N\id$ and $\cC$ a compact set of $E_p$. Set $\gK^p_N = \{ \cK\in \gK_N : \mbox{ for all } K\in \cK, |K|\le p\}$. We show that there exists $c_\cC$ a positive finite constant such that
$$
\cC \subset \bigcup_{\cK \in \gK^{p-1}_N} \Big\{ x\in (\rr^2)^N : \min_{K \in \cK} \min_{i\in K, j\notin K} \|x^i-x^j\|\ge c_\cC \Big\}.
$$
However this is clear since the application
$$
f : x\in \cC \mapsto \sum_{\cK \in \gK^{p-1}_N} \min_{K \in \cK} \min_{i\in K, j\notin K} \|x^i-x^j\|
$$
does not vanish on $\cC$. Since $f$
is continuous on the compact set $\cC$, there exists a constant $\tilde c_\cC$ such that $f(x)\ge \tilde c_\cC$ for all $x\in \cC$. Thus for all $x\in \cC$, there exists $\cK \in \gK^p$ such that 
$$
\min_{K \in \cK} \min_{i\in K, j\notin K} \|x^i-x^j\| \ge \frac{\tilde c_\cC}{|\gK^p|} := c_\cC.
$$

\vip

{\it Step 2.} We now prove that for all $\cC$ compact set of $E_{k+1}\cap E_N$, 
$$
\tilde I := \int_\cC \|x^1-x^2\|^{-1} \mu^N_\theta (\dd x) <\infty.
$$
According to Step 1, there exists $c_\cK>0$ be such that
$$
\cC \subset \bigcup_{K\subset \ig 1,N\id : |K|\ne N} F_K, \quad \mbox{ with } \quad F_K =\Big\{ x\in (\rr^2)^N : \min_{i\in K,j\notin K} \|x^i-x^j\| \ge c_\cC  \Big\}.
$$
This implies that $\tilde I \le \sum_{K : |K|\notin \{0,N\}}\tilde I_K$ with
\begin{align*}
    \tilde I_K = \int_{\cC\cap F_K} \|x^1-x^2\|^{-1} \mu^N_\theta(\dd x).
\end{align*}
Let $K\subset \ig 1,N\id$ with $|K|\ne N$, it remains to show that $\tilde I_K <\infty$. If $1,2$ does not belong to the same set among $K$ and $K^c$, then  we easily conclude according to Lemma \ref{lem:radon}-(i). We can thus suppose without loss of generality that $1,2\in K$. We have,
\begin{align*}
    \tilde I_K \le C_\cC \Big( \int_{\cC} \|x^1-x^2\|^{-1} \prod_{i\ne j\in K} \|x^i-x^j\|^{-\tilde\theta/|K|} \dd x\Big) \Big( \int_{\cC} \prod_{i\ne j\in K^c} \|x^i-x^j\|^{-\tilde\theta^c/|K^c|}\dd x\Big), 
\end{align*}
where $\tilde \theta = \theta |K|/N$, $\tilde \theta^c = \theta |K^c|/N$ and $C_\cC$ is a positive and finite constant. The last term of the product is bounded according to Lemma \ref{lem:radon}-(i) and Remark \ref{rk:unicitepreKS}-(v) since $d_{\tilde \theta^c,|K^c|}(k) = d_{\theta,N}(k) >1$ and $\tilde \theta^c /|K^c| = \theta/N$. The first term is also bounded by induction hypothesis according to Remark \ref{rk:unicitepreKS}-(v) again since $\cC$ is a compact set of $E_{k+1}$ seen as a subset of $(\rr^2)^{|K|}$.

\vip

From now we fix $\theta \in (0,N)$, $k\ge 2$ such that $d_{\theta,N}(k)>1$ and $\cC$ a compact set of $E_{k+1}$.

\vip

{\it Step 3.} We prove that for all compact set $\cC_\sS$ of $E_{k+1}\cap \sS$,
$$
\int_{\cC_\sS} \|u^1-u^2\|^{-1} \beta^N_\theta (\dd u) <\infty,
$$
Let $\cC' := \{ \gamma (z)+ru,,z\in B_{\rr^2}(0,1), u\in \cC_\sS, r\in [1,2]\}$ which is a compact set of $E_{k+1}\cap E_N$. According to Step 2, 
$$
\int_{\cC'}  \|x^1-x^2\|^{-1} \mu^N_\theta (\dd x) <\infty.
$$
Moreover, according to \cite[Lemma A.2]{ft} we have
$$
\int_{\cC'}  \|x^1-x^2\|^{-1} \mu^N_\theta (\dd x) = \frac12\mathrm{Leb}(B_{\rr^2}(0,1)) \Big( \int_1^2 r^{\nu} \dd r \Big) \Big( \int_{\cC_\sS}  \|u^1-u^2\|^{-1} \beta^N_\theta (\dd u)\Big),
$$
where $\nu = d_{\theta,N}(N)/2-1$, which concludes. 
\vip

{\it Step 4.} We conclude the proof here by showing that 
$$
I := \int_\cC \|x^1-x^2\|^{-1} \mu^N_\theta (\dd x) <\infty.
$$
According to Step 1, we get $I\le \sum_{\cK \in \gK^{k}_N} I_\cK$ where
\begin{align*}
    I_\cK = I := \int_{\cC\cap F_\cK} \|x^1-x^2\|^{-1} \mu^N_\theta (\dd x), \quad \mbox{ with } \quad F_\cK =\Big\{ x\in (\rr^2)^N : \min_{K\in \cK}\min_{i\in K,j\notin K} \|x^i-x^j\| \ge c_\cC \Big\},
\end{align*}
with $c_\cC$ a positive and finite constant. We fix $\cK \in \gK^{k}_N$ and we conclude by showing that $I_\cK < \infty$. Here also, one can suppose without loss of generality that $1,2$ belong to some fixed $K_0\in \cK$ according to Lemma \ref{lem:radon}-(i). 

\vip

By definition of $F_\cK$ we get the existence of $C_\cC$ positive and finite such that
\begin{align*}
    I_\cK \le C_\cC \Big( \int_{\cC} \|x^i-x^j\|^{-1}\prod_{i\ne j\in K_0}  \|x^i-x^j\|^{-\theta/N}\dd x \Big) \prod_{K\ne K_0\in \cK} \Big( \int_{\cC} \prod_{i\ne j\in K} \|x^i-x^j\|^{-\theta/N}\dd x \Big).
\end{align*}
According to Lemma \ref{lem:radon}-(i) and Remark \ref{rk:unicitepreKS}-(v) again, the only non trivially finite integral in the above is the first one which we denote $I_{\cK,K_0}$. We use \cite[Lemma A.2]{ft} and get
\begin{align*}
    I_{\cK,K_0} = \frac12\int_{\rr^2\times \rr_+\times \sS_{|K_0|}} \indiq_{\Psi (z,r,u) \in \cC} \|u^1-u^2\|^{-1} r^\alpha \dd z \dd r \prod_{i\ne j\in K_0}\|u^i-u^j\|^{-\theta/N}\sigma (\dd u),
\end{align*}
where $\Psi (z,r,u) = \gamma (z) + \sqrt{r}u$, $\alpha = d_{|K_0|\theta/N,|K_0|}(|K_0|)/2-3/2 = d_{\theta,N}(|K_0|) /2 - 3/2$ according to Remark \ref{rk:unicitepreKS}-(v). Since $\cC$ is compact, one can find $M_\cC>0$ such that $\indiq_{\Psi (z,r,u) \in \cC} \le \indiq_{|z|\le M_\cC}\indiq_{r\le M_\cC}$. This implies that
\begin{align*}
    I_{\cK,K_0} \le \Big(\int_{\rr^2} \indiq_{|z|\le M_\cC}\dd z\Big) \Big(\int_0^{M_\cC} r^\alpha \dd r\Big) \Big( \int_{\sS_{|K_0|}} \|u^1-u^2\|^{-1} \beta^{|K_0|}_{\theta |K_0|/N} (\dd u) \Big).
\end{align*}
The first integral is clearly finite. According to Lemma \ref{lem:dimension}, since $d_{\theta,N}(k)>1$ and $|K_0|\le k$, we have $d_{\theta,N}(|K_0|)>1$ so that $\alpha >-1$, implying that the second integral is finite. Let us conclude by justifying the finiteness of the last integral according to two scenarios. If first $|K_0|<N\wedge (k+1)$, then the induction hypothesis concludes because $\sS_{|K_0|}$ is a compact set of $E_{k+1}$, $d_{\theta |K_0|/N,|K_0|}(|K_0|) = d_{\theta,N}(|K_0|) >1$ (see the reasoning above). Secondly, in the case where $|K_0| = N <k+1$, $\sS_{|K_0|}$ is a compact set of $E_{k+1}$ and we can use Step 3 and \cite[Lemma A.2]{ft} to conclude. 

\end{proof}
\subsection{Computation of the explicit density starting from an all-particle collision}\label{app:densité}

We now prove Proposition \ref{prop:densité}
\begin{proof}
   Fix $\theta$, $N\ge 2$ as in the statement and $m\in \rr^2$. According to Lemma \ref{lem:symdesangles}, we get for all $\varphi \in C((\rr^2)^N)$,
   \begin{align*}
        \E^N_{\gamma (m)}[\varphi(X_t)] =\int \E[\varphi (\gamma (M_t) + D_t^{1/2}u') ]Z^{-1}\beta^N_\theta(\dd u'),
    \end{align*}
    where $(M_t)_{t\ge 0}$ is a $2$-dimensional Brownian motion starting at $m$ and $(D_t)_{t\ge 0}$ is a $d_{\theta,N}(N)$-dimensional squared Bessel process starting at $0$. According to \cite{ry}, we get, setting $\nu = d_{\theta,N}(N)/2-1$, $\E^N_{\gamma (m)}[\varphi(X_t)] $ equals
    \begin{align*}
        \frac{2^{-\nu}Z^{-1}N}{2\pi t^{2+\nu} \Gamma (\nu+1)}\int_{\rr^2}\int_0^\infty\int_{\sS_N} \varphi (\gamma (m') + ru')e^{-N\|m'-m\|^2/(2t)}r^{2\nu+1}e^{-r^2/(2t)}\dd m' \dd r\beta^N_\theta(\dd u').
    \end{align*}
    Using \cite[Lemma A.2]{ft} we deduce that $\E^N_{\gamma (m)}[\varphi(X_t)] $ equals
    \begin{align*}
        \frac{2^{-\nu}Z^{-1}}{2\pi t^{2+\nu} \Gamma (\nu+1)}&\int_{(\rr^2)^N} \varphi (x) e^{-N\|S_{\ig 1,N\id}(x)-m\|^2/(2t)}e^{-R_{\ig 1,N\id}(x)/(2t)}\mu^N_\theta(\dd x)\\
        &=\frac{Z^{-1}}{2^{\nu+1}\pi t^{2+\nu} \Gamma (\nu+1)}\int_{(\rr^2)^N} \varphi (x)e^{-\|x - \gamma (m)\|^2/(2t)}\mu^N_\theta(\dd x),
    \end{align*}
   where we used Lemma \ref{lem:dispersionformula} in the last line.
\end{proof}
\section*{Acknowledgements}
The author would like to express his sincere gratitude to Nicolas Fournier for insightful discussions, to Loïc Béthencourt for fruitful exchanges focusing on Lévy processes and excursion theory and Louis-Pierre Chaintron for many conversations about entropy on path space which inspired Appendix \ref{sec:unicitefaible}. Moreover some AI has been used in the process to grasp some ideas concerning the Dirichlet form aspect of this work and to review the paper. 

\bibliographystyle{plainnat}
\bibliography{references}

\end{document}